%% file: top.tex
\documentclass[lecture]{pcms-l}
\input pcmslmod-modified.tex

\def\currentyear{2004}
\usepackage{amssymb,amscd,verbatim,ifthen,amsthm}
\usepackage{graphicx}
\input{CrayolaColors}

\newtheorem{theorem}{Theorem}[section]

\newtheorem{prop}[theorem]{Proposition}
\newtheorem{cor}[theorem]{Corollary}
\newtheorem{con}[theorem]{Conjecture}

\theoremstyle{definition}
\newtheorem{definition}[theorem]{Definition}
\newtheorem{example}[theorem]{Example}
\newtheorem{xca}[theorem]{Exercise}
\newtheorem{prob}[theorem]{Problem}

\theoremstyle{remark}
\newtheorem{remark}[theorem]{Remark}

\numberwithin{section}{chapter}
\numberwithin{equation}{section}
\numberwithin{figure}{section}

\begin{document}


\part*{\hbox{Poset Topology: Tools and Applications}}
\pauth{Michelle L. Wachs\\[5in]
{\small \sc IAS/Park City Mathematics Institute, Summer 2004}}

\vbox{
\tableofcontents
}

\mainmatter
\setcounter{page}{3}

\LogoOn

\lectureseries[Poset Topology]{\hbox{Poset Topology: 
Tools and Applications}}
\auth[Wachs]{Michelle L. Wachs}
\address{Department of Mathematics, 
University of Miami, 
Coral Gables, Fl 33124}
\email{wachs@math.miami.edu}
\thanks{This work was partially supported by NSF grant DMS 0302310}

\setaddress

\input macro

\lecture*{Introduction}

  The   theory of poset topology evolved from the
seminal 1964  paper of  Gian-Carlo Rota on the M\"obius function of a partially
ordered set.  This theory   provides
a deep and fundamental link between  combinatorics and other branches
of mathematics.  
Early impetus for this theory   came from diverse fields such as
\begin{itemize}
\item commutative algebra (Stanley's 1975 proof of the upper bound conjecture)
\item 
group theory  (the work of Brown (1974) and Quillen (1978)  on
$p$-subgroup posets)
\item combinatorics 
(Bj\"orner's 1980 paper on poset shellability)
\item representation theory (Stanley's 1982 paper on group actions on the homology of
posets)
\item topology
 (the Orlik-Solomon theory of hyperplane arrangements (1980))
\item complexity theory (the 1984 paper of 
 Kahn, Saks, and Sturtevant on the evasiveness
conjecture).
\end{itemize}  Later developments  have kept the theory vital. I
mention just a few examples: 
Goresky-MacPherson formula for  subspace arrangements, Bj\"orner-Lov\'asz-Yao
 complexity
theory  results,
Bj\"orner-Wachs extension of  shellability to nonpure complexes,
 Forman's  discrete version of Morse theory, and
 Vassiliev's work on knot invariants and graph connectivity.

So,  what is poset topology?  By the topology  of a partially ordered set (poset)
we mean the topology  of a certain simplicial complex associated with the poset, called
the order complex of the poset.  In these lectures I will present some  of the 
techniques  that have been developed over the years to study the topology of a poset,
and    discuss some of the  applications of poset topology to the fields mentioned above
as well as to  other fields.  In particular, I will  discuss  tools for computing 
homotopy type and (co)homology of posets, with an emphasis on group equivariant
(co)homology. 
 Although  posets and
simplicial complexes can be viewed as essentially the same topological object, we will
narrow our focus, for the most part, to tools that were developed specifically for
posets;   for example, lexicographical shellability, recursive atom orderings, Whitney
homology  techniques,
(co)homology bases/generating set techniques,  and fiber theorems.

Research in poset topology is very much driven by  the study of concrete
examples that arise in various contexts both inside and outside of combinatorics.  
These examples often turn
out to have a rich and interesting topological structure, whose analysis leads 
to the development of
new techniques in poset topology.
 These lecture notes are organized according 
to  techniques  rather than applications.  A recurring theme is  the use of original
examples in demonstrating a technique, where by original example I mean the example that
led to the development of the technique in the first place.   More recent examples  will 
be discussed as well.

With regard to the choice of topics, I was primarily motivated by my own research
interests and the desire to provide the students    at the PCMI graduate school with 
concrete skills in this subject.  Due to space and time constraints  and my decision to
focus on techniques specific to posets, there are a number of  very important tools for
general simplicial complexes that I have only been able to mention in passing (or not
 at all). I point out, in particular,
 discrete Morse theory (which is  a major part of the lecture
series of Robin Forman, its originator) and  basic techniques from algebraic topology
such as long exact sequences and  spectral sequences.    For
further techniques and applications, still of current interest, we strongly recommend
 the influential 1995 book chapter of Anders Bj\"orner \cite{bj95}.

 The exercises
vary in difficulty and are there to reinforce and supplement the material treated in these
notes.  There are many open problems (simply referred to as problems) and conjectures
sprinkled throughout the text.   

I would like to thank the organizers (Ezra Miller,  Vic Reiner and Bernd Sturmfels) of
 the 2004 PCMI Graduate Summer School for inviting me to deliver these lectures.     I am
very  grateful to Vic Reiner for  his   encouragement and support.  I would
also like to thank  Tricia Hersh for  the help and  support she provided as my
overqualified teaching assistant.  Finally, I would like to express my gratitude  to the
graduate students at the summer school for their interest and inspiration.

\input lect1

\input lect2

\input lect3

\input lect4

\input lect5

\input biblio

\end{document}

%% file: pcmslmod-modified.tex
\makeatletter

\newif\iffirstlecture\firstlecturefalse

\newcommand{\lectureseries}{\firstlecturetrue
              \secdef\@lectureseries\@slectureseries} 

\newcommand{\@lectureseries}[2][default]{\chapter*{#2}%
              \gdef\thelectureseries{#1}} 

\newcommand{\@slectureseries}[1]{\chapter*{#1}}

\renewcommand{\auth}{\secdef\@auth\@sauth}

\newcommand{\@auth}[2][default]{\vspace{-1pc}{\raggedleft
        \Large\bf\noindent
        #2\endgraf
        \vspace*{2pc}
        }
        \def\@author{#1}%
}

\newcommand{\@sauth}[1]{\vspace{-1pc}{\raggedleft
        \Large\bf\noindent
        #1\endgraf
        \vspace*{2pc}
        }
        \def\@author{#1}%
}

\def\lecture#1{\global\Lecturetrue\global\Monographfalse
\iffirstlecture\else\chapter*{}\fi\firstlecturefalse
  \global\let\sectionmark\@gobble 
  \addtocounter{lecture}1\relax
  \refstepcounter{chapter}%
\gdef\thelecturename{#1\unskip}
  {\Large\bfseries
   \raggedleft
   \@xp\uppercase\@xp{\thelecturelabel} {\LARGE\thelecturenum}\\
   \vspace*{3pt}%
   \thelecturename
   \endgraf}%
  \let\@secnumber=\thelecturenum
  \@xp\lecturemark\@xp{\thelecturename}%
  \addcontentsline{toc}{chapter}{%
    \thelecturelabel\ \thelecturenum.\ \thelecturename}%
  \vspace{34\p@}\noindent}
  
\def\lecture{\global\Lecturetrue\global\Monographfalse
\iffirstlecture\else\chapter*{}\fi%
  \global\let\sectionmark\@gobble 
\secdef\@lecture\@slecture}

\def\@lecture[#1]#2{%
  \addtocounter{lecture}1\relax
  \refstepcounter{chapter}%
\gdef\thelecturename{#1\unskip}\firstlecturefalse
  {\Large\bfseries
   \raggedleft
   \@xp\uppercase\@xp{\thelecturelabel} {\LARGE\thelecturenum}\\
   \vspace*{3pt}%
    #2\unskip
   \endgraf}%
  \let\@secnumber=\thelecturenum
  \@xp\lecturemark\@xp{\thelecturename}%
  \addcontentsline{toc}{chapter}{%
    \thelecturelabel\ \thelecturenum.\ #2}%
  \vspace{34\p@}\noindent}
  
\def\slecturerunhead#1#2#3{%
    \let\@tempa\chaptername
    \uppercasenonmath{\@tempa}%
    \def\@tempb{#3\unskip}%
    \uppercasenonmath{\@tempb}%
    {\normalfont\@tempb}
    }
\def\slecturemark{
    \@secmark\markright\slecturerunhead\chaptername}%

\def\@slecture#1{%
\iffirstlecture
\gdef\thelecturename{#1\unskip}\firstlecturefalse
  {\Large\bfseries
\noindent\thelecturename
   \endgraf}%
  \let\@secnumber=\thelecturenum
  \@xp\slecturemark\@xp{\thelecturename}%
  \addcontentsline{toc}{chapter}{%
    \thelecturename}%
 \vspace{-6\p@}\noindent
\else
\gdef\thelecturename{#1\unskip}\firstlecturefalse
  {\Large\bfseries
   \raggedleft
   \@xp\uppercase\@xp{\thelecturename}
   \endgraf}%
  \let\@secnumber=\thelecturenum
  \@xp\slecturemark\@xp{\thelecturename}%
  \addcontentsline{toc}{chapter}{%
    \thelecturename}%
  \vspace{34\p@}\noindent
\fi}

\ifLecture
  \def\chapterrunhead#1#2#3{%
    \let\@tempa\@author
    \uppercasenonmath{\@tempa}%
    \uppercasenonmath{\thelectureseries}%
    \textmd{\@tempa, \thelectureseries}
    }
  \def\lecturerunhead#1#2#3{%
    \let\@tempa\chaptername
    \uppercasenonmath{\@tempa}%
    \def\@tempb{#3\unskip}%
    \uppercasenonmath{\@tempb}%
    \textmd{\@tempa\ #2. \@tempb}
    }
\else
  \let\chapterrunhead\partrunhead
\fi

\newif\ifBibliographyIsASection\BibliographyIsASectionfalse

  \def\bibliomark{
    \@secmark\markright\bibliorunhead\chaptername}%

  \def\bibliorunhead#1#2#3{%
    \let\@tempa\chaptername
    \uppercasenonmath{\@tempa}%
    \def\@tempb{#3\unskip}%
    \uppercasenonmath{\@tempb}%
    \textmd{\@tempb}
    }

\def\thebibliography#1{%
  \ifBibliographyIsASection
    \section*\refname
    \if@backmatter
      \markboth{\refname}{\refname}%
    \fi
  \else
\chapter*{}
  {\Large\bfseries
   \raggedleft
   \@xp\uppercase\@xp{\bibname} \\
   \endgraf}%
  \let\@secnumber=\thelecturenum
  \@xp\bibliomark\@xp{\bibname}%
  \addcontentsline{toc}{chapter}{%
    \bibname}%
  \vspace{34\p@}\noindent
  \fi
  \normalsize\labelsep .5em\relax
  \list{\@arabic\c@enumi.}{\settowidth\labelwidth{\@biblabel{#1}}%
  \leftmargin\labelwidth
  \advance\leftmargin\labelsep
	\usecounter{enumi}}\sloppy
  \clubpenalty9999 \widowpenalty\clubpenalty  \sfcode`\.\@m}

  \def\indexmark{
    \@secmark\markright\indexrunhead\chaptername}%

  \def\indexrunhead#1#2#3{%
    \let\@tempa\chaptername
    \uppercasenonmath{\@tempa}%
    \def\@tempb{#3\unskip}%
    \uppercasenonmath{\@tempb}%
    \textmd{\@tempb}
    }

\ifLecture
\def\theindex{\cleardoublepage
\@restonecoltrue\if@twocolumn\@restonecolfalse\fi
\columnseprule \z@ \columnsep 35pt
\def\indexchap{\@startsection
		{chapter}{1}{\z@}{8pc}{34pt}%
		{\raggedleft
		\Large\bfseries}}%
 \twocolumn[\indexchap[{\indexname}]{\@xp\uppercase\@xp{\indexname}}]
  \@xp\indexmark\@xp{\indexname}%
	\thispagestyle{plain}\let\item\@idxitem\parindent\z@
	 \footnotesize\parskip\z@ plus .3pt\relax\let\item\@idxitem}
\fi

\def\@makefntext{\noindent\@makefnmark}

\def\setaddress{%
  {\let\@makefnmark\relax  \let\@thefnmark\relax
        \nobreak
        \addressnum@=\z@
        \loop\ifnum\addressnum@<\addresscount@\advance\addressnum@\@ne
           \footnote{
           \csname @address\number\addressnum@\endcsname
           \csname @curraddr\number\addressnum@\endcsname
           \csname @email\number\addressnum@\endcsname}\repeat
  \ifx\@empty\@date\else \@footnotetext{\@setdate}\fi
  \ifx\@empty\@subjclass\else \@footnotetext{\@setsubjclass}\fi
  \ifx\@empty\@keywords\else \@footnotetext{\@setkeywords}\fi
  \ifx\@empty\thankses\else \@footnotetext{%
    \def\par{\let\par\@par}\@setthanks}\fi
    }%
  \@setcopyright
}

\def\@tmpevenhead{\relax}

\def\cleardoublepage{\clearpage\if@twoside \ifodd\c@page\else
    \let\@tmpevenhead\@evenhead \let\@evenhead\relax\hbox{}\eject 
    \let\@evenhead\@tmpevenhead\if@twocolumn\hbox{}\newpage\fi\fi\fi}

\def\@setcopyright{%
  \let\copyrightyear\currentyear             
  \insert\copyins{\hsize\textwidth
    \parfillskip\z@ \leftskip\z@\@plus.9\textwidth
    \fontsize{6}{7\p@}\normalfont\upshape
    \everypar{}%
    \vskip-\skip\copyins \nointerlineskip
    \noindent\vrule\@width\z@\@height\skip\copyins
    \copyright\copyrightyear\ 
Michelle L. Wachs\par
    \kern\z@}%
}

\renewcommand{\@auth}[2][default]{{\raggedleft
        \begingroup
  \fontsize{\@xivpt}{18}\bfseries
  #2\par \endgroup
        \vspace*{2pc}
        }
        \def\@author{#1}%
}

\renewcommand{\@sauth}[1]{{\raggedleft
        \begingroup
  \fontsize{\@xivpt}{18}\bfseries
  #1\par \endgroup
        \vspace*{2pc}
        }
        \def\@author{#1}%
}

\def\@lecture[#1]#2{%
  \addtocounter{lecture}1\relax
  \refstepcounter{chapter}%
\gdef\thelecturename{#1\unskip}\firstlecturefalse
  {\Large\bfseries
   \raggedleft
   \@xp\uppercase\@xp{\thelecturelabel} {\LARGE\thelecturenum}\\
   \vspace*{3pt}%
    #2\unskip
   \endgraf}%
  \let\@secnumber=\thelecturenum
  \@xp\lecturemark\@xp{\thelecturename}%
  \addcontentsline{toc}{chapter}{%
    \thelecturelabel\ \thelecturenum.\ #2}%
  \vspace{10\p@}\noindent}
  
\def\@slecture#1{%
\iffirstlecture
\gdef\thelecturename{#1\unskip}\firstlecturefalse
  {\Large\bfseries
\noindent\thelecturename
   \endgraf}%
  \let\@secnumber=\thelecturenum
  \@xp\slecturemark\@xp{\thelecturename}%
  \addcontentsline{toc}{chapter}{%
    \thelecturename}%
 \vspace{-6\p@}\noindent
\else
\gdef\thelecturename{#1\unskip}\firstlecturefalse
  {\Large\bfseries
   \raggedleft
   \@xp\uppercase\@xp{\thelecturename}
   \endgraf}%
  \let\@secnumber=\thelecturenum
  \@xp\slecturemark\@xp{\thelecturename}%
  \addcontentsline{toc}{chapter}{%
    \thelecturename}%
  \vspace{10\p@}\noindent
\fi}

\makeatother

%% file: CrayolaColors
\def\Red{} 
\def\Black{} 
\def\White{} 

%% file: macro.tex
\def\<{<\!\!\!\cdot\,}
\def\k{{\mathbf k}}
\def\C{{\mathbb{C}}}
\def\Q{{\mathbb{Q}}}
\def\R{{\mathbb R}}
\def\Z{{\mathbb{Z}}}
\def\N{{\mathbb{N}}}
\def\p{{\mathbb{P}}}
\def\S{{\mathbb{S}}}
\def\lk{{\rm lk}}
\def\ch{{\rm ch}}
\def\im{{\rm im}}
\def\sh{{\rm sh}}
\def\sgn{{\rm sgn}}
\def\wh{{\rm WH}}
\def\des{{\rm des}}
\def\s{{\mathfrak S}}
\def\ltimes{\mathop{\Large \textsf{X}}}

\newcommand\begth{\begin{theorem}}
\newcommand\enth{\end{theorem}}
\newcommand\beq{\begin{eqnarray*}}
\newcommand\eeq{\end{eqnarray*}}
\newcommand\bq{\begin{eqnarray}}
\newcommand\eq{\end{eqnarray}}
\newcommand\ben{\begin{enumerate}}
\newcommand\een{\end{enumerate}}
\newcommand\bit{\begin{itemize}}
\newcommand\eit{\end{itemize}}
\newcommand\itm{\vspace{-.4in}\item}
\newcommand\bitm{\vspace{-.5in}\item}

\def\multi#1{\vbox{\baselineskip=0pt\halign{\hfil$\scriptstyle\vphantom{(_)}##$\hfil\cr#1\crcr}}}
\def\picture #1 by #2 (#3){
  \vbox to #2{
    \hrule width #1 height 0pt depth 0pt
    \vfill
    \special{picture #3} 
    }
  }

\def\scaledpicture #1 by #2 (#3 scaled #4){{
  \dimen0=#1 \dimen1=#2
  \divide\dimen0 by 1000 \multiply\dimen0 by #4
  \divide\dimen1 by 1000 \multiply\dimen1 by #4
  \picture \dimen0 by \dimen1 (#3 scaled #4)}
  }

\edef\savecatcodeat{\the\catcode`@}
\catcode`\@=11

\def\tb@ifSpecChars#1#2{#1}
\def\tb@ifNoSpecChars#1#2{#2}

\def\tableau{%
  \bgroup
  \@ifstar{\let\Tif\tb@ifNoSpecChars\tb@tableauB}
          {\let\Tif\tb@ifSpecChars\tb@tableauB}}

\def\tb@tableauB{
  \@ifnextchar[{\tb@tableauC}{\tb@tableauC[]}}

\def\tb@tableauC[#1]{\hbox\bgroup%
    \let\\=\cr
    \def\bl{\global\let\tbcellF\tb@cellNF}%
    \def\tf{\global\let\tbcellF\tb@cellH}
%
    \dimen2=\ht\strutbox \advance\dimen2 by\dp\strutbox%
    \ifx\baselinestretch\undefined\relax%
    \else%
       \dimen0=100sp \dimen0=\baselinestretch\dimen0%
       \dimen2=100\dimen2 \divide\dimen2 by\dimen0%
    \fi%
    \let\tpos\tb@vcenter
    \tb@initYoung
    \tb@options#1\eoo
    \let\arrow\tb@arrow%
    \dimen0=\Tscale\dimen2%
    \dimen1=\dimen0 \advance\dimen1 by \tb@fframe%
    \lineskip=0pt\baselineskip=0pt
%
    \def\tb@nothing{}%
    \def\endcellno{$\rss\egroup\bss\egroup}
    \def\endcell{\endcellno\kern-\dimen0}
    \def\begincell{\vbox to\dimen0\bgroup\vss\hbox to\dimen0\bgroup\hss$}%
    \let\overlay\tb@overlay%
    \let\fl\tb@fl%
    \let\lss\hss\let\rss\hss\let\tss\vss\let\bss\vss
    \def\mkcell##1{
        \let\tbcellF\tb@cellD
        \def\tb@cellarg{##1}
        \ifx\tb@cellarg\tb@nothing\let\tb@cellarg\tb@cellE\fi%
        \begincell\tb@cellarg\endcellno
        \tbcellF}
    \let\savecellF\tbcellF
     \Tif{\catcode`,=4\catcode`|=\active}{}\tb@tableauD}%

\let\tb@savetableauD\tableauD
{
    \catcode`|=\active \catcode`*=\active \catcode`~=\active%
    \catcode`@=\active
\gdef\tableauD#1{%
  \Tif{
    \mathcode`|="8000 \mathcode`*="8000%
    \mathcode`~="8000 \mathcode`@="8000%
    \def@{\bullet}%
    \let|\cr
    \let*\tf
    \let~\sk
  }{}%
  \tpos{\tabskip=0pt\halign{&\mkcell{##}\cr#1\crcr}}%
  \global\let\tbcellF\savecellF
  \egroup
  \egroup}
}
\let\tb@tableauD\tableauD
\let\tableauD\tb@savetableauD
\let\tb@savetableauD\undefined


\def\tb@options#1{\ifx#1\eoo\relax\else\tb@option#1\expandafter\tb@options\fi}

\def\tb@option#1{%
  \if#1t\let\tpos\tb@vtop\fi
  \if#1c\let\tpos\tb@vcenter\fi
  \if#1b\let\tpos\vbox\fi
  \if#1F\tb@initFerrers\fi
  \if#1Y\tb@initYoung\fi
  \if#1s\tb@initSmall\fi
  \if#1m\tb@initMedium\fi
  \if#1l\tb@initLarge\fi
  \if#1p\tb@initPartition\fi
  \if#1a\tb@initArrow\fi
}

\def\tb@vcenter#1{\ifmmode\vcenter{#1}\else$\vcenter{#1}$\fi}

\def\tb@vtop#1{\hbox{\raise\ht\strutbox\hbox{\lower\dimen0\vtop{#1}}}}

\def\tb@initPartition{\def\Tscale{.3}}
\def\tb@initSmall{\def\Tscale{1}}
\def\tb@initMedium{\def\Tscale{2}}
\def\tb@initLarge{\def\Tscale{3}}

\def\tb@initArrow{\dimen2=1.25em}

\def\tb@initYoung{%
  \def\tb@cellE{}
  \let\tb@cellD\tb@cellN
  \def\sk{\global\let\tbcellF\tb@cellNF}}
\def\tb@initFerrers{%
  \def\tb@cellE{\bullet}
  \let\tb@cellD\tb@cellNF
  \def\sk{\bullet}}

\tb@initMedium

\def\tb@sframe#1{%
  \vbox to0pt{
    \vss
    \hbox to0pt{%
      \hss
      \vbox to\dimen1{
        \hrule depth #1 height0pt
        \vss
        \hbox to\dimen1{
          \vrule width #1 height\dimen1
          \hss
          \vrule width #1
          }%
        \vss
        \hrule height #1 depth 0in
        }%
      \kern-\tb@hframe
      }%
    \kern-\tb@hframe}}

\def\tb@hframe{.2pt}\def\tb@fframe{.4pt}\def\tb@bframe{2pt}
\def\tb@cellH{\tb@sframe{\tb@bframe}}       
\def\tb@cellNF{}                            
\def\tb@cellN{\tb@sframe{\tb@fframe}}       
\let\tbcellF\tb@cellN                       

\def\tb@rpad{1pt}
\def\tb@lpad{1pt}
\def\tb@tpad{1.8pt}
\def\tb@bpad{1.8pt}

\def\tb@overlay{\endcell\@ifnextchar[{\tb@overlaya}{\begincell}}
\def\tb@overlaya[#1]{\vbox to\dimen0\bgroup%
  \tb@overlayoptions#1\eoo%
  \tss\hbox to\dimen0\bgroup\lss}
\def\tb@overlayoptions#1{\ifx#1\eoo\relax\else\tb@overlayoption#1\expandafter\tb@overlayoptions\fi}

\def\tb@overlayoption#1{
  \if#1t\def\tss{\vskip\tb@tpad}\let\bss\vss\fi
  \if#1c\let\tss\vss\let\bss\vss\fi
  \if#1b\def\bss{\vskip\tb@bpad}\let\tss\vss\fi
  \if#1l\def\lss{\hskip\tb@lpad}\let\rss\hss\fi
  \if#1m\let\lss\hss\let\rss\hss\fi
  \if#1r\def\rss{\hskip\tb@rpad}\let\lss\hss\fi
}

\def\tb@fl{\endcell\begincell\vrule depth 0pt width \dimen0 height \dimen0 \endcell\begincell}



\@ifundefined{diagram}{}{
\def\tb@arrowpad{.5}

\newoptcommand{\tb@arrow}{\@ne}[2]{%
  \endcell
   \begingroup%
   \let\dg@getnodesize\tb@getnodesize
   \dg@USERSIZE=#1\relax%
   \ifnum\dg@USERSIZE<\@ne \dg@USERSIZE=\@ne \fi%
   \dg@parse{#2}%
   \dg@label{\tb@draw{#1}{#2}}}

\def\tb@getnodesize#1#2#3#4#5{\dimen3=\tb@arrowpad\dimen2 #4=\dimen3 #5=\dimen3\relax}
\def\tb@getnodesize#1#2#3#4#5{\ifnum#2=0\ifnum#3=0\tb@getnodesizetail{#4}{#5}\else\tb@getnodesizehead{#4}{#5}\fi\else\tb@getnodesizehead{#4}{#5}\fi}
\def\tb@getnodesizetail#1#2{\dimen3=.5\dimen2 #1=\dimen3 #2=\dimen3}
\def\tb@getnodesizehead#1#2{\dimen3=.5\dimen2 #1=\dimen3 #2=\dimen3}

\def\tb@draw#1#2#3#4{%
        \dg@X=0\dg@Y=0\dg@XGRID=1\dg@YGRID=1\unitlength=.001\dimen0%
        \dg@LBLOFF=\dgLABELOFFSET \divide\dg@LBLOFF\unitlength%
        \dg@drawcalc
        \begincell
        \let\lams@arrow\tb@lams@arrow
        \begin{picture}(0,0)\begingroup\dg@draw{#1}{#2}{#3}{#4}\end{picture}%
        \endcell
        \endgroup
        \begincell}
}

%
%
%
\def\tb@lams@arrow#1#2{%
 \lams@firstx\z@\lams@firsty\z@
 \lams@lastx#1\relax\lams@lasty#2\relax
 \lams@center\z@
 %
 \N@false\E@false\H@false\V@false
 \ifdim\lams@lastx>\z@\E@true\fi
 \ifdim\lams@lastx=\z@\V@true\fi
 \ifdim\lams@lasty>\z@\N@true\fi
 \ifdim\lams@lasty=\z@\H@true\fi
 \NESW@false
 \ifN@\ifE@\NESW@true\fi\else\ifE@\else\NESW@true\fi\fi
 %
 \ifH@\else\ifV@\else
  \lams@slope
  \ifnum\lams@tani>\lams@tanii
   \lams@ht\ten@\p@\lams@wd\ten@\p@
   \multiply\lams@wd\lams@tanii\divide\lams@wd\lams@tani
  \else
   \lams@wd\ten@\p@\lams@ht\ten@\p@
   \divide\lams@ht\lams@tanii\multiply\lams@ht\lams@tani
  \fi
 \fi\fi
 %
 \ifH@  \lams@harrow
 \else\ifV@ \lams@varrow
 \else \lams@darrow
 \fi\fi
}

\catcode`\@=\savecatcodeat
\let\savecatcodeat\undefined

%% file: lect1.tex



\lecture{Basic definitions, results,  and examples}



\section{Order complexes and face posets}
We begin by defining the order complex of a poset and the face poset of a
simplicial complex.  These constructions  enable us to view posets and simplicial
complexes as essentially the same topological object.  We shall assume throughout these
lectures that all posets and simplicial complexes are finite, unless otherwise stated.

 An {\em abstract simplicial complex}
$\Delta$ on finite vertex set
$V$ is a nonempty collection of subsets of
$V$ such that 
\begin{itemize}
\item $\{v\} \in \Delta$ for all $v \in V$ 
\item if $G \in \Delta$  and $F \subseteq G$ then $F \in \Delta$.
\end{itemize}  
The elements of
$\Delta$ are called {\em faces} (or  simplices) of $\Delta$ and the maximal faces are
called {\em facets}.  We say that a face $F$ has dimension $d$ and write $\dim F = d$  if
$d = |F|-1$.  Faces of dimension $d$ are referred to as $d$-faces.
 The {\em dimension} $\dim \Delta$ of $\Delta$ is defined to be $\max_{F \in
\Delta }\dim F$.  We also allow the  (-1)-dimensional complex
$\{\emptyset\}$, which we refer to as the {\em empty simplicial complex}.  It will be
convenient to refer to the empty set
$\emptyset$, as the   {\em degenerate empty  complex} and say that it has dimension
$-2$, even though we don't really consider it to be a simplicial complex.  If all
facets  of $\Delta$ have the same dimension then
$\Delta$ is said to be {\em pure}.

A $d$-dimensional {\em geometric simplex} in $\R^n$ is defined to be the   convex hull of
$d+1$  affinely  independent  points in $\R^n$ called {\em vertices}.  The convex hull of
any subset of the vertices is called a {\em face} of the geometric simplex.   A {\em
geometric simplicial complex} $K$ in
$\R^n$ is a nonempty collection of geometric simplices in 
$\R^n$ such that
\begin{itemize}
\item Every face of a simplex in $K$ is in $K$. 
\item The intersection of any two simplices of $K$ is a face of both of them.
\end{itemize}

From a geometric simplicial complex $K$, one gets an abstract simplicial complex
$\Delta(K)$ by letting  the faces of $\Delta(K)$ be the vertex sets of the 
simplices of
$K$.  Every abstract simplicial complex $\Delta$ can be obtained in this way, i.e., there
is a geometric simplicial complex $K$ such that $\Delta(K) = \Delta$.  Although $K$ is not
unique, the underlying topological space, obtained by taking the union of the simplices of
$K$ under the usual topology on $\R^n$,  is unique up to homeomorphism.  We refer to
this space as the {\em geometric realization} of
$\Delta$ and denote it by
$ \|\Delta\|$. We will usually drop the $\| \, \|$ and let
$\Delta$ denote an abstract simplicial complex as well as its geometric realization.

\begin{figure}\vspace{-.8in}
\begin{center}
\includegraphics[width=9cm]{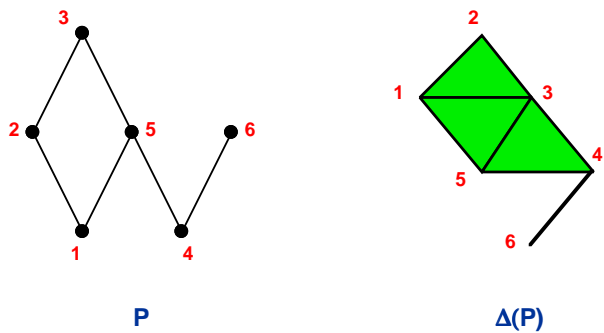}
\end{center}\begin{center}\caption{Order complex of a poset} \label{figord}
\end{center} 
\end{figure}

To every  poset
$P$, one can associate an abstract  simplicial complex $\Delta(P)$ called the {\em order
complex} of $P$.  The vertices of $\Delta(P)$ are  the elements of
$P$ and the faces of $\Delta(P)$ are the chains (i.e., totally ordered subsets) of $P$. 
(The order complex of the empty poset is the empty simplicial complex $\{\emptyset\}$.)  
For example, the Hasse diagram of a poset
$P$ and the geometric realization of  its order complex are given in Figure~\ref{figord}.  

\begin{figure} \vspace{-.3in}
\begin{center}
\includegraphics[width=10cm]{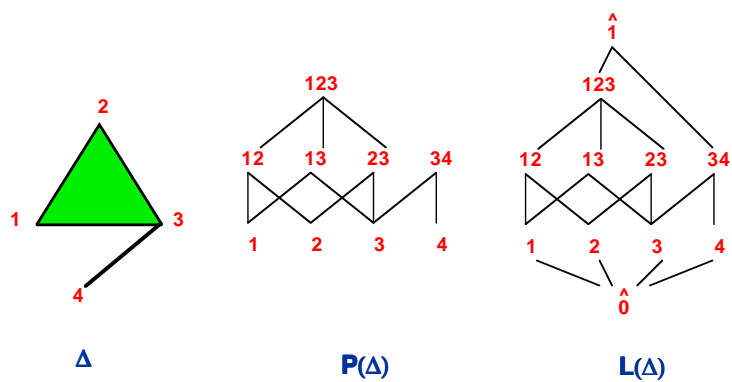}
\end{center}\begin{center}\caption{Face poset and face lattice of a simplicial complex}
\label{figfaceposet}
\end{center}
\end{figure}

To every simplicial complex $\Delta$, one can associate a poset $P(\Delta)$
called the {\em face poset} of $\Delta$, which is defined to be the poset of nonempty
faces ordered by inclusion.  The {\em face lattice} $L(\Delta)$ is $P(\Delta) $ with a
smallest element
$\hat 0$ and a largest element $\hat 1$ attached.  An example  is given in
Figure~\ref{figfaceposet}. 

\begin{figure} \vspace{-.3in}
\begin{center}
\includegraphics[width=10cm]{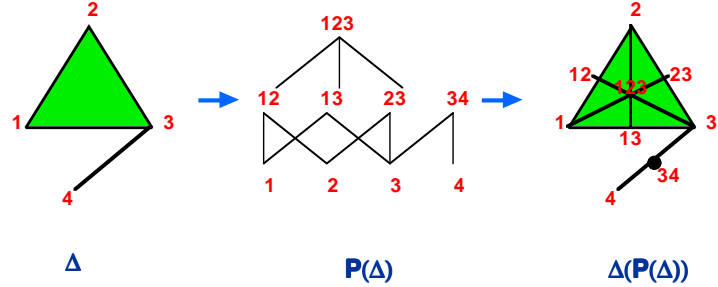}
\end{center}\begin{center}\caption{ Barycentric subdivision}
\label{figbary}
\end{center}
\end{figure}

If we start with a simplicial complex $\Delta$, take its face poset $P(\Delta)$, and
then take the order complex $\Delta(P(\Delta))$, we get a simplicial complex  known as
the {\em barycentric subdivision} of $\Delta$; see Figure~\ref{figbary}.  The geometric
realizations  are always homeomorphic,
$$\Delta \cong \Delta(P(\Delta)).$$ 

When we attribute a topological property to a poset, we mean that the geometric
realization of the order complex of the poset has that property.   For instance, if we
say that the poset
$P$ is homeomorphic to the $n$-sphere $\S^n$ we mean that $\| \Delta(P)\|$ is homeomorphic
to
$\S^n$. 

\begin{figure} 
\begin{center}
\includegraphics[width=6cm]{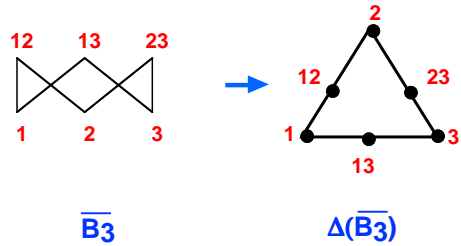}
\end{center}\begin{center}\caption{Order complex of the subset lattice (Boolean algebra)}
\label{figbool}
\end{center}
\end{figure}

\begin{example} The Boolean algebra. Let $B_n$ denote the lattice of subsets of
$[n]:=\{1,2,\dots,n\}$ ordered by containment, and let
$\bar B_n := B_n -
\{\emptyset, [n]\}$.  Then 
$$\bar B_n
\cong \S^{n-2}$$ because $\Delta(\bar B_n) $ is the barycentric subdivision of the boundary
of the
$(n-1)$-simplex.   See Figure~\ref{figbool}.
\end{example}

We now review some basic poset terminology.  An $m$-chain
of a poset $P$ is a totally ordered subset $c= \{ x_1 < x_2 <\dots < x_{m+1}\}$ of  $P$.
We say the length $l(c)$ of $c$ is $m$.    We
 consider the empty chain to be a $(-1)$-chain. The length
$l(P)$ of 
$P$ is defined to be
$$l(P) := \max\{l(c): c \mbox{ is a chain of } P\}.$$   Thus, $l(P) = \dim \Delta(P)$ and
 $l(P(\Delta)) = \dim \Delta$. 

A chain of $P$ is said to be {\em
maximal} if it is inclusionwise maximal.  Thus, the set $\mathcal M(P)$ of maximal chains
of $P$ is the set of facets of $\Delta(P)$.  
A poset $P$ is said to be {\em pure}
(also known as ranked or graded) if all maximal chains have the same length.  Thus, $P$ is
pure if and only if $\Delta(P)$ is pure.   Also a simplicial complex $\Delta$ is pure if
and only if its face poset
$P(\Delta)$ is pure.  The posets  and  simplicial complexes of
Figures~\ref{figord} and~\ref{figfaceposet} are all nonpure, while the poset  and 
simplicial complex of Figure~\ref{figbool} are both pure. 

 For
$x
\le y$ in
$P$, let 
$(x,y)$ denote the open interval $\{z \in P : x < z < y\}$ and let $[x,y]$
denote the closed interval $\{z \in P : x \le z \le y\}$. Half open 
intervals $(x,y]$ and $[x,y)$ are defined  similarly.

If $P$ has a unique minimum element, it is usual to denote it by
$\hat 0$ and refer to it as the {\em bottom} element. Similarly, the unique maximum
element, if it exists, is denoted 
$\hat 1$ and is referred to as the {\em top} element.  Note that if
$P$ has a bottom element $\hat 0$ or top element $\hat 1$ then $\Delta(P)$ is contractible
since it is a cone.  We  usually  remove the top and bottom elements and study
the more interesting topology of the remaining poset.  Define the {\em proper part} of a
poset
$P$, for which $|P| >1$, to be
$$\bar P := P - \{\hat 0,\hat 1\}. $$ In  the case that $|P| = 1$,
it will be convenient to define  
$\Delta(\bar P)$ to be the degenerate empty complex $ \emptyset$.
  We will  also say
$\Delta((x,y)) =
\emptyset$ and $l((x,y)) = -2$ if
$x = y$.  

For posets with a bottom element $\hat 0$, the  elements that cover $\hat 0$ are called
{\em atoms}.   For posets with a top element $\hat 1$, the elements that are covered by
$\hat 1$ are called {\em coatoms}.

 A poset $P$ is said
to be {\em bounded} if it has a top element $\hat 1$  and a bottom element $\hat 0$.
 Given a poset $P$, we
define the bounded extension 
$$\hat P := P \cup \{\hat 0,\hat 1\},$$
where new elements $\hat 0$ and $\hat 1$ are adjoined (even if $P$ already has a bottom or
top element).

A poset $P$ is said to be a {\em meet semilattice} if every pair of elements $x,y \in P$
has a meet
$x
\land y$, i.e. an element  less than or equal to both $x$ and $y$ that is greater
than all other such elements.  A poset $P$ is said to be a {\em join semilattice} if every
pair of elements $x,y \in P$ has a join
$x
\lor y$, i.e. a unique element  greater than or equal to both $x$ and $y$ that is less
than all other such elements.  If $P$ is both a join semilattice and a meet semilattice
then
$P$ is said to be a {\em lattice}.   It is a basic  fact of lattice theory that any
finite meet (join)  semilattice with  a top (bottom) element is a lattice.

 The {\em dual} of a 
poset $P$ is the poset $P^*$ on the same underlying set with the order relation reversed. 
Topologically there is no difference between a poset and its dual since  $\Delta(P)$ and
$\Delta(P^*)$ are identical simplicial complexes. The direct product $P \times Q$
of two posets
$P$ and $Q$ is the poset whose underlying set is the cartesian product $\{(p,q) : p
\in P,
\,  q
\in Q\}$ and whose order relation is given by   $$(p_1,q_1) \le_{P\times Q} (p_2,
q_2)\,\,\mbox{  if
} \,\, p_1 \le_P p_2\,\, \mbox{ and } \,\, q_1 \le_Q q_2.$$

Define the {\em join}   of two simplicial complexes
$\Delta$ and
$\Gamma$ on  disjoint vertex sets to be the simplicial complex given by
\bq \label{joineq} \Delta * \Gamma :=\{A \cup B : A \in \Delta , B \in
\Gamma\}.\eq  The {\em join} (or ordinal sum) $P*Q$ of   posets $P$ and $Q$ is the poset whose
underlying set is the disjoint union of $P$ and $Q$  and whose order relation is given by
$x <y$ if either (i) $x  <_P y$, (ii) $x<_Q y$, or (iii) $x \in P$ and $y \in Q$. Clearly
$$\Delta(P*Q) = \Delta(P) * \Delta(Q).$$  There are topological relationships between the
join and product of posets, which are discussed in Section~\ref{prodsec}.

\section{The M\"obius function} 
\begin{figure}
\begin{center}
\includegraphics[width=5cm]{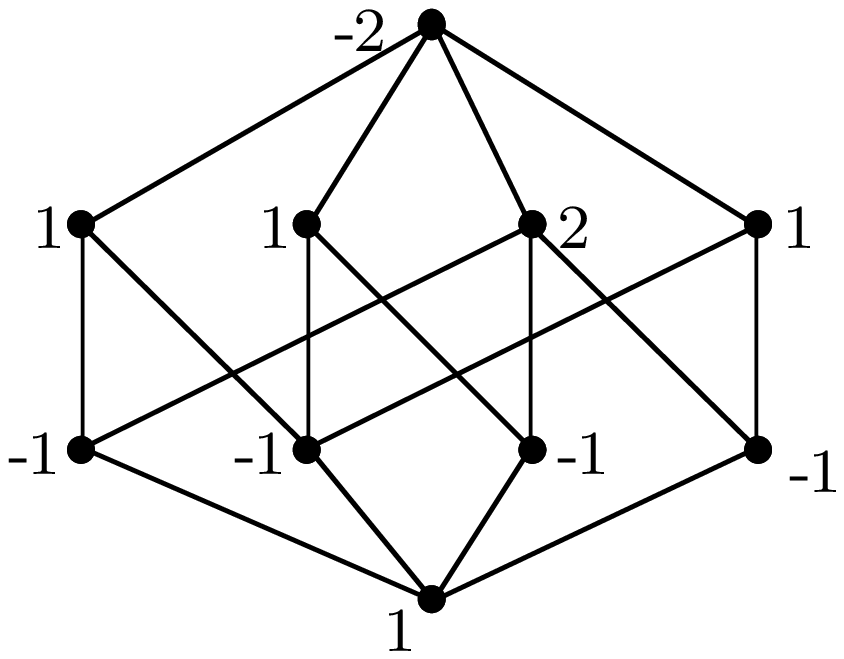}
\end{center}\begin{center}\caption{$\mu(\hat 0,x)$} \label{figmob}
\end{center}
\end{figure}

The story of poset topology  begins with the
 M\"obius function 
$\mu (=
\mu_P)$  of a poset $P$ defined recursively on closed intervals of  $P$ as follows:
\begin{eqnarray*} \mu(x,x) &=& 1, \quad \mbox{ for all } x \in P \\
\mu(x,y) &=& - \sum_{x \le z < y} \mu(x,z),  \quad \mbox{ for all } x <y \in P.
\end{eqnarray*}
For  a bounded poset $P$,  define the M\"obius invariant 
$$\mu(P) := \mu_P(\hat 0, \hat 1).$$
In  Figure~\ref{figmob}, the values of $\mu(\hat 0,x)$ are shown for each
element
$x$ of the poset.

There are various techniques for computing the M\"obius function of a poset; see
\cite{st1}. Perhaps  the most basic technique is given by the product formula.  

\begin{prop} \label{prodrule} Let $P$ and $Q$ be posets.  Then for $(p_1,q_1) \le (p_2,q_2)
\in      P\times Q$,
$$\mu_{P\times Q}((p_1,q_1),(p_2,q_2)) = \mu_P(p_1,p_2) \mu_Q(q_1,q_2).$$
\end{prop}

\begin{xca} Prove Proposition~\ref{prodrule}.
\end{xca}  

\begin{xca} Use the product formula to show that the M\"obius function for the subset
lattice
$B_n$ is given by

$$\mu(X,Y) = (-1)^{|Y-X|} \quad \mbox{ for all } X \subseteq Y.$$
\end{xca}

\begin{xca} For positive integer $n$, the lattice
$D_n$ of divisors of
$n$ is the set of positive  divisors of $n$ ordered by $a \le b$ if $a$ divides
$b$.  Show that the M\"obius function for 
$D_n$ is 
 given by
$$\mu(d,m) = \mu(m/d),$$ where $\mu(\cdot)$ is the classical M\"obius function of number
theory, which is defined on the set of positive integers by 
$$\mu(n) = \begin{cases} (-1)^k &\mbox{ if } n \mbox { is the product of $k$ distinct
primes}\\ 0 &\mbox{ if $n$ is divisible by a square}.
\end{cases}$$
This example is the reason for the name {\em M\"obius function of a poset}.
\end{xca}  

The combinatorial significance of the M\"obius function was first demonstrated
by Rota in 1964 in his Steele-prize winning paper \cite{ro64}.  The
M\"obius function of a poset is used in enumerative combinatorics to obtain inversion
formulas.  

\begin{prop}[M\"obius inversion] Let $P$ be a poset and let $f,g:P \to \C$.  Then 
$$g(y) = \sum_{x \le y} f(x)$$
if and only if 
$$f(y) = \sum_{x \le y} \mu(x,y)\, g(x).$$
\end{prop}

Three examples of
M\"obius inversion  are  classical M\"obius inversion ($P= D_n$), inclusion-exclusion ($P
= B_n$),  and  Gaussian inversion ($P= B_n(q)$, the lattice of subspaces of an
$n$-dimensional vector space over the field with $q$ elements); see \cite{st1} for
details.

Our interest in the M\"obius function stems from its connection to the Euler
characteristic.
The reduced Euler characteristic $\tilde\chi(\Delta)$ of a simplicial complex
$\Delta$ is defined to be $$ \tilde\chi(\Delta) := \sum_{i = -1}^{\dim \Delta} (-1)^i \,
f_i(\Delta),$$
where $f_i(\Delta)$ is the number of $i$-faces of $\Delta$. 

\begin{prop}[Philip Hall Theorem] \label{hall} For any  poset $P$,
$$\mu(\hat P) = \tilde\chi(\Delta( P)).$$
\end{prop}

\begin{xca} Prove Proposition~\ref{hall}.
\end{xca}

It follows from the Euler-Poincar\'e formula below that the Euler characteristic is a
topological invariant. Hence by Proposition~\ref{hall}, $\mu_P(x,y)$ depends only on the
topology of the open interval $(x,y)$ of $P$. 

\begin{theorem}[Euler-Poincar\'e formula]\label{eupo}
For any simplicial complex $\Delta$,
$$\tilde \chi(\Delta) := \sum_{i=-1}^{\dim \Delta} (-1)^i \, \tilde \beta_i(\Delta),$$
where $\tilde \beta_i(\Delta)$ is the $i$th reduced Betti number of $\Delta$, i.e., the
rank, as an abelian group, of the $i$th reduced homology of $\Delta$ over $\Z$.
\end{theorem}

The M\"obius function of a poset plays a fundamental role in the theory of hyperplane
arrangements and the homology of a poset plays a fundamental role in the theory of
subspace arrangements.  We discuss the connection with arrangements in the next section.

\section{Hyperplane and subspace arrangements}

A  {\em hyperplane arrangement} ${\mathcal A}$ is a finite collection of (affine) 
hyperplanes in
some vector space $V$.  We will consider only real hyperplane arrangements ($V= \R^n$)
and complex hyperplane arrangements ($V = \C^n$) here.  

Real hyperplane arrangements 
divide $\R^n$ into regions. 
A remarkable formula for the number of regions was given by Zaslavsky \cite{za} in 1975. 
This formula involves  the notion  of intersection
semilattice  of a hyperplane arrangement.  

\begin{figure}
\begin{center}
\includegraphics[width=12cm]{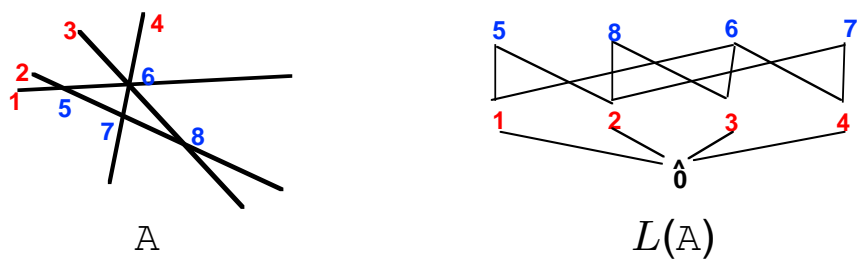}
\end{center}\begin{center}\caption{Intersection semilattice of a hyperplane arrangement}
\label{fighyper}
\end{center}
\end{figure}

The {\em intersection semilattice} $L({\mathcal A})$ of a hyperplane arrangement
${\mathcal A}$ is defined to be 
 the meet semilattice  of  nonempty intersections of hyperplanes in
${\mathcal A}$ ordered by reverse inclusion.  Note that we include the  intersection over
the empty set which is the bottom element $\hat 0$ of  $L({\mathcal A})$.  Note also that 
$L({\mathcal A})$  has a top element if and only if    $\cap \mathcal A \ne \emptyset$. 
Such an arrangement is called a {\em central arrangement}.  Hence for central arrangements
$\mathcal A$, the intersection semilattice $L(\mathcal A)$ is actually a lattice.    An
example of a hyperplane arrangement in
$\R^2$ and its intersection semilattice are given in Figure~\ref{fighyper}. 

Before stating Zaslavsky's formula, we discuss four fundamental examples of  real
hyperplane arrangements and their intersection lattices, to which we refer  throughout
these lectures.

\begin{example} \label{boolarr}  {\em The (type A) coordinate hyperplane arrangement and
the Boolean algebra
$B_n$.} The coordinate hyperplane arrangement is the central hyperplane arrangement
consisting of the coordinate hyperplanes
$x_i = 0$ in
$\R^n$.  It is easy to see that the intersection lattice of this arrangement is isomorphic
to the subset lattice $B_n$.  Indeed, the intersection  $$\{ \bold x \in
\R^n :x_{i_1} = x_{i_2} =\dots = x_{i_k} =0\},$$ where $1 \le i_1 <i_2 <\dots<i_k \le n$, 
corresponds to the subset
$\{ i_1,i_2, \dots, i_k\}$.  This correspondence is an isomorphism from the intersection
lattice to  $B_n$.   
\end{example}

\begin{example} \label{bboolarr} {\em The type B coordinate hyperplane arrangement and the
face lattice of
 the $n$-cross-polytope
$C_n$.} The type B coordinate hyperplane arrangement is the affine hyperplane arrangement
consisting of the hyperplanes $x_i = \pm 1$ in $\R^n$.    One can see  that if we attach a
top element $\hat 1$ to the intersection semilattice of this arrangement we have a lattice
that is isomorphic to the lattice of  faces of the $n$-cross-polytope, which we
denote by
$C_n$.
(This is  dual  to the lattice of
 faces of the
$n$-cube.)   Indeed, the   intersection  $$\{ \bold x \in
\R^n :\epsilon_1 x_{i_1} = \epsilon_2 x_{i_2} = \dots =  \epsilon_kx_{i_k} =1 \},$$
where $1 \le i_1 <i_2 <\dots<i_k \le n$ and $\epsilon_i \in \{-1,1\}$, maps to the 
$(n-k)$-face $$\{ \bold x \in
[-1,1]^n :\epsilon_1 x_{i_1} = \epsilon_2 x_{i_2} = \dots =  \epsilon_kx_{i_k} =1 \}$$ of
the
$n$-cube.  This correspondence is an isomorphism from the intersection lattice to the dual
of the face lattice of the cube.
\end{example}

\begin{example}\label{braidarr} {\em The (type A) braid arrangement and the partition
lattice
$\Pi_n$}.  For
$1
\le i < j
\le n$, let
$$H_{i,j}=
\{\bold x
\in 
\R^n : x_i=x_j\}.$$ The hyperplane arrangement $$\mathcal A_{n-1} := \{H_{i,j} :1 \le i < j
\le n\} $$ is known as the {\em braid arrangement} or the {\em type A Coxeter
arrangement}.  The intersection lattice
$L(\mathcal A_{n-1})$ is isomorphic to
$\Pi_n$, the lattice of partitions of the set $[n]$ ordered by refinement.  Indeed, for
each partition
$\pi
\in
\Pi_n$, let $\ell_\pi$ be the linear subspace of $\R^n$ consisting of all points $(x_1,\dots,
x_n)$ such that $x_i = x_j$ whenever $i$ and $j$ are in the same block of $\pi$.  The map
$\pi \mapsto \ell_\pi$ is a poset isomorphism from $\Pi_n$ to $L(\mathcal A_{n-1})$.  The
braid arrangement $\mathcal A_{2} $ intersected with the plane $x_1+x_2+x_3 = 0$ and the
partition lattice
$\Pi_3$ are shown in Figure~\ref{figbraid}.

\begin{figure} 
\begin{center}
\includegraphics[width=9cm]{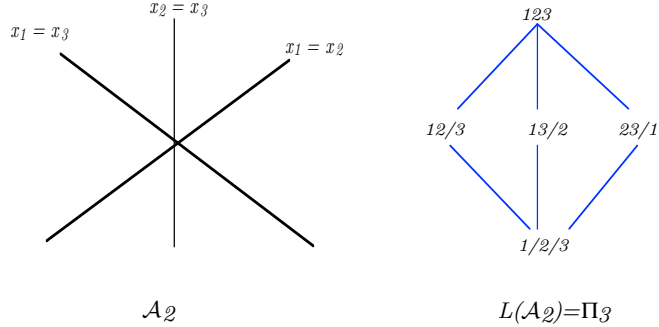}
\end{center}\vspace{-.3in}\begin{center}\caption{Intersection lattice of braid
arrangement}
\label{figbraid}
\end{center}
\end{figure}

\end{example}

\begin{example} \label{bbraidex}{\em The type B braid arrangement and the type B
partition lattice $\Pi_n^B$}.   For $1 \le i < j \le n$, let $$H^+_{i,j}=
\{\bold x
\in 
\R^n : x_i=x_j\} \quad \mbox{ and } \quad H^-_{i,j}=
\{\bold x
\in 
\R^n : x_i= -x_j\}.$$
For $i = 1,\dots n$, let 
$$H_{i}=
\{\bold x
\in 
\R^n : x_i=0\}.$$
The hyperplane arrangement $$\mathcal B_n := \{H^+_{i,j} :1 \le i < j
\le n\} \cup \{H^-_{i,j} :1 \le i < j
\le n\} \cup \{H_i :1 \le i 
\le n\} $$ is called the {\em type B braid arrangement} or the {\em type B Coxeter
arrangement}.  The intersection lattice
$L(\mathcal B_{n})$ is isomorphic to  the  {\em type B (or signed)  partition lattice}. The
elements of the type B partition  lattice $\Pi^B_n$ are partitions of $\{0,\dots,n\}$ for
which any of  the elements of $[n]$ can have a bar except for the elements of the block
that contains $0$ and the smallest element of each block.  For example, $025/ 1 \bar7 9 /
3 4
\bar 6 \bar 8 $ is an element of $\Pi^B_9$, while  $02\bar5 / 1 \bar7 9 / \bar 3 4
\bar 6 \bar 8 $ is not because $5$ and $3$ are not allowed to be barred.  The covering relation is given  by $\pi_1 <\hspace{-.1in}\cdot \,\pi_2$ if
$\pi_2$ is obtained from $\pi_1$ by merging two blocks $B_1$ and $B_2$  into a single
block $B$ in the following manner:  Suppose $\min B_1 < \min B_2$. Then
\begin{itemize}
\item if $0 \in B_1$, let $B$ be the union of $B_1$ and $B_2$ with all bars 
removed,
\item if  $ 0 \notin B_1$, let $B$ be the
union of either \begin{itemize}
\item $B_1$ and $B_2$ with all bars intact or
\item  $B_1$ and $\bar B_2$, where $\bar B_2$  is obtained from $B_2$ by barring  all
unbarred elements and unbarring  all barred elements.
\end{itemize}
 \end{itemize}
For example, the type B partitions that cover $025 / 1 \bar7 9 / 3 4
\bar 6 \bar 8 $  in $\Pi^B_9$ are 
$$0125 7 9 / 3 4\bar 6 \bar 8,\qquad 
0234568 / 1 \bar7 9, \qquad 025 / 134\bar6 \bar7 \bar 89, \qquad 025 / 1\bar3\bar46 \bar7
89.$$ The isomorphism from $\Pi_n^B$ to $L(\mathcal B_n)$ is quite natural. Take a
typical type  B partition 
$$025/ 1 \bar7 9 /
3 4
\bar 6 \bar 8 .$$  It maps to the subspace
$$\{ \bold x \in \R^9:  x_2 = x_5 = 0, \, x_1 = -x_7 = x_9,  \, x_3 = x_4 = -x_6 =
-x_8\},$$ in $L(\mathcal B_n)$.
   The type B braid
arrangement
$\mathcal B_{2} $ and the type $B$ partition lattice
$\Pi^B_2$ are shown in Figure~\ref{figbbraid}.

\begin{figure} 
\begin{center}
\includegraphics[width=9cm]{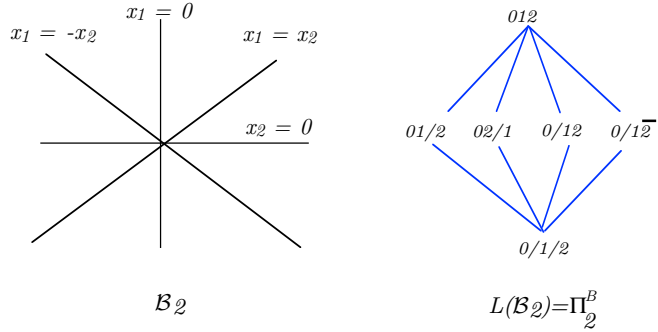}
\end{center}\vspace{-.3in}\begin{center}\caption{Intersection lattice of type B braid
arrangement}
\label{figbbraid}
\end{center}
\end{figure}

\end{example}

 Examples \ref{boolarr} and \ref{braidarr} are referred to as type A
examples, and Examples~\ref{bboolarr} and~\ref{bbraidex} are referred to as type B
examples because of their connection with Coxeter groups.   
Indeed, 
  associated
with every finite Coxeter group (i.e., finite group generated by Euclidean reflections) is
a simplicial complex called its Coxeter complex.    The order complex of the Boolean
algebra
$B_n$ is  the Coxeter complex   of the symmetric group $\s_n$,  which is  the type A
Coxeter group, and the order complex of the face lattice of the
cross-polytope
$C_n$ is the Coxeter complex of the  hyperoctahedral group, which is  the type B
Coxeter group. (The use of notation is unfortunate here; $B_n$ is  type A and $C_n$ is type
$B$.)  

Also associated with every finite Coxeter group is a hyperplane arrangement, called its
Coxeter arrangement, which
 consists of all its reflecting hyperplanes.    The group generated by the reflections
about hyperplanes in the Coxeter arrangement is the Coxeter group.  The braid arrangement
is the Coxeter arrangement of the symmetric group (type A Coxeter group) and type B braid
arrangement is the Coxeter arrangement of the hyperoctahedral group (type B Coxeter group). 
Coxeter groups are discussed further in Section~\ref{CLsec}.  See the chapters in this
volume  by 
Fomin and Reading  \cite{fr04} and   Stanley 
\cite{st04}  for further discussion of Coxeter arrangements.

The types A and B partition lattices 
belong to another  family of well-studied lattices, namely the     Dowling lattices.    
We will not define Dowling lattices, but we will occasionally refer to them; see
\cite{gw} for the  definition.  A broad class of Dowling lattices arise as 
intersection lattices of  complex hyperplane arrangements  $\mathcal A_{m,n}$ consisting of
hyperplanes of the forms
$z_j=0$, where $j=1,\dots,n$, and
$z_j =
\omega^h z_i$, where
$\omega$ is the $m$th primitive root of unity $e^{2\pi i \over m}$,  $ 1 \le i <j\le n$,
and
$h
\in [m]$.  This class includes the types A and B partition lattices.

We now state Zaslavsky's seminal result.  
\begin{theorem}[Zaslavsky \cite{za}] Suppose $\mathcal A$ is a hyperplane
arrangement in
$\R^n$. Let $r(\mathcal A)$ be the number of regions into which $\mathcal A$ divides
$\R^n$ and let $b(\mathcal A)$ be the number of these regions that are bounded.  Then 
\begin{equation} \label{reg} r(\mathcal A) = \sum_{x \in L(\mathcal A)} |\mu(\hat 0,
x)|\end{equation} and 
\begin{equation} \label{bounded} b(\mathcal A) = |\mu(L(\mathcal A) \cup \hat
1)|.\end{equation}
\end{theorem}

The arrangement of Figure~\ref{fighyper} has a total of 10 regions with 2 of them 
bounded.  One can use the values of the M\"obius function given in Figure~\ref{figmob} to
confirm (\ref{reg}) and (\ref{bounded}) for the arrangement of Figure~\ref{fighyper}. 
Note that if
$\mathcal A$ is a central arrangement,
$  L(\mathcal A)\cup
\hat 1$ has an artificial top element above the top element of $L(\mathcal A)$.  In other
words
$  L(\mathcal A)\cup \hat 1$ has exactly one coatom.  It is easy to see that posets with
only one coatom have M\"obius invariant 0.  Since central arrangements clearly have no
bounded regions, (\ref{bounded}) is trivial for central arrangements.

\begin{xca}  Suppose we have a hyperplane $H$ of $\R^n$ which is generic with
respect to a central hyperplane arrangement $\mathcal A$ in $\R^n$.   This means that
$\dim(H \cap X) =
\dim(X) -1$ for all $X \in L(\mathcal A)$.   Let $\mathcal A_H = \{H \cap K : K \in
\mathcal A\}$.  This is a  hyperplane arrangement induced in
$H \cong R^{n-1}$.  Show that the number of bounded regions of $\mathcal A_H$ is
independent of the choice of generic hyperplane $H$ (see \cite{bw1}).
\end{xca}

The next major development in the combinatorial theory of hyperplane arrangements is a
1980 formula of Orlik and Solomon
\cite{os}, which can be viewed as a complex analog of (\ref{reg}).    The number of regions
in a real hyperplane arrangement $\mathcal A$ is equal to the sum of all the Betti numbers
of the complement  $\R^n - \cup \mathcal A$.  Indeed, since   each  region is
contractible, all the Betti numbers are 0 except for the degree $0$ Betti number, which
equals the number of regions.  Hence (\ref{reg}) can be interpreted as a formula for the
sum of the Betti numbers of the complement
$\R^n -
\cup
\mathcal A$.  The analog  for complex arrangements is given by 
following result.

\begin{theorem} [Orlik and Solomon \cite{os}] Let $\mathcal A$ be a hyperplane
arrangement in
$\C^n$.  The complement $M_{\mathcal A } := \C^n - \cup \mathcal A$ has torsion-free
integral cohomology and has Betti numbers given by, 
$$\beta_i(M_{\mathcal A}) =\sum_{\scriptsize\begin{array}{c}x \in L({\mathcal
A})\\ \dim^\C(x) = n-i \end{array}} |\mu(\hat 0, x)|,$$
for all $i$.
\end{theorem}

There is a striking common generalization of  the Zaslavsky formula (\ref{reg}) 
and the Orlik-Solomon formula, obtained by  Goresky and MacPherson in 1988, which
involves subspace arrangements.  A {\em real subspace arrangement} is a finite collection
of (affine) subspaces in
$\R^n$.  Real hyperplane arrangements and complex hyperplane arrangements are both examples
of real subspace arrangements.  Indeed, hyperplanes in $\C^n$ can be viewed as codimension
$2$ subspaces of $\R^{2n}$.   Again the intersection semilattice $L(\mathcal A)$ is defined
to be the semilattice of nonempty intersections of subspaces in the subspace arrangement
$\mathcal A$.  

\begin{theorem} [Goresky and MacPherson \cite{gm}] \label{GMform} Let $\mathcal A$
be a subspace arrangement in
$\R^n$.  The reduced integral cohomology of the complement $M_{\mathcal A } := \R^n - \cup
\mathcal A$ is given by the group isomorphism
$$\tilde H^i(M_{\mathcal A};\Z) \cong \bigoplus_{x \in L({\mathcal
A}) \setminus \{\hat 0\}} \tilde H_{n-\dim x -2-i} ((\hat 0, x);\Z),$$
for all $i$.  
\end{theorem}

To see that the Goresky-MacPherson formula
reduces to the Zaslavsky formula and to the Orlik-Solomon formula, one needs to understand
the homology of the intersection lattice of a  central hyperplane arrangement.  The
intersection lattice belongs to a well-understood class of lattices called geometric
lattices.  A  fundamental result due to Folkman
\cite{fol} states that the proper part of any
geometric lattice $L$ has vanishing reduced homology in every dimension except the top
dimension (i.e. dimension equal to $l(L)-2$).  In fact,  the
homotopy type is that of a wedge of spheres of top dimension. 
The intersection lattice of an affine hyperplane arrangement  belongs to a more general
class of lattices called geometric semilattices, which were introduced and studied by
Wachs and Walker
\cite{ww}.   The proper part of a geometric semilattice also has
the homotopy type of a wedge of spheres of  top dimension.  Topology of geometric
(semi)lattices is discussed further  in Sections~\ref{geosec} and~\ref{RAO}.

\begin{xca}
Use Folkman's result to show that the Goresky-MacPherson formula reduces to both the 
Zaslavsky formula and  the
Orlik-Solomon formula.
\end{xca}

The intersection lattice of a hyperplane arrangement determines more
than the additive group structure of the integral cohomology of the complement.  Orlik and
Solomon show that it determines the ring structure as well.  Ziegler \cite{z93} showed
that,  in general, for 
 subspace arrangements the combinatorial data
(intersection lattice and dimension information) does not determine ring structure. 
However in certain special cases the combinatorial data  does  determine the cohomology
algebra, see
\cite{fz}, \cite{y}, \cite{dp}.   In
Section~\ref{arrangsec} we discuss some stronger versions of the Goresky-MacPherson
formula, namely  a homotopy version due to Ziegler and 
\v Zivaljevi\'c \cite{zz}, and an equivariant  version due to Sundaram and
Welker \cite{suwe}.   For further reading on hyperplane arrangements,  see the chapter by
Stanley in this volume \cite{st04} and the text by Orlik and Terao \cite{ot}.  Further
information on subspace arrangements can be found in  Bj\"orner  \cite{bj92} and Ziegler
\cite{zhab}. 

\section{Some connections with graphs, groups and lattices} \label{grth} 

In this section we briefly discuss some results and questions  in which
poset topology plays a  role. We start with a  old conjecture of Karp in  graph complexity
theory.   An algorithm for deciding whether a graph
 with
$n$ nodes  has a certain property  checks the entries of the graph's adjacency
matrix  until a determination can be made.  A graph property  is said to be
 evasive  if the best algorithm needs to check all $n \choose 2$ entries (in the worst
case).  Here are some examples of  evasive graph properties:

\begin{itemize}
 \item property of being  connected
 \item property of containing a  perfect matching 
 \item property of having  degree at most $b$ for some fixed $b$. 
\end{itemize}

 A monotone  graph property is a property of graphs that is
isomorphism invariant and closed under addition of edges or closed under removal of edges. 
The graph properties listed above are clearly monotone graph properties.   We say that a
graph property is trivial if every graph has the property or every graph lacks the
property.

\begin{con} [Karp's Evasiveness Conjecture]  Every  nontrivial monotone graph
property is evasive.
\end{con}

  Kahn, Saks, and Sturtevant \cite{kss} proved the evasiveness conjecture for $n$ 
a prime power by using topological techniques and group actions. 
 Since determining whether a graph has a
certain property is equivalent to determining whether the graph lacks the property,
one can require without loss of generality that a  monotone graph property be closed under
removal of edges. 
 Given 
such 
a monotone graph property,
$\mathcal P$, let
$\Delta^n_{\mathcal P}$ be the simplicial complex whose vertex set is ${[n] \choose 2}$
and whose faces are the edge sets of graphs on node set $[n]$ that have the property. 
Alternatively, $\Delta^n_{\mathcal P}$ is the simplicial complex whose face poset is the
poset of graphs on node set $[n]$ that have property $\mathcal P$, ordered by
edge set inclusion.
   Kahn, Saks, and Sturtevant show that
\begin{itemize} \item nonvanishing reduced simplicial homology of $\Delta^n_{\mathcal P}$
 implies $\mathcal P$ is
evasive  
\item $\Delta^n_{\mathcal P}$ has nonvanishing reduced simplicial homology when $n$ is a
prime power, using  a topological fixed point theorem.
\end{itemize}

   Although this connection between
evasiveness  and topology doesn't really involve posets directly,   we mention it here
because posets and  simplicial complexes can be viewed as
the same object, and  as we will see in later in these lectures, the tools of poset
topology are useful in the study of the topology of graph complexes.     For other 
significant results on evasiveness and topology of graph
complexes, see eg.,  
\cite{cks}, \cite{yao}, \cite{for00}. This topic is  discussed in greater
depth in  Forman's chapter of this volume \cite{for05}.   Applications of graph
complexes in knot theory and group theory are discussed in Section~\ref{quilfibsec}. 
There are also connections between the
topology of graph complexes and commutative algebra, which are explored in the work of
Reiner and Roberts \cite{rr} and Dong \cite{dongth}.  A direct
 application of poset topology in a different complexity theory problem is discussed in
Section~\ref{kesec}.

  Representability  questions in lattice theory deal with
whether  an arbitrary lattice can be represented as a sublattice, subposet or interval in
a given class of lattices.  We briefly discuss three examples that have 
connections to poset topology.   

A  result of Pudl\'ak and Tuma \cite{pt}  states that every  lattice is isomorphic
to a sublattice of some  partition lattice
$\Pi_n$.  This implies that every  lattice can be  represented as the
intersection lattice of a subspace arrangement embedded in the braid arrangement. 
There is another representability result that is much easier to prove; namely that every
meet semilattice can be represented as the intersection semilattice of some subspace
arrangement, see \cite{zhab}. From either of these representability results, we see
that, in contrast to the situation with hyperplane arrangements, where the topology of the
proper part of the intersection semilattice is rather special (a wedge of spheres), any
topology is possible for the intersection semilattice of a general subspace arrangement. 
Indeed, given any simplicial complex $\Delta$, there is a linear subspace arrangement
$\mathcal A$ such that $\overline{ L(\mathcal A)} $ is homeomorphic to $\Delta$;
namely
$\mathcal A$ is the linear subspace arrangement whose intersection lattice $L(\mathcal
A)$ is isomorphic to the face lattice 
$L(\Delta)$.

An open  representability question is whether  every
lattice can be represented as an interval in the lattice of subgroups of some group
ordered by inclusion.  An approach to obtaining a negative answer to this question, 
proposed by  Shareshian \cite{sh03}, is to establish restrictions on the topology of
intervals in the subgroup lattice.

\begin{con}[Shareshian \cite{sh03}]  Let  $G$ be a finite group.  Then every open interval
in the lattice of subgroups of $G$ has the homotopy type of a wedge of spheres. 
\end{con}

This conjecture was shown to hold for  solvable groups by Kratzer and Th\'evenaz
\cite{kt} (see Theorem~\ref{shth} which strengthens the Kratzer-Th\'evenaz result).
  Further discussion of
connections between poset topology and group theory can be found  in
Section~\ref{quilfibsec}

Our last example 
deals with the {\em order dimension} of a poset $P$, which is defined to be the
smallest integer $n$ such that $P$ can be represented as an induced subposet of a product
of
$n$ chains.  Order dimension is  an important and extensively studied poset invariant, see
\cite{trot}. Reiner and Welker give a lower bound on  order dimension of a lattice in
terms of its homology.  

\begth[Reiner and Welker \cite{rwel}]  
Let $L$ be a lattice and let $d$ be the largest dimension for which the reduced 
integral simplicial homology  of  the proper part of $L$ is nonvanishing.  Then 
the order dimension of $L$ is at least  $d+2$.
\enth

\section{Poset homology and cohomology}

By (co)homology of a poset, we usually mean  the reduced simplicial (co)homology of its
order complex.  On rare occasions, we will deal with nonreduced simplicial homology. 
Although it is presumed that the reader is familiar with simplicial homology and
cohomology, we review these concepts for  posets
 in terms of chains  of the poset. For each poset
$P$ and integer
$j$,  define the chain space 
$$ C_j(P;{\mathbf k}) := { \mathbf k}\mbox{-module freely generated by $j$-chains of
} P,
$$ where ${\mathbf k}$ is a field or the ring of integers.

The boundary map $\partial_j:  C_j(P;{\mathbf k}) \to   C_{j-1}(P;{\mathbf k})$
is defined by
$$\partial_j(x_1 < \dots < x_{j+1}) = \sum_{i=1}^{j+1} (-1)^i\,  (x_1 < \dots < \hat x_i <
\dots < x_{j+1}),$$
where the $\hat \cdot $ denotes deletion.
We have that $\partial_{j-1}
\partial_j = 0$, which makes $( C_j(P;{\mathbf k}) ,\partial_j)$  an algebraic
complex.  Define the cycle space
$Z_j(P;{\mathbf k}):= \ker \partial_j$ and the boundary space
$B_j(P;{\mathbf k}) := \im \partial_{j+1}$.  Homology of the poset $P$ in dimension $j$
is defined by
$$\tilde H_j(P;{\mathbf k}) : = Z_j(P;{\mathbf k}) / B_j(P;{\mathbf k}).$$

 The coboundary map $\delta_j:  C_j(P;{\mathbf k}) \to   C_{j+1}(P;{\mathbf
k})$ is defined by
\begin{equation}\label{form} \langle \delta_j(\alpha),\beta \rangle = \langle
\alpha,\partial_{j+1}(\beta)
\rangle\end{equation} where $\alpha \in  C_j(P;{\mathbf k})$, $\beta \in
 C_{j+1}(P;{\mathbf k})$, and
$\langle
\cdot,
\cdot
\rangle$ is the  bilinear
form on $\oplus_{j \ge -1}  C_j(P;{\mathbf k})$ for which the chains of $P$ form an
orthonormal basis. 
This is
equivalent to saying
\begin{eqnarray*}  \delta_j(x_1 < \dots < x_{j}) = \qquad \qquad\\
\sum_{i=1}^{j+1} (-1)^i \sum_{x \in (x_{i-1},x_i)} 
 (x_1  < \dots < x_{i-1} < x < x_i <\dots < x_{j}),\end{eqnarray*}
for all chains $x_1 < \dots < x_{j}$, where $x_0$ is the bottom element of $\hat P$
and $x_{j+1}$ is the top element of
$\hat P$. 
Define the
cocycle space to  be $Z^j(P;{\mathbf k}):= \ker \delta_j$ and the coboundary space to
be $B^j(P;{\mathbf k})  := \im
\delta_{j-1}$.
Cohomology of the poset $P$ in dimension $j$
is defined to be
$$\tilde H^j(P;{\mathbf k}) : = Z^j(P;{\mathbf k}) / B^j(P;{\mathbf k}).$$

When $\k$ is  a field, $\tilde H^j(P;\k)$ and $\tilde H_j(P;\k)$ are
isomorphic vector spaces. The
$j$th  (reduced) Betti number  of
$P$ is given by
$$\tilde\beta_j(P) := \dim \tilde H_j(P;\C),$$
which   is the same as the rank of the free part of $\tilde
H_j(P;\Z)$.

We will work primarily with homology over $\C$ and $\Z$. 
For $x <y$ in $P$, we  write $\tilde
H_j(x,y)$ for the complex homology of the open interval
$(x,y)$  of
$P$, and 
$\tilde\beta_j(x,y)$ for the $j$th Betti number of the open interval
$(x,y)$.  When $x=y$,  define   $\tilde H_j(x,y)$ to be $\C$ and $\tilde \beta_j(x,y)$ to
be
$1$ if $j = -2$, and to be $0$ for all other $j$. 

Many of the posets that arise have the homotopy type of a wedge of spheres.  We review
a basic  fact   pertaining to wedges of spheres and a partial converse.

\begth Suppose $\Delta$ has the homotopy type of a
wedge of spheres of various dimensions,  where 
$r_i $  is the number of spheres of dimension $i$.  Then for each  $i = 0,1,\dots, \dim
\Delta$,
\begin{equation}\label{homshell}\tilde H_i(\Delta;\Z ) \cong \tilde H^i(\Delta;\Z) \cong
\Z^{r_i}.\end{equation} 
\enth

\begth If $\Delta$ is simply connected and has
vanishing reduced integral homology in all dimensions but  dimension $n$, where  homology
is free of rank $r$, then
$\Delta$ has the homotopy type of a wedge of $r$ spheres of dimension $n$.
\enth

The first tool that we mention for computing homology of posets and simplicial complexes is
 a very efficient computer software package called
``SimplicialHomology'', developed by  Dumas, Heckenbach, Sauders, and Welker \cite{DHSW}. 
 One can run it interactively or download the source file at the web site:
\begin{center}{\tt
http://www.cis.udel.edu/$\sim$dumas/Homology},\end{center} 
 where a manual can also be
found.   This package has been responsible for many of the more recent conjectures in the
field.  Its output was also part of the proofs of  (at least) three results on integral
homology appearing in the literature; see
\cite{sh04, shwa2, wa3}.  

\section{Top cohomology of the partition lattice} \label{liesec} 
\begin{figure}
\begin{center}
\includegraphics[width=6cm]{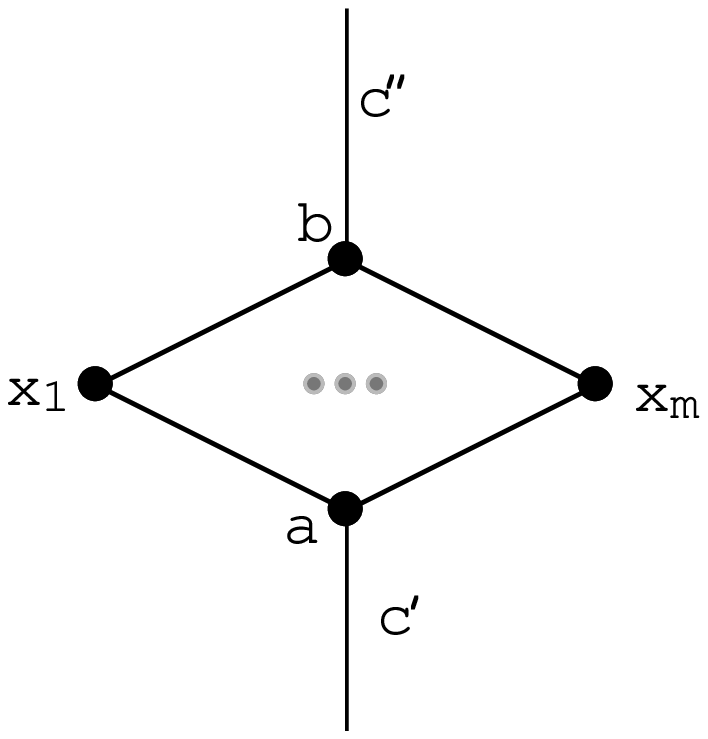}
\end{center}\begin{center}\caption{Top coboundary relations: $\delta(c^\prime
c^{\prime\prime})$}
\label{figcorel}
\end{center}
\end{figure}

The top dimensional
cohomology of a poset has a particularly simple description.  For the sake of simplicity
assume
$P$ is pure of length $ d$. Let
$\mathcal M(P)$ be the set of maximal chains of
$P$ and let $\mathcal M^\prime(P)$ be the set of chains of length $d-1$.  Since $\ker
\delta_d =  C_d(P;\k)$, we have the following presentation of top cohomology as a
quotient of
$ C_d(P;\k)$:
 $$\tilde H^d(P;\k) = \langle
\mathcal M(P) \mid \mbox{coboundary relations}\rangle,$$
where the coboundary relations have the form $\delta_{d-1}(c)$ for $c \in \mathcal
M^\prime(P)$.  Each chain $c$ in  $\mathcal
M^\prime(P)$ is the concatenation $c^\prime c^{\prime\prime}$ of two unrefinable chains
$c^\prime$ and
$c^{\prime\prime}$.  If $c^{\prime}$ is not empty, let $a$ be the maximum
element of
$c_1$, and if
$c^{\prime}$ is empty, let $a$ be $\hat 0$ of $\hat P$.  If $c^{\prime\prime}$ is not
empty, let
$b$ be the minimum element of
$c^{\prime\prime}$, and if
$c^{\prime\prime}$ is empty, let  $b$ be $\hat 1 $ of $\hat P$.  Clearly, $[a,b]$ is an
interval of length 2 in $\hat P$.  Let $\{x_1,\dots,x_m\}$ be the set of elements in
the open interval $(a,b)$.  We have
$$\delta_{d-1} (c) = \pm(c^{\prime} x_1c^{\prime\prime} + \dots + c^{\prime}x_m
c^{\prime\prime} ).$$
Hence  the cohomology
relations can be associated with the intervals of length 2 in $\hat P$.  See
Figure~\ref{figcorel}.

\begin{figure}\begin{center}
\hspace{-.2in}\mbox{\includegraphics[width=6 cm]{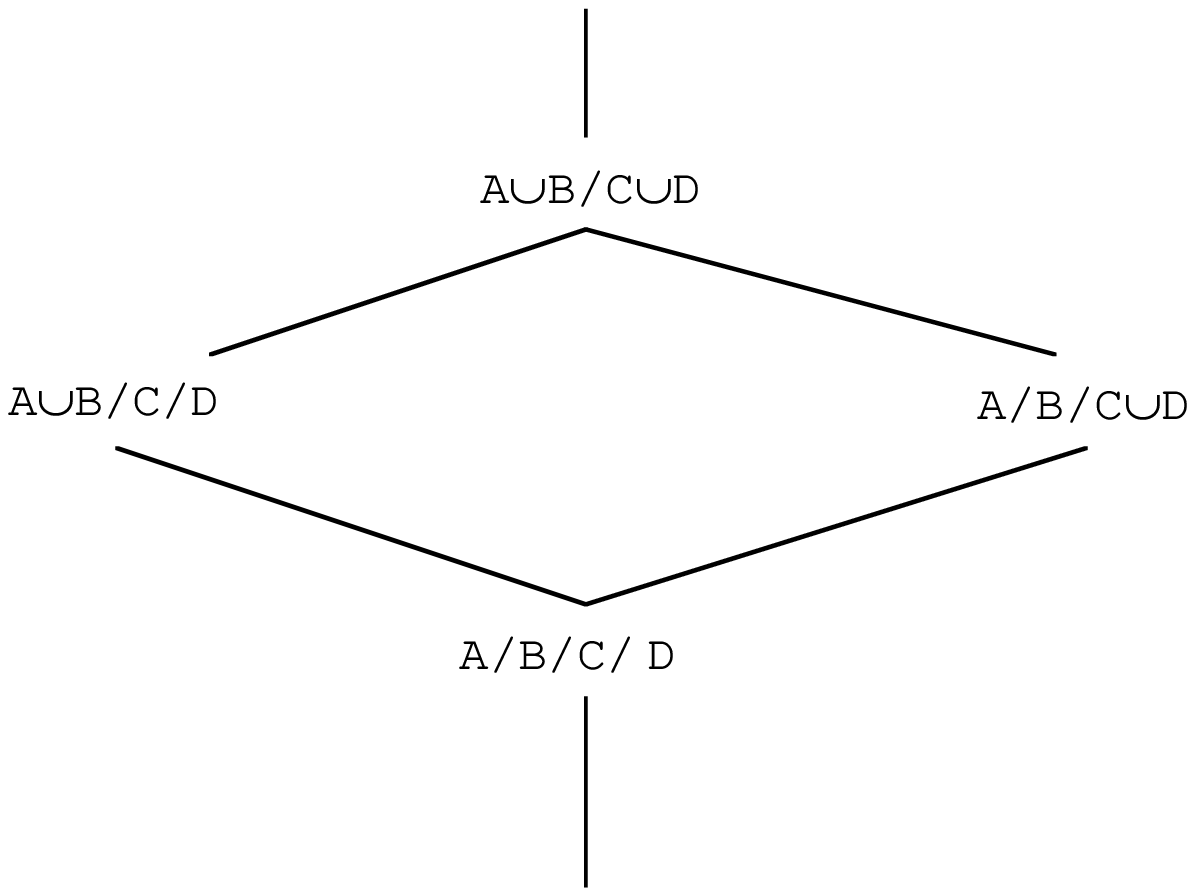}} \hspace{.3in}
\mbox{\includegraphics[width=6cm]{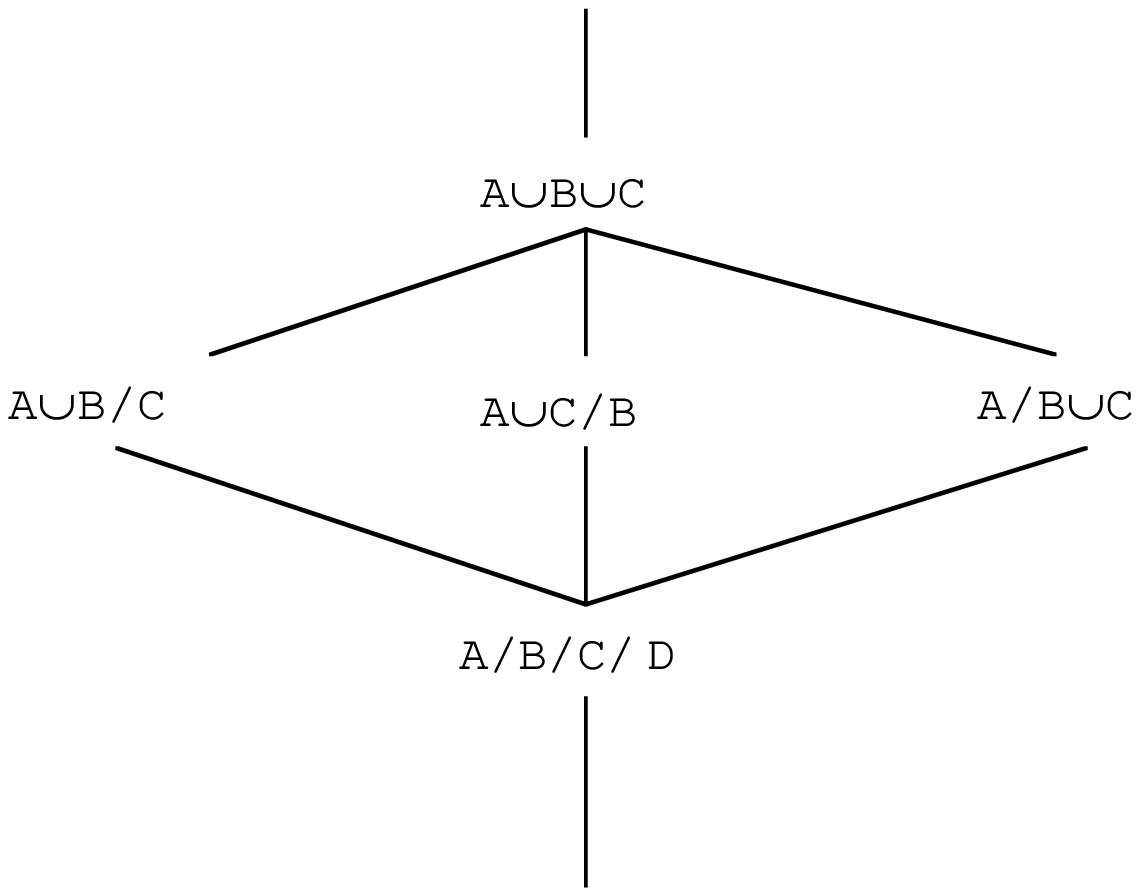}}
\end{center}
\hspace{-.3in}Type I interval \hspace{1.5in} Type II interval
\caption{\,}
\label{figtoppart}
\end{figure}

We demonstrate the use of intervals of length 2 by deriving a presentation for the top
cohomology of  the proper part of the partition lattice
$\bar \Pi_n$.  There are two types of length 2 closed intervals in
$\Pi_n$; see Figure~\ref{figtoppart}.  In Type I intervals, there are 2 pairs of blocks
$\{A,B\}$ and
$\{C,D\}$ which  are   separately merged resulting in  blocks $A\cup B$ and $C\cup D$.  In
Type II intervals, there are 3 blocks
$A,B,C$, which are merged into one block $A\cup B\cup C$.  Type I
intervals have 4 elements and type II intervals have 5 elements. 

\begin{figure}\begin{center}
\mbox{\includegraphics[width=3cm]{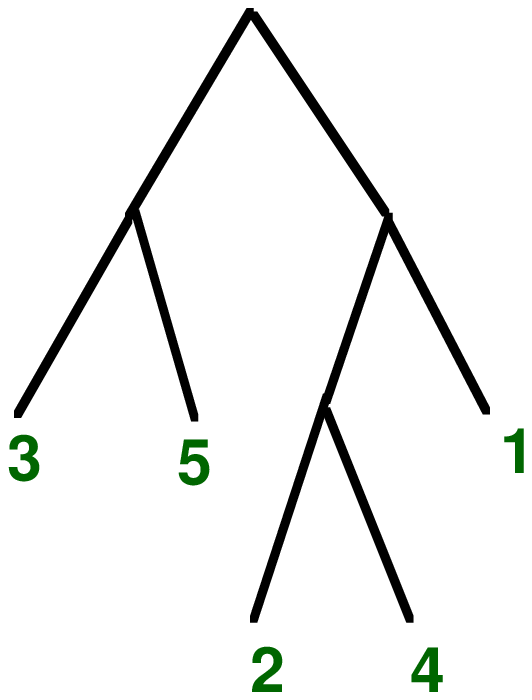}}
\end{center}
\caption{\,}
\label{figtreepart}
\end{figure}

 The two types of
intervals induce two types of cohomology relations, Type I and Type II
cohomology relations on maximal chains.  It is convenient to use 
binary trees on leaf set $[n]$ to describe these relations.   A maximal chain of
$\bar\Pi_n$ is  just a sequence of merges of  pairs of
blocks.  The  binary tree given in Figure~\ref{figtreepart}  corresponds 
to the sequence of merges:
\begin{enumerate}
\item 
merge blocks $\{3\}$ and $\{5\}$
\item  merge blocks $\{2\}$ and $\{4\}$
\item merge blocks $\{2,4\}$
and $\{1\}$.
\end{enumerate}
This corresponds to the maximal chain $$ 1/2/35/4 \,\, <\!\!\!\!\cdot \,\,
1/24/35 \,\,<\!\!\!\!\cdot \,\,\,124/35 $$ of $\bar\Pi_5$  The internal nodes of the tree
represent  the  merges, and the leaf sets of the left and right subtrees of the internal
nodes are the blocks that are merged.  The sequence of merges follows the postorder
traversal  of the internal nodes, i.e. first traverse the left subtree in postorder, then
the right subtree in postorder, then the root.  

Given a  binary tree $T$ on leaf set $[n]$, let $c(T)$ be the maximal chain of $\Pi_n$
obtained by the procedure described above.  Although not all maximal chains can be
obtained in this way, it can be seen that every maximal chain is equal, modulo the
cohomology relations of Type I, to $\pm c(T)$ for some $T$.  So the set 
$$\{c(T) : T \mbox{ is a  binary tree on leaf set } [n]\}$$
generates top cohomology $\tilde H^{n-3}(\bar \Pi_n;\k)$.
The Type I cohomology relations induce the following relations 
\begin{equation} \label{type1} c(\cdots (A \land B)\cdots) = (-1)^{|A||B|} \,\,c(\cdots(B
\land A)\cdots),\end{equation} where $X\land Y$ denotes the binary tree whose left subtree
is
$X$ and whose right subtree is $Y$, and $|X|$ denotes the number of internal nodes of
$X$.  The Type II cohomology relations induce the following relations 
\begin{eqnarray} \label{type2}   c(\cdots (A \land (B \land C))\cdots)
&+&(-1)^{|C|} c(\cdots ((A
\land B)\land C)\cdots)\\ \nonumber  
&+&(-1)^{|A|B|}  c(\cdots(B \land(A \land C))
\cdots) = 0.\end{eqnarray}

\begin{xca} Show that the Type I cohomology relations yield (\ref{type1}) and the Type II
cohomology relations yield   (\ref{type2}).
\end{xca} 

The relations (\ref{type1}) and (\ref{type2}) resemble   the relations satisfied
by  the bracket operation of a Lie algebra.  The Type I relation
(\ref{type1}) corresponds to the anticommuting relation and the Type II relation
(\ref{type2}) corresponds to the Jacobi relation.  Indeed there is a well-known connection
between the top homology of the partition lattice and the free Lie algebra which involves 
representations of the symmetric group.  In the next lecture, we discuss representation
theory.

\begin{theorem}[Stanley \cite{st82}, Klyachko \cite{kl}, Joyal \cite{jo}] \label{SKJ} The
representation of the symmetric group
$\mathfrak S_n$ on $\tilde H_{n-3}(\bar \Pi_n;\C)$ is isomorphic to the
representation of $\mathfrak S_n$ on the multilinear component of the free Lie algebra
over $\C$ on
$n$ generators tensored with the sign representation.
\end{theorem}

This result follows from a formula of Stanley for the representation of the
symmetric group on homology of the partition lattice (Theorem~\ref{thlie}) and an earlier
similar formula of Klyachko for the free Lie algebra.  The first purely combinatorial
proof  was obtained by Barcelo
\cite{bar}.   The presentation of top cohomology discussed above  appeared in an
alternative combinatorial proof of  Wachs
\cite{wa1}.  It also appeared in the proof of a  
superalgebra version of this result obtained by Hanlon and Wachs \cite{hw}.  A $k$-analog
of the Lie superalgebra result was also obtained by Hanlon and Wachs \cite{hw}.  A type B
version (Example~\ref{bbraidex}) was obtained by  Bergeron
\cite{ber} and a generalization to Dowling lattices was obtained by Gottlieb and Wachs
\cite{gw}.

%% file: lect2.tex

\lecture{Group actions on posets}

In this lecture we  give a crash course on the representation theory of the symmetric
group  and then discuss some representations on homology that are induced by symmetric
group actions on posets.  For further details on  the representation theory of the
symmetric group and symmetric functions, we refer the reader to the following 
excellent standard references 
\cite{ful,mac,sa,st2}.

There are various  reasons that we are interested in understanding  how a group acts
on the homology of a poset.  One is that this can be a useful tool in computing the
homology of the poset.  Another is that interesting representations often arise.  We 
limit our discussion to the symmetric group, but point out there are often interesting
analogous  results for other groups such as the hyperoctahedral group, wreath product
groups, and the general linear group.

\section{Group representations} We restrict our discussion to finite groups $G$ and
finite dimensional vector spaces over the field $\C$.  A
 finite dimensional vector space $V$ over $\C$ is said to be a {\em representation} of $G$
if there is a group homomorphism
$$\phi:G \to GL(V).$$  
For $g \in G$ and $v \in V$, we write $gv$ instead of $\phi(g)(v)$ and view $V$ as a
module over the ring $\C G$ ($G$-module for short).
The {\em dimension} of the representation $V$ is defined to be the dimension of $V$ as a
vector space.

There are two particular  representations of every group that are very important;  the
trivial representation and the regular representation.  The {\em trivial representation},
denoted $1_G$, is  the 1-dimensional representation $V=\C$, where $gz = z$ for all $g \in
G$ and
$z
\in
\C$. The {\em (left) regular representation} is the $G$-module $\C G$ where $G$ acts on
itself by left multiplication, i.e., the action of $g \in G$  on generator $h \in G$ is
$gh$.

We say that $V_1$ and $V_2$ are {\it isomorphic}
representations of
$G$ and write
$V_1
\cong_G V_2$, if there is a vector space isomorphism $\psi:V_1 \to V_2$ such that 
$$\psi(gv) = g \psi(v)$$
for all $g \in G$ and $v \in V_1$.  In other words, $V_1$ and $V_2$ are isomorphic
representations of $G$ means that they are isomorphic $G$-modules.

The {\em character} of a $G$-module $V$  is a function $\chi^{V}:G \to \C$ defined by
$$\chi^{V}(g) = \mbox{trace}(\phi(g)).$$
One basic fact of representation theory is that the character of a representation
completely determines the representation.   Another is that $\chi^V(g)$ depends only on
the conjugacy class of $g$.  

A $G$-module $V$ is said to be {\em irreducible} if its only submodules are the trivial
submodule $0$ and $V$ itself.  A basic result of representation theory is that the number
of irreducible representations of $G$ is the same as the number of conjugacy classes of
$G$. Another very important fact is that  every
$G$-module decomposes into a direct sum of irreducible submodules,
$$V \cong_G V_1 \oplus \dots \oplus V_m.$$
The  decomposition is unique (up to order and up to isomorphism).  Hence it
makes sense to talk about the multiplicity of an irreducible in a representation.
We have the following  fundamental fact.

\begth \label{regdecomp} The multiplicity of any irreducible representation of $G$ in the
regular representation of $G$ is equal to the dimension of the irreducible.
\enth

There are two operations on representations that are quite useful.  The first is called
{\em restriction}.    For $H$  a subgroup of $G$ and $V$ a representation of $G$, the
restriction of
$V$ to
$H$, denoted
$V\downarrow^G_H$, is the representation of $H$ obtained by restricting  $\phi$ to
$H$. Thus the restriction has the same underlying vector space with a smaller 
group action.
 The other operation, which is   called induction, is  a bit more complicated.
For $H$  a subgroup of $G$ and $V$ a representation of $H$,  the induction of $V$ to $G$
is given by $$V\uparrow^G_H := \C G \otimes_{\C H} V,$$
where the tensor product $A \otimes_{S} B$ denotes the usual tensor product of an
$(R,S)$-bimodule $A$ and a left $S$-module $B$ resulting in a left $R$-module.  Now
the underlying vector space of the induction  is larger than $V$.  

\begin{xca} Show 
$$\dim V\uparrow^G_H = {|G|
\over |H|}
\dim V.$$

\end{xca}  Although restriction and induction are not inverse operations, they are related
by a formula called Frobenius reciprocity.   We will state an
important special case,  which is, in fact, equivalent to  Frobenius reciprocity.

\begin{theorem}  Let $U$ be an irreducible representation of $H$ and let $V$ be an
irreducible representation of $G$, where $H$ is a subgroup of $G$.  Then the
multiplicity of $U$ in $V\downarrow^G_H$ is equal to the multiplicity of $V$ in
$U\uparrow^G_H$.
\end{theorem}

There are  two types of tensor products of representations.  Given a representation $U$ of 
$G$ and a representation  $V$ of $H$, the (outer) tensor product $U \otimes V$ is a
representation of $G\times H$ defined by 
$(g,h)(u,v) = (gu,hv)$.  Given two representations $U$ and $V$ of $G$, the (inner) tensor
product, also denoted $U \otimes V$, is the representation of $G$ defined by
$g(u,v) = (gu,gv)$.

In these lectures we will describe  representations in any of the following ways:
\begin{itemize}\item giving the character
\item giving an isomorphic representation
\item giving the multiplicity of each irreducible
\item using  operations such as restriction, induction and tensor product.
\end{itemize}

\section{Representations of the symmetric group} \label{secrepsym}

In this section, we construct the irreducible representations of the symmetric group
$\mathfrak S_n$, which are called Specht modules and are denoted by $S^\lambda$, where
 $\lambda$ is a partition of
$n$.   We also discuss skew shaped Specht modules.  

Let $\lambda$ be a partition of $n$, i.e., 
 a weakly decreasing sequence of positive integers $\lambda=(\lambda_1 \ge \dots  \ge
\lambda_k)$ whose sum is
$n$.  We write
$\lambda \vdash n$ (or $|\lambda| = n$) and say that the length $l(\lambda)$  is $k$.  We
will also write
$\lambda = 1^{m_1}2^{m_2}\cdots n^{m_n}$, if $\lambda$ has $m_i$ parts of size $i$ for each
$i$. Each partition
$\lambda$ is identified with  a {\em Young (or Ferrers) diagram} whose
$i$th row has
$\lambda_i$ cells.  For example, the partition $(4,2,2,1) \vdash 9$
is identified with the Young diagram

$$ \qquad\tableau[scY]{ , , ,|, | , ||}
$$

A {\em Young tableau} of shape $\lambda \vdash n$  is a filling of the Young diagram
corresponding to $\lambda$, with distinct positive integers in [n].   A Young tableau is
said to be {\em standard } if the  entries   increase along each
row and column.  For example, the  Young tableau on the left is not standard, while the
one on the right is.  

$$  \tableau[scY]{8,2,4, 1|7,5| 9,3|6|} 
 \hspace{1.5in} \tableau[scY]{1 ,2 ,4 ,7
|3,6 | 5, 8|9|}$$

 Let $\mathcal  T_\lambda$ be the set of Young tableaux of shape $\lambda$ and let
$M^\lambda$ be the complex vector space generated by elements of 
$\mathcal T_\lambda$. The symmetric group $\mathfrak S_n$ acts on $M^\lambda$ by permuting entries
of the Young tableaux. That is, for transposition $\sigma = (i,j) \in \mathfrak S_n$, the
tableau
$\sigma T$ is obtained from $T$ by switching entries $i $ and $j$. For example, 
$$(2,3)\,\,\, \tableau[scY]{8,\Red \mathbf 2\Black,4, 1|7,5| 9,\Red \mathbf
3\Black|6|}\quad  =\quad
\tableau[scY]{8,\Red \mathbf 3\Black,4, 1|7,5| 9,\Red \mathbf 2\Black|6|}.
$$
The representation  that we have described is clearly the left
regular representation of $\mathfrak S_n$.   One can also let $\mathfrak S_n$ act as
the right regular representation on $M^\lambda$.  That is for transposition $\sigma =
(i,j) \in
\mathfrak S_n$ and $T \in  \mathcal  T_\lambda$, the tableau $T \sigma$ is obtained from
$T$ by switching the contents of the $i$th and $j$th cell under some fixed ordering of the
cells of $\lambda$.

 We will say that two tableaux in $\mathcal T_\lambda$ are {\em row-equivalent} if they
have the same sequence of row sets.  For example,  the tableaux
$$  \tableau[scY]{8,2,4, 1|7,5| 9,3|6|}  \hspace{.5in}
 \mbox{ and }  \hspace{.5in} \tableau[scY]{1 ,2 ,4 ,8
|5,7 | 3,9|6|}$$ 
are row-equivalent.  {\em Column-equivalent} is defined similarly.
For shape
$\lambda$, the row stabilizer
$R_\lambda$  is defined to be the subgroup
$$R_\lambda := \{\sigma \in \mathfrak S_n :  T \sigma \mbox{ and $T$ are row-equivalent for
all
$T\in
\mathcal   T_\lambda$}\}. $$ Similarly, the column stabilizer $C_\lambda$ is defined to be
the subgroup
$$C_\lambda := \{\sigma \in \mathfrak S_n :  T \sigma \mbox{ and $T$ are column-equivalent
for all
$T\in
\mathcal  T_\lambda$}\}.$$

We now give two  characterizations of the Specht module $S^\lambda$; one as a 
 subspace of $M^\lambda$ generated by certain signed sums of Young
tableaux called polytabloids; and the other as a quotient of $M^\lambda$ by certain
relations called row relations and Garnir relations.  We caution the reader that our
notions of polytabloids and Garnir relations are dual to the usual notions  given
in  standard texts such as
\cite{sa}. 

We begin with the submodule characterization.  For each $T  \in \mathcal T_\lambda$, define the
{\it polytabloid} of shape $\lambda$,
$$e_T := \sum_{\alpha \in R_\lambda} \sum_{\beta \in C_\lambda}\sgn(\beta)\, \,
T\alpha\beta.$$
For example if
$$T= \tableau[scY]{1,2 |3|}
\quad \mbox{ then }\quad  e_T = \left(\,\, \tableau[scY]{1,2 |3|} - \tableau[scY]{3,2
|1|}\,\,\right) +
\left(\,\,\tableau[scY]{2,1|3|}-
\tableau[scY]{3,1 |2|}\,\,\right ).$$
Since the left and right action of $\mathfrak S_n$ on $\mathcal T_\lambda$ commute, we
have  
\begin{equation} \label{invariant} \pi e_T = e_{\pi T}, \end{equation}  
 for all $T \in \mathcal  T_\lambda$ and $\pi \in \mathfrak S_n$.
We can now define the {\em Specht module} $S^{\lambda}$ to be the subspace of $M^\lambda$
given by
$$S^{\lambda} := \langle e_T : T \in \mathcal  T_\lambda \rangle.$$  It follows from
(\ref {invariant}) that
 $S^\lambda$ is an $\mathfrak S_n$-submodule
of
$M^\lambda$ (under the left action).

\begin{theorem} \label{irred} The Specht modules $S^\lambda$ for all $\lambda \vdash n$
form a complete set of irreducible
$\mathfrak S_n$-modules.
\end{theorem}

A polytabloid $e_T$ is said to be a {\em standard polytabloid} if $T$ is
a standard Young tableau.  We will see shortly that the standard
polytabloids of shape $\lambda$ form a basis for the Specht module $S^\lambda$.  

Now we give the  quotient characterization.  The {\em row relations}  are defined for all 
$T
\in
\mathcal \mathcal T_\lambda$ and  $\sigma
\in R_\lambda$ by 
\begin{equation} \label{rowrel} r_\sigma(T) :=  T \sigma - T.\end{equation}

For  all
$i,j$ such that
$1
\le j
\le
\lambda_{i}$, let 
$C_{i,j}(\lambda)$ be the set of cells  in columns $j$ through $\lambda_{i}$ of
row
$i$ and in columns $1$ through $j$ of row
$i+1$.  Let $G_{i,j}(\lambda)$ be the subgroup of $\mathfrak S_n$ consisting of
permutations $\sigma$
that fix all entries of the cells  that are not in $C_{i,j}(\lambda)$ under the right
action of $\sigma $ on
 tableaux of shape $\lambda$. 
 The {\em Garnir relations} are defined for all $i,j$ such that $1 \le j \le
\lambda_{i}$ and for all $T
\in
\mathcal   T_\lambda$ by
\begin{equation}\label{garrel} g_{i,j}(T) := \sum_{\sigma \in G_{i,j}(\lambda)} \, 
T\sigma.
\end{equation}  For example if 
$$T = \tableau[scY]{7,\Red  \mathbf 1, \mathbf 5  ,\mathbf {10} |\mathbf 3,\Red \mathbf 
4  ,
\Black 2\Black |9, 8\Black |11 , 6\Black} \quad \mbox{ and } \quad  (i,j)  = (1,2)
$$
then  the entries $1,5,10,3,4$ are permuted while the remaining entries are fixed. 
So $$ g_{1,2}(T) = \tableau[scY]{7,\Red \mathbf 1, \mathbf 3 ,  \mathbf{4} |\mathbf 5,\Red
\mathbf {10 } ,
\Black 2\Black |9, 8\Black |11 , 6\Black}
+ \tableau[scY]{7,\Red \mathbf 1, \mathbf 3  , \mathbf{4 }| \mathbf {10},\Red \mathbf {5} 
,
\Black 2\Black |9, 8\Black |11 , 6\Black}
 + \tableau[scY]{7,\Red \mathbf 1, \mathbf  3  ,\mathbf 5|\mathbf{4} ,\Red \mathbf {10}  ,
\Black 2\Black |9, 8\Black |11 , 6\Black}
+ \tableau[scY]{7,\Red \mathbf 1, \mathbf 3  , \mathbf 5 |\mathbf {10},\Red \mathbf{4}  ,
\Black 2\Black |9, 8\Black |11 , 6\Black} +
\dots$$

Again, since the left and right action of $\mathfrak S_n$ on $\mathcal 
T_\lambda$ commute, we have  
\begin{eqnarray*}  \pi r_{\sigma}(T) &=& r_{\sigma}(\pi T) \\ 
\pi  g_{i,j}(T) &=& g_{i,j}(\pi T),
\end{eqnarray*}
 for all $\pi \in \mathfrak S_n$. 
Consequently, the subspace $U^\lambda$ of $M^\lambda$ generated by the row
relations~(\ref{rowrel}) and the Garnir relations~(\ref{garrel}) is an $\mathfrak
S_n$-submodule of
$M^\lambda$.

\begin{theorem} \label{specht} For all $\lambda \vdash n$,
$$
 S^\lambda \cong_{\mathfrak S_n} M^\lambda/U^\lambda.$$
\end{theorem}

Now we can view the Specht module  $S^\lambda$ as the module generated by
tableaux of shape $\lambda$ subject to the row and Garnir relations.     

\begin{xca}\label{exbasis} Prove Theorem~\ref{specht} by first showing that,
\newline (a)  $U^\lambda \subseteq \ker \psi,$ where $\psi:M^\lambda \to S^\lambda$ is
defined by $\psi(T) = e_T$,
\newline (b)  the standard polytabloids $e_T $ are
linearly independent,
\newline(c)  the  standard tableaux span $M^\lambda/ U^\lambda$.
\end{xca}

We have the following consequence of Exercise~\ref{exbasis}.

\begin{cor} \label{skewcor} The standard polytabloids of shape $\lambda$ form a basis for
$S^\lambda$. The standard tableaux of shape $\lambda$ form a basis for
$M^\lambda/U^\lambda$.  Consequently $\dim S^\lambda$ is equal to the number of standard
tableaux of shape $\lambda$.
\end{cor}

There is a remarkable  formula for the number of standard tableaux of a fixed shape
$\lambda$.
 \begin{theorem}[Frame-Robinson-Thrall hook length formula] \label{frt} For all $\lambda
\vdash n$, 
$$ \dim S^\lambda = {n! \over\prod_{x \in \lambda} h_x},$$
where the product is taken over all cells $x$ in the Young diagram $\lambda$, and $h_x$ is
the number of cells in the hook formed by
$x$, which consists of $x$, the cells that are  below $x$ in the same column, and the cells
to the right of
$x$ in the same row. 
\end{theorem} 

One can generalize  Specht modules to skew shapes.  By removing a smaller skew
diagram $\mu$ from the northwest corner of a skew diagram $\lambda$, one gets a {\em skew}
diagram denoted by $\lambda/\mu$.  For example if $\lambda = (4,3,3)$ and $\mu = (2,1)$
then
$$\lambda/ \mu = 
\tiny{\tableau[scY]{
\White  , ,\Black , |\White ,\Black ,|,,|}}.$$ 
Skew Specht modules $S^{\lambda/\mu}$ are defined analogously to \lq\lq straight'' 
Specht modules. There is a submodule characterization and a quotient characterization. 
Theorem~\ref{specht} and Corollary~\ref{skewcor} hold in the skew setting.  There is a
classical combinatorial rule for decomposing Specht modules of skew shape into irreducible
straight shape Specht modules called the Littlewood-Richardson rule, which we will not
present here.

\begin{example}\label{impex} Some important classes of skew and straight Specht modules are
listed below.

\vspace{.2in} $\lambda/\mu =\tiny{\tableau[scY]{\White , , ,,,\Black |\White , 
,,,\Black\cdot^\cdot\White|\White , , ,\Black|\White,
,\Black|\White,
\Black|}}$\qquad $S^{\lambda/\mu} = $ regular representation.

\vspace{.2in} $\lambda \hspace{.18in}= \,\,\tiny{\tableau[scY]{,
,|}}$ $\cdots \tiny{\tableau[scY]{|}}$
\qquad
$S^\lambda =$ trivial representation

\vspace{.2in} $\lambda \hspace{.18in}= \,\,
\begin{array}{c}\tiny{\tableau[scY]{||}}\\\vdots \\\tiny{\tableau[scY]{|}} \end{array}
$
\qquad\qquad
$S^\lambda =$ sign representation

\vspace{.1in}
\noindent where the sign representation is the 1-dimensional representation $V=\C$ whose
character is $\sgn(\sigma)$, that is
$\sigma z =
\sgn(\sigma) z$ for all $\sigma \in \mathfrak S_n $ and $z \in
\C$.  We denote the sign representation by $\sgn_n$ or $S^{(1^n)}$ and we denote the
trivial representation by $1_{\s_n}$ or $S^{(n)}$.  
\end{example}

\begin{xca} For  skew or straight shape $\lambda/\mu$, let $\chi^{\lambda/\mu}$ denote the
character of the representation $S^{\lambda/\mu}$, and for $\sigma \in \mathfrak S_n$, let
$f(\sigma)$ denote the number of fixed points of $\sigma$.
\begin{enumerate}
\item[(a)] Show $\chi^{(n-1,1)}(\sigma) = f(\sigma) -1$ for all $\sigma \in \mathfrak S_n$.
\item[(b)] Show $\chi^{(n,1) /(1)} (\sigma) = f(\sigma)$ for all $\sigma \in \mathfrak
S_n$.
\item[(c)] Find the character of each of the representations in Example~\ref{impex}.
\end{enumerate}
\end{xca}

A {\em skew hook}
is a connected skew diagram that does not contain the subdiagram $(2,2)$. Each cell of a
skew hook, except for 
 the southwestern most and northeastern most  end cells, has exactly two cells adjacent to
it.  Each of the end cells has only one cell adjacent to it.  Let $H$ be a skew hook with
$n$ cells.  We label the cells of $H$ with numbers
$1$ through $n$, starting 
 at the  southwestern  end cell, moving through the adjacent cells,   and ending at the
northeastern  end cell.  For example, we have the labeled skew hook 
$${\tableau[scY]{\White , , ,,,\Black 11|\White , ,
,\Black 8,9,10|\White, ,,\Black 7|\White, ,,\Black 6|\White,,,
\Black 5|1,2,3,4}}.$$ If cell
$i+1$ is above cell $i$ in
$H$ then we  say that the skew hook
$H$ has a {\em descent} at
$i$.  Let $\mbox{des}(H)$ denote the set of descents of $H$.  For each subset $S$ of
$[n-1]$, there is exactly one skew hook with
$n$ cells and descent set
$S$. For example, the skew hook $\tiny{\tableau[scY]{\White , , ,,,\Black |\White , ,
,\Black,,|\White, ,,\Black|\White, ,,\Black|\White,,,
\Black|,,,}}$ is the only skew hook with 11 cells and descent set $\{4,5,6,7,10\}$. 
 
\vspace{.1in}The Specht modules of skew hook shape are called  {\em Foulkes
representations}.   Note
that for any skew hook $H$ with $n$ cells, the set of standard
tableaux of  shape
$H$ corresponds bijectively to the set of permutations in $\mathfrak S_n$ with descent set
$\mbox{des}(H)$. (The   descent set $\mbox{des}(\sigma)$ of a permutation $\sigma \in
\mathfrak S_n$ is the set of all
$i
\in [n-1]$ such that $\sigma(i) > \sigma(i+1)$.)  Indeed, by listing the entries of cells
$1$ through $n$,  one gets a permutation with
descent set $\mbox{des}(H)$.  Hence by Corollary~\ref{skewcor} for skew shapes, the
dimension of the Foulkes representation
$S^H$ is the number of permutations in
$\mathfrak S_n$ with descent set  $\mbox{des}(H)$. 

 A {\em descent} of a standard Young tableau is
an entry
$i$ that   is in a higher row than
$i+1$. By applying the
Littlewood-Richardson rule mentioned above, one gets the following decomposition of the
Foulkes representation into irreducibles, 
\bq \label{foulkes} S^H = \bigoplus_{\lambda \vdash n} c_{H,\lambda} S^\lambda,\eq
where $c_{H,\lambda} $ is the number of standard Young tableaux of shape $\lambda$ and
descent set $\mbox{des}(H)$.

\begin{xca} \label{foudec} Use (\ref{foulkes}) to show that the regular representation of
$\mathfrak S_n$ decomposes into Foulkes representations as follows:
$$ \C \s_n \cong_{\s_n} \bigoplus_{H \in \mbox{SH}_n} S^H,$$ where $\mbox{SH}_n $ is the
set of  skew hooks with
$n$ cells.
\end{xca}

The {\em induction product} of an
$\mathfrak S_j$-module
$U$ and an
$\mathfrak S_k$-module $V$ is the $\mathfrak S_{j+k}$-module
$$U \bullet V := (U \otimes V)\uparrow_{\mathfrak S_j \times \mathfrak S_k}^{\mathfrak
S_{j+k} }.$$
(We are viewing $\s_j\times\s_k$ as the subgroup of $\s_{j+k}$ consisting of permutations
that stabilize the sets
$\{1,2\dots,j\}$ and $\{j+1,j+2,\dots,j+k\}$.)

\begin{xca}  If a skew shape $D$ consists of two  shapes $\lambda$ and $\mu$, where
$\lambda$ and $\mu$ have no rows or columns in common, we say that $D$ is the {\em
disjoint union} of $\lambda$ and $\mu$.  Show that  $S^D = S^\lambda \bullet S^\mu$ if $D$
is the disjoint union of $\lambda$ and $\mu$.   For example  

$$ S^{^{\tiny{\tableau[scY]{\White,,,\Black, ,
|\White,,,\Black ,| ,,||}}}} = S^{^{\tiny{\tableau[scY]{,,||}}}} \bullet
S^{^{\tiny{\tableau[scY]{,,|,|}}}}.$$ 
\end{xca}

\begin{xca} \label{pieri} Let $\lambda\vdash n$.  
\begin{itemize}
\item[(a)]Show that $$S^\lambda \downarrow_{\s_{n-1}}^{\s_n} \cong_{\s_{n-1}}\,\,\,
\bigoplus_{\mu } S^{\mu}$$ summed over all Young diagrams $\mu$ obtained from $\lambda$ by
removing a cell from the end of one of the rows of $\lambda$.  (We are viewing $\s_{n-1}$
as the subgroup of $\s_n$ consisting of permutations that fix $n$.) For example, 
$$ S^{^{\tiny{\tableau[scY]{,,|,|,||}}}}
\downarrow_{\s_7}^{\s_8}\,\,\,\,\,\,\, \cong_{\mathfrak S_7}
\,\,\,\,\, S^{^{\tiny{\tableau[scY]{,|,|,||}}}}\oplus S^{^{\tiny{\tableau[scY]{,,|,|||}}}}
\oplus S^{^{\tiny{\tableau[scY]{,,|,|,|}}}}.$$
\item[(b)]Show that $$S^\lambda \uparrow_{\s_{n}}^{\s_{n+1} }= S^\lambda\bullet S^{(1)}  
\,\, \cong_{\s{n+1}} \,\,
\bigoplus_{\mu } S^{\mu}$$ summed over all Young diagrams $\mu$ obtained from $\lambda$ by
adding a cell to the end of one of the rows of $\lambda$.  For example, 
$$S^{^{\tiny{\tableau[scY]{,,|,|,|}}}}\bullet S^{{^{\tiny{\tableau[scY]{|}}}}
}\,\,\,\,\cong_{\s_8}\,\,\,\, S^{^{\tiny{\tableau[scY]{,,,|,|,|}}}}\oplus
S^{^{\tiny{\tableau[scY]{,,|,,|,|}}}} \oplus S^{^{\tiny{\tableau[scY]{,,|,|,||}}}}.$$
 \end{itemize}
\end{xca}

There is an important generalization of Exercise~\ref{pieri} (b) known as Pieri's rule.
\begth[Pieri's rule]\label{pierith} Let $m,n \in \Z^{+}$.  If  $\lambda \vdash n$ then
$$S^\lambda \bullet S^{(m)} \cong_{\s_{m+n}} \bigoplus_{\mu} S^\mu,$$
summed over all partitions $\mu$ of $m+n$ such that $\mu$ contains $\lambda$ and the skew
shape $\mu/\lambda$ has at most one cell in each column.  Similarly
$$S^\lambda \bullet S^{(1^m)} \cong_{\s_{m+n}} \bigoplus_{\mu} S^\mu,$$
summed over all partitions $\mu$ of $m+n$ such that $\mu$ contains $\lambda$ and the skew
shape
$\mu/\lambda$ has at most one cell in each row.
\enth

The conjugate of a partition $\lambda$ is the partition  $\lambda^\prime$ whose Young
diagram is the transpose of that of $\lambda$.  
 
\begin{theorem} For all partitions
$\lambda\vdash n$,
$$ S^\lambda \otimes \sgn_n \cong_{\s_n} S^{\lambda^\prime},$$
where  the tensor product is an
inner tensor product. This also holds for skew diagrams.
\end{theorem}

\begin{xca}  Let $V$ be a representation of $\s_n$.  Show that $$(V \otimes
\sgn_n)\uparrow_{\s_n}^{\s_{n+1}}\,\,\, \cong_{\s_{n+1}} \,\,\,
V\uparrow_{\s_n}^{\s_{n+1}}
\otimes\,\, \sgn_{n+1},$$ and 
$$(V \otimes
\sgn_n)\downarrow^{\s_n}_{\s_{n-1}}\,\,\, \cong_{\s_{n-1}} \,\,\,
V\downarrow^{\s_n}_{\s_{n-1}}
\otimes\,\, \sgn_{n-1}.$$
\end{xca}

\section{Group actions on poset (co)homology}

Let $G$ be a finite group.  A  $G$-simplicial complex  is a simplicial complex
together with an
action of $G$ on its vertices that takes faces to faces.  
  A 
$G$-poset  is a poset together with a
$G$-action on its elements that preserves the partial order; i.e.,  $x <y \Rightarrow gx
<gy$.  So if $P$ is a $G$-poset then its order complex $\Delta(P)$ is a
$G$-simplicial complex and if
$\Delta$ is a
$G$-simplicial complex then its face poset $P(\Delta)$ is a $G$-poset. 

A $G$-{\it space} is a topological space on which $G$ acts as a group of homeomorphisms.  
If $\Delta$ is a $G$-simplicial complex then    
 the geometric realization
$\|\Delta\|$ is a
$G$-space under the natural induced action of $G$.

\begin{example} The subset lattice $B_n$ is an $\mathfrak S_n$-poset.  The action of a
permutation
$\sigma \in \mathfrak S_n$ on a subset $\{a_1,\dots,a_k\}$ is given by
\bq \label{snaction} \sigma \{a_1,\dots,a_k\} = \{\sigma(a_1),\dots,\sigma(a_k)\}.\eq   
\end{example}

\begin{example}\label{gposetsym} The partition lattice $\Pi_n$ is an $\mathfrak
S_n$-poset.   The action of a permutation
$\sigma \in \mathfrak S_n$ on a partition $\{B_1,\dots, B_k\}$ is given by
\bq\label{partact} \sigma \{B_1,\dots, B_k\} = \{\sigma B_1, \dots, \sigma B_k\},\eq 
where $\sigma B_i$ is defined in (\ref{snaction}). The
symmetric group is isomorphic to the group generated by reflections about hyperplanes in
the braid arrangement.  The action described here is simply the action of the reflection
group on intersections of hyperplanes in the braid arrangement.
\end{example}

\begin{example} \label{crossaction} The face lattice $C_n$ of the $n$-cross-polytope (or
the  face lattice of the
$n$-cube) is an $\s_n[\Z_2]$-poset.  
 The  wreath product group
$\s_n[\Z_2]$ is also known as the hyperoctahedral group or the type B Coxeter group  (see
Section~\ref{symm} for the definition of wreath product).  It
is  the group generated by reflections about the hyperplanes in the type B braid
arrangement. 

By viewing the $n$-cross-polytope as the convex hull of the points $\pm e_i$,
$i = 1,\dots, n$,  one obtains the action of  $\s_n[\Z_2]$ on $C_n$.  We describe this
action in combinatorial terms. The   
$(k-1)$-dimensional faces  of the cross-polytope are convex hulls of certain $k$
element subsets of 
$\{\pm e_i:i  \in [n]\}$.  These are the ones that don't contain both $e_i$ and $-e_i$ for
any $i$.  Thus the $(k-1)$-faces can be  identified with 
$k$-subsets
$T$  of
$[n]
\cup
\{\bar i :i \in [n]\}$ such that $\{i,\bar i\} \nsubseteq T$ for all $i$.  By ordering 
these   sets by containment, one gets the face poset of the $n$-cross-polytope. For
example, the
$4$-subset
$\{3,\bar 5,6,
\bar 8\}$ is identified with the convex hull of the points $e_3, -e_5, e_6, -e_8$.    The
elements of
$\s_n[\Z_2]$ are identified with permutations 
$\sigma$ of
$[n]
\cup
\{\bar i :i \in [n]\}$ for which $\overline{\sigma(i)} = \sigma(\bar i)$ for all $i$
(where $\bar{\bar a} = a$).  
  Then the action of a permutation $\sigma \in \s_n[\Z_2]$  on a $k$-subset
$\{a_1,\dots,a_k\}$ of $[n]
\cup
\{\bar i :i \in [n]\}$ is given  by (\ref{snaction}).

\end{example}

\begin{example} \label{bbraidact} The type B partition lattice $\Pi^B_n$ is also an
$\s_n[\Z_2]$-poset.  Recall from Example~\ref{bbraidex} that the type B partition lattice
$\Pi^B_n$ is the intersection lattice of the type
$B$ braid arrangement. Since this  arrangement is  invariant under  reflection  about any
hyperplane in the arrangement, elements of the reflection group $\s_n[\Z_2]$ map
intersections of hyperplanes to intersections of hyperplanes. This gives the
action of $\s_n[\Z_2]$ on $\Pi^B_n$. There is a combinatorial description of the action
analogous to (\ref{partact}) (see e.g. \cite{gw}).
 \end{example}

\begin{example} The lattice of subspaces of an $n$-dimensional vector space over a finite
field
$F$ is a $GL_n(F)$-poset, where the general linear group $GL_n(F)$ acts in the obvious
way. 
\end{example}

\begin{example} The lattice of subgroups of a finite group $G$ ordered by inclusion is a
$G$-poset, where $G$ acts by conjugation.
\end{example}

\begin{example}
The  semilattice of 
$p$-subgroups of a group
$G$ ordered by inclusion is a $G$-poset, where   $G$ acts by conjugation.
\end{example}

Let $P$ be a $G$-poset.  Since  $g \in G$ takes $j$-chains to $j$-chains,
$g$ acts as a linear map on
$C_j(P;\C)$.  It is easy to  see that $$g \partial (c) = \partial(g c) \mbox{ and } g
\delta (c) = \delta(g c).$$  Hence $g$ acts as a linear map on $\tilde H_j(P;\C)$ and 
on $\tilde H^j(P;\C)$.  This means that whenever $P$ is a $G$-poset, $\tilde H_j(P;\C)$ and 
 $\tilde H^j(P;\C)$ are $G$-modules.  The bilinear form (\ref{form}), induces a pairing
between
$\tilde H_j(P;\C)$ and 
 $\tilde H^j(P;\C)$, which allows one to view them as dual $G$-modules.
For $G =  \mathfrak S_n$ or $G= \s_n[\Z_2]$, 
$$\tilde H_j(P,\C) \cong_G \tilde H^j(P,\C)$$ since dual $\s_n$-modules (resp., dual
$\s_n[\Z_2]$-modules) are isomorphic.

Given a $G$-simplicial complex $\Delta$,  the natural homeomorphism from
$\Delta$ to its barycentric subdivision $\Delta(P(\Delta))$ commutes with the
$G$-action.  Consequently, for all $j \in \Z$
$$
\tilde H_j(\Delta;\C) \cong_G \tilde H_j(P(\Delta);\C),$$
and 
$$
\tilde H^j(\Delta;\C) \cong_G \tilde H^j(P(\Delta);\C).$$

\begin{xca} The maximal chains of $\bar B_n$ correspond bijectively to tableaux of shape
$1^n$ via the map
$$\tableau[scY]{t_1| t_2|^{\vdots}  |t_n } \quad \mapsto \quad (\{t_1\}\subset
\{t_1,t_2\}\subset\dots\subset\{t_1,\dots,t_{n-1}\}).$$ 
\begin{itemize} \item[(a)] Show that the Garnir relations
map to the  coboundary relations.
Consequently, the representation of
$\mathfrak S_n$ on the top cohomology $\tilde H^{n-2}(\bar B_n;\C)$  is the sign
representation.
\item[(b)]   Show that the polytabloids map to cycles in top homology. 
\end{itemize}
\end{xca}

There is an  equivariant version of the Euler-Poincar\'e formula (Theorem~\ref{eupo}),
known as the Hopf-trace formula.  It is convenient to state this formula  in
terms of virtual representations,
which we  define first.
The  {\em representation
group} $\mathcal G(G)$ of a group
$G$  is the free abelian group on the set of  all isomorphism classes $[V]$ of 
$G$-modules
$V$ modulo the subgroup generated by all $[V\oplus W] -[ V]-[W]$.  Elements of the 
representation group are called {\em virtual representations}.  If $V$ is an actual
representation, we denote the virtual representation $[V]$ by $V$.  Note that  two
virtual representations
$A-B$ and $C-D$, where
$A,B,C,D$ are actual representations of $G$, are  equal in the representation group if and
only if 
$A\oplus D
\cong_G B
\oplus C$ in the usual sense. We will  write  $A-B
\cong_G C-D$.

\begth[Hopf trace formula] \label{hopf} For any $G$-simplicial complex $\Delta$,
$$\bigoplus_{i=-1}^{\dim \Delta} (-1)^i C_i(\Delta; \C)  \cong_G \bigoplus_{i=-1}^{\dim
\Delta} (-1)^i \tilde H_i(\Delta; \C) $$
\enth

We will usually suppress the $\C$ from our notation $\tilde H_i(P;\C)$ (resp., $\tilde
H^i(P;\C)$) and write $\tilde H_i(P)$ (resp., $\tilde
H^i(P)$) instead when viewing (co)homology as a  $G$-module.

\section{Symmetric functions, plethysm, and wreath product modules} \label{symm}

 Symmetric
functions  provide a convenient way of describing and computing
representations of the symmetric group.   In this section we  give the basics of
symmetric function theory.  Then we demonstrate its use in computing homology of an
interesting example known as the matching complex.

\subsection{Symmetric functions}

Let $ x =(x_1,x_2,\dots)$ be an infinite sequence of indeterminates.  
 A {\em
homogeneous symmetric function of degree
$n$}  is a formal power series $f(x) \in \Q[[x]]$ in which each term has degree $n$ and
$ f(x_{\sigma(1)},x_{\sigma(2)}, \dots) = f(x_1,x_2,\dots)$ for all permutations $\sigma$
of
$\Z^+$.

Let $\Lambda^n$ denote the set of homogeneous symmetric functions of degree $n$ and let 
$\Lambda = \bigoplus_{n \ge 0} \Lambda^n$.  Then $\Lambda$ is a graded $\Q$-algebra, since
 $\alpha f(x) + \beta g(x) \in \Lambda^n$ if  $f(x),g(x) \in \Lambda^n$ and
$\alpha,\beta \in \Q$; and $f(x)g(x)
\in
\Lambda_{m+n}$ if $f(x) \in \Lambda^m$ and $g(x) \in \Lambda^n$.

There are several important bases for the vector space $\Lambda^n$.  We  mention
just two of them here; the basis of power sum symmetric functions 
and the  basis of Schur functions.   These bases are indexed by partitions of
$n$.  For $n \ge 1$, let
$$p_n =
\sum_{i \ge 1} x_i^n$$ and let $p_0 =1$.  The {\em power sum symmetric function} indexed by
$\lambda = (\lambda_1\ge\lambda_2\ge\dots\ge\lambda_k)$ is defined by 
$$ p_\lambda := p_{\lambda_1} p_{\lambda_2} \dots p_{\lambda_k}.$$

Given a Young diagram  $\lambda$, a {\em semistandard  tableau} of shape $\lambda$
is a filling of $\lambda$ with positive integers so that the rows weakly increase and the
columns strictly increase.  Let $SS_\lambda$ be the set of semistandard tableaux of shape
$\lambda$.  Given a semistandard tableau $T$, define 
$w(T):= x_1^{m_1}x_2^{m_2}\cdots$, where for each $i$, $m_i$ is the number of times $i$
appears in
$T$.  Now define the
Schur function indexed by
$\lambda$ to be
$$s_\lambda: = \sum_{T \in SS_\lambda} w(T).$$
Skew shaped Schur functions $s_D$ are defined analogously for all skew diagrams $D$.  

While it  is obvious that the power sum symmetric functions are symmetric functions,  it is
not obvious that the Schur functions defined this way are. 

\begth The sets $\{p_\lambda : \lambda \vdash n\}$ and $\{s_\lambda : \lambda \vdash n\}$
form  bases for $\Lambda^n$.  Moreover, $\{s_\lambda : \lambda \vdash n\}$ is an integral
basis, i.e., a basis for the $\Z$-module $\Lambda^n_\Z$ of homogeneous symmetric functions
of degree
$n$ with integer coefficients. 
\enth

The Schur function $s_{(n)}$ is known as the complete homogeneous symmetric function of
degree $n$ and is denoted by $h_n$.  The Schur function $s_{(1^n)}$ is known as the
elementary symmetric function of degree $n$ and is denoted by $e_n$. There is an important involution
$\omega:\Lambda_n
\to
\Lambda_n$ defined by 
$$ \omega (s_\lambda) = s_{\lambda^\prime},$$ where $\lambda^\prime$ is the conjugate of
$\lambda$.  Clearly
$$ \omega(h_n) = e_n.$$ For the power sum symmetric functions we have
\bq \label{omega}\omega(p_\lambda) = (-1)^{|\lambda|-l(\lambda)} p_\lambda.\eq

\begin{xca} \label{elemhom}
Prove
$$\sum_{n \ge 0} h_n = \prod_{i\ge 1} (1-x_i)^{-1},$$
and
$$ \sum_{n \ge 0} e_n = \prod_{i\ge 1} (1+x_i),$$
where $h_0 = e_0 =1$.
\end{xca}

The {\em Frobenius 
characteristic} $\ch(V)$ of a representation $V$ of $\mathfrak S_n$  is
the symmetric function given by
$$\ch(V) := \sum_{\mu \vdash n}{1 \over z_\mu } \chi^V(\mu)\, p_\mu,$$ where
$z_\mu := 1^{m_1} m_1! 2^{m_2} m_2!\dots n^{m_n}m_n!$ for $\mu =1^{m_1} 2^{m_2} \dots
n^{m_n}$ and $\chi^V(\mu)$ is the character $\chi^V(\sigma)$ for $\sigma \in \s_n$ of
conjugacy type
$\mu$.  Some basic facts on Frobenius characteristic are compiled in the next result.

\vspace{.2in}\begth \label{omegasgn} \begin{enumerate}
\item[]
\item[(a)] For all (skew or straight) shapes $\lambda$, $$\ch(S^\lambda) = s_\lambda.$$
 \item[(b)] For all representations $V$ of $\mathfrak S_n$, 
$$  \omega(\ch
V) =
\ch(V \otimes
\sgn_n).
$$
\item[(c)]  For all representations $U, V$ of $\mathfrak S_n$,
$$\ch(U \oplus V) = \ch(U) + \ch(V)$$
\item[(d)] For all representations $U$ of $\mathfrak S_m$ and  $V$ of $\mathfrak S_n$,
$$\ch(U \bullet V) = \ch(U)\,\, \ch(V).$$
\end{enumerate}
\enth

The direct sum $\bigoplus_{n
\ge 0}
\mathcal G(\mathfrak S_n)$ of representation groups  is a ring under the induction
product.   It follows from Theorem~\ref{omegasgn} that the Frobenius characteristic map is
an isomorphism from  the  ring $\bigoplus_{n
\ge 0}
\mathcal G(\mathfrak S_n)$ 
 to 
the ring of  symmetric functions over $\Z$.

\begin{definition}  Let $f \in \Lambda$ and let $g$ be a    formal power series   with
positive integer coefficients.   Choose any ordering of the
monomials of $g$, where a monomial appears in the ordering  $m_i$ times if its
coefficient is
$m_i$.
   For example, if
$g = 3 y_1 y_2^2 + 2 y_2y_3 + \dots$ then  the monomials can be arranged as $$(y_1 y_2^2,
y_1 y_2^2, y_1 y_2^2, y_2y_3, y_2y_3,\dots).$$  
Pad the sequence of monomials with zero's if $g$ has a finite number of terms.  Define the
{\em plethysm} of
$f$ and
$g$, denoted
$f[g]$, to be the formal power series obtained from $f$ by replacing the indeterminate
$x_i$ with the $i$th monomial of $g$ for each $i$.  Since $f$ is a symmetric function, the
chosen order of the monomials doesn't matter.  For example, if $f= \sum_{n\ge 0} e_n =
\prod_{i\ge 1} (1+x_i)$ and $g$ is as above then 
$$f[g] = (1+y_1 y_2^2)(1+y_1 y_2^2)(1+y_1 y_2^2)(1+y_2y_3)(1+y_2y_3) \cdots.$$  
\end{definition}

 The following proposition is immediate.
\begin{prop}\label{pleth}
Suppose $f,g \in \Lambda$ and   $h$ is a formal power series
with positive integer coefficients.  Then
\begin{itemize} \item
If $f$ has positive integer coefficients then
$  f[p_n] = p_n[f]$ is obtained by replacing each indeterminate $x_i$ of $f$ by $x_i^n$.
\item $(af +bg)[h] = af[h] + bg[h]$, where $a,b \in \Q$
\item $fg[h] = f[h]g[h]$.  
\end{itemize}
\end{prop}

Note that if $g \in \Lambda$ has positive integer coefficients then $f[g] \in \Lambda$. 
One can extend the definition of plethysm to  all $g \in \Lambda $, by using
Proposition~\ref{pleth} and the fact that the power sum symmetric functions form a basis
for
$\Lambda$.  Hence plethysm is a binary operation  on $\Lambda$, which is clearly
associative, but not commutative.   Note that the plethystic identity
 is $p_1 = h_1$.   We say that $f,g \in \Lambda$ are {\em plethystic inverses} of each
other, and write $g = f^{[-1]}$, if
$f[g] = g[f] = h_1$.

\subsection{Composition product and wreath product}  Our purpose for  introducing plethysm
in these lectures is that plethysm encodes a product operation on symmetric group
representations called composition product, which is described below.  (This description is
based on an  exposition given in
\cite{su96}.)

Let $G$ be a finite group. The {\em wreath product} of $\s_m$ and $G$,  denoted
by $\s_m[G]$, is defined to be the set of $(m+1)$-tuples $(g_1,g_2,
\dots,g_m; \tau)$ such that $g_i \in G$ and  $\tau
\in
\s_m$ with multiplication given by 
$$(g_1,
\dots,g_m; \tau) (h_1,
\dots,h_m; \gamma) = (g_1h_{\tau^{-1}(1)},
\dots,g_m h_{\tau^{-1}(m)}; \tau \gamma) .$$
 The following proposition is immediate.

\begin{prop} \label{wreathhomo} The map $(\alpha_1,\alpha_2,
\dots,\alpha_m; \tau) \mapsto \tau$ is a homomorphism from $\s_m[\s_n]$ onto $\s_m$.
\end{prop}

\begin{definition} Let
$V$ be an
$\s_m$-module and
$W$ be a
$G$-module. Then the {\em wreath product} of $V$ with $W$, denoted
$V[W]$, is the inner tensor product of two $\s_m[G]$-modules:
$$V[W] = \widetilde{W^{\otimes m}} \otimes \hat V,$$
where $\widetilde{W^{\otimes m}}$ is the vector space $W^{ \otimes
m}$
with $\s_m[G]$ action given by 
\bq \label{wreathac}(\alpha_1,
\dots,\alpha_m; \tau) (w_1 \otimes \cdots \otimes w_m) =
\alpha_1 w_{\tau^{-1}(1)} \otimes \cdots \otimes  \alpha_m
w_{\tau^{-1}(m)}\eq
and $\hat V$ is the pullback of the representation of $\s_m$
on $V$ to $\s_m[G]$ through the homomorphism given in Proposition~\ref{wreathhomo}. 
That is,
$\hat V$ is the representation of $\s_m[G]$  on
$V$ defined by
$$(\alpha_1,
\dots,\alpha_m; \tau)v = \tau v.$$
\end{definition}

Given a finite set $A= \{a_1 < a_2 < \dots < a_n\}$ of positive integers, let $\s_A$
be the set of permutations of the set $A$.  We shall view a
permutation in
$\s_A$ as a word whose letters come from $A$. For $\sigma \in
S_n$,  let $\sigma^A$ denote the word
$a_{\sigma(1)}a_{\sigma(2)}
\cdots a_{\sigma(n)}.$
We shall view an element of the Young subgroup $\s_k \times
\s_{n-k}$ of $\s_n$ as the concatenation $\alpha \star
\beta$ of words $\alpha \in \s_k$ and $\beta \in
\s_{\{k+1,\dots,n\}}$. 
The wreath product $\s_m[\s_n]$ is isomorphic to  the normalizer of the
Young subgroup $ \underbrace{\s_n \times \cdots
\times \s_n}_{\mbox{$m$
times}}$ in $\s_{mn}$.  The isomorphism is given by
$$(\alpha_1,\dots,\alpha_m;\tau) \quad\mapsto\quad
\alpha_{\tau(1)}^{A_{\tau(1)}}\star
\cdots
\star
\alpha_{\tau(m)}^{A_{\tau(m)}},$$
where $A_i =[in]\setminus [(i-1)n]$.

Define the {\em composition product} of an $\s_m$-module $V$ and an $\s_n$-module $W$ by
$$V\circ W:= V[W]\uparrow_{\s_m[\s_n]}^{\s_{mn}}$$
The following result relates plethysm to composition product.

 \begth \label{compprod} Let $V$ be an 
$\s_m$-module and $W$  be an $\s_n$-module.    Then
$$\ch (V\circ W) = \ch V [\ch W].$$
\enth

\begin{example}\label{exd}
Given  $\lambda \vdash n$, let $\Pi(\lambda)$ be the set of  partitions  of $[n]$
whose block sizes form the partition $\lambda$.  In this example, we set  $\lambda =
d^{b}$, where $d$ and $b$ are positive integers.

(a)  The symmetric group $\s_{bd}$ acts on $\Pi(d^b)$ as in
Example~\ref{gposetsym}, and this action induces a representation of $\s_{bd}$ on the
complex vector space  $\C\Pi(d^b)$
 generated
by elements of  $\Pi(d^b)$. 
It is not difficult to see that the $\s_{bd}$-module $\C\Pi(d^b)$ is the induction of the
 trivial representation of the stabilizer of the partition $$\pi(d^b) :=
1,\dots,d
\,\,/
\,\,d+1,\dots, 2d\,\,/
\,
\dots
\,\,/
\, (b-1)d+1,\dots,bd$$
to $\s_{bd}$.  The stabilizer of $\pi(d^b)$ is $\s_b[\s_d]$ and the trivial representation
of $\s_b[\s_d]$ is $S^{(b)}[S^{(d)}]$.  So $\C\Pi(d^b)$ is the composition product
$S^{(b)}\circ S^{(d)}$.  By Theorem~\ref{compprod}
$$\ch\,\, \C\Pi(d^b) = h_b[h_d]. $$

(b)  The stabilizer $\s_b[\s_d]$ of $\pi(d^b)$ acts on the interval
$(\pi(d^b),\hat 1)$ of $\Pi_{bd}$.  This induces a representation of $\s_b[\s_d]$ on the
top homology $\tilde H_{b-3}(\pi(d^b),\hat 1)$ (recall $\tilde H_i(x,y)$ denotes complex
homology of the open interval $(x,y)$).  By observing that  the interval
$(\pi(d^b),\hat 1)$ is isomorphic to the poset  $\bar \Pi_b$, one can see that
$$\tilde H_{b-3}(\pi(d^b),\hat 1) \cong_{\s_b[\s_d]} \tilde H_{b-3}(\bar \Pi_b)[S^{(d)}].$$
Next observe that $\bigoplus_{x \in \Pi(d^b)} \tilde H_{b-3}(x,\hat
1)$ is the induction of $\tilde H_{b-3}(\pi(d^b),\hat 1)$ from $\s_b[\s_d]$ to $\s_{bd}$.
So  $\bigoplus_{x \in \Pi(d^b)} \tilde H_{b-3}(x,\hat
1)$ is the composition product $\tilde H_{b-3}(\bar \Pi_b) \circ S^{(d)}$. By
Theorem~\ref{compprod},
\begin{eqnarray*}\ch\,\, \bigoplus_{x \in \Pi(d^b)} \tilde H_{b-3}(x,\hat
1) &=& \ch \,\,\tilde H_{b-3}(\pi(d^b),\hat
1)\uparrow_{\s_b[\s_d]}^{\s_{bd}}
\\&=&
(\ch
\,\,\tilde H_{b-3}(\bar
\Pi_b))[h_d].
\end{eqnarray*}
 Theorem~\ref{compprod} is inadequate when
$\lambda$ is not of the form
$d^b$.  In \cite{su,wa2}, more general results are given, which
enable one to derive  the following generating function,  
\begin{equation} \label{plethT} \sum_{\lambda \in
\mbox{Par}(T,b)}\,\,\,{\Big(}\ch
\,\,
\bigoplus_{x
\in
\Pi(\lambda)}
\tilde H_{b-3}(x,\hat 1){\Big)}\,\,z_{\lambda_1}\cdots z_{\lambda_b} = (\ch
\,\,\tilde H_{b-3}(\bar
\Pi_b))\,\,{\Big[}\sum_{i\in T}z_i h_i{\Big]},\end{equation}
where $T$ is any set 
of positive integers, $\mbox{Par}(T,b)$ is the set of  partitions of length $b$, 
all of whose parts are in
$T$, the $\lambda_i$ are the parts of $\lambda$, and the $z_i$ are
(commuting) indeterminates.
\end{example}

\subsection{The matching complex} \label{matchex}  In this subsection, we further
demonstrate the power of symmetric function theory and plethysm in computing  homology. 
Our example is a well-studied simplicial complex known as the matching complex, which is
defined to be the  
 simplicial complex $M_n$ whose vertices are the 2-element
subsets of
$[n]$ and whose faces are collections of mutually  disjoint 2-element subsets of $[n]$.  
Alternatively,
$M_n$ is the simplicial complex of graphs on node set $[n]$ whose degree is at most
$2$.   Its face poset is the proper part of the poset of partitions of $[n]$ whose block
sizes are at most $2$.  It is not difficult to see that 
\bq \label{compsym} \ch\, C_{k-1}(M_n)  = e_k[h_2] h_{n-2k} .\eq It follows from
this and the Hopf trace formula (Theorem~\ref{hopf}) that 
$$\sum_{k\ge -1} (-1)^{k-1} \ch \tilde H_{k}(M_n) = \sum_{k\ge 0} (-1)^k e_k[h_2]
h_{n-2k},$$ which implies by Exercise~\ref{elemhom} that 
$$\sum_{n \ge 0} \sum_{k\ge -1} (-1)^{k-1} \ch \tilde H_k(M_n) = \prod_{i \le j} (1-x_ix_j)
\prod_{i\ge 1} (1-x_i)^{-1} .$$   The right hand side can be decomposed into Schur
functions by using the following symmetric function identity of Littlewood 
\cite[p.238]{lit}:
$$\prod_{i \le j} (1-x_i x_j)
\prod_{i\ge 1} (1-x_i)^{-1} = \sum_{\lambda = \lambda^\prime} (-1)^{{|\lambda| - r(\lambda)
\over 2}} s_\lambda,$$
where $r(\lambda)$ is the rank of $\lambda$, i.e., the size of the main diagonal (or Durfee
square) of the Young diagram for $\lambda$. From this we conclude that 
$$\bigoplus_{k\ge -1} (-1)^{k-1} \tilde H_k(M_n)  \cong_{\s_n}
\bigoplus_{\scriptsize\begin{array}{c}
\lambda:
\,\lambda
\vdash n
\\
\lambda  =
\lambda^\prime  \end{array}} (-1)^{{|\lambda| - r(\lambda)
\over 2}}S^{\lambda}.$$  With additional work involving long exact sequences of relative
homology, Bouc
 obtains the following beautiful refinement.
\begth[Bouc \cite{bo}] \label{boucrepth} For all $n \ge 1 $ and $k \in \Z$, 
\bq \label{boucresult}\tilde H_{k-1}(M_n) \cong_{\s_n}
\bigoplus_{\scriptsize\begin{array}{c}
\lambda:
\,\lambda
\vdash n
\\
\lambda =
\lambda^\prime \\ r(\lambda)= |\lambda| -2k \end{array} } S^{\lambda}.\eq
\enth
 From Bouc's formula one  obtains a formula for the  Betti number in  dimension $ k-1$  as
the number of standard Young tableaux of self-conjugate shape and rank $n -2k$ (which can
be computed from the hook-length formula, Theorem~\ref{frt}).  This result provides an
excellent illustration of  the use of representation theory in the computation of  Betti
numbers.     

We sketch a proof of Theorem \ref{boucrepth} due to Dong and Wachs
\cite{dw}, which involves a technique called discrete
Hodge theory.  Let $\Delta$ be a $G$-simplicial complex.  The {\it combinatorial Laplacian}
$\Lambda_k:C_k(\Delta)
\to C_k(\Delta)$ is defined by
$$\Lambda_k = \delta_{k-1}\partial_k + \partial_{k+1}\delta_k.$$   A basic result of
discrete  Hodge theory  is that
\bq \label{dhodge}\tilde H_k(\Delta) \, \cong_G \, \mbox{ker} \Lambda_k.  \eq

The key observation of \cite{dw} is that  when one applies the Laplacian $\Lambda_{k-1}$
to  an oriented simplex
$\gamma
\in C_{k-1}(M_n)$, one gets
$$\Lambda_{k-1}(\gamma) = T_n \cdot \gamma,$$ where $T_n = \sum_{1 \le i < j\le n} (i,j)
\in 
\C\s_n$ and $(i,j)$ denotes a
transposition. It is then shown that $T_n$ acts on the Specht module
$S^\lambda$ as multiplication by the scalar $c_\lambda$ defined by 
\beq \label{contenteq} c_\lambda =\sum_{i=1}^{r} \left (
\binom{\alpha_i +1}2 -
\binom {\beta_i+1} 2\right ),\eeq where $r=r(\lambda)$, $\alpha_i = \lambda_i-i$, and
$\beta_i = \lambda^\prime_i-i$.  That is, $\alpha_i$ is the number of cells to the right
of and in the same row as the $i$th cell of the diagonal of $\lambda$ and $\beta_i$ is
the number of cells below and in the same column as the $i$th cell of the
diagonal.  The array $(\alpha_1,\dots,\alpha_r\mid
\beta_1,
\dots,\beta_r)$ is known as  Frobenius notation for $\lambda$.  Note that $\lambda$ is
uniquely determined by its Frobenius notation.  

Next  we   decompose $C_{k-1}(M_n)$ into Specht modules by using (\ref{compsym}) and 
another symmetric function identity of Littlewood, cf., \cite[I
5 Ex. 9b]{mac}, namely,
\beq \label{L1}\prod_{i \le j}(1-x_ix_j) = \sum_{\nu 
\in {\mathcal B}} (-1)^{|\nu|/2}s_\nu,\eeq where ${\mathcal B}$ is the
set of all partitions of the form $(\alpha_1+1,\dots,\alpha_r+1\mid
\alpha_1, \dots, \alpha_r )$ for some $r$. This 
 and Pieri's rule (Theorem~\ref{pierith}) yield  the  decomposition into irreducibles:
$$C_{k-1}(M_n) \cong \bigoplus_{\lambda \in A_n} a_\lambda^k S^\lambda,$$ where $$A_n =
\{ (\alpha_1,\dots,\alpha_r\mid\beta_1,\dots,\beta_r) \vdash n : r \ge 1, \, 
\alpha_i \ge 
\beta_i \,\,\,\,\forall\, i\in [r]\}$$ and $a^k_\lambda$ is a nonnegative integer. If 
$\lambda$ is self-conjugate then  
$$a_{\lambda}^k = \begin{cases} 1 & \text {if } r(\lambda) = n-2k \\ 0 & \text{otherwise}
\end{cases} .$$

Clearly for $\lambda \in A_n$, we have $c_\lambda = 0 $ if and only if $\lambda$ is
self-conjugate.   It follows that
$$\ker \Lambda_{p-1} \cong \bigoplus_{\scriptsize\begin{array}{c}
\lambda:
\,\lambda
\vdash n
\\
\lambda  =
\lambda^\prime \end{array} }  a^k_\lambda \,\, S^{\lambda}
\cong\bigoplus_{\scriptsize\begin{array}{c}
\lambda:
\,\lambda
\vdash n
\\
\lambda  =
\lambda^\prime \\  r(\lambda) = n -2k \end{array} }S^{\lambda}.$$ Theorem~\ref{boucrepth}
now follows by discrete Hodge theory (\ref{dhodge}). 

There is a bipartite analog of the matching complex called the $m \times n$
chessboard complex.  This is the  simplicial complex of
subgraphs of the complete bipartite graph on node sets of size $m$ and $n$, whose nodes
have degree at most 1.  Alternatively, it is the simplicial complex of nonattacking rook
placements on an $m \times n$ chessboard.    Friedman and Hanlon
\cite{fh} use discrete Hodge theory to obtain a chessboard complex analog of Bouc's
formula (\ref{boucresult}).  (In fact the  proof of 
(\ref{boucresult}) described above was patterned on the Friedman-Hanlon proof.)  Not only
can the Betti numbers be computed from this formula, but  Shareshian and Wachs
\cite{shwa2}  use the formula to compute
 torsion in the integral homology of the chessboard complex; see Theorem~\ref{chess}. So
representation theory has proved to be a useful tool in computing torsion in integral
homology as well as Betti numbers.  Karaguezian, Reiner and Wachs \cite{krw}  use
formula~(\ref{boucresult}) and its chessboard complex analog to derive a more
general result of Reiner and Roberts \cite{rr}, which computes the   homology of general
bounded degree graph and bipartite graph complexes.

The matching complex and chessboard complex have arisen in various  areas of mathematics,
such as group theory (see Section~\ref{quilfibsec}), commutative algebra  \cite{rr}
and discrete geometry \cite{zv}.  These complexes are discussed further in
Lecture~\ref{opsmaps}.  See
\cite{wa3} for a survey article on the topology of matching complexes, chessboard
complexes and   general bounded degree graph complexes.

%% file: lect3.tex

\lecture{Shellability and edge labelings} \label{lexshell}
\section{Shellable simplicial complexes}   Shellability is a combinatorial property of
simplicial and more general cell complexes, with
strong topological and algebraic consequences.   
Shellability first appeared in the middle of the nineteenth century  in 
Schl\"afli's computation of the   Euler characteristic of a convex
polytope
\cite{sch}.  Schl\"afli made the assumption, without proof, that the boundary complex of a
convex polytope is shellable.  This assumption was eventually proved in  1970  by
Brugesser and Mani \cite{bm} and was used in McMullen's
 proof \cite{mc}  of the famous upper bound conjecture for convex polytopes (the upper
bound conjecture is stated in Section~\ref{CMsec}).  

The original theory of shellability  applied only  to pure  complexes.  In the
early 1990's, a nonpure simplicial complex arose in the  complexity theory work
of Bj\"orner, Lov\'asz and Yao \cite{bly} with  topological
properties somewhat similar to those  of  pure shellable complexes.  This led Bj\"orner
and Wachs
\cite{bw96,bw97}  to extend  the theory  of shellability to nonpure  complexes.  The 
Goresky-MacPherson formula created a need for such an extension, since, unlike for
hyperplane arrangements,  the intersection semilattice of a subspace arrangement is not
necessarily pure.  As it turned out,  there were many other uses for nonpure
shellability.

 For each face
$F$ of a simplicial complex
$\Delta$, let $\langle F\rangle$ denote the subcomplex generated by
$F$, i.e.,
$\langle F\rangle = \{G : G \subseteq F\}$.   A simplicial complex $\Delta$ is said to be
{\em shellable} if its facets can be arranged in linear order $F_1, F_2, \,\ldots, F_t$ in
such a way that the subcomplex $\left(\bigcup^{k-1}_{i=1}\langle
F_i\rangle\right)\cap\langle F_k\rangle$ is pure and $(\dim F_k-1)$-dimensional for all
$k=2,\, \ldots, t$. Such an ordering of facets is called a {\em shelling}.  We emphasize
that we are not assuming purity; here shellability refers to what is commonly
called nonpure shellability.   
 
\begin{figure}\begin{center}
\mbox{\includegraphics[width=12cm]{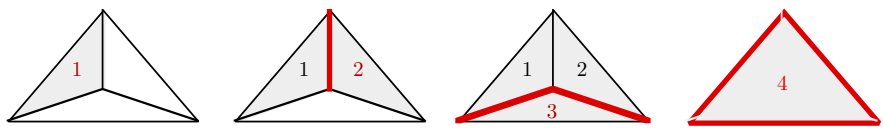}}
\caption{Shelling of the boundary of $3$-simplex}
\label{figshell}\end{center}\end{figure}

\begin{figure} \begin{center}
\mbox{\includegraphics[width=10cm]{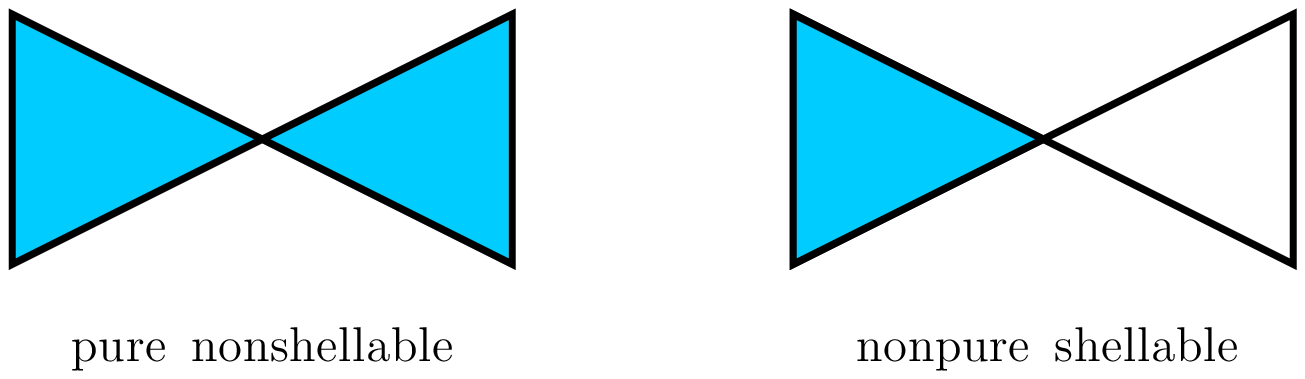}}
\caption{} \label{fignoshell}
\end{center}
\end{figure}

 A shelling of the
boundary complex of the $3$-simplex is given in Figure~\ref{figshell}.  
 A
pure nonshellable  complex   and  a  nonpure shellable complex are given in
Figure~\ref{fignoshell}, where the shading indicates that the triangles are filled in.  

\newpage\begin{xca} \label{35sc} Verify that the third and fifth simplicial complexes below
are not shellable, while the others are.

\vspace{-.2in}
\includegraphics[width=10cm]{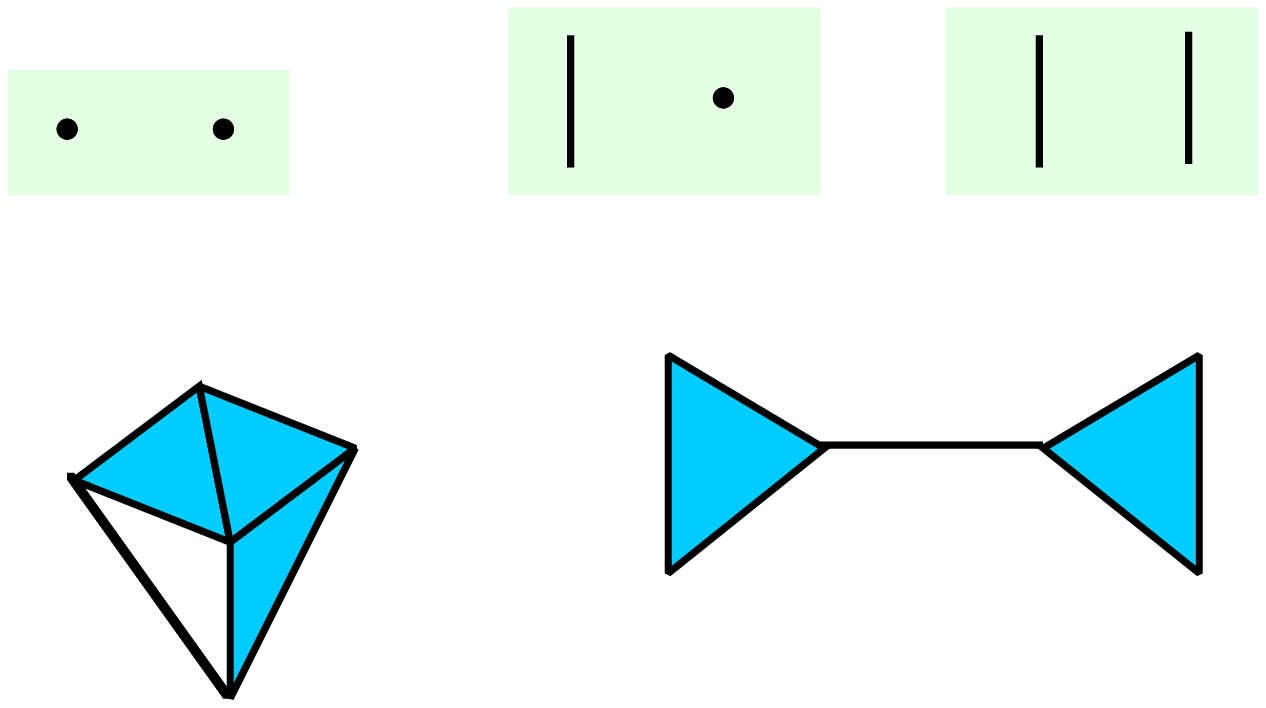}
\end{xca}

\vspace{.2in} 

\begin{theorem}[Brugesser and Mani \cite{bm}]  The boundary  complex of a convex
polytope is shellable.
\end{theorem}

A geometric construction called a line shelling is used to prove this result.  See
Ziegler's book \cite{z95} on polytopes  for a description of this basic construction.

Shellability is not a topological property. By this we mean that shellability of a complex
is not determined by the topology  of its geometric realization; a shellable  simplicial
complex can be homeomorphic to a nonshellable simplicial complex.  In fact, there exist
nonshellable triangulations of the 3-ball and the 3-sphere; see
\cite{li},\cite[Chapter 8]{z95}, \cite{z98}, \cite{lu}.  

Shellability does have strong topological consequences, however, as is shown by the
following result.  By a {\it wedge} $\lor_{i=1}^n X_i$ of $n$ mutually disjoint
connected topological spaces
$X_i$, we mean the  space obtained by selecting  a base point for each $X_i$ and
then identifying all the base points with each other. 

\begin{theorem}[Bj\"orner and Wachs \cite{bw96}] \label{wedge} A shellable simplicial
complex  has the homotopy type of a wedge of spheres (in varying dimensions), where for
each $i$, the number of $i$-spheres is the number of $i$-facets whose entire boundary is
contained in the union of the earlier facets.  Such facets are usually called {\em homology
facets}.
\end{theorem}

\begin{proof}[Proof idea]  Let $\Delta$ be a shellable simplicial complex.  We first
observe that any shelling of 
$\Delta$ can be rearranged to produce a shelling in which  the homology
facets come last.    So $\Delta$ has a
shelling
$F_1,F_2,
\dots F_k$  in which $F_1,\dots, F_j$ are not homology facets and $F_{j+1}, \dots, F_k$
are, where $1 \le j \le k$. The basic idea of the proof is that as we
attach the first $j$ facets, we  construct a contractible simplicial complex at each
step.  The remaining homology facets  attach as spheres since the entire boundary of
the facet is identified with a point.   
\end{proof}

\begin{cor} \label{shellhom} If $\Delta$ is shellable then for all $i$,
\begin{equation}\label{homshell}\tilde H_i(\Delta;\Z ) \cong \tilde H^i(\Delta;\Z) \cong
\Z^{r_i},\end{equation} where $r_i $ is the number of homology $i$-facets of $\Delta$.
\end{cor}

 The  homology facets yield more than just  the Betti numbers;
they form a basis for  cohomology.    We discuss the connection with cohomology further
 in the context of lexicographic shellability in the  next section.

There are certain operations on complexes and posets that preserve shellability. 
 For a simplicial complex $\Delta$ and $F \in \Delta$, define the {\em link} to be the
subcomplex given by
$$\lk_{\Delta} F := \{G \in \Delta : G \cup F \in \Delta, \, G \cap F = \emptyset\}.$$

\begth[\cite{bw96}] \label{linkth} The link of every face of a shellable complex is
shellable.  
\enth

Define the {\em suspension} of a simplicial complex $\Delta$ to be
$$\mbox{susp}(\Delta) := \Delta * \{\{a\}, \{b\}\},$$
where $*$ denotes the join, and $a\ne b$ are not vertices of $\Delta$.

\begth[\cite{bw96}] \label{join} The
join of two  simplicial complexes is shellable if and only if each complex is shellable.
In particular,  $\mbox{susp}(\Delta)$  is shellable if and only if  $\Delta$ is shellable.
\enth

Define the $k$-skeleton of a simplicial complex $\Delta$ to be the subcomplex consisting
of all faces of dimension $k$ or less.

\begth[\cite{bw97}] \label{skelshell} The $k$-skeleton of a shellable simplicial complex is
shellable for all
$k
\ge 0$.
\enth

 A simplicial complex   $\Delta$ is said to be
$r$-{\it connected} (for $r\ge 0$) if it is nonempty and
connected and its
$j$th homotopy group $\pi_j(\Delta)$ is trivial for all
$j=1,\dots,r$. So $0$-connected is the same as connected and $1$-connected is the same as
simply connected. A nonempty simplicial complex
$\Delta$ is said to be
$r$-{\it acyclic} if its
$j$th reduced integral homology group
$\tilde H_j(\Delta)$ is trivial for all
$j=0,1,\dots,r$.  We say that $X$ is
$(-1)$-connected and $(-1)$-acyclic  when $\Delta$ is nonempty. It is also convenient to
say that every simplicial complex is $r$-connected and $r$-acyclic for all $r\le -2$. 

It is a basic fact of homotopy theory that  $r$-connected implies $r$-acyclic and that the
converse holds only for simply connected complexes.  Another basic fact is that a
a simplicial complex is $r$-connected ($r$-acyclic) if and only if its $(r+1)$-skeleton has
the  homotopy type (homology) of a wedge of $(r+1)$-spheres.  This makes shellability  a
useful tool in  establishing
$r$-connectivity and
$r$-acyclicity.

\begin{theorem}\label{shellcon} A simplicial complex is $r$-connected (and therefore
$r$-acyclic)  if its
$(r+1)$-skeleton is pure shellable. 
\end{theorem}

A poset $P$ is said to be shellable if its order complex $\Delta(P)$ is shellable. 
The following is an immediate consequence of Theorems~\ref{linkth} and~\ref{join}.

\begin{cor}\label{boundshell}\begin{enumerate}
\item[] 
\item[(a)] A bounded poset $P$ is shellable if and only
if its proper part
$\bar P$ is shellable.  
\item[(b)] Every (open or closed) interval of a shellable poset is
shellable.  
\item[(c)] The join of two posets is
shellable if and only if each of the posets is shellable.
\end{enumerate}
\end{cor}
 
\begth[Bj\"orner and Wachs \cite{bw97}] \label{prodshell} 
The product of  bounded posets is shellable if and only if each of the posets is
shellable.  
\enth

Sometimes shellability has stronger topological consequences than homotopy type.

\begin{theorem}[Danaraj and Klee \cite{dk}] \label{danklee} Let $\Delta$ be a pure 
shellable
$d$-dimensional simplicial complex in which every codimension $1$ face is contained in at
most
$2$ facets.  Then
$\Delta$ is homeomorphic to a  $d$-sphere or a $d$-ball.  
Moreover, $\Delta$ is
homeomorphic to a $d$-sphere if and only if every codimension $1$ face is contained in
exactly $2$ facets. 
\end{theorem} 

Note that  the condition  on
codimension
$1$ faces in the Danaraj and Klee theorem can be expressed as follows:  Closed length $2$
intervals of  $L(\Delta)$ have at most 4 elements. Pure posets in which every length $2$
interval has exactly $4$ elements are said to be {\em thin}.  We have the following poset
version of the Danaraj and Klee result.

\begin{theorem}[Bj\"orner
\cite{bj84}] \label{bjthin}
If
$\hat P$ is pure, thin and shellable then
$P$ is isomorphic to the face poset of a regular cell decomposition 
(a generalization of simplicial
decomposition; see \cite{hat}) of an $l(P)$-sphere.
\end{theorem}

Many natural classes of simplicial complexes and posets, which have arisen in various
fields of mathematics, have turned out to be shellable (pure or nonpure).  A striking
illustration of the ubiquity of shellability is given by the following result.

\begin{theorem}[Shareshian \cite{sh}] \label{shth} The lattice of subgroups of a finite
group
$G$ is shellable if and only if $G$ is solvable. 
\end{theorem}

The  pure version of this result, which states that {\em the lattice of subgroups of a
finite group $G$ is pure shellable if and only if $G$ is supersolvable}, was proved
by Bj\"orner \cite{bj80} in the late 70's  by introducing a technique called
lexicographic shellability. Shareshian uses a  general (nonpure) version of
lexicographic shellability  more recently introduced by Bj\"orner and Wachs \cite{bw96}. 
We discuss lexicographic shellability in the remaining sections of this lecture.  Further
discussion of connections between group theory and poset topology can be found in
Section~\ref{quilfibsec}.

We end this section with two interesting conjectures. A pure simplicial complex
$\Delta$ is said to be {\it extendably shellable} if  every partial shelling of $\Delta$ 
extends to a shelling of
$\Delta$.  Tverberg (see \cite{dk}) 
conjectured that the boundary complex of every polytope is extendably shellable.
  The conjecture was settled in the negative by Ziegler \cite{z98} who showed that
there are simple and simplicial polytopes whose boundary complex is not extendably
shellable.  The conjecture is, however, true for the simplex since every ordering of the
facets of the boundary complex is a shelling.  Hence the $k = d-1$ case of the following
conjecture holds.

\begin{con}[Simon \cite{si94}] For all $k \le d$, the
$k$-skeleton of the
$d$-simplex is extendably shellable. 
\end{con} 

Simon's conjecture was shown to be true for $k \le 2$ by Bj\"orner and Eriksson
\cite{bjer}, who  extended the conjecture to all matroid complexes.  Recently Hall
\cite{h04} showed that the boundary complex of the $12$-dimensional cross-polytope is a
counterexample to the extended conjecture of Bj\"orner and Eriksson.

A simplicial complex $\Delta$ is said to be  {\it minimally nonshellable} (called an 
obstruction to shellability in \cite{wa99}) if $\Delta$ is not shellable, but every proper
induced subcomplex of $\Delta$ is   shellable.  For example, the complex consisting of
two disjoint $1$-simplexes is a minimally nonshellable complex (see 
Exercise~\ref{35sc}).  Billera and Myers \cite{bm99} showed that this is the only
minimally nonshellable order complex.   It is also the only $1$-dimensional minimally
nonshellable simplicial complex.   It was shown by Wachs
\cite{wa99} that the number of vertices in any $2$-dimensional minimally nonshellable   
simplicial complex is either $5,6$ or $7$; hence there are only finitely many such
complexes.

\begin{con}[Wachs]  There are only finitely many
$d$-dimensional minimally nonshellable simplicial complexes (obstructions to
shellability) for each
$d$.
\end{con}

\section{Lexicographic shellability} \label{ELsec} In the early 1970's, Stanley
\cite{st72,st74} introduced a technique for showing that the M\"obius function of
rank-selected subposets of  certain posets 
alternates in sign.  This technique involved labeling the edges of
the Hasse diagram of the poset in a certain way.  Stanley conjectured that the posets that
he was considering were Cohen-Macaulay, a topological (and algebraic) property of
simplicial complexes  implied by shellability (cf., Section~\ref{CMsec}).  Bj\"orner
\cite{bj80} proved this conjecture by finding a condition
on edge labelings which implies shellability of the poset. From this emerged the theory
of  lexicographic shellability, which was further developed 
in a series of papers 
by 
Bj\"orner and Wachs, first in the pure case
\cite{bw82,bw83} and later in the general (nonpure) case \cite{ bw96,bw97}.

There are two basic versions of lexicographic shellability, EL-shellability and
CL-shellability.  In this section we begin with the  simpler but less powerful version,
EL-shellability, and discuss some of its consequences.

An {\em edge labeling} of a bounded poset $P$ is a map $\lambda: \mathcal E(P) \to
\Lambda$, where
$\mathcal E(P)$ is the set of edges of the Hasse diagram of $P$, i.e., the covering
relations $x <\!\!\!\!\cdot \,\, y$ of $P$, and $\Lambda$ is some poset (usually the
integers
$\Z$ with its natural total order relation).  Given an edge labeling $\lambda:
\mathcal E(P) \to \Lambda$, one can associate  a word $$\lambda(c) = \lambda(\hat 0, x_1)
\lambda(x_1, x_2) \cdots \lambda(x_{t}, \hat 1)$$  with each maximal chain $c = (\hat 0
<\!\!\!\!\cdot
\,\,x_1 <\!\!\!\!\cdot \,\,
\cdots<\!\!\!\!\cdot\,\, x_{t} <\!\!\!\!\cdot\,\, \hat 1)$.   We say that  $c
$ is  {\em increasing} if the associated word $\lambda(c)$ is
{\em strictly}  increasing.   That is, $c$ is  increasing if 
$$ \lambda(\hat 0, x_1) <
\lambda(x_1, x_2)<  \cdots < \lambda(x_{t}, \hat 1).$$  We say that  $c
$ is  {\em decreasing} if the associated word $\lambda(c)$ is
{\em weakly}  increasing.    We can order the
maximal chains lexicographically by using the lexicographic order on the corresponding
words.  Any
edge labeling
$\lambda$ of
$P$ restricts to an edge labeling of any closed interval $[x,y]$ of $P$.  So we may refer
to increasing and decreasing maximal chains of $[x,y]$, and lexicographic order of maximal
chains of
$[x,y]$.

 \begin{definition}  Let $P$ be a bounded poset. {\em  An edge-lexicographical
labeling} (EL-labeling, for short)  of
$P$ is an edge labeling such that in each closed
interval $[x,y]$ of $P$, there is a unique  increasing maximal chain, which
lexicographically precedes all other maximal chains of $[x,y]$.
\end{definition}

\begin{figure}  \begin{center}
\mbox{\includegraphics[width=6cm]{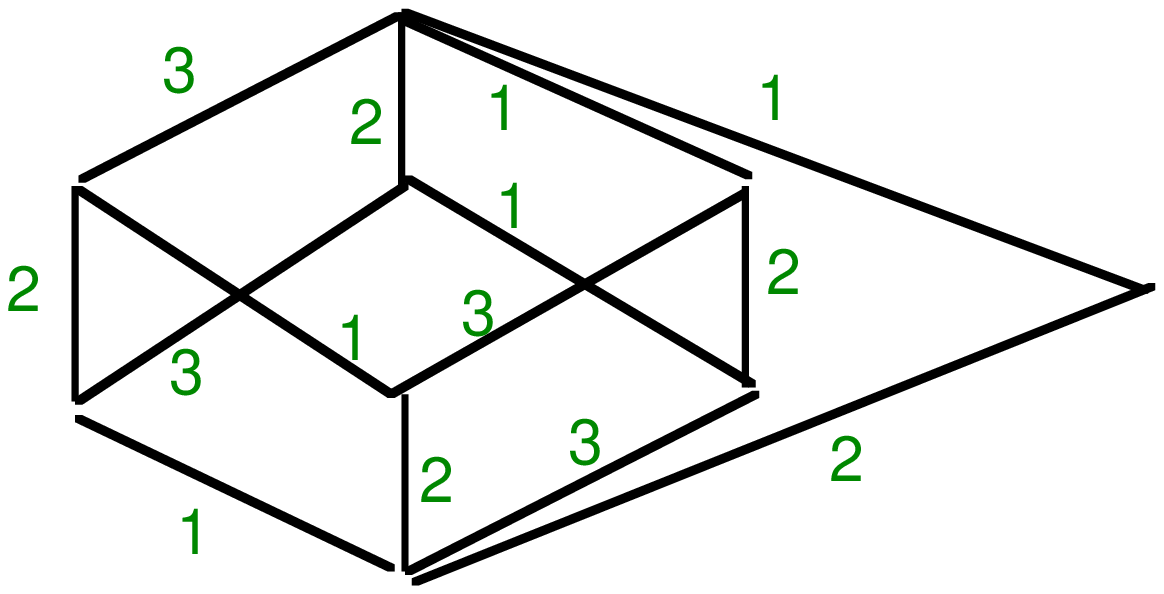}}
\caption{EL-labeling} \label{figel}\end{center}
\end{figure}

An example of an EL-labeling of a poset is given in Figure~\ref{figel}. The leftmost
chain, which has associated word  $123$, is the only  increasing maximal chain of the
interval
$[\hat 0,
\hat 1]$.  It  is also  lexicographically less than all other  maximal chains.  One
needs to check  each interval to verify that
the labeling is indeed an EL-labeling.

A bounded poset that admits an EL-labeling is said to be {\em
edge-lexicographic shellable} (EL-shellable, for short).  The following theorem justifies
the name.
\begin{theorem}[Bj\"orner \cite{bj80}, Bj\"orner and Wachs \cite{bw96}]\label{ELshell}
Suppose
$ P$ is a bounded poset with  an EL-labeling.  Then the lexicographic order of the maximal
chains of
$ P$ is a shelling of
$\Delta(P)$.  Moreover, the corresponding order of the maximal chains of $\bar P$  is a
shelling of
$\Delta(\bar P)$.
\end{theorem}

\begin{xca} Prove Theorem~\ref{ELshell}.
\end{xca}

It is for nonbounded posets that shellability has interesting topological consequences
since the order complex of a  bounded poset is a just a cone.

\begth[Bj\"orner and Wachs \cite{bw96}] \label{lshellhom} Suppose $P$ is a poset for which
$\hat P$ admits an EL-labeling.  Then
$P$ has the homotopy type of a wedge of spheres, where the number of $i$-spheres is the
number of  decreasing maximal $(i+2)$-chains of $\hat P$.  The   decreasing
maximal $(i+2)$-chains, with $\hat 0$ and $\hat 1$ removed, form a basis for
cohomology
$\tilde H^i(P;\Z)$.
\enth

\begin{figure}\begin{center}
\includegraphics[width=6cm]{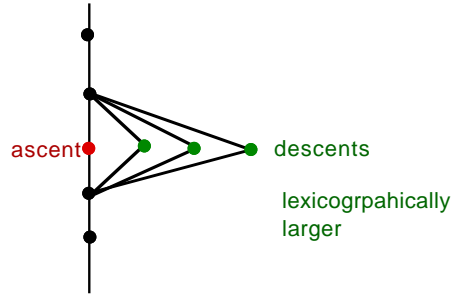}
\caption{Straightening step}\label{figdecrease}
\end{center}
\end{figure}

The first part of Theorem~\ref{lshellhom} is a consequence of
Theorems~\ref{wedge} and~\ref{ELshell}.  Indeed, the 
proper parts of the decreasing chains are the homology facets of the shelling of
$\Delta(P)$ induced by the lexicographic order of maximal chains of $\hat P$.  To 
establish the
second part, one needs only show that the proper parts of the decreasing chains span
cohomology.  This is proved by showing that the cohomology relations enable one to
express a maximal chain with an ascent as the negative of a  sum of lexicographically
larger maximal chains.  This provides a step in a ``straightening'' algorithm for
expressing maximum chains as linear combinations of decreasing maximum chains.  See 
Figure~\ref{figdecrease}.

\begin{remark}  In Bj\"orner's original version of EL-shellability, the unique increasing
maximal chain was required  to be  weakly increasing and the decreasing chains  were
required to be strictly decreasing. The two versions have the same topological and
algebraic consequences, but it is unknown whether they are equivalent.
\end{remark}

\begin{example}  In the EL-labeling given in Figure~\ref{figel}, the two rightmost maximal
chains are the only decreasing maximal chains.  One has length 3 and the other has length
2. It follows form Theorem~\ref{lshellhom} that the order complex of the proper part of the
poset has the homotopy type of a wedge of a $1$-sphere and a $0$-sphere.  This
is consistent with the fact that the order complex of the proper part of the poset
consists of the barycentric subdivision of the boundary of a $2$-simplex 
and an isolated point.
\end{example} 

\subsection{The Boolean algebra  $B_n$.} \label{ELboolex}  
\begin{figure}\begin{center}
\includegraphics[width=4cm]{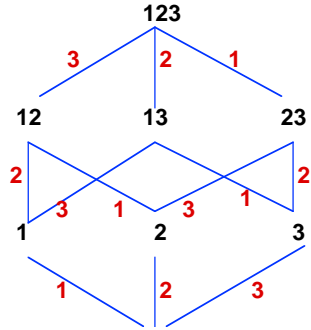}
\caption{EL-labeling of $B_3$} \label{figelbool}
\end{center}
\end{figure}

There is a very natural
EL-labeling of the Boolean algebra $B_n$; simply label the covering relation $A_1
\subset\!\!\!\!\cdot\,\, A_2$ with the  unique element in the singleton set $A_2-A_1$.  The
maximal chains correspond to the permutations in $\mathfrak S_n$.  It is easy to see that
each interval has a unique increasing chain that is lexicographically first.  There is
only one decreasing chain which is consistent with the fact that $\Delta(\bar B_n)$ is a
sphere.  See Figure~\ref{figelbool}.

For each $k \le n$, define the truncated Boolean algebra
$B_n^k$ to be the subposet of
$B_n$ given by
$$B_n^k = \{A \subseteq [n] : |A| \ge k \}.$$  Define an edge labeling $\lambda$
of
$B_n^k \cup \{\hat 0\}$ as follows: 
$$\lambda(A_1, A_2) =\begin{cases} \max A_2 & \mbox{ if }A_1 = \hat 0 \mbox{ and }|A_2| = k
\\
 a  &\mbox{ if }
A_2-A_1 = \{a\}. \end{cases}$$
  It is easy to check that this is an EL-labeling.  The
decreasing chains correspond to permutations with descent set
$\{k,k+1,\dots,n-1\}$.   Hence $\dim
\tilde H^{n-k-1}(\bar B_n^k) $ equals the number of permutations in $\mathfrak S_n$ with
descent set
$\{k,\dots,n-1\}$, which  equals the number of standard Young tableaux of hook
shape
$k 1^{n-k}$.  An equivariant version of
this result is given by the following special case of a result of Solomon (the general
result is given in Section~\ref{ranksec}).

\begth[Solomon \cite{sol68}] \label{solth} For all $k \le n$,
$$\tilde H_{n-k-1}(\bar B_n^k) \cong_{\mathfrak S_n} S^{(k1^{n-k})}.$$
\enth

\begin{proof}[Combinatorial proof idea]  We use a surjection from the set of  tableaux of
shape
$k1^{n-k}
$ to the set of  maximal chains  of $\bar B_n^k$, illustrated by the following example in
which $n =5$ and $k=2$.
$$\tableau[scY]{3,1|2|4|5} \mapsto ( \{3,1\}  <\!\!\!\cdot \, \{3,1,2\} <\!\!\!\cdot \,
\{3,1,2,4\} ).$$
The surjection  determines a surjective  $\mathfrak S_n$-homomorphism from the tableaux
module
$M^\lambda$ to the chain space $C_{n-k-2}(\bar B_n^k;\C)$.  To show that this homomorphism
induces a surjective homomorphism from the quotient space $S^\lambda$ to the quotient space
$\tilde H^{n-k-2}(\bar B_n^k) = C_{n-k-2}(\bar B_n^k;\C)/B^{n-k-2}(\bar B_n^k;\C)$, we
observe that the row relations map to
$0$ in
$C_{n-k-2}(\bar B_n^k;\C)$ and  show that   the  Garnir relations map to coboundary
relations.  Since the dimensions of the two vector spaces are equal, the homomorphism is
an isomorphism.    We demonstrate
 the fact that Garnir relations map to cohomology relations on two examples and leave the
general proof as an exercise. The Garnir relation
$$\tableau[scY]{3,1| 2|4|} + \tableau[scY]{3,1| 4|2|} $$

\vspace{.1in}
 \noindent maps to the sum of maximal chains of the proper part of the poset

\vspace{.1in}\begin{center}
\includegraphics[width=2cm]{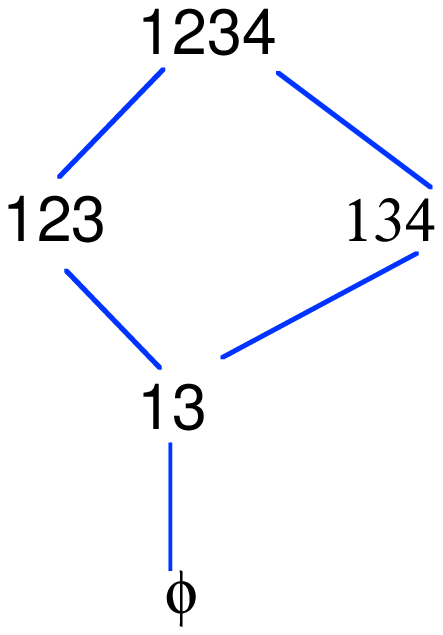}
\end{center}

\noindent The Garnir relation
$$\tableau[scY]{ 3,1|2| 4|} + \tableau[scY]{ 2,1|3| 4|} +\tableau[scY]{
2,3|1| 4|}
 $$ 

\vspace{.1in} 
\noindent maps to the sum of maximal chains of the proper part of the poset
\vspace{.1in}\begin{center}
\includegraphics[width=2cm]{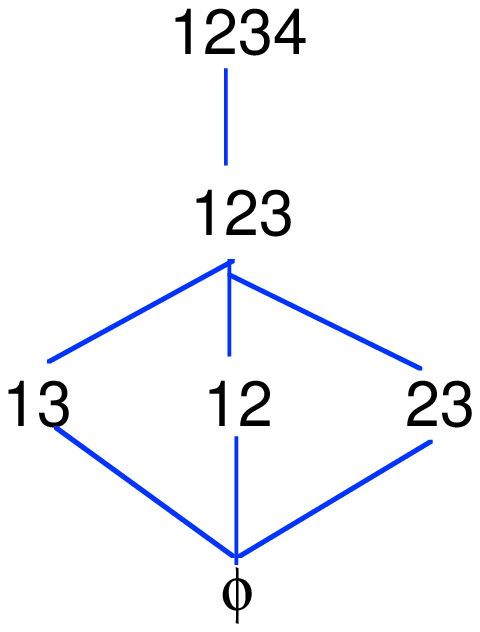}
\end{center}

\end{proof}

\subsection{The  partition lattice $\Pi_n$.} \label{ELpart}   We
give two different EL-labelings of the partition lattice.  The first
one, due to Gessel, appears in
\cite{bj80}, and the second one is due to Wachs \cite{wa96}. If
$x <\!\!\!
\cdot\, y$ in
$\Pi_n$ then
$y$ is  obtained from
$x$ by merging  two blocks, say $B_1$ and $B_2$.  
For the first edge labeling $\lambda_1$, let
 $$\lambda_1(x,y) = \max \{\min B_1,
\min B_2\}$$  and for the second edge labeling $\lambda_2$, let $$\lambda_2(x,y)= \max B_1
\cup B_2.$$ 
The increasing  chain from $\hat 0$ to $\hat 1$ is the same  for both labelings; it  consists of
partitions with only one nonsingleton block. More precisely, the chain is given by 
$$ \hat 0 <\!\!\! \cdot\, \{1,2\} <\!\!\! \cdot\, \{1,2,3\} <\!\!\! \cdot\, \cdots <\!\!\!
\cdot\, \hat 1 ,$$  where we have written only the nonsingleton block of each partition in
the chain.  We leave it as an exercise to show that these labelings are EL-labelings of
the partition lattice.

The decreasing  maximal chains for $\lambda_1$ and $\lambda_2$ are not the same.  We
describe those for  $\lambda_2$ first. They  also consist of partitions
with only one nonsingleton block and  are of the form 
$$c_{\sigma}:= (\hat 0 <\!\!\! \cdot\, \{\sigma(n),\sigma(n-1)\} <\!\!\! \cdot\,
\{\sigma(n),\sigma(n-1),\sigma(n-2)\} <\!\!\!
\cdot\,
\cdots$$ $$ <\!\!\!
\cdot\, \{\sigma(n),\sigma(n-1),\sigma(n-2),\dots, \sigma(1)\} =\hat 1), $$
where $\sigma \in \mathfrak S_n$ and $\sigma(n) = n$. 
We conclude that the homotopy type of $\bar \Pi_n$ is given by 
\begin{equation} \label{parhom} \bar \Pi_n \simeq \bigvee_{(n-1)!} \S^{n-3} \end{equation} 
and that the chains $\bar c_\sigma$, where
$\sigma \in
\mathfrak S_{n}$ and
$\sigma(n)= n$, form a basis for
$\tilde H^{n-3}(\bar \Pi_n, \Z)$.  From this basis one can immediately see that
\begin{equation}\tilde H_{n-3}(\Pi_n) \downarrow^{\s_n}_{\s_{n-1}} \,\, \cong_{\s_{n-1}}
\,\,  \tilde H^{n-3}(\Pi_n)
\downarrow^{\s_n}_{\s_{n-1}}\,\, \cong_{\s_{n-1}}
\,\, \C
\s_{n-1},
\end{equation} 
a result obtained by
Stanley
\cite{st82} as a consequence of his formula for the full representation of $\s_n$ on
$\tilde H_{n-3}(\bar \Pi_n)$ (Theorem~\ref{thlie}).

\begin{figure} \begin{center}
\includegraphics[width=6 cm]{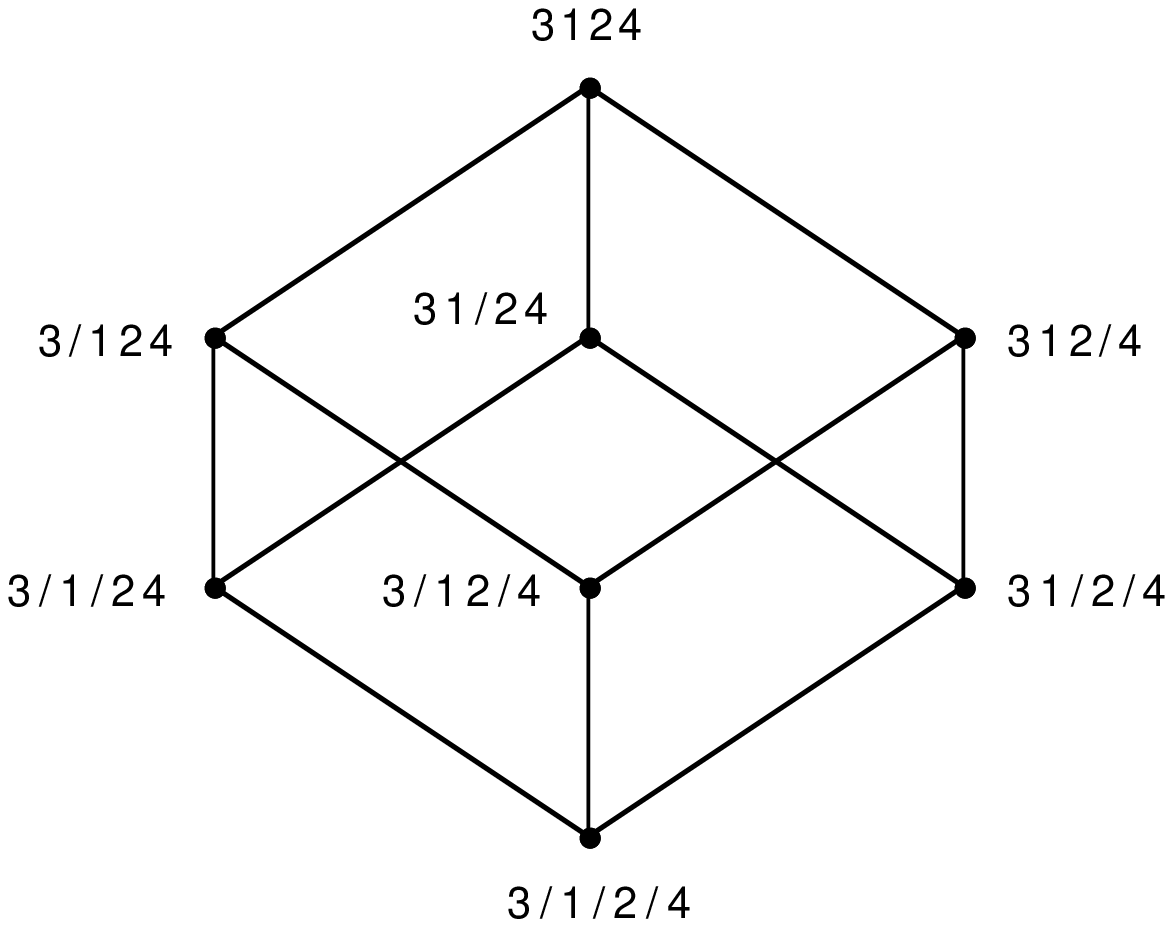}\\
\end{center}
\begin{center}\caption{ $\Pi_{3124}$ }\label{figposet1}
\end{center}
\end{figure}

Next we describe a nice basis for homology of $\bar\Pi_n$ that is
dual to the decreasing chain basis for $\lambda_2$.    To {\em split} a permutation $\sigma
\in
\s_n$ at positions
$j_1 < j_2 <\dots < j_k$ in $[n-1]$ is to form  the partition $$\sigma(1),
\sigma(2),
\dots,
\sigma(j_1)\,/ \, \sigma(j_1+1), \sigma(j_1+ 2), \dots,
\sigma(j_2)\,/  \dots / \,\sigma(j_{k}+1), \sigma(j_k+ 2), \dots,
\sigma(n)$$ of $[n]$.  For each
$\sigma \in
\mathfrak S_n$, let
$\Pi_\sigma$ be the induced subposet of the partition lattice $\Pi_n$ consisting of
partitions obtained by splitting the permutation $\sigma$ at any set of positions in
$[n-1]$  .   The subposet $\Pi_{3124}$ of $\Pi_4$ is shown in  Figure~\ref{figposet1}. 
Each poset
$\Pi_\sigma$ is isomorphic to the subset lattice $B_{n-1}$.  Therefore $\Delta(\bar
\Pi_\sigma)$ is an $(n-3)$-sphere 
embedded
in 
$\Delta(\bar
\Pi_n)$, and hence it determines a fundamental cycle $\rho_\sigma \in \tilde H_{n-3}(\bar
\Pi_n;\Z)$. 

\begth[Wachs \cite{wa96}] \label{splitperm} The set $\{\rho_\sigma : \sigma \in \mathfrak S_n,
\sigma(n) = n\}$ forms a basis for 
$\tilde H_{n-3}(\bar
\Pi_n;\Z)$ dual to the decreasing chain basis $\{\bar c_\sigma : \sigma \in \mathfrak S_n,
\sigma(n) = n\}$ for cohomology.  Call the homology basis, the
{\em splitting basis}.
\enth

\begin{figure}\begin{center}
\includegraphics[width=7 cm]{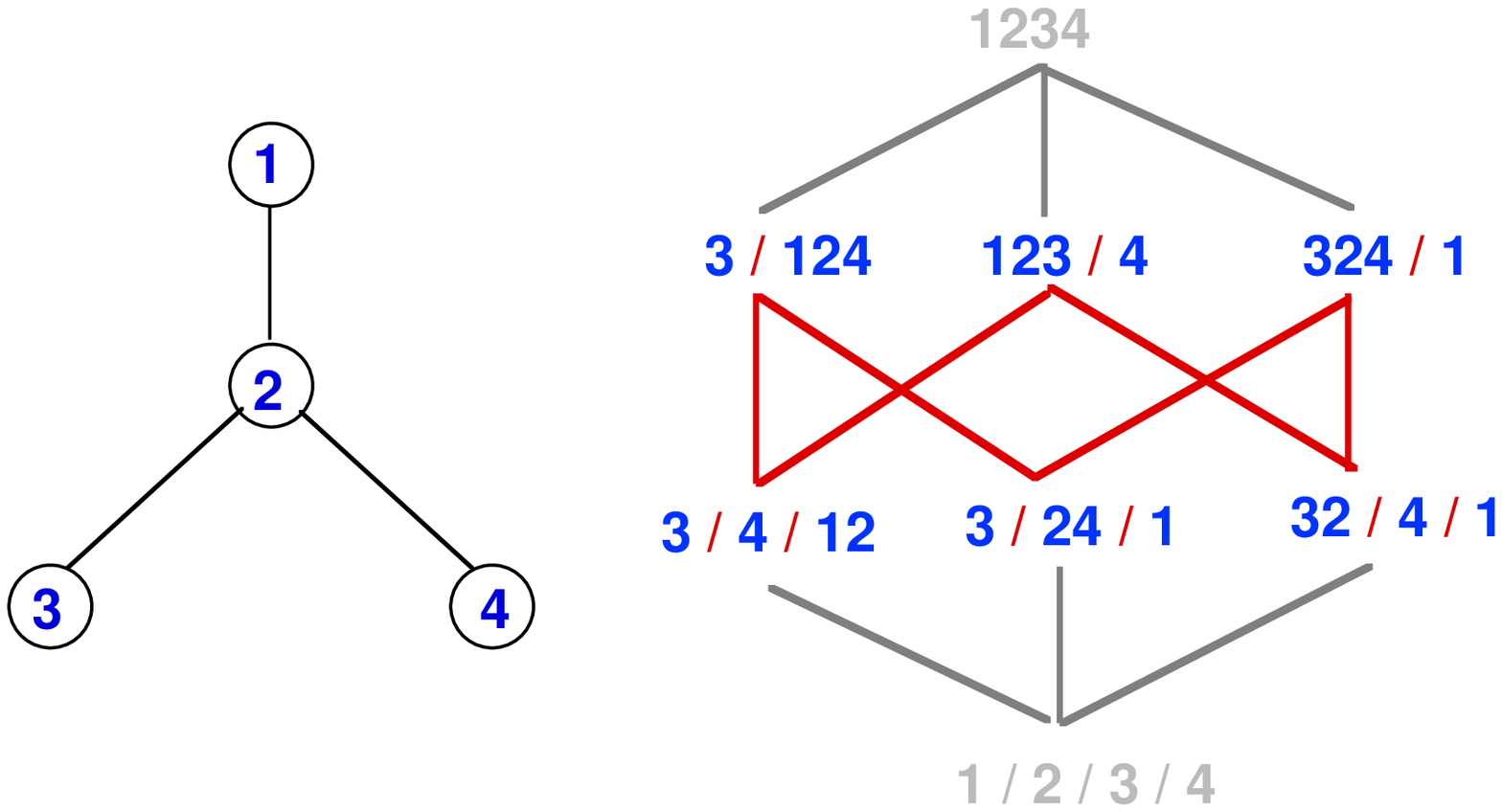}\\
\end{center}
\begin{center}\caption{ $T$ and $\Pi_T$ }\label{nbcfig}
\end{center}
\end{figure}

Now we describe the decreasing chain basis for cohomology for  the  EL-labeling
$\lambda_1$ and its dual basis for homology.  Given any rooted
nonplanar (i.e. children of a node are unordered) tree $T$ on node set $[n]$, by removing
any  set of edges of
$T$, one forms a partition of $[n]$ whose blocks are the node sets of the connected
components of the resulting graph.   Let $\Pi_T$ be the induced subposet of the partition
lattice
$\Pi_n$ consisting of partitions obtained by removing edges of  $T$; see
Figure~\ref{nbcfig}.  (If
$T$ is a linear tree then $\Pi_T$ is the same as $\Pi_\sigma$, where $\sigma$ is the
permutation obtained by reading the nodes of the tree from the root down.)   Each poset
$\Pi_T$ is isomorphic to the subset lattice $B_{n-1}$.  We let $\rho_T$ be the fundamental cycle
of the spherical complex $\Delta(\bar \Pi_T)$.     Let $T$ be an increasing tree on node set
$[n]$, i.e., a rooted nonplanar tree on node set $[n]$ in which each node $i$ is greater
than its parent $p(i)$.  We form the chain $ c_T$ in $\Pi_T$, from top down,  by removing
the edges
$\{i,p(i)\}$, one at a time, in increasing order of $i$. For the increasing tree $T$ in
Figure~\ref{nbcfig}, $$c_T = (1/2/3/4 \,\,<\!\!\!\cdot \,\,\, 3/24/1 \,\, <\!\!\!\cdot
\,\,\, 324/1 \,\,<\!\!\!\cdot\,\, 1234).$$ We claim that the
$c_T$, where
$T$ is an increasing tree on node set $[n]$, are the decreasing chains of $\lambda_1$. 

\begth \label{splittree}  Let $\mathcal T_n$ be the set of increasing trees on node set $[n]$. The
set
$\{\rho_T : T
\in
\mathcal T_n\}$ forms a basis for 
$\tilde H_{n-3}(\bar
\Pi_n;\Z)$ dual to the decreasing chain basis $\{\bar c_T: T
\in
\mathcal T_n\}$ for cohomology.  Call the homology basis, the
{\em tree splitting basis}.
\enth

\begin{xca}
\begin{enumerate} \item[]
\item[(a)] Show $\lambda_1$ is an EL-labeling whose decreasing maximal chains are of the
form $c_T$, where
$T$ is an increasing tree.
\item[(b)] Prove Theorem~\ref{splittree}
\item[(c)] Show $\lambda_2$ is an EL-labeling whose decreasing maximal chains are of
the form $c_\sigma$,
 where
$\sigma \in \s_n$ is such that $\sigma(n) = n$.
\item[(d)] Prove Theorem~\ref{splitperm}. 
\end{enumerate}
\end{xca}

Both of the  EL-labelings and their corresponding bases are special cases of more general
constructions.  The first EL-labeling
$\lambda_1$ and its bases are specializations of constructions for geometric
lattices due to
Bj\"orner \cite{bj80,bj82}; this is discussed in the next section.  The second EL-labeling
$\lambda_2$ generalizes to all Dowling lattices; see the next exercise.   A geometric
interpretation of the splitting basis,
 in which the fundamental cycles correspond to
bounded regions in an affine slice of the real braid arrangement, is given by Bj\"orner
and Wachs
\cite{bw1}.  This leads to analogs of the splitting basis for  intersection lattices of
other Coxeter arrangements, in particular the type B  partition lattice
$\Pi_n^B$.  An analog of the splitting basis for all Dowling lattices is given 
by Gottlieb and Wachs \cite{gw}.

\begin{xca}[Gottlieb {\cite[Section 7.3]{got98}}] \label{bnel}  
\begin{enumerate} \item[]
\item[(a)] Find an EL-labeling of the
type
$B$ partition lattice $\Pi_n^B$ analogous to the EL-labeling $\lambda_2$ of $\Pi_n$.
\item [(b)]  
Describe the decreasing chains in terms of signed permutations in $
\s_n[\Z_2]$ whose right-to-left maxima are positive, where for $i = 1,2,\dots n$, we
say that   $\sigma(i)$ is a right-to-left maxima of $\sigma$ if
$|\sigma(i)| > |\sigma(j)|$ for all
$j=i\dots n$. 
\item [(c)] Show that the number of decreasing chains is
$(2n-1)!! := 1\cdot 3 \cdots (2n-1)$, thereby recovering the well-known fact that $\bar
\Pi_n^B$ has the homotopy type of a wedge of $(2n-1)!!$-spheres of dimension $n-2$.
\end{enumerate}
\end{xca}

A partition  $\pi \in \Pi_n$ is said to be
{\em noncrossing} if for all $a<b<c<d$, whenever $a,c$ are in a block $B$ of $\pi$ and
$b,d$ are in a block $B^\prime$ of $\pi$ then $B = B^\prime$.
  Let 
$\mbox{NC}_n$ be the induced subposet of $\Pi_n$ consisting of noncrossing
partitions.   This poset, known as the {\em noncrossing partition lattice}, was
first introduced by Kreweras \cite{kr72} who showed  that it is a graded
lattice with M\"obius invariant equal to  the signed Catalan number $(-1)^{n-1}{1 \over
n} {2n-2 \choose n-1}$.  It has since undergone extensive study  due to its many
fascinating properties and its connections to various diverse fields of mathematics such
as combinatorics, discrete geometry, mathematical biology,  geometric group theory, low
dimensional topology,  and free probability; see the 
survey papers of Simion
\cite{sim00} and McCammond
\cite{mc05}.

 Bj\"orner and Edelman
 showed that  $\mbox{NC}_n$ is EL-shellable (see 
Exercise~\ref{crossELex}).  Reiner
\cite{r97} introduced   an analog of the noncrossing partition lattice for the type B
Coxeter group, and established many of the same properties held by $\mbox{NC}_n$, including
EL-shellability.  A generalization of the types A and B noncrossing partition
lattices to  all finite Coxeter groups was recently introduced in the geometric group
theory /low dimensional topology  work of Bessis
\cite{b03} and Brady and Watt
\cite{braw}.  This generalization is discussed in Section~\ref{CLsec}.   The
M\"obius function for the type D noncrossing partition lattice was computed by
Athanasiadis and Reiner and the question of EL-shellability was raised
\cite{ar04}.
  The general Coxeter group noncrossing partition lattices were subsequently shown to be
EL-shellable by Athanasiadis, Brady and Watt \cite{abw05}. 

\begin{xca}[Bj\"orner and Edelman, cf. \cite{bj80}]  \label{crossELex}
\begin{enumerate}
\item[] \item[(a)] Show that  the restriction to $\mbox{NC}_n$ of one of the EL-labelings
of $\Pi_n$ described above is an EL-labeling.
\item[(b)] Describe the  decreasing chains for the EL-labeling and show that their number
is the Catalan number 
${1 \over n} {2n-2 \choose n-1}$.
\end{enumerate}
\end{xca}

\begin{xca}[Stanley \cite{st97}] Define an edge labeling for $\mbox{NC}_n$ by
$$\lambda(x,y) = \max\{i \in B_1 : i < \min B_2\},$$
where 
  $B_1$ and $B_2$ are the blocks of $x$ that are merged to obtain $y$ 
 and $\min
B_1 <
\min B_2$.

\begin{enumerate}
\item[(a)] Show that $\lambda$ is an EL-labeling of $\mbox{NC}_n$.
\item[(b)] Show that     the decreasing chains correspond bijectively, via the labeling,
to the
 sequences of the form $(a_1\ge a_2 \ge \dots \ge a_{n-1})$ where $ a_{n-i} \in [i]
$ for all $i
\in [n-1]$.  Show that  the number of such sequences is a  Catalan
number.
\end{enumerate}
This labeling is particularly interesting because the maximal chains correspond
bijectively, via this labeling, to a well-studied class of sequences called parking
functions; see
\cite{st97}.
\end{xca}

\subsection{ Geometric lattices} \label{geosec} A {\em geometric lattice} is
a lattice
$L$ that is
 semimodular (which means that for all $x,y \in L$, the
join
$x
\lor y$ covers $y$ whenever $x$ covers the meet $x
\land y$) and atomic (which means every element of $L$ is the join of atoms).  Some
fundamental examples are the subset lattice $B_n$, the subspace lattice $B_n(q)$, the
partition lattice $\Pi_n$, the type $B$ partition lattice $\Pi_n^B$, and more generally,
the intersection lattice
$L(\mathcal A)$ of any central hyperplane arrangement.   Geometric lattices are
fundamental structures of matroid theory.    It is easy to show that geometric lattices
are pure.

We describe  an edge labeling for  geometric lattices, which comes from Stanley's
early  work \cite{st74} on rank-selected M\"obius invariants  and was one of the main
motivating examples for Bj\"orner's original work \cite{bj80} on EL-shellability.   
Fix an ordering
$a_1,a_2,\dots, a_k$ of the atoms of the geometric lattice $L$.  Then label each edge $x
<\!\!\!
\cdot\, y$ of the Hasse diagram with the smallest $i$ for which $x \lor a_i = y$.         Note
that if the atoms of the subset lattice
$B_n$ are ordered
$\{1\},\{2\}, \dots,
\{n\}$, then the geometric lattice EL-labeling is precisely the labeling given in
Section~\ref{ELboolex}.

\begin{xca}\label{nbc} \begin{enumerate} 
\item[]
\item[(a)] Show that the edge labeling for geometric lattices described above is an EL-labeling.
\item[(b)] Find an ordering of the atoms of $\Pi_n$ such that the induced geometric
lattice EL-labeling is equivalent to the  EL-labeling
$\lambda_1$ of Section~\ref{ELpart}. 
\item[(c)] Is $\lambda_2$ of Section~\ref{ELpart} equivalent to a    geometric lattice
EL-labeling for some ordering of the atoms?
\item[(d)] Show that every semimodular lattice has an EL-labeling.
 \end{enumerate}
\end{xca}

The decreasing chains of the geometric lattice EL-labeling have a very nice characterization, due
to Bj\"orner \cite{bj82}, which is described in the language of matroid theory.    A set
$A$ of atoms in a geometric lattice $L$ is said to be independent if $r(\lor A) = |A|$.  A set of
atoms that is minimally dependent (i.e., every proper subset is independent) is
called a {\em circuit}.  An independent set of atoms 
  that  can be obtained from a circuit by removing 
its smallest element (with respect to the  fixed  ordering
$a_1,a_2,\dots, a_k$ of the atoms of $L$) is called a {\em broken circuit}.  A maximal  
independent set of atoms is said to be an NBC base if   contains no broken circuits.  There is a
natural bijection from the NBC bases of $L$ to the decreasing chains of $L$.  Indeed,  
given any NBC base
$A=\{a_{i_1},\dots,a_{i_r}\}$, where $1 \le i_1 < i_2 < \dots < i_r \le k$, construct the maximal
chain 
$$c_A:= (\hat 0 < a_{i_r} < a_{i_r} \lor a_{i_{r-1}} < \dots < a_{i_r}
\lor a_{i_{r-1}} \lor \dots \lor a_{i_1} = \hat 1).$$ 
It is not difficult  to check that the label sequence of $c_A$ is  $(i_r,i_{r-1},\dots, i_1)$,
which is decreasing, and that the map $A \mapsto c_A$ is a bijection from the NBC bases of $L$ to
the decreasing chains of $L$. 
We conclude that the set 
$\{\bar c_A : A \mbox{ is an NBC base of } L\}$ 
is a basis for top cohomology of $L$.  

Bj\"orner \cite{bj82} also constructs a  basis for homology of $L$ indexed by the NBC
bases, which is dual to the cohomology basis described above. Any independent set of atoms
in a geometric lattice generates (by taking joins) a Boolean algebra  embedded in the
geometric lattice.  Given any  independent  set of atoms
$A$, let
$L_A$ be the Boolean algebra generated by
$A$ and let  $\rho_A$ be the fundamental cycle of $\bar L_A$.  

\begth[Bj\"orner \cite{bj82}] \label{geombases} Fix an ordering
of the set of atoms of a geometric lattice $L$.  The set
$$\{
\rho_A : A \mbox{ is an NBC base of } L\}$$ is a basis for top homology of $\bar L$, which
is dual to the decreasing chain basis $$\{\bar c_A : A \mbox{ is an NBC base of } L\}$$ 
 for top cohomology.  
\enth

\begin{xca} Prove Theorem~\ref{geombases}.
\end{xca}

\begin{xca} Show that the tree splitting basis for homology of the partition lattice is an
example of a Bj\"orner NBC basis.  
\end{xca}

For further reading on the homology of geometric lattices see Bj\"orner's book
chapter \cite{bj921}.

\subsection{The $k$-equal partition lattice} \label{kesec} In this section we  present the
original example that motivated  the Bj\"orner-Wachs extension of the  theory of 
shellability from pure to nonpure complexes.   This interesting example  played a pivotal
role in two other important developments, namely the use of poset topology in complexity
theory (Bj\"orner, Lov\'asz and Yao \cite{bly, bl}) and the derivation  of the equivariant
version of the Goresky-MacPherson formula (Sundaram and Welker \cite{suwe}).

\begin{figure}\begin{center}
\includegraphics[width=9cm]{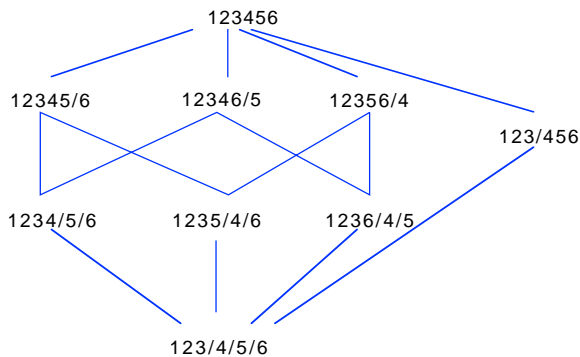}
\end{center}
\begin{center} \caption{ Interval of $\Pi_{6,3}$}\label{figkequal}
\end{center}
\end{figure}

The $k$-equal partition lattice
$\Pi_{n,k}$ is the subposet of $\Pi_n$ consisting of partitions that have no blocks of size
$\{2,3,\dots, k-1\}$.  This lattice is not a geometric lattice when $k >2$. It's not even
pure; consider the interval
 of $\Pi_{6,3}$ shown in Figure~\ref{figkequal}.

The $k$-equal partition lattice arose in the early 1990's from the ``$k$-equal problem'' in
complexity theory: Given a sequence of $n$ real numbers $x_1,\dots,x_n$, determine whether
or not some
$k$ of them are equal.  Bj\"orner, Lov\'asz and Yao \cite{bly,bl} obtained the following
lower bound on the complexity of the $k$-equal problem by studying the topology of $\bar
\Pi_{n,k}$,
$$ \max\{n-1,n\log_3{n \over 3k}\}.$$ This lower bound differs by a factor of only 16 from
the best upper bound.
The
$k$-equal problem translates into a  subspace arrangement problem; namely determine
whether or not the point
$(x_1,\dots, x_n)$ is in the complement of the real  subspace arrangement $\mathcal
A_{n,k}$, called the $k$-equal arrangement, consisting of subspaces of the form 
$$\{(x_1,\dots,x_n) \in \R^n : x_{i_1} = x_{i_2} = \dots =x_{i_k}
\},$$  where $1 \le i_1 < \cdots < i_k \le n$.   Bj\"orner and Lov\'asz \cite{bl}
show that the sum of the Betti numbers of the complement of
the
$k$-equal arrangement gives a lower bound for the depth of the best decision tree for the
$k$-equal problem.  They then use the Goresky-MacPherson formula
(Theorem~\ref{GMform}) to reduce the complexity problem to that of  studying the topology
of lower intervals in
 the intersection lattice of $\mathcal A_{n,k}$, which is isomorphic to $\Pi_{n,k}$.

Bj\"orner and Welker \cite{bjwe95} prove that all open intervals of $\Pi_{n,k}$ have the
homotopy type of a wedge of spheres of varying dimensions.  Since homotopy type of a wedge
of spheres (of top dimension)  is  the main topological consequence of pure shellability,
 this result led Bj\"orner and Wachs to  consider 
shellability  for nonpure complexes. 

We now present the nonpure EL-labeling of
$\Pi_{n,k}$ obtained in \cite{bw96} .
First linearly order the label set $\{\bar 1< \bar 2 < \cdots < \bar n < 1 < 2 <
\cdots < n\}$.  Now label the edge $\pi \,<\!\!\!\cdot \,\,\tau$ as follows:

$$ \lambda(\pi,\tau) = \begin{cases}\max B \quad &\mbox{if a new block $B$ is formed from
singleton blocks} \\ a \quad &\mbox{if a nonsingleton block is merged with a
singleton
$\{a\}$}\\ \overline{\max B_1\cup B_2} \quad &\mbox{if two nonsingletons $ B_1$ and $B_2$
are merged.}
\end{cases}
$$

\begin{xca}
 Show that  $\lambda$ is an EL-labeling of $\Pi_{n,k}$.
\end{xca}

We next describe the decreasing maximal chains.  Given a hook
shape Young diagram  $k 1^m$, by the {\em corner} of the hook we mean
the cell in the first row and first column, and by
the leg of
$T$ we mean the first column of $T$ minus the corner.  We  refer to a skew diagram as a
$k$-broken skew hook diagram if it is the disjoint union of  hooks $H_1, H_2, \dots, H_t$
of the form $k1^m$, where $m$ is arbitrary.  The  tableaux in Figure~\ref{figtab} are
reverse standard tableaux (i.e. column and row entries are decreasing rather than
increasing) of 
$3$-broken skew hook shape.  For each reverse standard  tableau
$T$  of
$k$-broken skew hook shape, with entries in $[n]$, let  $\pi_{_T}$ be the
partition in
$\Pi_{n,k}$ whose nonsingleton blocks are  
$B_1,
\dots, B_t$, where
$B_i$ is the set of  
 entries of $T$ in the hook
$H_i$.  For the tableaux in Figure~\ref{figtab}, 
$$\pi_{_T} = 15,14,3,9,8,5/ 10,7,1,2/ 13,11,4,12,6$$ and
$$\pi_{_{T^\prime}} = 15,14,3,9,8,5/ 10,7,1/ 13,11,4,12,6/2.$$ 
Given a  reverse standard  tableau
$T$  of
$k$-broken skew hook shape, we  
let $a$ be the smallest entry  among all entries of $T$ that are either in the  leg of a
hook
 or in the corner of a hook that has no leg.
In the former case we let 
$T^\prime$ be the standard  tableau of
$k$-broken skew hook shape obtained by removing $a$ from $T$.   In the latter case 
we let
$T^\prime$ be the standard tableau of
$k$-broken skew hook shape obtained by removing the entire hook (which is a row) containing
$a$.   In either case, we  refer to the tableau
$T^\prime$ as the {\em predecessor} of
$T$.  In Figure~\ref{figtab}, 
$T^\prime$ is the
 predecessor of $T$. For each reverse standard  tableau
$T$  of
$k$-broken skew hook shape, with $n$ cells and entries in $[n]$, let
$c_T$ be the maximal chain of $\Pi_{n,k}$ such that
\begin{itemize}
\item $\pi_{_T} \in c_T$
\item   the upper segment $c_T\cap [\pi_{_T},\hat 1]$ is the chain
$$\pi_{_T} = B_1/B_2/\dots/B_t < B_1 \cup B_2 /B_3/\dots/B_t < \dots < B_1 \cup \dots \cup
B_t = \hat 1$$
\item the lower segment $c_T\cap [\hat 0, \pi_{_T}]$ is the chain
$$\hat 0 = \pi_{_{T_0}} < \pi_{_{T_1}} < \dots < \pi_{_{T_{n-t(k-1)}}} = \pi_{_T},$$
where $T_i$ is the predecessor  of $T_{i+1}$.
\end{itemize}
For the tableau $T$ of Figure~\ref{figtab}, the maximal chain $c_T$ is given by

\begin{eqnarray*}\hat 1 
&>&  15,14,3,9,8,5, 10,7,1,2/13,11,4,12,6 
\\ &>&  15,14,3,9,8,5/ 10,7,1,2/ 13,11,4,12,6 
\\ &>&  15,14,3,9,8,5/ 10,7,1/13,11,4,12,6
\\&>&  15,14,3,9,8/ 10,7,1/ 13,11,4,12,6
\\ &>&  15,14,3,9,8/ 10,7,1/13,11,4,12
\\ &>& 15,14,3,9/ 10,7,1/13,11,4,12
\\ &>&  15,14,3/ 10,7,1/ 13,11,4,12
\\ &>&  15,14,3/  13,11,4,12
\\ &>& 15,14,3/ 13,11,4
\\ &>& 15,14,3
\\ &>& \hat 0
\end{eqnarray*}
where only the nonsingleton blocks are shown.

\begin{figure}\begin{center}\hspace{.4in} $T=\tableau[scY]{\White
,,,,,,\Black 13,11,4|\White ,,,,,,\Black 12|\White ,,,,,,\Black 6|\White ,,,\Black
10,7,1|\White ,,,\Black 2| 15,\Black 14,3
|9|8|5}$ \hspace{.4in}
 $T^\prime= \tableau[scY]{\White ,,,,,,\Black 13,11,4|\White ,,,,,,\Black 12|\White
,,,,,,\Black 6|\White ,,,\Black 10,7,1| 15,\Black 14,3
|9|8|5}$
\caption{\,} \label{figtab} \end{center}\end{figure}

\begin{xca}\label{exkequal} \begin{enumerate} 
\item[] 
\item[(a)] \cite{bw96} Show that the set of decreasing chains with respect to
$\lambda$ is the set
$$\{c_T : T \mbox{ is a reverse standard tableau of $k$-broken skew hook shape, }$$ 
\vspace{-.3in}\begin{center} with $n$ cells
and $n$ in the  corner of the leftmost  hook$\}$.\end{center}
\item[(b)] \cite{bw96} Show that the Betti numbers are given by 
 $$
\tilde\beta_{n-3-t(k-2)}(\bar \Pi_{n,k})= \sum_{\scriptsize\begin{array}{c}
j_1+j_2+\ldots+j_t = n\\  j_i
\ge k\end{array}}
 \binom {n-1} {j_1-1, j_2, \ldots, j_t} \prod_{i=1}^t 
\binom {j_i-1} {k-1}
$$ where $t \ge 1$,
and 
$$\tilde\beta_{i}(\bar \Pi_{n,k})= 0$$ if $i$ is not of the form
$n-3-t(k-2)$ for any $t \ge 0$.
\item[(c)] \cite{suwa} Prove  that 
  $$\tilde  H_{n-3-t(k-2)}(\bar \Pi_{n,k})\downarrow_{\mathfrak S_{n-1}}^{\mathfrak S_n}
\cong_{\s_{n-1}} 
\bigoplus_D S^D $$
summed over  all $k$-broken skew hook diagrams $D$ of size $n$ in which the corner of the
leftmost  hook is removed.
This is equivalent to
$$\tilde  H_{n-3-t(k-2)}(\bar \Pi_{n,k})\downarrow_{\mathfrak S_{n-1}}^{\mathfrak S_n}
\cong \bigoplus_{\scriptsize\begin{array}{c} j_1+j_2+\ldots+j_t = n\\  j_i
\ge k\end{array}}S^{1^{j_t-k}} \bullet S^{(k-1)} \bullet \prod_{i=1}^{t-1} S^{(k1^{j_i-k})}
,
$$ where $\bullet$ and $\prod$ denote induction product.
\end{enumerate}
\end{xca} 

The decreasing chains given in Exercise~\ref{exkequal}(a) were also used by Sundaram and
Wachs \cite{suwa} to obtain a  formula for the  unrestricted representation of $\mathfrak
S_n$ on
$\tilde  H_{n-3-t(k-2)}(\bar
\Pi_{n,k})$,
involving   composition product of representations.  In order to transfer the
representation of the symmetric group on homology of the $k$-equal partition poset to 
representations  of the symmetric group on the homology of the complement $M_{\mathcal
A_{n,k}}:=
\R^n -
\cup
\mathcal A_{n,k}$ of the real
$k$-equal arrangement $\mathcal A_{n,k}$ and the complement $M_{\mathcal A^\C_{n,k}}:=
\C^n -
\cup
\mathcal A^\C_{n,k}$ of the complex $k$-equal arrangement $\mathcal A^{\C}_{n,k}$, 
Sundaram and Welker derived an  equivariant version of the
Goresky-MacPherson  formula; see Theorem~\ref{swthm}.  By
computing the multiplicity of the trivial representation in the homology of the
complement, they   obtain Betti numbers for the orbit spaces $M_{\mathcal
A_{n,k}}/\mathfrak S_n$ and $M^\C_{\mathcal
A_{n,k}}/\mathfrak S_n$.  For instance, they recover
the following result of Arnol'd \cite{a}: 
\begin{equation} \label{arn} \tilde\beta_i(M^\C_{\mathcal A_{n,k}}/\mathfrak
S_n) =
\begin{cases} 1 &\mbox{if
$i =2k-3$} \\
0& \mbox{otherwise}.\end{cases}
\end{equation}  The orbit space $M^\C_{\mathcal A_{n,k}}/\mathfrak
S_n$ is  homeomorphic to the space of monic polynomials of degree $n$ 
whose roots have multiplicity at most
$k-1$.

Type $B$ and $D$ analogs of the $k$-equal arrangement and the $k$-equal partition lattice
were studied by Bj\"orner and Sagan \cite{bs96}.   Gottlieb \cite{got03} extended this
study to Dowling lattices.

\section{CL-shellability and Coxeter groups} \label{CLsec}

We now consider a more powerful version of lexicographic shellability called 
 chain-lexicographic shellability (CL-shellability for short).  This tool was
introduced by Bj\"orner and Wachs in the early 1980's in order to establish shellability
of Bruhat order on a Coxeter group \cite{bw82}. 

For a bounded poset $P$, let $\mathcal M\mathcal E(P)$ be the set of pairs
$(c,x<\!\!\!\cdot \, y)$ consisting of a maximal chain
$c$ and an edge
$x<\!\!\!\cdot \,  y$ along that chain.
 A {\em chain-edge labeling} of $P$ is a map
$\lambda:\mathcal M\mathcal E(P)\to\Lambda$, where $\Lambda$ is some poset, satisfying:
 If two maximal chains 
coincide along their bottom $d$ edges, then their labels
also coincide along these edges.

\begin{figure}\begin{center}
\includegraphics[width=4cm]{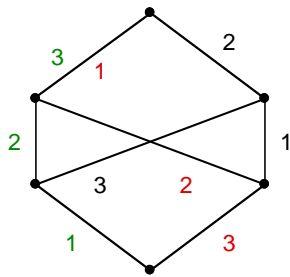}
\caption{\bf chain edge labeling} \label{figbruhat}
\end{center}
\end{figure}

 An example of a chain-edge labeling is given in Figure~\ref{figbruhat}.  Note that one of
the edges has two labels. This edge receives label 3 if it is paired with the leftmost
maximal chain and receives label 1 if it is paired with the other maximal chain that
contains the edge.

Just as  for edge labelings, we need to restrict chain-edge labelings $\lambda:\mathcal
M\mathcal E(P)\to\Lambda$ to intervals
$[x,y]$.  The complication is that since there can be more than one way to label an edge
depending upon which maximal chain it is paired with, there can  be more than
one way to restrict
$\lambda$ to
$\mathcal M\mathcal E([x,y])$.  It follows from the definition of chain-edge labeling
that each maximal chain
$r$ of
$[\hat 0,x]$ determines a unique restriction of $\lambda$ to $\mathcal M\mathcal
E([x,y])$.  
 This enables
one to talk about increasing and decreasing maximal chains and  lexicographic order
of maximal chains in the rooted interval $[x,y]_r$.

\begin{definition}  Let $P$ be a bounded poset. {\em  A chain-lexicographic
labeling} (CL-labeling, for short)  of
$P$ is a chain-edge labeling such that in each closed
rooted interval $[x,y]_r$ of $P$, there is a unique strictly increasing maximal chain,
which lexicographically precedes all other maximal chains of $[x,y]_r$.  A poset that
admits a CL-labeling is said to be CL-shellable.
\end{definition}

The chain-edge labeling of Figure~\ref{figbruhat} is a CL-labeling.  The unique increasing
maximal chain of the poset is the leftmost maximal chain.

It is easy to see that EL-shellability implies CL-shellability.  All the
consequences of EL-shellability discussed in Section~\ref{ELsec} are also consequences of
the more general CL-shellability.  It is unknown whether CL-shellability and
EL-shellability are equivalent notions.

We now present the original example that motivated this more technical version of
lexicographic shellability.

\begin{definition}
 A Coxeter system $(W,S)$ consists of a
a group $W$ together with a set of generators $S$ such that the following relations form a
presentation of $W$:
\begin{itemize} 
\item $s^2 = e$, for all $s \in S$
\item $(st)^{m_{s,t}} = e$, where $m_{s,t} \ge 2$, for certain $s \ne t \in S$.
\end{itemize}
The group $W$ is said to be a Coxeter group.
\end{definition}

Finite Coxeter groups can be characterized as finite reflection groups, i.e, finite groups
generated by linear reflections in Euclidean space.  Coxeter groups are an important class
of groups, which have fascinating connections to many areas of mathematics, including
combinatorics; see the chapter by  Fomin and Reading in this volume \cite{fr04} and the
recent book of Bj\"orner and Brenti \cite{bb05}.   The finite irreducible Coxeter groups
have been completely classified.  There are four infinite families and six exceptional
irreducible Coxeter groups.       The most basic family consists of the type A Coxeter
groups, which are the symmetric groups
$\mathfrak S_n$ with the adjacent transpositions $(i,i+1),\, i=1,\dots, n-1$, forming the
generating set.  This is the reflection group of the (type A) braid arrangement. It is
also the group of symmetries  of the $n$-simplex.  The hyperoctahedral groups
$\s_n[\Z_2]$ with generators given by signed adjacent transpositions $(1,-1)$ and
$(i,i+1)$,
$i=1,\dots ,n-1$, form the type B family.  This is the reflection
group of the type B braid arrangement.  It is also the symmetry group of the $n$-cube and
the $n$-cross-polytope.

Let $(W,S)$ be a Coxeter system. Every $\sigma \in W$ can be expressed as a product,
$$\sigma = s_1
\dots s_k,$$ where
$s_i
\in S$.  The word $s_1 \dots s_k$ is said to be a {\em reduced expression} for
$\sigma$ if its length
$k$ is minimum among all words whose product is $\sigma$.  The {\em length} $l(\sigma)$ of
$\sigma$ is defined to be the length of a reduced expression for $\sigma$.

{\em Bruhat order} on $W$ is a partial order relation on $W$ that is defined via the
covering relation,
$\sigma <
\!\!\! \cdot \tau$ if
\begin{itemize} 
\item   $\tau = t \sigma$, for some $t \in T:=\{\alpha s \alpha^{-1}: s \in S, \alpha \in
W\}$
\item
$l(\tau) = l(\sigma)+1$.
\end{itemize}

For any subset $J$ of the generating set $S$, there is an induced  subposet of Bruhat
order on $W$, called the {\em quotient} by $J$,  defined as follows
$$W^J: = \{\sigma \in W : s\sigma > \sigma \mbox{ for all } s \in J\}.$$
Note $ W = W^{\emptyset} $.

Bruhat order describes the inclusion relationships of the Schubert subvarieties of a flag
manifold.
It was conjectured by de Concini and Stanley in the late 1970's that any open
interval of
 of Bruhat order on a quotient of a Coxeter group is homeomorphic to
a $d$-sphere or a $d$-ball, where
$d$ is the length of the interval.   This  result was needed by  de Concini and Lakshmibai
\cite{dl} in their work  on Cohen-Macaulayness of homogeneous coordinate rings of certain
generalized Schubert varieties.  In an attempt to  prove the conjecture by establishing
EL-shellability, Bj\"orner and Wachs instead came up with the notion of
CL-shellability and constructed a CL-labeling of the dual Bruhat poset (call this a dual
CL-labeling).  This labeling  relied on the following well-known characterization of
Bruhat order.  

\begth[Subword characterization of Bruhat order] Let $(W,S)$ be a Coxeter system.  Then
$\sigma <
\tau$ in Bruhat order on
$W$ if and only if for any reduced expression $w$ for $\tau$ there is a reduced expression
for
$\sigma$ that is a subword of $w$.
\enth

\begin{figure}\begin{center}
\includegraphics[width=6cm]{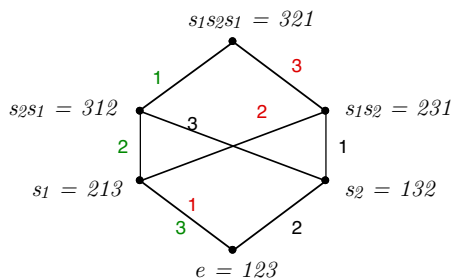}
\end{center}
\begin{center}\caption{Dual CL-labeling of $\mathfrak S_3$}\label{figbruhat2}
\end{center}
\end{figure}

We  describe the dual CL-labeling of intervals
$[\sigma,
\tau]$  of
$W^J$ (this works for infinite $W$ as well as finite $W$).  Fix a reduced expression $w$
for
$\tau$. It follows easily from the subword characterization that as we travel  down a
maximal chain, we delete (unique) letters of $w$ one at a time until we reach a
reduced word for
$\sigma$.  Label the edges of each maximal chain from top down with the position  of the
letter in $w$ that is crossed out.  This is illustrated in Figure~\ref{figbruhat2} on the
full interval $[e, 321]$ of Bruhat order on
$\mathfrak S_3$, where $s_i$ denotes the adjacent transposition $(i,i+1)$. 

\begth[Bj\"orner and Wachs \cite{bw82}]  \label{dualCL} The dual chain-edge labeling
described above is a dual CL-labeling of closed intervals of Bruhat order on $W^J$.  The
number of decreasing chains (from top to bottom) is
$1$ if $J = \emptyset$ and is at most $1$ otherwise. 
\enth

It follows from Theorem~\ref{dualCL} that every open interval $(\sigma,\tau)$ of $W^J$ has
the homotopy type of a $(l(\tau) -l(\sigma)-2)$-sphere or is contractible.  A stronger
topological consequence can be obtained by using the Danaraj and Klee result
(Theorem~\ref{danklee}).   Indeed, it can be shown that every
closed interval of length 2 in
$W$ has exactly
$4$ elements and that every closed interval of length 2 in $W^J$ has at most $4$
elements.   

\begin{cor} Every open interval $(\sigma,\tau)$ of Bruhat order on a quotient $W^J$ of  a
Coxeter group
$W$ is homeomorphic to a $(l(\tau) -l(\sigma)-2 )$-sphere or a $(l(\tau)
-l(\sigma)-2)$-ball.  If $J=\emptyset$ then $(\sigma,\tau)$ is homeomorphic to a $(l(\tau)
-l(\sigma)-2)$-sphere.
\end{cor}

Besides for quotients, there are other interesting classes of 
induced subposets of Bruhat order  such as descent classes, studied by Bj\"orner
and Wachs \cite{bw88}, and the subposet of involutions,
 which arose  in algebraic geometry work of Richardson and
Springer \cite{ricsp}.  There is a natural  notion of 
descent set for
general Coxeter groups which generalizes the descent set of a permutation.  A descent
class of a Coxeter system $(W,S)$ is the set of all elements whose descent set is some
fixed subset of $S$.  Bj\"orner and Wachs used the chain-edge labeling described above to
show  that every finite interval of  any descent class of any Coxeter group is
dual CL-shellable. This chain-edge labeling does not work for the subposet of involutions
because the maximal chains of the subposet are not maximal in the full poset. 
Recently Incitti found an EL-labeling that does work
for the classical Weyl groups (i.e.,  the type A, B or D Coxeter groups) and conjectured
that his result extends to all  intervals of all (finite or infinite)
Coxeter groups.  More recently, Hultman \cite{hul04} proved the main consequence of this
conjecture, namely  that all intervals of Bruhat order on the set of involutions of a
Coxeter group are   homeomorphic to spheres.

An alternative way to label the edges of Hasse diagram of $W^J$ is by labeling the edge
$\sigma < \!\!\!\!\cdot\,\, \tau$ with the element $t \in T$ such that $\tau = t \sigma$. 
By imposing a certain linear order on $T$,
Edelman \cite{e} showed that  this edge labeling
 is an EL-labeling for the symmetric group.  Proctor \cite{p} did the same for
 the classical Weyl groups.  Several years after   the
introduction of CL-shellability, Dyer \cite{dy} found a way to linearly order $T$ so that
the edge labeling by elements of $T$ is  an EL-labeling for all Coxeter groups and all
quotients.   Dyer's linear order (called reflection order) was recently used by Williams
\cite{wi05} to obtain an EL-labeling of the poset of cells in a certain cell decomposition
of  Rietsch
\cite{rie} of the totally nonnegative part of an arbitrary flag variety (for any reductive
linear algebraic group). This and the fact that the poset is thin (cf.,
Theorem~\ref{bjthin}) led Williams to conjecture that Rietsch's cell decomposition is a
regular cell decomposition of a ball (if true this would improve a result of Lusztig
\cite{lus} asserting that the totally nonnegative part of the flag variety is
contractible).

There are two other natural partial order relations on a Coxeter group $(W,S)$   that are
important and interesting. By replacing $T$ by $S$ in the
definition of Bruhat order, one gets the definition of {\it weak order}.  By replacing 
ordinary length $l$ by absolute
length  $al$ in the definition of Bruhat order, where
$al(\sigma)$ is the length of the shortest factorization of $\sigma$ in elements of 
$T$, one gets the definition of {\em absolute length
order}.

 The Hasse diagram  of
weak order is the same as the Cayley graph of the
 group
$W$ with respect to  $S$, directed away from the identity.  It is pure and bounded with the
same rank function as Bruhat order, ordinary length
$l$. Its topology was  first studied
by Bj\"orner \cite{bj842} who showed that although it is not lexicographically  shellable,
it  has the homotopy type of a sphere.

 The Hasse
diagram of absolute length order $W^{al}$ is the same as the Cayley graph of the
group
$W$ with  respect to
$T$, directed away from the identity.   It is pure with rank function, absolute length
$al$; its  bottom element is $e$,   but it has no unique top element.  In the symmetric
group
$\s_n$, the maximal elements are the $n$-cycles.   Interest in the  absolute length order
is a fairly recent development, which  arose in work on the braid group \cite{br01},
\cite{braw}, \cite{b03}.  It was shown by Brady
\cite{br01} that if
$W$ is the symmetric group
$\s_n$, then every  interval $[e,c]$ of $W^{al}$, where $c$ in an $n$-cycle, 
is isomorphic to the noncrossing partition lattice ${NC}_n$ discussed in
Section~\ref{ELpart}.   This observation (and connections to finite type Artin
groups) led Bessis \cite{b03} and Brady and Watt \cite{braw} to  define  the
noncrossing partition  poset for any finite Coxeter group to be the interval
$[e,c]$ of $W^{al}$,  where $c$ is a Coxeter
element of $(W,S)$.  A Coxeter element is  a product of all the elements of $S$ in some
order.  In $\s_n$ the Coxeter elements are the $n$-cycles.    Reiner's type B
noncrossing partition lattice as well as the classical (type A) noncrossing partition
lattice are recovered from this definition.   

\begin{xca}[Brady \cite{br01}] Show that the interval $[e,c]$, where $c$ is an $n$-cycle,
in absolute length order of the symmetric group is isomorphic to the noncrossing partition
lattice
$\mbox{NC}_n$.
\end{xca}

Although the topology of the interval $(e,c)$, where $c$ is a Coxeter element, is known to
have the homotopy type of a wedge of spheres via the Athanasiadis-Brady-Watt proof of
EL-shellability 
\cite{abw05}, little is known about the topology of the   full absolute
length poset, $W^{al}$, even for the symmetric group.

\begin{prob}[Reiner {\cite[Problem 3.1]{arm05}}]  What can be said about the topology of
$W^{al} -\{e\}$, where $W$ is an arbitrary finite Coxeter group?   Is 
$W^{al}$ lexicographically shellable, for types
$A$ and $B$?  It is known that for type $D$, the poset $W^{al}$ is not shellable.
\end{prob}

In the next exercise we see other examples of  posets  that admit CL-labelings.

\begin{xca} \label{normal} Let $\mathcal W_n$ be the poset of finite words over alphabet
$[n]$, ordered by the subword relation.   We have $34 <
23244 $ in $\mathcal W_n$.  Let $
\mathcal N_n$
 be the induced subposet  of normal words, where a word is said to be {\em
normal} if no two adjacent letters are equal.  For example, $2324$ is normal while $23244$
is not. 
\begin{enumerate} 
\item[(a)] (Bj\"orner and Wachs \cite{bw83}) Find a dual CL-labeling of each interval $[u,
v]$ in $
\mathcal N_n$, 
and show that $(u,v)$ is homeomorphic to an $(l(v) - l(u)-2)$-sphere.
\item[(b)] (Bj\"orner \cite{bj90}) Find a dual CL-labeling of each interval $[u,
v]$ in $
\mathcal W_n$.
\item[(c)] (Bj\"orner \cite{bj90, bj93})  Given a word $v=a_1\dots a_n$ in $\mathcal
W_n$, define its repetition set $R(v) = \{i: a_i = a_{i-1}\}$.  A {\it normal embedding} of
a  word
$u $ in $v= a_1 \dots a_n$   is a sequence $1 \le i_i < \dots < i_k\le n$ such that $u =
a_{i_1}\dots a_{i_k}$ and 
$R(v) \subseteq  \{i_1,\dots,i_k\}$. Show that the interval
$(u,v)$ in
$\mathcal W_n$ is homotopy equivalent to a wedge of 
$n_{u,v}$
spheres of dimension $l(v) - l(u)-2$, where $n_{u,v}$ is the number of normal
embeddings of $u$ in
$v$.
\end{enumerate}
\end{xca}

 In the next lecture, a formulation of CL-shellability,  called recursive atom ordering,
  which does not involve edge
labelings, will be presented.  This formulation has proved
to be  quite useful in  many applications, indeed, even more useful than the
original chain labeling version.    A more general form of lexicographic
shellability which does involve edge labelings was introduced by Kozlov
\cite{k97}. A
generalization of lexicographic shellability to  balanced  complexes was
formulated by Hersh
\cite{her} and was further studied  by Hultman
\cite{hul}.  A  relative version of lexicographic shellability was introduced by Stanley
\cite{st96}.  A Morse theory version of  lexicographical shellability, which is even more
general than Kozlov's version was formulated by Babson and Hersh \cite{bahe}.

\section{Rank selection} \label{ranksec}
Edge labelings of posets  were first introduced by Stanley \cite{st71, st72,st74} for the
purpose of studying the M\"obius function of rank-selected subposets of certain pure
lattices.  If one drops the requirement that the unique increasing maximal
chain be lexicographically first in the definition of EL-shellability then one has the kind
of labeling that  Stanley considered. Such a labeling was  called an R-labeling (we now
call it an ER-labling and we call the chain-edge version  a CR-labeling). 
 Here we discuss how 
CL-labelings are used to determine homotopy type and homology of rank-selected subposets of
pure CL-shellable posets. 

 Let
$P$ be a   bounded poset. For  $x \in P$, define the {\em rank},
$$r(x) := l([\hat 0, x]).$$  For $R \subseteq [l(P)-1]$, define the rank-selected subposet
$$P_R := \{ x \in P : r(x) \in R\}.$$  
Given a CL-labeling of a bounded poset $P$, the {\em descent set} of a maximal
chain $c: = (\hat 0 = x_0 \< x_1 \< x_2 \< \dots \< x_t \< x_{t+1} = \hat 1)$ of $P$ 
 is defined to be $$ \mbox{des}(c) := \{i \in [t] : \lambda(c,x_{i-1}\< x_i) \ge
\lambda(c,x_i\<x_{i+1})\}.$$ 

\begth[Bj\"orner and Wachs \cite{bw83}] \label{rankthm} Suppose $P$ is pure and
CL-shellable.
  Let
$R
\subseteq [l(P)-1]$.  Then
$\hat P_R$ is also CL-shellable, and hence $P_R$ has the homotopy type of a wedge of
$(|R|-1)$-spheres.  The number of spheres is the number of maximal chains of $P$ with
descent set
$R$.  Moreover, the set 
$$\{c_R : c \in \mathcal M(P) \,\,\& \,\, \des(c) = R \}$$
forms a basis for $\tilde H^{|R|-1}(P_R;\Z)$.
\enth

\begin{remark} Theorem~\ref{rankthm}  holds  for nonpure posets only when the rank set $R$
is rather special, see
\cite{bw97}.    It was shown by Bj\"orner \cite{bj80} that pure shellability is preserved
by rank selection. It is not known, however, whether pure EL-shellability  is preserved by
rank-selection.  Stanley \cite{st72,st74} proved the precursor to  Theorem~\ref{rankthm} 
that for pure posets $P$ that admit ER-labelings, the rank selected M\"obius invariant
$\mu(\hat P_R)$ is equal to $(-1)^{|R|-1}$ times 
 the number of maximal chains of $P$ with
descent set
$R$.  
\end{remark}

\begin{example}\label{rankboolex} {\em The rank-selected Boolean algebra $(B_n)_R$}.
Consider the EL-labeling of $B_n$ given in Example~\ref{ELboolex}.  Since the label
sequences of  the maximal chains are the permutations in
$\mathfrak S_n$, the number of maximal chains with descent set 
$R
\subseteq [n-1]$ is the number of permutations with descent set $R$.  So $(B_n)_R$ has
the homotopy type of a wedge of $d_{n,R}$ spheres of dimension $|R|-1$, where 
$d_{n,R}$ is the number of permutations in $\mathfrak S_n$ with descent set $R$.  (Note
that the truncated Boolean algebra given in Section~\ref{ELboolex} is an example of a
rank-selected Boolean algebra.)  We have the following equivariant homology version of the
homotopy result.
\end{example}

\begth[Solomon \cite{sol68}] \label{rankboolth} Let $R \subseteq [n-1]$ and let $H$ be the
skew hook with
$n$ cells and descent set $R$ (the definition of descent set of a skew hook was
given in Section~\ref{secrepsym}). Then 
$$\tilde H_{|R|-1} ((B_n)_R) \cong_{\mathfrak S_n} S^H.$$
\enth

\begin{xca} Prove Theorem~\ref{rankboolth} by using a map from the set of
tableaux of shape H to the set of maximal chains of $(B_n)_R$ thereby generalizing the 
proof of Theorem~\ref{solth}.
\end{xca}

We consider a type B-analog and a q-analog of
Example~\ref{rankboolex} in the next two exercises.

\begin{xca} \label{rankcrossex} {\em The rank-selected lattice of faces of the
$n$-cross-polytope
$(C_n)_R$}.   Recall form Example~\ref{crossaction} that we can identify the   
$(k-1)$-faces of $C_n$ with  $k$-subsets
$T$  of
$[n]
\cup
\{- i :i \in [n]\}$ such that $\{i,- i\} \nsubseteq T$ for all $i$.   Note that the
maximal chains of
$C_n$   correspond bijectively to elements of the hyperoctahedral group
$\s_n[\Z_2]$, i.e. the signed permutations.  As was mentioned earlier, there is a natural 
notion of  descent set for
general Coxeter groups which generalizes the descent set of a permutation.  For the
hyperoctahedral group,  descent can be characterized as follows: 
$i\in \{0,1,\dots,n-1\}$ is a  {\em descent} of a signed permutation $\sigma \in
\s_n[\Z_n]$ if either
\begin{itemize}
\item[(a)] $\sigma(i) > \sigma(i+1)$ (viewed as elements of $\Z$) or
\item[(b)] $i=0$ and $\sigma(1) < 0$.
\end{itemize}
 Find an
$EL$-labeling of
$C_n$ such that for all
$R \subseteq [n]$, the number of maximal chains with descent set $R$ equals $\bar d_{n,R}$,
the number of signed permutations with descent set $R$.  Consequently, 
$(C_n)_R$ has the homotopy type of a wedge of $\bar d_{n,R}$ spheres of dimension 
$|R|-1$.
\end{xca}

\begin{xca}[Simion \cite{sim95}] \label{rankqex} {\em The rank-selected subspace lattice
$(B_n(q))_R$}.  Find an EL-labeling
$\lambda$ of
$B_n(q)$ with label set $[n]$ such that  
\begin{itemize}
\item for each maximal chain $c$, the label sequence $\lambda(c)$ is a permutation in
$\s_n$,
\item for each permutation $\sigma \in \s_n$,
the number of maximal chains $c$ such that $\lambda(c) = \sigma$ is $q^{{\rm
inv}(\sigma)}$, where ${\rm inv}(\sigma)$ is the number of inversions of $\sigma$.
\end{itemize} 
Consequently, for all $R \subseteq [n-1]$, the rank-selected subspace lattice $(B_n(q))_R$
has the homotopy type of a wedge of
$d_{n,R}(q)
$
 spheres of dimension $|R|-1$, where $$d_{n,R}(q) = \sum_{\scriptsize{\begin{array}{c}
\sigma \in
\s_n
\\
\des(\sigma) = R\end{array}}} q^{{\rm inv}(\sigma)}.$$
\end{xca}

There are also interesting results on $\s_n$-equivariant rank-selected homology of the
partition lattice with ties to permutation enumeration due to Sundaram \cite{su,su95}. 
See also the related recent work of Hanlon and Hersh \cite{hanher} on rank-selected
partition lattices
\cite{hanher}

\begin{xca} Let $L$ be a geometric lattice with the EL-labeling given in
Section~\ref{geosec}.  Also let $R \subseteq [l(L)-1]$. Show that there is a basis for top
homology of $L_R$ consisting of fundamental cycles of subposets of $L_R$ that are all
isomorphic to the join  
$\bar B_{s_1} * \bar B_{s_2} * \dots * \bar B_{s_k}$, where
$\{s_1,s_1+s_2,\dots,s_1+s_2 + \dots + s_k\} = [l(L)] - R$.  
\end{xca}

Stanley \cite{st82}  gives a  type B analog of Theorem~\ref{rankboolth} for
the lattice of Exercise~\ref{rankcrossex}  and a $q$-analog of
Theorem~\ref{rankboolth} for the lattice of Exercise~\ref{rankqex}.  These  results involve
the   representation theories of the   hyperoctahedral group and the general linear
group, which are analogous to that of the symmetric group.  He proves
 Theorem~\ref{rankboolth} and the analogous results by invoking the following
general theorem.

\begth[\cite{st82}]\label{stanrank} Let $P$ be a bounded pure
shellable
$G$-poset of length $\l$.  If $R \subseteq [\l-1]$ then
\bq \label{rkhom1}\tilde H_{|R|-1}( P_R) \cong_G \bigoplus_{T\subseteq R}(-1)^{|R-T|}
 C_{|T|-1}( P_T).
\eq
Equivalently,
\bq \label{rkhom2} C_{|R|-1}( P_R) \cong_G \bigoplus_{T\subseteq R} \tilde
H_{|T|-1}( P_T).
\eq

\enth

In \cite{st82} Theorem~\ref{stanrank} is proved for a more general class
of posets called Cohen-Macaulay posets.  These posets are discussed in the
next section.  Equation~(\ref{rkhom1}) is a consequence of the Hopf trace formula and the
fact that pure shellability, or more generally Cohen-Macaulayness, is preserved by
rank-selection.  Equation~(\ref{rkhom2}) follows by the principle of inclusion-exclusion
or M\"obius inversion on the subset lattice.

Note that Theorems~\ref{rankboolth} and
\ref{stanrank} imply the  fact that the regular representation of $\s_n$ decomposes
into a direct sum of Foulkes modules, cf., Exercise~\ref{foudec}.

\begin{xca} Prove Theorem~\ref{rankboolth} by using Theorem~\ref{stanrank}.
\end{xca}

As was mentioned in the introduction to this lecture, shellability theory is intimately 
connected with the enumeration of faces of  polyhedral complexes,  a  central   topic in
geometric combinatorics; see the books of Stanley
\cite{st83} and Ziegler~\cite{z95}.
The   descent set of a maximal chain is a specialization of a more general concept in
shellability theory known as the restriction  of a facet, which is just the smallest new
face that is added to the complex when the facet is added.  When facets of pure shellable
simplicial complexes are enumerated according to the size of their restriction, an
important combinatorial invariant of the simplicial complex, known as the
$h$-vector,  is computed. For any
$(d-1)$-dimensional simplicial  complex $\Delta$ (shellable or not),
the
$h$-vector $(h_{0}(\Delta),h_1(\Delta),\dots,h_d(\Delta))$ and the $f$-vector
$(f_{-1}(\Delta),f_0(\Delta),\dots,f_{d-1}(\Delta))$ determine each other,
$$\sum_{i=0}^d h_i x^{d-i} = \sum_{i=0}^d f_{i-1} (x-1)^{d-i}, $$
 but the $h$-vector
is usually more convenient for expressing relations such as the upper bound conjecture, a
symmetry relation known as the Dehn-Sommerville equations, the  celebrated Billera-Lee
\cite{bl81} and Stanley \cite{st80}  characterization of the f-vector of a simplicial
polytope, and the conjectured extension to convex polytopes and homology spheres. 
Although  shellability has played an important role in the study of $f$-vectors and
$h$-vectors,  much fancier tools from commutative algebra and algebraic geometry have come
into play.   See the books of Stanley
\cite{st83} and Ziegler
 \cite{z95} for basic treatments of this material, the survey
article of Stanley \cite{st04} for important recent developments, and the 
paper of Swartz \cite{swa} for even more recent developments.

 Refinements of the
$f$-vector and
$h$-vector called the flag
$f$-vector and the flag
$h$-vector, respectively,   (defined for all pure posets and the more general  balanced 
complexes) have been extensively studied, beginning with the Bayer-Billera \cite{baybil}
analog of the Dehn-Sommerville equations;
see
\cite{st83} for further information.  For pure lexicographically
shellable posets the entries of the flag
$h$-vector have a simple combinatorial interpretation as the number of maximal chains
with fixed descent set.

For nonpure complexes,   two-parameter generalizations of the f-vector and
h-vector are defined in
\cite{bw96}.  The f-triangle and h-triangle also determine each other.  The f-triangle
counts faces of a simplicial complex according to the maximum size of a facet containing
the face and  the size of the face.  For  shellable complexes, it is shown in \cite{bw96}
that  the 
$h$-triangle counts facets   according to the size of the
restriction set and the size of the facet.  Duval \cite{duv} shows that the
entries of the
$h$-triangle are nonnegative for a more general class of complexes than the shellable
complexes; namely the sequentially Cohen-Macaulay complexes, which are discussed in the
next lecture.     It is pointed out  by Herzog, Reiner, and
Welker
\cite{hrw} that    the h-triangle of a sequentially Cohen-Macaulay complex $\Delta$
encodes the multigraded Betti numbers  appearing in the minimal free resolution of the
Stanley-Reisner ideal of the Alexander dual of $\Delta$.  Stanley-Reisner rings and ideals
are discussed briefly in the next lecture and Alexander duality is discussed in the last
lecture.

Descent sets and the more general restriction sets also play an important role in
  direct sum decompositions  of   Stanley-Reisner
rings of shellable complexes; due to 
 Kind and Kleinschmidt \cite{kk79}, Garsia \cite{garsia}, and Bj\"orner and Wachs
\cite[Section 12]{bw97}.

%% file: lect4.tex

\lecture{Recursive techniques}

The recursive definition of the M\"obius function of a poset provides a recursive
technique for computing the reduced Euler characteristic of the order complex of a
poset.  More  refined recursive techniques for computing the homology of a poset
are discussed in this lecture.    A general class of   posets to which
these techniques can be applied, the  Cohen-Macaulay posets or more generally
the  sequentially Cohen-Macaulay posets, are discussed in Section~\ref{CMsec}.  
A recursive formulation of
CL-shellability is presented in Section~\ref{RAO}.   In Section~\ref{ex} 
the recursive techniques for computing Betti numbers are demonstrated on various examples,
and in Section~\ref{whit} equivariant versions of these techniques are also demonstrated. 
In Section~\ref{basissec}   bases and generating sets for  homology and cohomology of some
of these examples are presented together with their use in computing $\s_n$-equivariant
homology.

\section{Cohen-Macaulay complexes} \label{CMsec}

Recall that the  link of a face $F$ of a simplicial complex $\Delta$ is the
subcomplex 
$$\lk_{\Delta} F := \{G \in \Delta : G \cup F \in \Delta, \,\,G \cap F = \emptyset\}.$$ 

\begin{definition} \label{CMdef} A
simplicial complex $\Delta$ is said to be {\em Cohen-Macaulay} over  $\mathbf k$ if 
$$\tilde H_i(\lk_{\Delta} F;\mathbf k) = 0,$$
for all $F \in \Delta$ and  $i < \dim\lk_{\Delta} F$.
\end{definition}

The following  result 
shows that Cohen-Macaulayness is a topological property.

\begth[Munkres \cite{munk84}] \label{Munk} The simplicial complex $\Delta$ is
Cohen-Macaulay over $\mathbf k$ if and only if
 for all $p \in \|\Delta\|$ and all $i < \dim \Delta$,
$$\tilde H_i(\|\Delta\|;\mathbf k) =  H_i(\|\Delta\|, \|\Delta\| - p; \mathbf k) =
0,$$
where the homology is reduced singular homology and relative singular homology,
respectively.
\enth

The Cohen-Macaulay property has its
origins in commutative algebra, in the theory of
Cohen-Macaulay rings.  Associated with every simplicial complex $\Delta$ on vertex set
$[n]$ is a ring
${\bf k}[\Delta]$ called the  {\em Stanley-Reisner ring} of the simplicial complex, which
is defined to be the quotient of the polynomial ring 
${\bf k}[x_1,\dots, x_n]$ by the {\it Stanley-Reisner ideal} $I_\Delta$ generated by
monomials of the form
$x_{i_1}x_{i_2}\dots x_{i_j}
 $ where $\{{i_1}, {i_2},\dots,i_j \} \notin \Delta$.  The
Stanley-Reisner construction is 
 a two-way bridge used to obtain topological and enumerative properties of
the simplicial complex from properties of the  ring and vice versa.
  A simplicial complex
is Cohen-Macaulay if and only if its Stanley-Reisner ring is a Cohen-Macaulay ring (see
\cite{st83} for the definition of Cohen-Macaulay ring).   The equivalence of the 
characterization given in Definition~\ref{CMdef} and the ring theoretic characterization
is due to Reisner
\cite{r76}.

It follows from Theorem~\ref{Munk} that any triangulation of a $d$-sphere is
Cohen-Macaulay.  The ring theoretic consequence of this fact played an essential role in
Stanley's celebrated proof of  the Upper Bound Conjecture for spheres.  
For $n > d$, define the cyclic polytope $C(n,d)$ to be the convex hull of any
$n$ distinct points on the moment curve $\{(t,t^2, \dots, t^d) \in \R^d : t \in \R \}$ (the
face poset of the polytope is independent of the choice of points).  The boundary complex
of the cyclic polytope is a simplicial complex.  The upper bound conjecture for spheres 
asserts that the boundary complex of the cyclic polytope achieves the maximum number of
faces of each dimension, over all   simplicial complexes on $n$ vertices that triangulate a
$d$-sphere.  See \cite{st83} for Stanley's proof of this conjecture and other very
important uses of the Stanley-Reisner ring and commutative algebra in   the enumeration of
faces of simplicial complexes.

The   Stanley-Reisner bridge can also be crossed in the opposite
direction obtaining ring theoretic information from the topology and combinatorics of the
simplicial complex,  as exemplified by the commutative algebra results mentioned
at the end of the last lecture.  Another example is a fundamental result of Eagon and
Reiner
\cite{er} which states that a square free monomial ideal  has a linear resolution if and
only if it is the Stanley-Reisner ideal $I_\Delta$ of a simplicial complex $\Delta$ whose
Alexander dual is Cohen-Macaulay.   We will not define linear resolution, but Alexander
duality is discussed in Section~\ref{prodsec}.   A nonpure generalization of this result
involving sequential Cohen-Macaulayness (discussed below) appears in papers of Herzog and 
Hibi
\cite{herhib} and Herzog, Reiner and Welker \cite{hrw}.   For further reading on the
extensive connections between simplicial topology and  commutative algebra, see Stanley's
classic book 
\cite{st83} and the recent book of Miller and Sturmfels
\cite{ms}.

\begin{xca}[Walker \cite{wa81}] \label{walker}
 Show that if $\lk_\Delta F $ is empty, $0$-dimensional, or connected for all $F \in
\Delta$, then
$\Delta$ is pure.  Consequently, Cohen-Macaulay simplicial complexes are pure. 
\end{xca}

The main tool for establishing Cohen-Macaulayness is shellability.  Indeed,  it follows
from Theorem~\ref{linkth} and Corollary ~\ref{shellhom} that pure shellability  implies the
Cohen-Macaulay property over any
$\mathbf k$.  From Exercise~\ref{walker} we see that this implication  does not hold for
nonpure shellability.  In order to extend the implication to  the nonpure setting,  
Stanley
\cite{st83} introduced   nonpure versions of  Cohen-Macaulay
 for complexes and rings, called sequentially Cohen-Macaulay, and showed that all
  shellable complexes are sequentially Cohen-Macaulay.   Duval
\cite{duv} and Wachs
\cite{wa2} found similar simpler characterizations of Stanley's sequential
Cohen-Macaulayness for simplicial complexes.  Here we take Wachs' characterization as the
definition.

\begin{definition} \label{seqCMdef} Let $\Delta$ be a simplicial complex.
For each
$m = 1,2,\dots,\dim(\Delta)$, let
$\Delta^{\langle m\rangle}$ be the subcomplex of $\Delta$ generated by facets of
dimension at least $m$. The complex $\Delta$ is said to be {\em sequentially acyclic over}
$\mathbf k$ if $\tilde H_i(\Delta^{\langle m \rangle} ; \k) = 0 $ for all $i <m$.  We say
that
$\Delta$ is {\em sequentially Cohen-Macaulay over} $\k$ if $\lk_{\Delta} F$ is sequentially
acyclic over $\k$ for all $F \in \Delta$.
\end{definition}

\begin{xca}  The {\em pure $m$-skeleton} $\Delta^{[m]}$ of a simplicial complex $\Delta$ is
defined to be the subcomplex generated by all faces of dimension $m$.   Show that the
formulation  of sequentially Cohen-Macaulay complex given in Definition~\ref{seqCMdef} is
equivalent to the following  formulation  of Duval \cite{duv}:
$\Delta$ is sequentially Cohen-Macaulay if and only if
$\Delta^{[ m]}$ is Cohen-Macaulay for all  $m =
1,2,\dots,\dim(\Delta)$.
\end{xca}
It is clear  that  pure  sequential
Cohen-Macaulayness is the same thing as  Cohen-Macaulayness.  The following generalization
of  Theorem~\ref{Munk} shows that sequential Cohen-Macaulayness is a topological
property, just as Cohen-Macaulayness is.   Given a nonnegative integer $m$ and a
topological space $X$, define
$X^{\langle m \rangle }$ to be the topological closure of the set 
$$\{ p \in X : p \mbox{ has a
neighborhood  homeomorphic to an open $d$-ball where $d \ge m$}\}.$$ 

\begth \label{genmunk} The simplicial complex $\Delta$ is sequentially Cohen-Macaulay over
$\k$ if and only if  for all $i <m$ and $p \in \|\Delta\|^{\langle m\rangle} $,
$$\tilde H_i(\|\Delta\|^{\langle m\rangle};\k)=H_i(\|\Delta\|^{\langle m\rangle},
\|\Delta\|^{\langle m\rangle} - p; \mathbf k) = 0.$$
\enth

\begin{xca}  Munkres proof of Theorem~\ref{Munk}
is based on the following lemma:  For any face $F $ of a simplicial complex $\Delta$, any
point $p$ in the interior of $F$, and any integer $i$,
$$ \tilde H_i(\lk F;\k) \cong H_{i+\dim F+1}(\|\Delta\|, \|\Delta\| -p;\k).$$
 Use this lemma to prove Theorem~\ref{genmunk}.
\end{xca}

\begin{xca}[Wachs \cite{wa2}]\label{free4} Show that if $\Delta$ is sequentially acyclic
then
$\tilde H_i(\Delta;\Z)$ is free for all $i$ and vanishes whenever there is no facet of
dimension
$i$.
\end{xca}

\begth[Stanley \cite{st83}] \label{shellCM} Every shellable simplicial complex is
sequentially Cohen-Macaulay over $\mathbf k$ for all $\mathbf k$.
\enth

\begin{xca}  Prove Theorem~\ref{shellCM}
by first showing  that if $\Delta$ is shellable then there
is  a shelling order of the  facets of $\Delta$ in which the dimensions of the facets
weakly decrease (see Bj\"orner and Wachs \cite{bw96}). 
\end{xca}

We say that a poset is (sequentially) Cohen-Macaulay if its order complex is (sequentially)
Cohen-Macaulay.   Early  work on Cohen-Macaulay posets can be found in the paper
of Baclawski
\cite{bac80} and in the survey article of Bj\"orner, Garsia, and Stanley \cite{bgs}.     
We have the following nice   characterization of sequentially Cohen-Macaulay posets, which,
in the pure case, is well-known and follows easily from Definition~\ref{CMdef}.  

\begth[Bj\"orner, Wachs, and Welker \cite{bww2}] A  poset
$P$ is sequentially Cohen-Macaulay if and only if every open   interval of $\hat
P$ is sequentially acyclic.
\enth

\noindent One direction of the proof follows
immediately from the fact that open intervals are links of chains.  The other is a
consequence of  the fact that  joins of sequentially acyclic simplicial complexes are
sequentially acyclic, which is proved in
\cite{bww2} by using a fiber theorem of Quillen.  Fiber theorems are discussed in Lecture
5.   It is also shown in
\cite{bww2} that other poset operations such as product preserve sequential
Cohen-Macaulayness, a fact that had been known for some time in the pure case \cite{wa88};
see Exercise~\ref{prodcm}.

To compute the unique nonvanishing Betti number of a  Cohen-Macaulay poset (or any poset
in which homology is concentrated in a single dimension),   one needs only to compute its
M\"obius invariant.  Indeed, it follows from the Phillip Hall Theorem
(Proposition~\ref{hall}) and the Euler-Poincar\'e formula (Theorem~\ref{eupo}) that  if the
homology of
$ P$ vanishes below the top dimension then 
\begin{equation} \label{mubetti} \tilde\beta_{l( P)}( P) =
(-1)^{l(P)}\mu(\hat P).\end{equation} By applying the recursive definition of
M\"obius function we get the  recursive formula observed by  Sundaram in \cite{su},
\begin{equation}\label{bettipure}\tilde \beta_{l(P)-1}( P\setminus \{\hat 0\}) = 
\sum_{\scriptsize{\begin{array}{c} x
\in P
\\ \end{array}}}(-1)^{l(P)+r(x)} \tilde \beta_{r(x)-2}(\hat 0, x),\end{equation} 
where $P$ is a  Cohen-Macaulay poset with a bottom element $\hat 0$ and $r(x)$ is the
rank of $x$.

Since general sequentially Cohen-Macaulay posets can
have homology in multiple dimensions, (\ref{mubetti}) and (\ref{bettipure}) do not hold in
the nonpure setting.  However, if the  sequentially Cohen-Macaulay poset has a property
known as semipure, then the following
generalization of (\ref{bettipure}) shows that  M\"obius function can still be used to
compute its Betti numbers.  A  poset
$P$ is said to be {\em semipure } if $P_{\le y} := \{x \in P : x \le y\}$ is pure for all
$y \in P$.    The proper part of the poset  given in Figure~\ref{figel}  is semipure. 
Also the face poset of any simplicial complex is semipure. 

\begth[Wachs \cite{wa2}] \label{betti4} Let $P$  be a  semipure sequentially
Cohen-Macaulay poset with a bottom element $\hat 0$. Then for all $m$,
\begin{equation}\label{bettisemipure}\tilde \beta_{m-1}( P\setminus \{\hat 0\}) = 
\sum_{\scriptsize{\begin{array}{c} x
\in P
\\ m(x) = m \end{array}}}(-1)^{m+r(x)} \tilde \beta_{r(x)-2}(\hat 0, x),\end{equation} 
where $m(x)$ is the length of the longest chain of $P$ containing $x$. 
\enth

\begin{xca} Prove Theorem~\ref{betti4}.
\end{xca}

We will demonstrate the
effectiveness
of the formulas (\ref{bettipure}) and (\ref{bettisemipure}) in computing Betti numbers on
concrete examples in Section~\ref{ex}.   Equivariant versions of these formulas will be
discussed in Section~\ref{whit}.  For more general versions see
\cite{wa2}.

There is a homotopy version of Cohen-Macaulay complexes due to Quillen \cite{q}  and a
homotopy version of sequentially Cohen-Macaulay complexes studied by  Bj\"orner, Wachs
and Welker \cite{bww1,bww2}. The requirement that homology vanish below a certain dimension
is replaced by the requirement that the homotopy groups vanish below that dimension. The
homotopy versions are stronger than the homology versions, but are still still weaker than
shellability.  They  are not topological properties, however.  For instance there is a
triangulation of the
$5$-sphere that is not homotopy Cohen-Macaulay, see \cite[Section 8]{q}.

We summarize the implications for bounded posets in the following diagram.

$$\begin{array}{ccc} \mbox{pure EL-shellable} & \Longrightarrow &  \mbox{EL-shellable}\\
\big{\Downarrow} & & \big{\Downarrow}\\
\mbox{pure CL-shellable} & \Longrightarrow &  \mbox{CL-shellable} \\
\big{\Downarrow} & & \big{\Downarrow}\\
\mbox{pure shellable} & \Longrightarrow &  \mbox{shellable}\\
\big{\Downarrow} & & \big{\Downarrow}\\
\mbox{homotopy Cohen-Macaulay} & \Longrightarrow &  \mbox{sequentially homotopy
Cohen-Macaulay}
\\\big{\Downarrow} & & \big{\Downarrow}\\
\mbox{Cohen-Macaulay over } \Z & \Longrightarrow &  \mbox{sequentially Cohen-Macaulay over
} \Z\\
\big{\Downarrow} & & \big{\Downarrow}\\
\mbox{Cohen-Macaulay over } \k & \Longrightarrow &  \mbox{sequentially Cohen-Macaulay over
} \k
\end{array}
$$

\noindent It is known that all the implications but $$\mbox{(pure) EL-shellable} 
\Longrightarrow \mbox{(pure) CL-shellable}$$
are strict.  Whether or not there are CL-shellable posets that are not EL-shellable is an
open question.  It is also unknown whether or not CL-shellability and dual CL-shellability
are equivalent. Examples of shellable posets that are not CL-shellable were obtained by
Vince and Wachs
\cite{vw} and Walker \cite{wal85}. 

 There are other important properties of simplicial
complexes which have recursive formulations such as vertex decomposability (introduced by
Provan and Billera
\cite{bp} for pure simplicial complexes and extended to  nonpure simplicial complexes by
Bj\"orner and Wachs \cite{bw97}) and constructible complexes (introduced by Hochster
\cite{hoc} and extended to the nonpure case by Jonsson \cite{jon05}).  There are no
special poset versions of these tools, however, as there are for shellability.   

\section{Recursive atom orderings} \label{RAO} In this section we present a rather
technical, but very useful tool for establishing poset shellability (and in turn
sequential Cohen-Macaulayness).  This technique, called recursive atom ordering, was
introduced by Bj\"orner and Wachs   in the early 1980's.   It 
 is equivalent to CL-shellability, but does not involve edge labelings.

\begin{definition} \label{recatdef} A bounded poset $P$ is said to {\em admit a
recursive atom ordering} if its length $l(P)$ is $1$, or if $l(P)>1$
and there is an ordering $a_1,a_2,\,\ldots,a_t$ of the atoms of $P$
that satisfies:
\begin{itemize}
\item[(i)] For all $j=1,2,\,\ldots,t$ the interval $[a_j,\hat1]$ admits
a recursive atom ordering in which the atoms of $[a_j,\hat1]$ that
belong to $[a_i,\hat1]$ for some $i<j$ come first.
\item[(ii)] For all $i<j$, if $a_i,a_j<y$ then there is a $k<j$ and an
atom $z$ of $[a_j,\hat1]$ such that $a_k<z\le y$.
\end{itemize}
 A {\em recursive coatom ordering} is a recursive atom
ordering of the dual poset $P^*$. 
\end{definition} 

\begin{figure}\begin{center}
\includegraphics[width=6 cm]{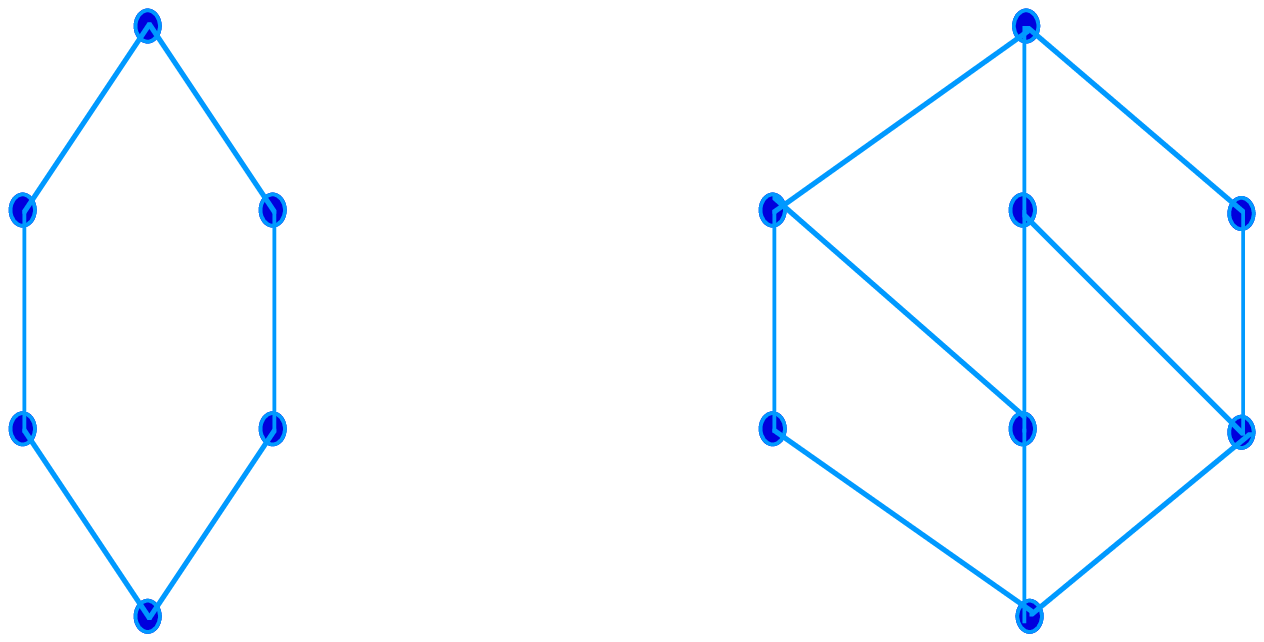}\\
\hspace{-.1in}(a) \hspace{1.5in}(b)
\end{center}
\begin{center}\caption{ }\label{figrecatom}
\end{center}
\end{figure}

Figure~\ref{figrecatom} (a) gives an example of a poset that does not admit a recursive
atom ordering since condition (ii) fails for every ordering of the atoms. The poset in
Figure~\ref{figrecatom} (b)  does admit a recursive atom ordering. The left to right order
in the  drawing of the Hasse diagram gives an  ordering of the atoms
 that satisfies (i) and (ii), for {\em each } interval.

\begth[Bj\"orner and Wachs \cite{bw83, bw96}] \label{CLrao} A bounded poset $P$ is
CL-shellable if and only if
$P$ admits a recursive atom ordering. 
\enth

\begin{proof}[Proof idea] We use a recursive procedure to obtain a recursive atom
ordering from a CL-labeling and vice versa.  
Given a CL-labeling
$\lambda$, order the atoms of $P$  in increasing order of the labels
$\lambda(c,\hat 0 <\!\!\!\cdot \, a)$, where $c$ is any maximal chain containing the atom
$a$ (the label is independent of the maximal chain
$c$ since
$\hat 0
\<a$ is a bottom edge).  Then recursively use the restriction of the CL-labeling to each 
interval
$[a,\hat 1]$ to obtain a recursive atom ordering of  
$[a,\hat 1]$.

Conversely, given a recursive atom ordering of $P$, label the bottom edge of each maximum
chain with the position of the atom in the atom ordering.  Then  recursively use the
recursive atom ordering of $[a,\hat 1]$ to obtain an appropriate chain-edge labeling of
$[a,\hat 1]$, for each atom $a$. 
\end{proof}

A bounded poset $P$ is said to be {\em  semimodular} if for all
$u,v \in P$ that cover some $x \in P$, there is an
element
$y
\in P$  that covers both
$u$ and
$v$.   Semimodular posets are  not necessarily shellable.  A poset is said to
be {\em totally semimodular} if every closed interval is semimodular.  Note that
semimodular lattices are totally semimodular.  In particular, the boolean algebra and the
partition lattice are totally semimodular. 

\begth[Bj\"orner and Wachs \cite{bw83}] \label{totsemi} Every ordering of the atoms of
a totally semimodular poset is a recursive atom ordering.  
\enth

\begth[Bj\"orner and Wachs \cite{bw83}] \label{raoshell} An ordering of the facets of
a simplicial complex $\Delta$ is a shelling if and only if the ordering is a recursive
coatom ordering of the face lattice $ L(\Delta)$.  
\enth

\begin{xca} Prove \begin{itemize}
\item[(a)]  Theorem \ref{totsemi}.
\item[(b)] Theorem \ref{raoshell}
\end{itemize}
\end{xca}

\begin{definition}[Wachs and Walker \cite{ww}]
 A  meet semilattice $P$ with rank function
$r$ is said to be a {\em geometric semilattice} if
\begin{itemize}
\item[(i)] every interval is a geometric lattice
\item[(ii)] for all $x \in P$ and subset $A$ of atoms whose join exists, if $r(x) < r
(\bigvee A)$ then there is an $a \in A$ such that $a \not\leq x$ and $a \lor x$ exists.
\end{itemize}
\end{definition}

 Examples of geometric
semilattices include 
\begin{itemize}
\item  geometric lattices
\item intersection semilattices of  affine hyperplane arrangements
\item face posets of matroid complexes
\item top truncated geometric lattices; i.e., $L_S$ where $L$ is a geometric lattice and
$$S=
\{0,1,2,\dots,j\}$$ for some $j \le l(L)$
\item $L-[x,\hat 1]$, where $L$ is  a geometric lattice
 and $x \in L$ (this is a characterization). 
\end{itemize}

A set $A$ of atoms in a geometric semilattice $P$ is said to be {\em independent} if
$\bigvee A$ exists and $r(\bigvee A) = |A|$.  A {\em basic } set of atoms is a maximal
independent set of atoms.

\begth[Wachs and Walker \cite{ww}] \label{geosemi}
Any  ordering of the atoms of  a geometric semilattice $P$ that begins with a
basic set of atoms is a recursive atom ordering of $P
\cup \{\hat 1\}$.  
\enth

\begin{xca} Let $P$ be a geometric semilattice.
\begin{itemize}\item[(a)] Show $P$ is pure.
\item[(b)] Show every principal upper order ideal $P_{\ge x} := \{y \in P : y \ge x\}$ is
a geometric semilattice.
\item[(c)] Prove Theorem~ \ref{geosemi}.
\end{itemize}
 
\end{xca}

From the proof of Theorem~\ref{CLrao}, we know that to every 
recursive atom ordering of a poset one can associate a  CL-labeling.   A  natural question
to ask, when  given a particular  recursive atom ordering,
 is whether there is a nice associated CL-labeling.  Then one can  obtain Betti numbers
and bases from the decreasing chains. Ziegler
\cite{z92} describes such a  CL-labeling for one of the recursive atom orderings of
Theorem~\ref{geosemi} and shows that it is in fact  an EL-labeling.  He then constructs 
bases for homology and cohomology of geometric semilattices which generalize the Bj\"orner
NBC bases for geometric lattices  discussed in Section~\ref{geosec}.

\section{More examples} \label{ex}

In this section we present examples of semipure posets that admit recursive atom orderings
and we demonstrate the use of the recursive formulas
(\ref{bettipure}) and (\ref{bettisemipure}) in computing  their Betti numbers. Let's begin
with some {\em pure} examples.

\begin{example}\label{injex} Let $k\le n$.  The injective word poset $\mathcal I_{n,k}$ is
defined to be the poset of words of length at most $k$ over alphabet $[n]$ with no repeated
letters.  The order relation is the subword relation. This  poset is pure and is  a lower
order ideal  of the normal word poset $\mathcal N_{n}$
(discussed in Exercise~\ref{normal}) which in turn is a subposet of the word poset
$\mathcal W_n$ (also discussed in Exercise~\ref{normal}).     It is left as an exercise to
show that  lexicographic order on the injective words of length
$k$ is a recursive coatom ordering of
$\mathcal I_{n,k}
\cup
\{\hat 1 \}$.  One can use the fact that all closed intervals of
$\mathcal I_{n,k}$ are  isomorphic to  Boolean algebras.  Since the duals of the closed
intervals are totally semimodular, one  needs only to verify condition (i) of
Definition~\ref{recatdef}.    
\end{example}

\begin{xca} \label{injexer} \begin{itemize}
\item[]\item[(a)] Show that lexicographic order on the 
injective words of length $k$ is a recursive coatom ordering of $\mathcal I_{n,k}
\cup
\{\hat 1 \}$.
\item[(b)] Use formula (\ref{bettipure}) to show that the top
Betti number 
$\tilde \beta_{k-1}(\bar I_{n,n})$ is equal to
$d_n$, the number of derangements of $n$ elements.  Consequently, $\bar I_{n,n}$ has the
homotopy type of a wedge of
$d_n$ spheres of dimension $k-1$ .
\item[(c)]  Let $\mathcal N_{n,k}$ be the
induced subposet of $\mathcal N_n$ consisting of words of length at most $k$.  Use
Exercise~\ref{normal}, Theorem~\ref{rankthm}, and  formula (\ref{bettipure}) to show that
$\bar {\mathcal N}_{n,k}$ has the homotopy type of a wedge  of $(n-1)^k$ spheres of
dimension
$k-1$. 
\item[(d)]  Let $\mathcal W_{n,k}$ be the
induced subposet of $\mathcal W_n$ consisting of words of length at most $k$.  Show that
$\bar {\mathcal W}_{n,k}$ also has the homotopy type of a wedge  of $(n-1)^k$ spheres of
dimension
$k-1$.
\end{itemize}
\end{xca}

\begin{remark} The injective word poset and normal word poset of
Example~\ref{injex}  were introduced by Farmer \cite{far} who showed that  $\bar{\mathcal
I}_{n,k}$ and
$\bar{\mathcal N}_{n,k}$ have the homology of a wedge of
$(k-1)$-spheres. Bj\"orner and Wachs \cite{bw83} recovered these results and
strengthened them to homotopy by establishing shellability. The Betti
number computations of Exercise~\ref{injexer} (b) and (c) are  due to Reiner and Webb
\cite{rw} and Farmer \cite{far}, respectively.
\end{remark}

Recall that in Example~\ref{ELpart}, we used EL-shellability to show that the proper part
of the partition lattice $\Pi_n$ has the homotopy type of a wedge of $(n-1)!$ spheres of
dimension $n-3$. In the next two examples we use the techniques of Sections~\ref{CMsec}
and~\ref{RAO} to obtain analogous results for the even and odd block size partition
posets.

\begin{example} \label{evenex} Let $\Pi_{2n}^{\rm even}$ be the subposet of $\Pi_{2n}$
consisting of partitions  whose block sizes are  even.   The even
block size partition poset is pure but lacks a bottom element. The upper intervals $[x,
\hat 1]$
 are all partition
lattices, which are totally semimodular; so every atom ordering of $[x, \hat 1]$ is a
recursive atom ordering.  Associate a word with each
atom of $\Pi_{2n}^{\rm even}
\cup \{\hat 0\}$ by listing the elements of each block in increasing order and then listing
the blocks in lexicographic order.  We claim that lexicographic order on the words
corresponds to a recursive atom ordering of 
$\Pi_{2n}^{\rm even}
\cup \{\hat 0\}$.  The verification is left as an exercise.

We next use (\ref{bettipure}) to compute the Betti numbers of $\bar \Pi_{2n}^{\rm even}$.
   Since
$\Pi_{2n}^{\rm even}$ is pure and the upper intervals are isomorphic to partition
lattices, to compute the unique nonvanishing Betti number $\tilde 
\beta_{n-2}( 
\bar\Pi_{2n}^{\rm even})$, we  apply (\ref{bettipure}) to the dual of the
poset.    Let
$\beta_{2n}$ denote the  Betti number
$\tilde  \beta_{n-2}( 
\bar\Pi_{2n}^{\rm even})$. By (\ref{bettipure}) we have
\begin{eqnarray*}\beta_{2n} &=& \sum_{x }
(-1)^{b(x)+n }\tilde\beta_{b(x)-3}(x,\hat 1)\\ &=&
\sum_{x } (-1)^{b(x) +n} (b(x)-1)!\\ &=& 
\sum_{r=1}^n (-1)^{r+n} (r-1)! \,\,|\{x \in \Pi_{2n}^{\rm even}: b(x) =
r\}|,\end{eqnarray*} where $b(x)$ denotes the number of blocks of $x$.  
Note that 
$$|\{x \in \Pi_{2n}^{\rm even}: b(x) = r\}| = {1 \over r!} \sum_{(j_1,j_2,\dots, j_r)
\vDash n} \binom{2n}{2j_1,2j_2,\dots,2j_r},
$$ where $\vDash n$ denotes composition of $n$.
We now have
$$\beta_{2n} = \sum_{r=1}^n (-1)^{r+n}  {1 \over r} \sum_{(j_1,j_2,\dots, j_r)
\vDash n} \binom{2n}{2j_1,2j_2,\dots,2j_r}.
$$
The exponential generating function for the Betti numbers is thus given by, 
\begin{eqnarray}\nonumber \sum_{n \ge 1} \beta_{2n } {u^{2n}\over(2n)!} &=& \sum_{r \ge 1}
(-1)^r {1
\over r} \sum_{n\ge 1}
\sum_{(j_1,j_2,\dots, j_r)
\vDash n} \binom{2n}{2j_1,2j_2,\dots,2j_r} (-1)^n{u^{2n}\over (2n)!}\\ \nonumber
&=&  \sum_{r \ge 1} (-1)^r {1
\over r}\Big{(} \sum_{j\ge 1} (-1)^j{u^{2j} \over (2j)!}\Big{)}^r\\ \label{cos}
&=&- \ln(\cos u).\end{eqnarray} 
By taking derivatives we get,
\begin{eqnarray*}\sum_{n \ge 1}  \beta_{2n } {u^{2n-1}\over(2n-1)!}  &=& \mbox{tan }
u.
\end{eqnarray*} 
So $\beta_{2n}$ is equal to the coefficient of ${u^{2n-1}\over (2n-1)!}$ in the
Taylor series expansion of $\tan u$. It is well known that this coefficient  is  equal to
the Euler number
$E_{2n-1}$, where   
$E_m$ is defined to be the  number of  alternating  permutations in
$\mathfrak S_{m}$, i.e., permutations with descents
at all the even positions and ascents at all the odd positions.   We conclude that
$\bar \Pi_{2n}^{\rm even}$ has the homotopy type of a wedge of $E_{2n-1}$ spheres of
dimension $n-2$.

There is an alternative way to arrive at $\beta_{2n} = E_{2n-1}$.  Replace the last line of
(\ref{cos}) with
$$-\ln \Big{(} \sum_{j\ge 0} (-1)^j{u^{2j} \over (2j)!}\Big{)}$$ and then  take derivatives
of both sides of the equation.  This results in 
\begin{eqnarray*}\sum_{n \ge 1}  \beta_{2n } {u^{2n-1}\over(2n-1)!}  \sum_{j\ge 0}
(-1)^j{u^{2j} \over (2j)!} &=&  \sum_{j\ge 0} (-1)^{j-1}{u^{2j-1} \over (2j-1)!}.
\end{eqnarray*}
By equating coefficients we obtain the recurrence relation, 
\begin{eqnarray*}\sum_{r = 1} ^j \binom{2j-1}{2r-1} (-1)^{r} \beta_{2r} &=&-1
\end{eqnarray*}
It is not difficult to check that $E_{2r-1}$ satisfies the same  recurrence relation.
\end{example}

\begin{example} \label{oddex}
Let $\Pi_{2n+1}^{\rm odd}$ be the
subposet  of $\Pi_{2n+1}$ consisting of
partitions  whose block sizes are odd.  The odd block size partition poset is
totally semimodular; so any atom ordering is a recursive atom ordering.

Now let $\beta_{2n+1}$ denote the top Betti number
$\tilde  \beta_{n-2}( 
\bar\Pi_{2n+1}^{\rm odd})$.  An  argument similar to the one used in deriving (\ref{cos})
yields: 
$\sum_{n \ge 0} \beta_{2n+1 } {u^{2n+1}\over(2n+1)!}$ is the compositional
inverse of $\sum_{n \ge 0}  (-1)^n{u^{2n+1}\over(2n+1)!}$, i.e.,
\begin{eqnarray}\label{compinv}\sum_{n \ge 0} \beta_{2n+1 } {1\over(2n+1)!}\Big{(}\sum_{j
\ge 0}  (-1)^j{u^{2j+1}\over(2j+1)!}\Big{)}^{2n+1} = u.\end{eqnarray}
Since $\sum_{j \ge 0} 
(-1)^j{u^{2j+1}\over(2j+1)!} = \sin u$, we have 
$\sum_{n \ge 0} \beta_{2n+1 } {u^{2n+1}\over(2n+1)!} = \sin^{-1} u$.
The coefficient of $u^{2n+1}\over(2n+1)!$ in the Taylor series expansion of $\sin^{-1} u$
is $(2n+1)!!^2$, where $(2n+1)!! := 1\cdot 3 \cdot 5  \cdots (2n+1)$.
We conclude that $\bar \Pi_{2n+1}^{\rm odd}$ has the homotopy type of a wedge of
$(2n+1)!!^2$
spheres of dimension $n-2$. 
 \end{example}

\begin{xca}\label{modd}  For  integers $n,d\ge 1$ and  $k \ge 0$, let
$\Pi_{nd+k}^{k
\bmod
\,d}$ be the subposet of the partition lattice $\Pi_{nd+k}$ consisting of partitions whose
block sizes  are congruent to
$k \bmod d$.  This poset is pure  only for $k\equiv  0$ or $1 \mod d$, and is bounded only
for
$k\equiv 1 \mod d$. \begin{itemize}
\item[(a)] Show that lexicographical order on  the words associated with the atoms as in
Example~\ref{evenex}  (i.e., list the elements of each block in increasing order and then
list the blocks lexicographical order) yields a recursive atom ordering of $\Pi^{0\bmod
d}_{nd} \cup
\{\hat 0\}$. 
\item[(b)] 
 Show that $\Pi_{nd+1}^{1
\bmod
\,d}$ admits a recursive atom ordering.
\item[(c)] Find a recursive atom ordering of  $\Pi_{nd+k}^{k
\bmod
\,d} \cup \{\hat 0\} $ for $k \not \equiv 1 \mod d$.
\item[(d)] Show that $\bar \Pi^{0\bmod d}_{nd}$ has the homotopy type of a wedge of
$E^d_{dn-1}$
spheres of dimension $n-2$, where $$E^d_{dn-1} = |\{\sigma \in \mathfrak S_{nd-1} :
\mbox{des}(\sigma) = \{d,2d,\dots,(n-1)d\}|.$$
\item[(e)] Show that $\bar \Pi^{1\bmod d}_{nd+1}$ has the homotopy type of a wedge of $c_n$
spheres of dimension $n-2$, where $\sum_{n \ge 0} c_n {u^{nd+1} \over (nd+1)!}$ is the
compositional inverse of $\sum_{n \ge 0}  (-1)^n{u^{nd+1}\over(nd+1)!}$.
\end{itemize} 
\end{xca}

\begin{remark} The result on  the
M\"obius invariant  of the even block size partition lattice derived in
Example~\ref{evenex} first appeared in 
 in the 1976 MIT  thesis of Garrett Sylvestor
 on Ising ferromagnets \cite{syl}.  Stanley \cite{st78} extended this result to the
$d$-divisible partition lattice of Exercise~\ref{modd} (d).  The M\"obius invariant
result derived in   Example~\ref{oddex} and Exercise~\ref{modd} (e) are also due to
Stanley.  The recursive atom ordering for the $d$-divisible partition lattice
(Example~\ref{evenex} and Exercise~\ref{modd} (a)) is due to  Wachs (see \cite{sa86}), as
is an EL-labeling of this poset
\cite{wa96}.   The recursive atom ordering for the
$1
\bmod d$ partition lattice (Example~\ref{oddex}
and Exercise~\ref{modd} (b)) is due to Bj\"orner (see \cite{chr} and \cite{bw83}). 
The general  recursive 
atom ordering of $ \Pi^{k
\bmod
\,d}_{nd+k} \cup \{\hat 0\}$ of Exercise~\ref{modd} (c) is  due to Wachs~\cite{wa2}.
\end{remark}

In the next example, we demonstrate the full power of Theorem~\ref{betti4} by
finding the Betti numbers for a nonpure shellable poset.
\begin{example}\label{kex}  For $n \ge k \ge 3$, let $\Pi_n^{\ge k}$ be 
the subposet of $\Pi_{n}$ consisting of
partitions  whose block sizes are  at least $k$.  It was shown by Bj\"orner and Wachs
\cite{bw96} that the poset
$\Pi_n^{\ge k}
\cup
\{\hat 0\}$ admits  a recursive atom ordering similar to that of $\Pi_{2n}^{\rm even}
\cup
\{\hat 0\}$. The following computation of Betti numbers  appears in
\cite{wa2}.  The dual of $\Pi_n^{\ge k}$ is semipure and the upper intervals  are
partition lattices; so we  apply (\ref{bettisemipure}) to the dual.    First note that 
\bq \label{mk}m(x) = \sum_{i = 1}^{b(x)} \left\lfloor {|B_i|\over
k}
\right\rfloor -1,\eq
where $B_1,B_2,\dots,B_{b(x)}$ are the blocks of $x$.   By (\ref{bettisemipure}) we have,
\begin{eqnarray*}\tilde \beta_{m-2}(\bar
\Pi_n^{\ge k}) &=& \sum_{x: m(x)= m-1} (-1)^{b(x)+m }\tilde\beta_{b(x)-3}(x,\hat 1)\\  &=& 
\sum_{r\ge 1} (-1)^{r+m} (r-1)! \,\,|\{x \in \Pi_{n}^{\ge k}: b(x) =
r, \, m(x) =m-1\}|.\end{eqnarray*}   
We have, 
$$|\{x \in \Pi_{n}^{\ge n}: b(x) = r,\, m(x) = m-1\}| = {1 \over r!}
\sum_{\scriptsize{\begin{array}{c}(j_1,j_2,\dots, j_r)
\vDash n\\j_i \ge k \,\,\forall i \\ \sum_{i=1}^r \lfloor j_i/k\rfloor = m\end{array}
}}\binom{n}{j_1,j_2,\dots,j_r}.$$
The two parameter exponential generating function is thus,
\begin{eqnarray*}&&\hspace{-.5in}\sum_{m,n \ge 1} (-1)^m \tilde \beta_{m-2}(\bar
\Pi_n^{\ge k})\,t^m {u^{n}\over
n!}\\ &=&
\sum_{r
\ge 1} (-1)^r {1
\over r} \sum_{n,m\ge 1}
\sum_{\scriptsize{\begin{array}{c}(j_1,j_2,\dots, j_r)
\vDash n\\j_i \ge k \,\,\forall i \end{array}
}} \binom{n}{j_1,j_2,\dots,j_r} 
\,\,t^{\sum_{i=1}^r
\lfloor j_i/k\rfloor}{u^{n}\over  n!}\\ &=&  \sum_{r \ge 1} (-1)^r {1
\over r}\Big{(} \sum_{j\ge k}  \,\,\,t^{\lfloor j/k \rfloor}{u^{j} \over
 j!}\Big{)}^r\\
&=& -\ln \Big{(}1+\sum_{j\ge k} \,\,\,t^{\lfloor j/k \rfloor}{u^{j} \over
j!}\Big{)}.\end{eqnarray*} 
\end{example}

\begin{xca} Show that $\Pi_n^{\ge k} \cup \{\hat 0\}$ has a recursive atom ordering. 
\end{xca}

\begin{prob} Linusson
\cite{l} computed the M\"obius invariant of  the  poset $\Pi_n^{\ge k} \cup \{\hat 0\}$
 in order to compute lower bounds for the complexity of a problem similar to the
$k$-equal problem of Section~\ref{kesec}; namely that of determining whether   a given
list   of real numbers has the property that the number of occurrences  of each entry is
at least $k$. It was shown by Bj\"orner and Lov\'asz \cite{bl} that Betti number
computations give better bounds than  M\"obius function computations.   Can the
Betti number computation of Example~\ref{kex} be used to
 improve Linusson's lower bound for the complexity of the ``at least k problem''?  
\end{prob}

\begin{xca}[Wachs \cite{wa2}] \label{jmodd}  In this exercise, we generalize the results of
Examples~\ref{evenex}~and~\ref{oddex} to the nonpure case  of Exercise~\ref{modd} . For
positive integers
$n,d,k$, with
$ k
\le d$, let
 $k_0 = k / \mbox{gcd}(k,d)$ and  $d_0 = d / \mbox{gcd}(k,d)$.   Show that
$$\sum_{\scriptsize\begin{array}{c} m\ge 1 \\n \ge 0\end{array}}  (-1)^m
\beta_{m-2}(\bar
\Pi^{k
\bmod
\,d}_{nd+k})\,t^{md_0+1} {u^{nd+k}\over
(nd+k)!} = f\Big{(}\sum_{i\ge 0} \,\,\,t^{\lfloor i/k_0 \rfloor d_0+1}{u^{id+k} \over
(id+k)!}\Big{)},$$
where $f(y)$ is the compositional inverse of the formal power series
$$g(y) = \sum_{i\ge 0} {y^{i d_0+1} \over
(id_0+1)!}.$$
\end{xca}

\section{The Whitney homology technique} \label{whit}  In this section we discuss a
technique for computing group representations on the homology of sequentially
Cohen-Macaulay posets and
demonstrate its use on the examples of Section~\ref{ex}.  This technique was introduced by
Sundaram \cite{su} in the pure case and later generalized to semipure posets by Wachs
\cite{wa2}. The pure case  is based on an equivariant version of
(\ref{bettipure}) and the semipure version is based on the more
general Theorem~\ref{betti4}.

 For any
 Cohen-Macaulay $G$-poset $P$ with bottom element $\hat 0$,   {\em Whitney homology} of
$P$   is defined for each integer
$r$ as follows,
$$\wh_r(P)= \bigoplus_{x \in P_r}  \tilde H_{r-2}(\hat 0,x),$$
where $P_r := \{x \in P : r(x) = r\}$.   

The action of $G$ on $P$ induces a representation of $G$ on $\wh_r(P)$. Indeed $g\in G$
takes
$(r-2)$-chains of $(\hat 0,x)$ to $(r-2)$-chains of $(\hat 0,gx)$. More precisely, as a
$G$-module
$$\wh_r(P)= \bigoplus_{x \in P_r/G}  \tilde H_{r-2}(\hat
0,x)\uparrow_{G_x}^G,$$
where $P_r/G$ is a set of orbit representatives and $G_x$ is the
stabilizer of
$x$.

 Whitney
homology  for geometric lattices  was introduced by Baclawski \cite{bac75}  as the homology
of an algebraic complex
whose Betti numbers are the  signless   Whitney numbers of the first kind.  The
formulation given here is due to Bj\"orner  \cite{bj82}.  Whitney homology for
intersection lattices of complex hyperplane arrangements forms an algebra 
isomorphic to an algebra that   Orlik and Solomon used to give a combinatorial
presentation of the cohomology algebra of the complement of the arrangement; see
\cite{bj92}. Sundaram
\cite{su} recognized  that Whitney homology for any Cohen-Macaulay poset could be used as a
tool in computing group representations on homology.  Her technique is based on the 
following result, which she obtains as a consequence of the Hopf trace formula.

\begth[Sundaram \cite{su,su94}]\label{sundth}  
Suppose $P$ is a 
$G$-poset with a bottom element $\hat 0$.   Then
\begin{eqnarray} \label{sundhopf}\bigoplus_{r=-1}^{l(P)-1} (-1)^r \tilde H_{r}( P\setminus
\{\hat 0\})\,\,
\cong_G\,\,
\bigoplus_{r=0}^{l(P)} (-1)^{r-1}
\,\,\bigoplus_{x \in P/G}  \tilde H_{r-2}(\hat
0,x)\uparrow_{G_x}^G.
\end{eqnarray}
Consequently, if $P$ is also Cohen-Macaulay,
\begin{eqnarray} \label{sundeq}\tilde H_{l(P)-1}( P\setminus \{\hat 0\})\,\, \cong_G\,\,
\bigoplus_{r=0}^{l(P)} (-1)^{l(P)+r}
\,\,\wh_r(P).
\end{eqnarray}
\enth

\begin{xca} \begin{enumerate} \item[] \item[(a)] Prove Theorem~\ref{sundth}.
\item[(b)] Show that if $P$ is the face poset of a simplicial complex
then (\ref{sundhopf}) reduces to the Hopf trace formula.
\end{enumerate}
\end{xca}

\begin{example}[Reiner and Webb \cite{rw}] Consider the injective word poset $\mathcal
I_{n,k}$ of Example~\ref{injex}. Let $G$ be the  symmetric group $\mathfrak S_n$, which
acts on 
$\mathcal I_{n,k}$ in the obvious way.  Since $\mathcal I_{n,k}$ is 
Cohen-Macaulay, we can apply (\ref{sundeq}). Let $x \in  \mathcal I_{n,k}$ be a word of
length
$r$.  Clearly
$G_x $ is isomorphic to the Young subgroup 
$(\mathfrak S_1)^{\times r}
\times \mathfrak S_{n-r}$, since the letters of $x$  must be fixed and the letters
outside of
$x$ may be freely permuted.   $G_x$ acts trivially on the letters outside of $x$.  
The interval
$(\hat 0,x)$ is isomorphic to the proper part of the Boolean algebra
$B_{r}$; so its top homology is $1$-dimensional. Hence,
$$\tilde H_{r-2}(\hat 0,x) \cong_{G_x} \underbrace{S^{(1)}
\otimes
\dots \otimes S^{(1)}}_r\otimes S^{(n-r)}$$  
Since $\mathfrak S_n$ acts transitively on  each rank row of $\mathcal I_{n,k}$,
we have
\begin{eqnarray*} \wh_r(\mathcal I_{n,k}\cup\{\hat 1\}) &=\,\,\,\,& \tilde H_{r-2}(\hat
0,x)\uparrow_{G_x}^{\mathfrak S_n} \\ &\cong_{\mathfrak S_n}& (S^{(1)})^{\bullet r} \bullet
S^{(n-r)}.
\end{eqnarray*}

It follows from  (\ref{sundeq}) that
\begin{eqnarray*}\tilde H_{k-1}(\bar {\mathcal I}_{n,k}) &\cong_{\mathfrak S_n}&
\bigoplus_{r=0}^k (-1)^{k-r}\,
\,(S^{(1)})^{\bullet r} \bullet S^{(n-r)} \\ &\cong_{\mathfrak S_n}&  (-1)^k S^{(n)}
\,\,\oplus\,\, S^{(1)} \bullet 
\bigoplus_{r=0}^{k-1} (-1)^{k-1-r}\,
\,(S^{(1)})^{\bullet r} \bullet S^{(n-r-1)} \\ &\cong_{\mathfrak S_n}&  (-1)^k S^{(n)}
\,\,\oplus
\,\, S^{(1)} \bullet
 \tilde H_{k-2}(\bar {\mathcal I}_{n-1,k-1}). 
\end{eqnarray*}
This recurrence relation, for $n=k$,  is an equivariant version of the well-known
recurrence relation for derangement numbers: $d_n = (-1)^n + nd_{n-1}$.
From this recurrence relation, one can obtain the following decomposition  of $\tilde
H_{n-1}(\bar {\mathcal I}_{n,n})$ into irreducibles:
\begin{equation}\label{derange}\tilde H_{n-1}(\bar {\mathcal I}_{n,n}) \cong_{\s_n}
\bigoplus_{\lambda
\vdash n} c_\lambda S^\lambda,\end{equation}
where $c_\lambda$ is the number of standard Young tableaux of shape $\lambda$ whose first
descent is even.

\begin{xca} \label{derangedecomp} Use Exercise~\ref{pieri} and the recurrence relation to
prove (\ref{derange}).
\end{xca}

Note that by  taking dimensions on
both sides of (\ref{derange}), and applying the well-known Robinson-Schensted-Knuth
correspondence, one recovers an enumerative result of D\'esarm\'enien
\cite{des}  that the number of derangements in
$\mathfrak S_n$ is equal to the number of  permutations in $\mathfrak S_n$ with first
descent even.  D\'esarm\'enien gives an elegant direct combinatorial proof of this
result.  The Frobenius characteristic of the representation $\bigoplus_{\lambda
\vdash n} c_\lambda S^\lambda$ was used by D\'esarm\'enien and Wachs \cite{deswa} to obtain
deeper enumerative connections between the two classes of permutations. A refinement of
the Reiner-Webb decomposition was given by Hanlon and Hersh \cite{hanher2}. In \cite{rwa},
Reiner and Wachs use (\ref{derange}) to obtain a decomposition (into irreducible
representations of
$\s_n$) of the eigenspaces of the  so called  ``random to top'' operator in card  shuffling
theory.  Type B analogs of these enumerative and representation theoretic results can also
be found in \cite{rwa}.

\end{example}

\begin{xca} \label{equiinjex} Recall from Exercises~\ref{injexer} (c) and (d)  that the top
Betti number  of both 
$\bar {\mathcal N}_{n+1,k}$ and $\bar {\mathcal W}_{n+1,k}$ is $n^k$.
\begin{itemize}
 \item[(a)] (Shareshian and Wachs) Use~(\ref{sundeq}) to prove 
$$\ch (\tilde H_n(\bar {\mathcal N}_{n+1,k}) \downarrow^{\s_{n+1}}_{\s_n}) = 
\sum_{t\ge 1} S(k,t) h_1^t h_{n-t},$$
where $S(n,k)$ is the Stirling number of the second kind.  
\item[(b)]  Conclude that $\tilde H_n(\bar {\mathcal N}_{n+1,k})
\downarrow^{\s_{n+1}}_{\s_n}
$ is the $k$th tensor power of $S^{(n)} \oplus   S^{(n-1,1)}$ by comparing characters. 
 \item[(c)] Show (a) and (b) hold for $\mathcal W_{n+1,k}$.  (Part (b) for  ${\mathcal
N}_{n+1,k}$ and $\mathcal W_{n+1,k}$ was originally  observed by Stanley by means of  the
fixed point M\"obius invariant, see Section~\ref{fixedp}.)
\end{itemize}
\end{xca}

\begin{xca}[Sundaram \cite{su}]  Show  that by applying (\ref{sundeq})  to the Boolean
algebra
$B_n$, one obtains the  well-known symmetric function identity,
$$\sum_{i=0}^n (-1)^i e_i h_{n-i} = 0.$$
\end{xca}

Sundaram  developed her Whitney homology technique in order to study 
representations of the symmetric group on various Cohen-Macaulay subposets of the partition
lattice.  In fact, she applies it to the full partition  lattice and obtains a conceptual
representation theoretic proof of the following classical  result of Stanley.   The
original proof of Stanley used a computation, due to Hanlon \cite{han81}, of the fixed
point M\"obius invariant of the partition lattice.  This   technique is  discussed in
Section~\ref{fixedp}.

\begth[Stanley \cite{st82}] \label{thlie}  For all positive integers $n$,
\bq \label{lie}\tilde H_{n-3} (\bar \Pi_n) ) \cong_{\s_n} e^{2\pi i/n}\uparrow_{\mathfrak
C_n}^{\s_n}
\otimes
\,
\sgn_n
\eq
where $\mathfrak C_n$ is the cyclic subgroup of $\s_n$ generated by $\sigma:=(1,2,\dots,n)$
and
$e^{2\pi i/n}$ denotes the one dimensional representation of $\mathfrak C_n$ whose
character value at
$\sigma$ is $e^{2\pi i/n}$. 
\enth

\begin{proof} (Sundaram \cite{su}).   Theorem~\ref{sundth} is
applied to the dual of the partition lattice.  By setting $T= \Z^+$, $b = r+1$ and $z_i =
1$ in (\ref{plethT}), one obtains 
 \begin{equation} \label{chwh} \sum_{n> r}\ch \wh_r((\Pi_n)^*)=\ch
\tilde H_{r-2}(\bar \Pi_{r+1})
\Big{ [}\sum_{i \ge 1} h_i\Big{]},\end{equation}

Now  (\ref{sundeq}) yields,
\begin{eqnarray*}\sum_{n\ge 1} (-1)^{n-1} \ch\tilde H_{n-2}(\Pi_n \setminus \{\hat 1\})&=&
\sum_{n\ge 1}
\sum_{r=0}^{n-1}(-1)^{r}\ch
\wh_r((\Pi_n)^*)\\ 
&=&  \sum_{r\ge0}(-1)^{r}\ch
\tilde H_{r-2}(\bar\Pi_{r+1})
\Big{ [}\sum_{i \ge 1} h_i\Big{]}.
\end{eqnarray*}
Since $\Delta(\Pi_n \setminus \{\hat 1\})$ is contractible for all $n >1$ and is
$\{\emptyset\}$ when
$n = 1$, it follows that 
$$ h_1 = \sum_{r\ge 0}(-1)^{r}\ch
\tilde H_{r-2}(\bar\Pi_{r+1})
\Big{ [}\sum_{i \ge 1} h_i\Big{]}.$$
Since $h_1$ is the plethystic identity,   \begin{eqnarray}\label{plethinv}\sum_{r\ge
1}(-1)^{r-1}\ch
\tilde H_{r-3}(\bar\Pi_{r}) &=& \Big{(}\sum_{i
\ge 1} h_i\Big{)}^{[-1]} \\ \nonumber &=&\sum_{d \ge 1} {1 \over d} \mu(d) \log(1+p_d)
\end{eqnarray} 
where $[-1]$ denotes plethystic inverse,  $\mu $ is the number theoretic M\"obius function,
and
$p_d$  is the power sum symmetric function.  The last equation follows from a formula of
Cadogan \cite{cad}, which is also derived in \cite{su}. By extracting the degree
$n$ term, we have 
\bq \label{liepi}\ch
\tilde H_{n-3}(\bar\Pi_{n}) = {1\over n} \sum_{d|n}(-1)^{n-n/d}\mu(d) \, p_d^{n/d}  .\eq
A standard formula for the character of  an induced
representation yields
$$\ch \, e^{2\pi i/n}\uparrow_{\mathfrak C_n}^{\s_n} = {1\over n} \sum_{d|n}\mu(d) \,
p_d^{n/d}  . $$
By (\ref{omega}) and  Theorem~\ref{omegasgn} (b), this together with (\ref{liepi}) implies
(\ref{lie}), as desired.
\end{proof}

The representation $e^{2\pi i/n}\uparrow_{\mathfrak C_n}^{\s_n}$ is a well-studied
representation called  the {\em Lie representation} because
it is isomorphic to the representation of the symmetric group on the multilinear component
of the free Lie algebra on $n$ generators, cf., Theorem~\ref{SKJ}.  

\begin{xca}[Sundaram \cite{su}] \label{equidivex} Use the Whitney homology technique to
prove the following results of Calderbank, Hanlon and Robinson \cite{chr}.
\begin{enumerate}\item[(a)]  (equivariant version of Exercise~\ref{modd} (e))
\begin{eqnarray*} \sum_{n\ge 0}(-1)^{n}\ch
\tilde H_{n-2}(\bar\Pi^{1 \bmod d}_{nd+1}) = \Big{(}\sum_{i
\ge 0} h_{id+1}\Big{)}^{[-1]}. \end{eqnarray*}
\item[(b)] (equivariant version of Exercise~\ref{modd} (d))  
\begin{eqnarray} \label{equiddiv}\quad \qquad\sum_{n\ge 1} (-1)^{n-1} \ch\tilde
H_{n-2}(\bar
\Pi^{0\bmod d}_{nd}) = 
\Big{(}\sum_{i \ge 1} h_{i}\Big{)}^{[-1]}\Big{
[}\sum_{i \ge 1} h_{id}\Big{]}.
\end{eqnarray}\end{enumerate}
\end{xca}

\vspace{.1in}
We now  present an equivariant version of Theorem~\ref{betti4}, which extends
Theorem~\ref{sundth} to the nonpure setting. For any  semipure sequentially
Cohen-Macaulay
$G$-poset
$P$ with a bottom element $\hat 0$, define
$r,m$-Whitney homology to be the $G$-module
\bq\label{whrm} \wh_{r,m}(P): = \bigoplus_{\scriptsize\begin{array}{c}x \in P_r \\ m(x) =
m\end{array}}
\tilde H_{r-2}(\hat 0, x),\eq
where $m(x)$ is the length of the longest chain of $P$ containing $x$ and $P_r := \{x \in
P : r(x) = r\}$.

\begth[Wachs \cite{wa2}] \label{semiwhit} Let $P$  be a  semipure sequentially
Cohen-Macaulay $G$-poset. Then for all $m$,
\begin{equation}\label{whitsemipure}\tilde H_{m-1}( P \setminus\{\hat 0\}) \cong_{G} 
\bigoplus_{r=0}^{m}\,(-1)^{m+r}  \,\wh_{r,m}(P).\end{equation} 
\enth

\begin{example}[equivariant version of Example~\ref{kex} \cite{wa2}] \label{equiatleast} 
We
will use Theorem~\ref{semiwhit} and (\ref{plethT}) to obtain  the following formula for
the two parameter generating function for the homology of $ \bar
\Pi^{\ge k}_{n}$,
\bq \label{homk}\quad \qquad\sum_{\scriptsize \begin{array}{c} m\ge 1\\ n \ge k\end{array}}
(-1)^{m-1}
\ch\tilde H_{m-2}(\bar
\Pi^{\ge k}_{n})\,u^n\,t^m = 
 \Big{(}\sum_{i \ge 1} h_{i}  \Big{)}^{[-1]}
\Big{[}\sum_{i \ge k} h_{i}\, u^i\,\, t^{\lfloor {i \over k} \rfloor} \Big{]}.
\eq
We apply
(\ref{whrm}) to the dual of $\Pi^{\ge k}_{n} $, 
$$\wh_{r,m}((\Pi^{\ge k}_{n})^*) = \bigoplus_{\scriptsize{\begin{array}{c}\lambda \in
\mbox{Par}(T,r+1)\\
\lambda
\vdash n\\m(\lambda) = m
\end{array}}}\bigoplus_{x \in \Pi(\lambda)} \tilde H_{r-2}(x,\hat 1),$$
where $T= \{k,k+1, \dots\}$,   $m(\lambda) =
\sum_i
\lfloor
\lambda_i/k\rfloor-1$ (recall (\ref{mk}) here), and the remaining notation is defined in
Example~\ref{exd} . By setting
$z_i = u^i t^{\lfloor i/k
\rfloor}$ in (\ref{plethT}), we obtain
\begin{eqnarray*}\sum_{m,n} \ch \wh_{r,m-1}((\Pi^{\ge k}_{n} )^*)\,u^n\,t^m 
&=&\!\!\!\!\sum_{\lambda \in \mbox{Par}(T,r+1)}\Big{(}\ch
\bigoplus_{x
\in
\Pi(\lambda)}
\tilde H_{r-2}(x,\hat 1)\Big{)} u^{|\lambda|}\,\,t^{\sum_i \lfloor \lambda_i/k\rfloor}\\
&=&
\ch \tilde H_{r-2}(\Pi_{r+1})\Big{[}\sum_{i \ge k} h_{i}\, u^i\,\, t^{\lfloor {i \over k}
\rfloor}
\Big{]}.
\end{eqnarray*} Thus formula (\ref{homk}) follows from (\ref{whitsemipure}) and
(\ref{plethinv}).
\end{example}

The following result generalizes 
 (\ref{equiddiv}) and (\ref{homk}).  Its proof is similar to  that of
(\ref{homk}) described above.

\begth[Wachs \cite{wa2}] \label{genwac}
Suppose   $S\subseteq \{2,3,\dots\}$ is such that  $S$ and $\{s-\min S: s  \in S\}$ are
closed under addition. For $n \in S$, let $\Pi_n^S$ be the subposet of $\Pi_n$
consisting of partitions whose block sizes are in $S$. Then
$$ \quad \qquad\sum_{\scriptsize \begin{array}{c} m\ge 1\\ n \in S\end{array}}
(-1)^{m-1}
\ch\tilde H_{m-2}(\bar
\Pi^{S}_{n})\,u^n\,t^m = 
 \Big{(}\sum_{i \ge 1} h_{i}  \Big{)}^{[-1]}
\Big{[}\sum_{i \in S} h_{i}\, u^i\,\, t^ {\phi(i)} \Big{]},
$$ 
where $\phi(i) := \max\{j \in  \Z^+ : i-(j-1)\min S \in S\}$.
\enth

We remark that the restricted block size partition poset $\Pi_n^S \cup \{\hat 0\}$  of 
Theorem~\ref{genwac} is the intersection lattice of a subspace arrangement, which is
discussed further in Section~\ref{arrangsec}.

\begin{example}[equivariant version of Exercise~\ref{jmodd} \cite{wa2}]   
Theorem~\ref{semiwhit} and  a generalization of (\ref{plethT}) given in \cite{wa2} can be
used to obtain the following formula for the  two parameter generating function for the
homology of the $j \bmod d$ partition poset for $j =2,3,\dots,d $:

\begin{eqnarray*}\sum_{\scriptsize\begin{array}{c} m\ge 1 \\n \ge 0\end{array}}  
(-1)^m& &
\hspace{-.3in}
 \ch\tilde
 H_{m-1}(\bar \Pi_{nd +j}^{j \bmod d} )\,  u^{nd+j} t^{md_0+1} \\
& =&\Big{(}\sum_{i\ge 0} h_{id_0+1}\Big{)}^{[-1]} \Big{[}\sum_{i\ge 0} h_{id+j}
\, u^{id+j} t^{\lfloor \frac i {j_0}\rfloor d_0+1} \Big{]}, \end{eqnarray*}
where $j_0=
\frac j  {\mbox{ gcd}(j,d)}$ and $d_0=\frac  d {\mbox{ gcd}( j,d)}$. 
\end{example}

\begin{prob} Do the  results of this section on restricted block size partition
lattices have nice generalizations to Dowling lattices or intersection lattices of 
Coxeter arrangements?  
\end{prob}

For other restricted block size partition posets with very interesting equivariant 
homology, see the work of Sundaram \cite{su01}.

\section{Bases for the restricted block size partition posets} \label{basissec}

As we have seen in  previous lectures, the
construction of explicit bases for homology and cohomology is an effective tool
in studying group representations on homology. 
In this section we construct bases for the homology and cohomology of the restricted block
size partition posets studied in the previous section.   These bases are used to 
obtain further results on  the
representations of the symmetric group on homology of the posets and to relate these
representations to representations on  homology of certain interesting graph
complexes.

\subsection{The $d$-divisible partition lattice $\Pi_{nd}^{0
\bmod d}$} \label{ddivsec} 
\begin{figure}\begin{center}
\includegraphics[width=9 cm]{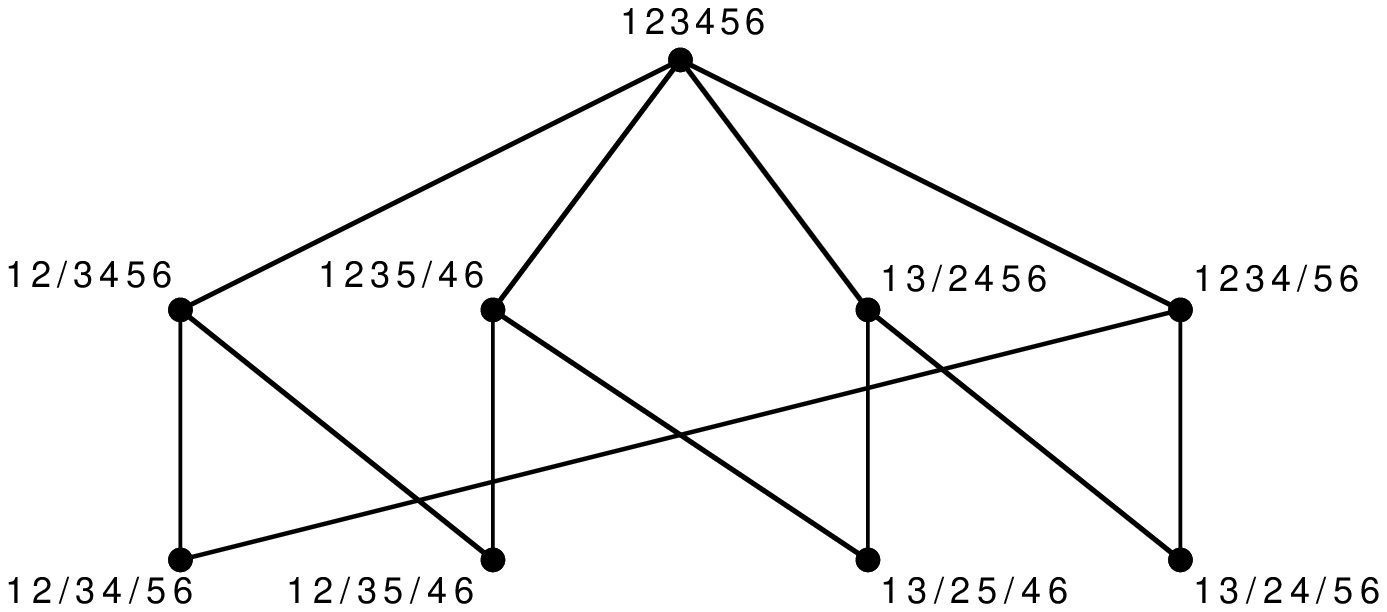}\\
\end{center}
\begin{center}\caption{ $\Pi^2_{123456}$ }\label{figposet2}
\end{center}
\end{figure}

We construct  analogs of the splitting basis for the homology
of the partition lattice given Section~\ref{ELpart} and its dual basis for cohomology.  To
{\em switch-split} a permutation
$\sigma$ at position
$j$ is to form the partition $$\sigma(1), \dots, \sigma(j-1),\sigma(j+1) \,\,\,|\,\,\,
\sigma(j),\sigma(j+2)
\dots,
\sigma(n)$$ of $[n]$. Switch-splitting a permutation at two or more
nonadjacent positions is defined similarly.  Let $d \ge 2$.  For each
$\sigma
\in
\mathfrak S_{nd}$, let $\Pi_\sigma^d$ be the induced subposet of $\Pi_{nd}^{0
\bmod d}$ consisting of partitions obtained by splitting or switch-splitting
$\sigma$ at any number of  positions in $\{d,2d, \dots, nd\}$. The subposet
$\Pi^2_{123456}$ of
$\Pi^{0 \bmod 2}_6$ is shown in  Figure~\ref{figposet2}.

Each poset $\Pi^d_\sigma \cup \{\hat 0\}$
is isomorphic to the face  lattice $C_{n-1}$ of the $(n-1)$-cross-polytope.  Therefore
$\Delta(\bar
\Pi^d_\sigma)$ is an $(n-2)$-sphere 
embedded
in 
$\Delta(\bar
\Pi^{0\bmod d}_{nd})$, and hence it determines a fundamental cycle $\rho^d_\sigma \in
\tilde H_{n-2}(\bar
\Pi^{0\bmod d}_{nd};\Z)$. 
In each poset $\bar \Pi^d_\sigma$, we select a distinguished maximal chain $c^d_\sigma$
whose $k$ block partition  is  obtained
by splitting $\sigma$ at positions $d, 2d,\dots,(k-1)d$, for $k = 2, 3, \dots, n$.   For
example,
$$c^2_{123456} = (12/34/56 < 12/3456).$$ 

 Let
$$A_n^d := \{\sigma \in \s_{nd} : \mbox{des}(\sigma) = \{d,2d,\dots, (n-1)d\}, \,\,
\sigma(nd) = nd\}.$$  Recall that in Exercise~\ref{modd} (d), it is stated that $\bar
\Pi_{nd}^{0\bmod d}$ has the homotopy type of a wedge of $|A_n^d|$ spheres of dimension
$n-2$. 

\begth[Wachs \cite{wa96}] \label{dsplit} The set $\{\rho^d_\sigma : \sigma \in A_n^d\}$
forms a basis for 
$\tilde H_{n-2}(\bar
\Pi^{0\bmod d}_{nd};\Z)$ and the set $\{c^d_\sigma : \sigma \in A_n^d\}$ forms a basis for 
$\tilde H^{n-2}(\bar
\Pi^{0\bmod d}_{nd};\Z)$.
\enth

The theorem is proved by first showing that for all $\alpha, \beta \in A_n^d$,
if $c_\alpha^d \in \Pi^d_\beta $ then $\alpha\le\beta$ in lexicographic order.   This is
used to establish linear independence of both $\{\rho^d_\sigma : \sigma
\in A_n^d\}$ and $\{c^d_\sigma : \sigma \in A_n^d\}$.  The result then follows from
Exercise~\ref{modd}~(d).

The ``splitting basis'' given in Theorem~\ref{dsplit} is used in \cite{wa96} to give a
combinatorial proof of the following result of Calderbank Hanlon and Robinson, which was
first conjectured by Stanley.  

\begth[Calderbank, Hanlon and Robinson \cite{chr}] \label{chrth} Let $
H_{n,d}$ be the skew hook of size $nd-1$ and  descent set $\{d,2d, \dots, (n-1)d\}$
(cf. Section~\ref{secrepsym}).  Then
$$\tilde H_{n-2}(\bar
\Pi^{0\bmod d}_{nd})\downarrow_{\s_{nd-1}}^{\s_{nd}} \cong_{\s_{nd-1}} S^{H_{n,d}}$$
\enth

Calderbank, Hanlon, Robinson obtain this result as a consequence of  (\ref{equiddiv}).   By
mapping  tableaux of  shape $H_{n,d}$ to maximal chains of $\bar
\Pi^{0\bmod d}_{nd}$, one gets a
combinatorial proof.  Indeed,  for each tableaux $T$ of shape $H_{n,d}$, let $\sigma_T$ be
the permutation obtained by reading  the entries of the skew hook tableaux
$T$ from the southwest end to the northeast end of $T$ and then attaching $nd$.  Now define
$$c_T:= c^d_{\sigma_T} \,\,\mbox{ and }\,\,\rho_T:= \rho^d_{\sigma_T}.$$

\begth[Wachs \cite{wa96}] The map $T \mapsto c_T$ induces a well-defined
$\s_{nd-1}$-isomorphism from  the skew hook Specht module
$S^{H_{n,d}}$ to  $\tilde H^{n-2}(\bar
\Pi^{0\bmod d}_{nd})\downarrow_{\s_{nd-1}}^{\s_{nd}}$.  
\enth

To prove this, one first observes that the row permutations leave $c_T$ invariant;
then one shows that the Garnir relations map to cohomology relations.

There is a dual version of  polytabloid defined for each tableaux $T$ by
$$e_T^*:= \sum_{\alpha \in R_\lambda} \sum_{\beta \in C_\lambda}\sgn(\beta)\, \,
T\beta\alpha.$$ (This is actually closer to the traditional notion of polytabloid than
the one we gave in Section~\ref{secrepsym}.)  For each skew or straight shape $\lambda
$, it is known that 
$$\langle e_T^*: T
\in
\mathcal T_\lambda\rangle  \cong_{\mathfrak S_n} S^\lambda.$$ 

\begth[Wachs \cite{wa96}] The map $e^*_T \mapsto \rho_T$ induces a well-defined
$\s_{nd-1}$-isomorphism from  $S^{H_{n,d}}$ to  $\tilde H_{n-2}(\bar
\Pi^{0\bmod d}_{nd})\downarrow_{\s_{nd-1}}^{\s_{nd}}$.  
\enth

All the results of this subsection were generalized to the restricted block size partition
posets
$\Pi_n^S$ of Theorem~\ref{genwac} by Browdy and Wachs \cite{br,brwa}. The nonpurity of
$\Pi_n^S$ in the general case significantly increases the complexity of the results.  The
``at least
$k$'' partition poset $\Pi_n^{\ge k}$ is an example of such a nonpure poset.

\subsection{The $1 \bmod d$ partition lattice $\Pi_{nd+1}^{1
\bmod d}$} \label{moddsec}

First we describe a basis for top cohomology of  $\bar \Pi_{nd+1}^{1
\bmod d}$, due to  Hanlon and Wachs \cite{hw}, which generalizes the 
decreasing chain basis $\{\bar c_\sigma : \sigma \in \mathfrak S_n,
\sigma(n) = n\}$ for cohomology of $\bar \Pi_n$ given in Section \ref{ELpart}.
This basis is used in \cite{hw} to prove a generalization of
Theorem~\ref{SKJ}, which relates the $\mathfrak S_{nd+1}$-module 
$\tilde H^{n-2}(\bar \Pi_{nd+1}^{1
\bmod d})$ to a $(d+1)$-ary version   of the free Lie algebra.  Then we describe a basis
for top homology found about ten  years later by  Shareshian and Wachs \cite{shwa1}, which
generalizes the tree-splitting basis for homology of $\bar \Pi_n$ given in Section
\ref{ELpart}.  This basis is used in \cite{shwa1} to relate the $\mathfrak
S_{nd+1}$-module 
$\tilde H^{n-2}(\bar \Pi_{nd+1}^{1
\bmod d})$ to the homology of   graph complexes studied by Linusson Shareshian
and Welker
\cite{lsw} and Jonsson \cite{jon05}. 

\begin{figure}\begin{center}
\includegraphics[width=7 cm]{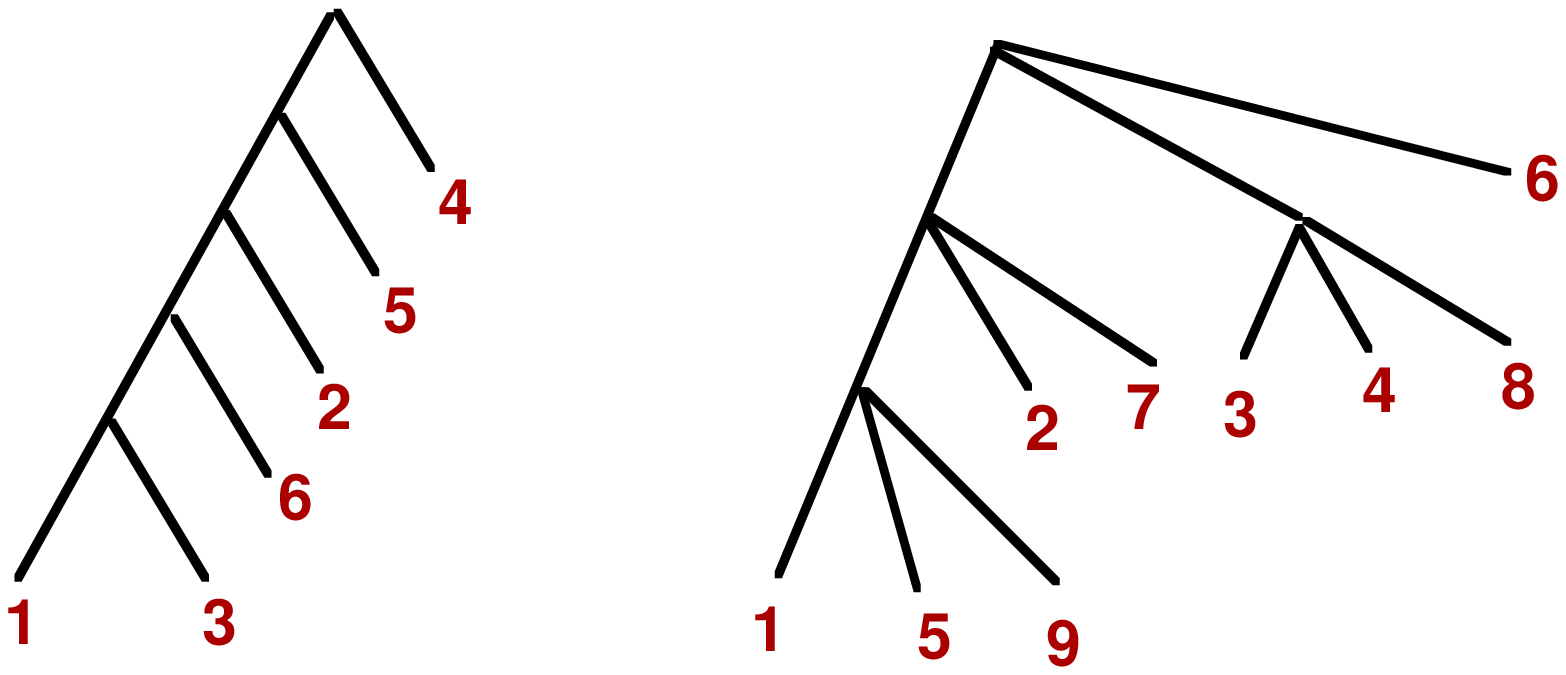}
\\
\hspace{-.6in}(a) \hspace{1.8in}(b)
\end{center}
\begin{center}\caption{}\label{figbrush}
\end{center}
\end{figure}

Let $\mathcal T_{nd+1}^{d}$ be the set of rooted planar $(d+1)$-ary  trees on leaf set
$[nd+1]$ (i.e., rooted trees in which each internal node has exactly $d+1$ children that
are  ordered from left to right).  For any node
$x$ of 
$T \in \mathcal T_{nd+1}^{d}$,  let
$m(x) $ be the smallest leaf in the tree rooted at $x$.  
A tree $T$ in $\mathcal T_{nd+1}^d$ is said to be a $d$-brush if for each node y of $T$,
the $m$-values of the children of $y$ increase from left to right, and the child with the
largest $m$-value is a leaf.  An example of a $1$-brush is given in Figure~\ref{figbrush}
(a) and of a $2$-brush is given in Figure~\ref{figbrush} (b).  Note that  every $1$-brush 
looks like  a comb, which is the reason for the terminology    ``brush".

\vspace{.1in}\begin{xca}\begin{itemize}
\item[]\item[(a)] Show that the number of $1$-brushes on leaf set $[n+1]$ is $n!$.
\item[(b)] Show that the number of $2$-brushes on leaf set $[2n+1]$ is $(2n-1)!!^2$. 
\end{itemize}
\end{xca}

Recall from Section~\ref{liesec} that the postorder traversal of  a binary tree on leaf set
$[n]$ yields a maximal chain of $\bar \Pi_n$.  Now we consider a more general
construction, which associates a maximal chain $c_T$ of
$\bar \Pi_{nd+1}^{1
\bmod d}$  to each tree
$T$ in
$\mathcal T_{nd+1}^d$. Each internal node $y$ of $T$ corresponds to a  merge of $d+1$
blocks that are the leaf sets of the trees rooted at the $d+1$ children of $y$. 
Postorder traversal of internal nodes of
$T$ yields a sequence of merges, which  corresponds to a maximal chain 
$c_T$ of
$\bar \Pi_{nd+1}^{1
\bmod d}$.  For example if $T$ is the tree of Figure~\ref{figbrush}
(b) then 
$$c_T = 159/2/7/3/4/8/6< \!\!\!\cdot \,\, 15927/3/4/8/6 <
\!\!\!\cdot
\,\, 15927/348/6 $$

\begth[Hanlon and Wachs \cite{hw}] \label{brushbasis} Let
$\mathcal B_{nd+1}^d$ be the set of $d$-brushes in 
$\mathcal T_{nd+1}^d$.  The set $\{c_T : T \in \mathcal
B_{nd+1}^d\}$ forms a basis for $\tilde H^{n-2}(\bar \Pi_{nd+1}^{1
\bmod d};\Z)$.
\enth

\begin{xca}[Hanlon and Wachs \cite{hw}] Prove Theorem~\ref{brushbasis} by showing
\begin{itemize} 
\item[(a)]  The set
$\{c_T : T \in \mathcal B_{nd+1}^d\}$ spans $\tilde H^{n-2}(\bar \Pi_{nd+1}^{1
\bmod d};\Z)$. 
\item[(b)]  $|\mathcal B_{nd+1}^d| = |\mu(\Pi_{nd+1}^{1
\bmod d})|$.
\end{itemize}
\end{xca}

We now construct the Shareshian-Wachs  basis for   homology of $\bar
\Pi_{nd+1}^{1
\bmod d}$.   Just like the tree-splitting basis, which it generalizes, it consists of 
fundamental cycles of Boolean algebras embedded in $\Pi_{nd+1}^{1
\bmod d}$.  

\begin{figure}\begin{center}
\includegraphics[width=8 cm]{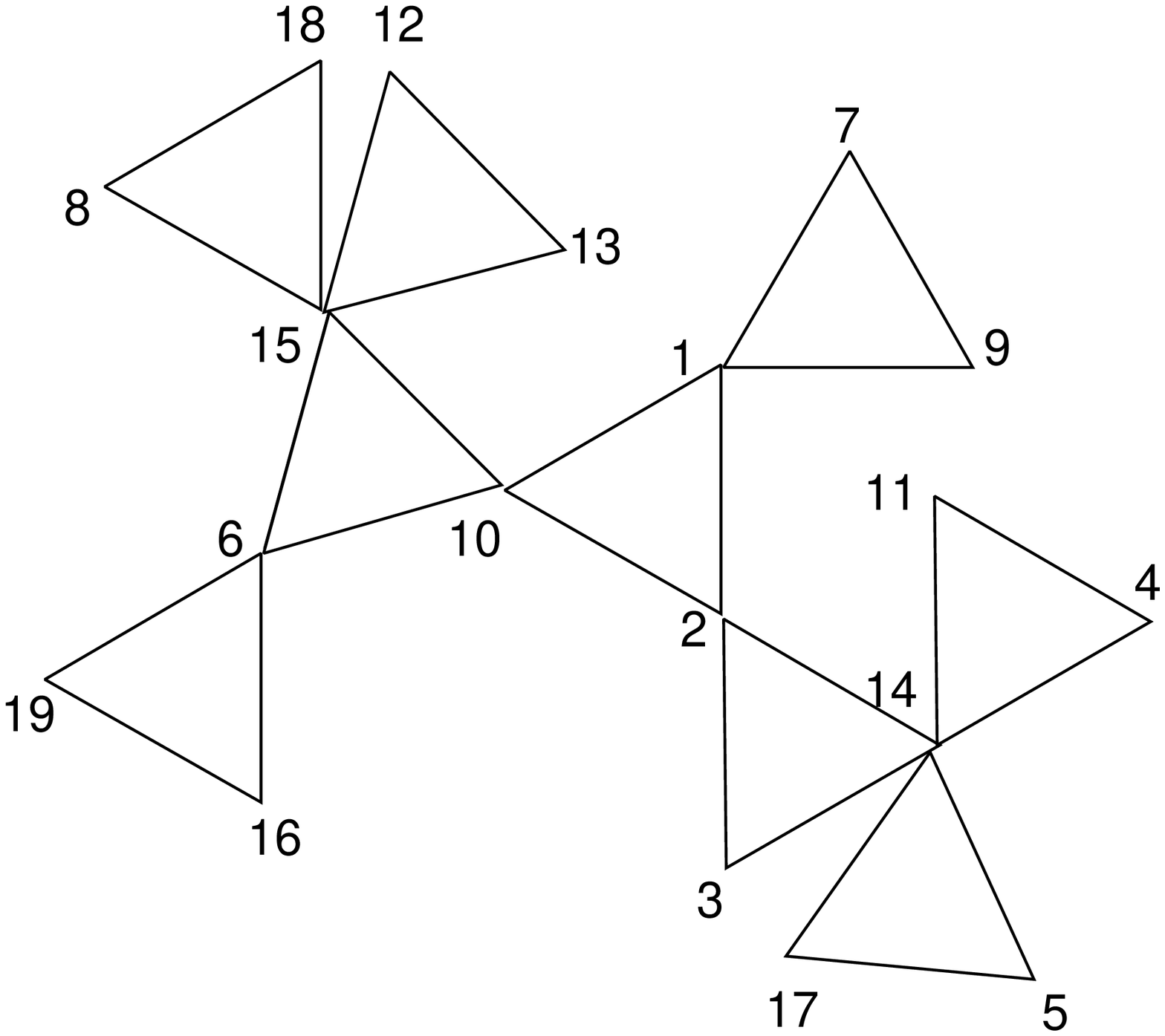}
\end{center}
\begin{center}\caption{}\label{figtree1}
\end{center}
\end{figure}

A connected graph $G$  is said to be a {\em $d$-clique tree} if either
$G$ consists of a single node, or $G$ contains a $d$-clique, the removal of whose edges
results in
a graph with $d$ connected components that
are all
$d$-clique trees.   Note that each edge of a $d$-clique tree is in a unique $d$-clique and
the removal of the edges of {\em any} $d$-clique from a $d$-clique tree results in
a graph with $d$ connected components that
are all
$d$-clique trees.
 Note also that a
$2$-clique tree is an ordinary tree.   An example of a $3$-clique tree, which we
refer to as a triangle tree, is given in Figure~\ref{figtree1}.

\begin{figure}\begin{center}
\includegraphics[width=8 cm]{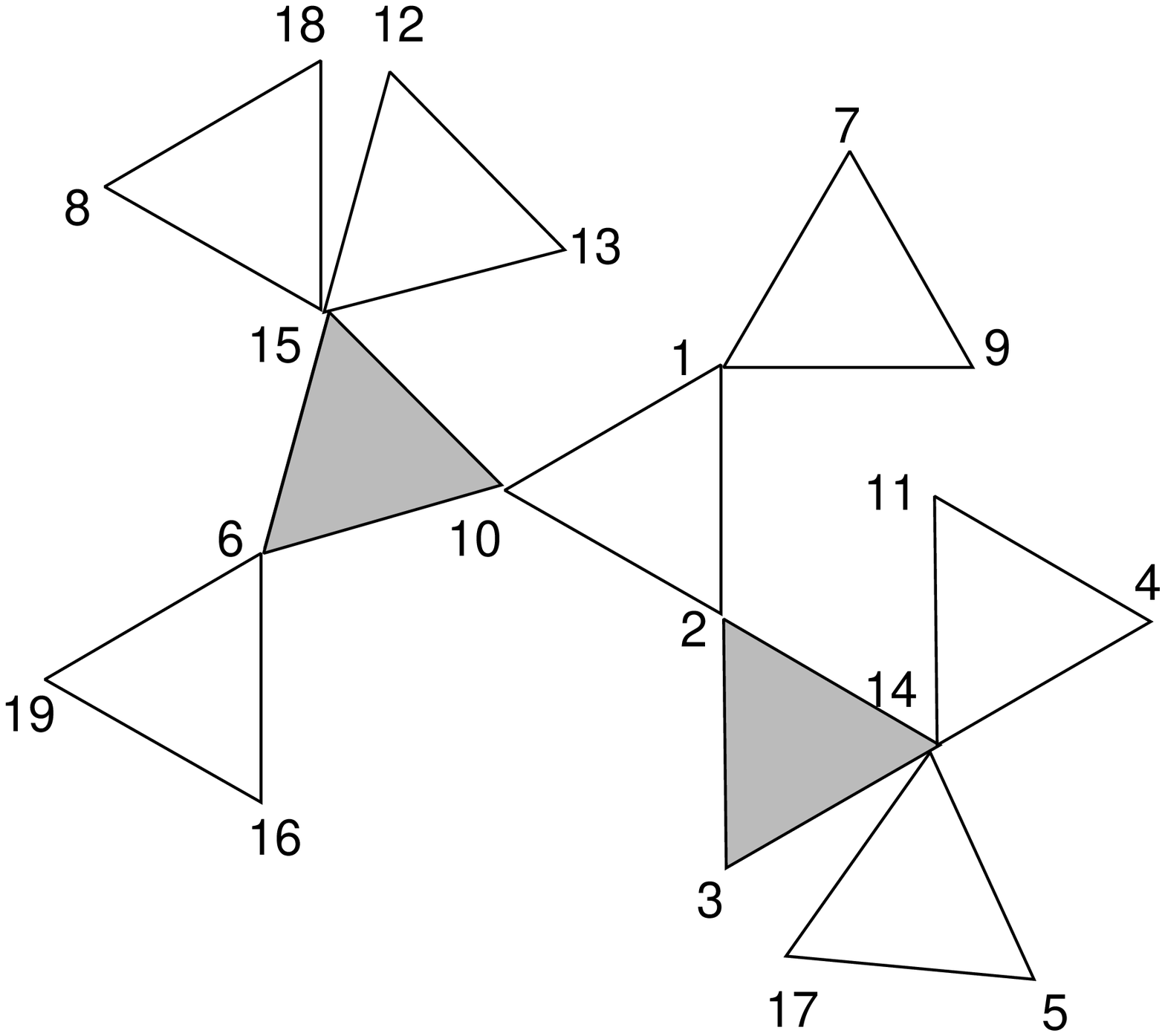}
\end{center}
\begin{center}\caption{}\label{figshtree2}
\end{center}
\end{figure}

Given a $(d+1)$-clique tree $T$ on node set $[nd+1]$, one obtains a partition in
$\Pi_{nd+1}^{1
\bmod d}$ by  choosing any set of $(d+1)$-cliques of $T$ and removing
the edges of each clique in the set.  The blocks of the partition are the node sets of the
connected components of the resulting graph.  We say that the partition is obtained by {\em
splitting} the
$(d+1)$-clique tree $T$ at the chosen set of $(d+1)$-cliques.  For example, the partition
obtained by splitting the triangle tree  at  the shaded triangles in
Figure~\ref{figshtree2} is
$$ 19,16,6 \,\,/ \,\,   8,15,18,12,13 \,\, / \,\,  10,2,1,9,7 \,\, / \,\, 3 \,\, / \,\,
17,5,4,11,14.$$  Now let
$\Pi_T$ be the subposet of
$\Pi_{nd+1}^{1 \bmod d}$ consisting of partitions obtained by splitting $T$. Clearly
$\Pi_T$ is isomorphic to the subset lattice $B_n$.  Therefore $\Delta(\bar \Pi_T)$ is an
$(n-2)$-sphere which determines a fundamental cycle $\rho_T$.  

\begth[Shareshian and Wachs \cite{shwa1}] \label{clique} The set of fundamental cycles
$\rho_T$ such that
$T$ is a
$(d+1)$-clique  tree on node set $[nd+1]$,  spans $\tilde H_{n-2}(\bar \Pi_{nd+1}^{1
\bmod d};\Z)$.
\enth

This is proved by identifying a set $S$ of $(d+1)$-clique  trees, called  
increasing
$(d+1)$-clique  trees, and establishing a bijection and  unitriangular relationship between
$(d+1)$-clique trees and $d$-brushes.  This shows  that 
$\{\rho_T: T \in S\}$ is  a basis for  $\tilde H_{n-2}(\bar \Pi_{nd+1}^{1
\bmod d};\Z)$.    The increasing $2$-clique trees are the  increasing trees
discussed in Section~\ref{ELpart} and the homology basis is the tree splitting basis.

Next we discuss an
application of Theorem~\ref{clique} that led to the discovery of the
clique tree splitting basis in the first place.  Let ${\rm NPM}_{2n}$ be the poset of
nonempty graphs on node set $[2n]$ that don't contain a perfect matching (i.e., 
a subgraph in which each of the $2n$ vertices has degree 1), ordered by inclusion of edge
sets.  Linusson, Shareshian, and Welker
\cite{lsw} show, using  discrete Morse theory,
that
${\rm NPM}_{2n}$  has the homotopy type of a wedge of $(2n-1)!!^2$ spheres of dimension
$3n-4$.   It  was in this work that  the increasing triangle trees first arose as the
critical elements of a Morse matching on 
${\rm NPM}_{2n}$ (increasing triangle trees are just called trees of triangles in
\cite{lsw}).   See
the chapter by Forman \cite{for05} in this volume to learn about discrete Morse
theory and critical elements.  Recall that the property of containing (or not containing) a
perfect matching is an example of a monotone graph property (see Section~\ref{grth}). 
When we discuss Alexander duality in the next lecture, we will see how the  
homology of the poset of graphs that have a monotone graph property is related to the
homology of the poset of graphs that don't have the property. 

Theorem~\ref{clique} and discrete Morse theory  play 
essential roles in the proof of the following equivariant version of the
Linusson-Shareshian-Welker result.  Indeed, discrete Morse theory is used to show
that there is a  sign twisted $\mathfrak S_{2n-1}$-isomorphism between the homology of
$\mbox{NPM}_{2n}$ and the cohomology of another poset called the {\em factor critical
graph poset}.  Discrete  Morse theory is also used to show that the cohomology of the
factor critical graph poset is  generated by certain maximal chains naturally indexed by
triangle trees.   This gives a natural map from the triangle tree generators $\rho_T$ of
homology of $\Pi_{2n-1}^{1
\bmod 2}$ to the maximal chains  of the factor
critical graph poset that are indexed by the triangle trees.  Relations on the triangle
tree generators of $\tilde H_{n-3}(\bar
\Pi_{2n-1}^{1
\bmod 2}) $ are derived, which are shown to map to coboundary relations in
the factor critical graph poset.

\begth[Shareshian and Wachs \cite{shwa1}] \label{perth}For all $n \ge 1$,
$$\tilde H_{3n-4}({\rm NPM}_{2n})\downarrow_{ \mathfrak 
S_{2n-1}}^{\mathfrak
S_{2n}} \cong_{\mathfrak S_{2n-1}} \tilde H_{n-3}(\bar
\Pi_{2n-1}^{1
\bmod 2})   \otimes\,\, \sgn_{2n-1}.$$
\enth

Another graph poset, recently studied by Jonsson \cite{jon05},  curiously has the same
bottom nonvanishing  homology  as
 ${\rm NPM}_{2n}$ and $\bar
\Pi_{2n-1}^{1
\bmod 2}$ (up to sign twists).  A
graph  is said to be  $d$-{\em edge-connected} if removal of any set of at most
$d-1$ edges leaves the graph connected.   So
a
$1$-edge-connected graph is just a connected graph.  The $d$-clique trees are examples of
graphs that are
$(d-1)$-edge-connected but not
$d$-edge-connected.  In fact, the $d$-clique trees are minimal elements of the poset of
$(d-1)$-edge-connected graphs.  Jonsson studied the integral homology of the poset
of graphs that are not $2$-edge-connected.  Let
${\rm NEC}^d_n$ be the  poset of nonempty graphs on node set
$[n]$ that are not $d$-edge-connected.  Jonsson discovered that
${\rm NEC}^2_{2n+1}$  has homology in multiple dimensions and that the bottom nonvanishing
reduced integral homology is in dimension $3n-2$ and is of rank
$(2n-1)!!^2$.   His proof  involves discrete Morse theory with the
increasing triangle trees as the critical elements. Shareshian and Wachs \cite{shwa1}
obtain the following equivariant version of  Jonsson's result by giving 
presentations of the two homology groups in terms of 
triangle-trees. 
 
\begth[Shareshian and Wachs \cite{shwa1}] \label{nec} For all $n$,
$$  \tilde H_{3 n- 2}({\rm
NEC}^2_{2n+1})  \cong_{\mathfrak S_{2n+1}} \tilde H_{n-2}(\bar
\Pi_{2n+1}^{1
\bmod 2}) \otimes \sgn_{2n+1}.$$
\enth

\begin{con}[Shareshian and Wachs \cite{shwa1}] For all $n,d \ge 1$, 
$$  \tilde H_{{d+1 \choose 2} n- 2}({\rm
NEC}^d_{dn+1})  \cong_{\mathfrak S_{dn+1}} \tilde H_{n-2}(\bar
\Pi_{dn+1}^{1
\bmod d})\otimes \sgn_{dn+1}^{\otimes d+1},$$ and 
$$  \tilde H_{i}({\rm
NEC}^d_{dn+1})  =  0$$
if $ i < {d+1 \choose 2} n- 2$.
\end{con}

The conjecture is true for $d=1,2$. Indeed, Jonsson's homology result and
Theorem~\ref{nec}  comprise the $d=2$ case.    The conjecture for $d=1$ says that  that
the poset of disconnected graphs on node set
$[n+1]$ has  homology 
isomorphic, as an
$\mathfrak S_{n+1}$-module, to that of the partition lattice.   
A proof of this well-known result is discussed in
Example~\ref{discex}.

There is another interesting poset with the same homotopy type and $\s_{nd+1}$-equivariant
homology  as that of
$\bar
\Pi_{nd+1}^{1
\bmod d}$,  worth mentioning here.  This poset is the proper part of the poset
$\mathcal T^{2 \bmod d}_{nd+2}$ of homeomorphically irreducible trees on leaf set
$[nd+2]$ in which each internal node has degree congruent to $2 \bmod d$.  By
{\em homeomorphically irreducible}  we mean  nonrooted and no  node has
degree 2.  The order relation is as follows: $T_1 < T_2$ if $T_1$ can be obtained from
$T_2$ by contracting internal edges.  So the bottom element of $\mathcal T^{2 \bmod
d}_{nd+2}$ is the  star tree (the tree with only one internal node), and the maximal
elements are trees in which each internal node has degree exactly $d+2$.  The $d=1$ case
of the tree poset $\mathcal T^{2 \bmod d}_{nd+2}$ has arisen in 
various areas such as  algebraic geometry, homotopy theory, geometric group theory,
 mathematical physics and mathematical biology; see eg., 
\cite{boa,vo90, robwh96, wh, rob, ak, bhv, read} and the references contained therein. 
Vogtmann
\cite{vo90} showed that the tree poset  in the
$d=1$ case is homotopy Cohen-Macaulay.  
 There
is a natural  action of
$\mathfrak S_{nd+2}$ on 
${
\mathcal T}^{2
\bmod d}_{nd+2}$ whose
representation on
 the homology of ${ \mathcal T}^{2 \bmod
d}_{nd+2} -\{\hat 0\}$ was computed by Robinson and  Whitehouse \cite{robwh96} in the $d=1$
case.  Hanlon
\cite{ha96} introduced the general tree poset
$\mathcal T^{2
\bmod d}_{nd+2}$, proved that it is Cohen-Macaulay, and   generalized the
Robinson-Whitehouse result. 

\begth[Robinson and Whitehouse, $d=1$ \cite{robwh96}, Hanlon \cite{ha96}] \label{lift}  For
all
$d$, $n
\ge 1$,
$$\tilde H_{n-2}(\bar{ \mathcal T}^{2 \bmod
d}_{nd+2}) \,\,\,\,\cong_{\s_{nd+2}}\,\,\,\, \tilde
H_{n-2}(\bar \Pi_{nd+1}^{1 \bmod d})\uparrow^{\s_{nd+2}}_{\s_{nd+1}}- \tilde
H_{n-1}(\bar \Pi_{nd+2}^{1 \bmod d}).$$
\enth

The $d=1$ case of the $\s_{nd+2}$-module given in Theorem~\ref{lift} has come to be known
as the Whitehouse module.  It has occurred in a variety of diverse contexts such as
homotopy theory \cite{robwh96,wh}, cyclic Lie operads      \cite{gk,kon}, homology of 
partition posets
\cite{su99, su01}, knot theory and graph complexes \cite{bblsw}, and hyperplane
arrangements and Lie algebra homology
\cite{hast}.  

From Theorem~\ref{lift}, one can show  that (see \cite{su012}),
\begin{equation}\label{restricttree}\tilde H_{n-2}(\bar{ \mathcal T}^{2 \bmod
d}_{nd+2})\downarrow^{\s_{nd+2}}_{\s_{nd+1}} \,\,\,\,\cong_{\s_{nd+1}}\,\,\, \tilde
H_{n-2}(\bar \Pi_{nd+1}^{1 \bmod d}).\end{equation}  A direct  combinatorial proof of
this is discussed in the following exercise.   

\begin{xca} \label{shelltree} By removing the leaf $nd+2$ from any tree $T$ in $ {
\mathcal T}^{2 \bmod d}_{nd+2}$ and designating the internal node that  had been adjacent
to the leaf $nd+2$ as the root, one turns $T$ into a rooted  nonplanar tree on leaf set
$[nd+1]$  in which each internal node has
$ id+1$ children, for some $i \ge 1$.  The maximal elements of ${ \mathcal T}^{2 \bmod
d}_{nd+2}$ are now rooted  nonplanar $(d+1)$-ary trees on leaf set $[nd+1]$.
\begin{itemize}
\item[(a)]  Show that any linear extension of the following partial ordering
of rooted nonplanar
$(d+1)$-ary trees on leaf set
$[nd+1]$ is a recursive coatom ordering of ${ \mathcal T}^{2 \bmod
d}_{nd+2} \cup \{\hat 1\}$.
Let $T_1 \land T_2 \land \dots \land T_{d+1}$ denote the $(d+1)$-tree
in which $T_1,T_2,\dots,T_{d+1}$ are the subtrees of the root.  The partial order is
defined to be the transitive closure of the relation 
given by $T < T^\prime$ if  $T^\prime$ can be  obtained from $T$ by replacing some subtree
$$(T_1 \land \dots \land T_{d+1}) \land T_{d+2}\land \dots \land T_{2d+1} $$  with the
subtree $$(T_1\land \dots \land T_{d}\land T_{d+2} )\land  T_{d+1} \land  T_{d+3}\land
T_{2d+1}
$$ where the  minimum leaf  of $T_{d+2}$ is less than the 
minimum leaf  of $ T_{d+1}$.
\item[(b)]   The poset $ { \mathcal T}^{2 \bmod
d}_{nd+2}$ can be viewed as the face poset of a simplicial complex whose
facets correspond to rooted nonplanar $(d+1)$-ary trees on leaf set $[nd+1]$.  Hence the
recursive coatom ordering of (a) is simply a shelling of the simplicial complex.  Show
that the homology facets of the shelling correspond to the $d$-brushes.  Consequently, 
\begin{equation}\label{topytree} \bar{ \mathcal
T}^{2
\bmod d}_{nd+2}  \simeq \bar \Pi_{nd+1}^{1\bmod d}.\end{equation}
\item[(c)] Prove (\ref{restricttree}) by first observing  that the set   of  the rooted
planar
$(d+1)$-ary trees on leaf set $[nd+1]$ indexes  respective  sets of maximal
chains of
$\bar \Pi_{nd+1}^{1
\bmod d}$ and $\bar{ \mathcal
T}^{2
\bmod d}_{nd+2}$ that generate top cohomology,    and then showing
that the two cohomology groups  have the same presentation in terms of rooted
planar
$(d+1)$-ary trees on leaf set $[nd+1]$. 
\end{itemize}
\end{xca}

Trappmann and Ziegler \cite{tz}  established  shellability of ${
\mathcal T}^{2
\bmod d}_{nd+2}$ in a different way from  Exercise~\ref{shelltree} (a), but  with the same
homology facets as in Exercise~\ref{shelltree} (b).
 Ardila and Klivans  \cite{ak} and  Robinson \cite{rob}
have recently independently  proved a  result stronger than the
 homotopy result (\ref{topytree}) and the homology result (\ref{restricttree}) in the $d=1$
case; namely that  
$\bar {\mathcal T}_{n+1}^{2
\bmod 1} $ and
$\bar
\Pi_n$ are $\s_n$-homeomorphic. Robinson's homeomorphism for the $d=1$ case 
restricts to a homeomorphism for the general case.  Hence    $\bar
{\mathcal T}_{nd+2}^{2
\bmod d} $ and
$\bar
\Pi_{nd+1}^{1\bmod d}$  are $\s_{nd+1}$-homeomorphic for all $n,d \ge 1$.

\begin{prob} It is well-known that for each permutation in $\s_n$, there is a distinct copy
of the face poset of the
$(n-1)$-dimensional associahedron embedded in $
{\mathcal T}_{n+2}^{2
\bmod 1} $ and that by taking   fundamental cycles,  one obtains  a basis for the
homology of
$\bar {\mathcal T}_{n+2}^{2
\bmod 1} $.  The associahedron is discussed in the chapter of Fomin and Reading in this
volumn
\cite{fr04}.  Is there a nice basis for homology of
$\bar {\mathcal T}_{nd+2}^{2
\bmod d} $ for general $d$, consisting of fundamental cycles of embedded face
posets of polytopes, which are copies of some sort of  $d$-analog of the associahedron?
\end{prob}

By Theorem~\ref{perth} and the isomorphism (\ref{restricttree}),
 the restriction  of the $\s_{2n}$-equivariant  homology of ${\rm NPM}_{2n}$ and  $\bar
{\mathcal T}_{2n}^{2 
\bmod 2}$ to
$\mathfrak S_{2n-1}$ are isomorphic (up to tensoring with the sign representation). For the
unrestricted homology modules, we have the following conjecture.

\begin{con}[Linusson, Shareshian and Welker \cite{lsw}]  For all $n \ge 1$,
$$\tilde H_{3n-4}({\rm NPM}_{2n}) \cong_{\mathfrak S_{2n}} \tilde H_{n-4}(\bar
{\mathcal T}_{2n}^{2 
\bmod 2}) \otimes \sgn_{2n} .$$ 
\end{con}

\section{Fixed point M\"obius invariant} \label{fixedp}  We conclude this lecture with a
very brief discussion of another commonly used recursive technique  for computing homology
representations.   Let $P$ be a
$G$-poset and for each $g
\in G$, let $P^g$ denote the induced subposet of
$P$ consisting of elements that are fixed by $g$.  Let $\chi^P$  be the character
of the virtual representation of $G$ on $\oplus_i (-1)^i\tilde H_i(P)$.  The following
consequence of the Hopf-Lefschetz fixed point theorem   first appeared in the
combinatorics literature in the work of Baclawski and  Bj\"orner
\cite{bacbj}.  

\begth \label{fixmob}  For any $G$-poset $P$ and $g \in G$,
$$\chi^P(g) = \mu(\hat P^g).$$  
\enth

There are many  interesting applications  of this theorem in
the literature.  We give a simple example here.  \begin{example}[Stanley,  see \cite{bj90}]
Consider the normal word poset
$\mathcal N_{n,k}$ under the  action of $\s_n$.  Since $P:=\bar {\mathcal N}_{n,k}$ is
Cohen-Macaulay,
$(-1)^{k-1}\chi^P $ is the character of the representation of $\s_n$ on $\tilde
H_{k-1}(\bar {\mathcal N}_{n,k})$.  The words that are fixed by $g$ are the words whose
letters are fixed points of $g$.  So $P^g$ is the poset of normal nonempty words of length
at most
$k$ on alphabet
$F(g)$, where $F(g)$ is the set of fixed points of $g$. By Exercise~\ref{injexer} (c),
$\mu(\hat P^g) = (-1)^{k-1} (|F(g)|-1)^k$.  Since the irreducible representation
$S^{(n-1,1)}$ has the character $(|F(g)|-1)$, its $k$th tensor power has the character
$(|F(g)|-1)^k$.  By Theorem~\ref{fixmob}, 
$$\tilde
H_{k-1}(\bar {\mathcal N}_{n,k}) \cong_{\s_n} (S^{(n-1,1)})^{\otimes k}.$$
Since the identical argument works for $\mathcal W_{n,k}$,
$$\tilde
H_{k-1}(\bar {\mathcal W}_{n,k}) \cong_{\s_n} (S^{(n-1,1)})^{\otimes k}.$$
\end{example}

\begin{xca} Show that if $\chi$ is the character of the representation of $\s_n$ on the
homology of the injective word poset $\bar{\mathcal I}_{n,n}$ then
$$\chi(g) = d_{|F(g)|},$$ where $d_j$ is the number of derangements of $j$ letters.
\end{xca}

%% file: lect5.tex

\lecture{Poset operations and maps} \label{opsmaps}

\section{Operations: Alexander duality and direct product} \label{prodsec}

In this section we consider some fundamental operations on posets and their affect on
homology.  The operations are  Alexander duality, join and direct
product. 
 \begth[Poset Alexander duality,
Stanley
\cite{st82,st1}]
\label{poalex} Let
$P$ be a 
$G$-poset  whose order complex triangulates an $n$-sphere.   If $Q$ is any induced
$G$-subposet of
$ P$ then for all
$i$,
$$\tilde H_i(Q) \cong_G \tilde H^{n-i-1}( P - Q) \otimes \tilde H_{n}( P). $$
If $Q$ is any induced subposet of $ P$ then
for all
$i$,
$$\tilde H_i(Q;\Z) \cong \tilde H^{n-i-1}( P - Q,\Z). $$
\enth

\begin{xca}\label{dcalex} Let $P$ be the poset of all graphs on node set $[n]$
ordered by inclusion of edge sets.  The symmetric group $\s_n$ acts on $P$ in the
obvious way, making $P$ an $\s_n$-poset.   Let $Q$ be an induced subposet  of  $P$
invariant under the action of
$\s_n$.
 Show that 
for all $i$,
$$\tilde H_i(\bar Q) \cong_{\mathfrak S_n} \tilde H_{{n \choose 2}-i-3}(\bar P - \bar Q)
\otimes 
\sgn_n^{\otimes n}.
$$
For example, if $Q$ is the poset of connected graphs on node set $[n]$ and $R$ is the
poset of disconnected graphs on node set $[n]$ then 
$$\tilde H_i(\bar Q) \cong_{\mathfrak S_n} \tilde H_{{n \choose 2}-i-3}(\bar R) \otimes 
\sgn_n^{\otimes n}.
$$
\end{xca}

\begin{xca}Given a simplicial complex $\Delta$ on vertex set $V$, its Alexander dual
$\Delta^{\vee}$ is defined to be the simplicial complex on $V$ consisting of
complements of nonfaces of
$\Delta$, i.e., $$ \Delta^{\vee} = \{V-F : F \subseteq V \mbox{ and }  F \notin \Delta\}.
$$ How are the homology of $\Delta$ and its Alexander dual related? 
\end{xca}

  Let $P$ be  a
$G$-poset and let $Q$ be an
$H$-poset.  Then the join 
 $P * Q$  and the product $P \times Q$ are  $(G \times H)$-posets with  respective $(G
\times H)$ actions given by  
$$(g,h) x = \begin{cases} gx \mbox{ if } x \in P \\ hx \mbox{ if } x \in Q \end{cases},$$
and  $$(g,h) (p,q) = (gp,hq).$$ We have the following poset version of the   
K\"unneth theorem of algebraic topology.  

\begth \label{jointh} Let $P$ be a $G$ poset and let  $Q$ be an $H$-poset. Then 
$$\tilde H_r(P*  Q)
\cong_{G\times H}\,\,
\bigoplus_{i} \tilde H_i( P) \otimes \tilde H_{r-i-1}(Q), $$
for all
$r$.
\enth

A result of Quillen \cite{q} states that if $P$ and $Q$  have
bottom elements
$\hat 0_P$ and  $\hat 0_Q$, respectively, then there is a $(G\times H)$-homeomorphism 
$$P\times  Q \setminus \{(\hat 0_P,\hat  0_Q)\} \cong_{G\times H} P\setminus
\{\hat 0_P\}
*
Q\setminus
\{\hat 0_Q\}.$$ 
 Walker \cite{wa81,wa88} proves the similar
result  that if
$P$ and
$Q$  have top elements as well as bottom elements then there is a $(G \times
H)$-homeomorphism
$$\overline{P\times Q} \cong_{G\times H} \bar P * \bar Q * A_2,$$ where
$A_2$ is a two element antichain on which the trivial group acts.    
Theorem~\ref{jointh} and Quillen's and Walker's results yield,

\begth \label{prodrep} Let $P$ be  a $G$-poset with a
bottom element
$\hat 0_P$ and let $Q$ be an
$H$-poset with a bottom element $\hat 0_Q$.  
 Then for all $r$,
$$\tilde H_r(P\times  Q \setminus \{(\hat 0_P,\hat  0_Q)\} \cong_{G\times H}\,\,
\bigoplus_{i} \tilde H_i(P\setminus \{\hat 0_P\}) \otimes \tilde H_{r-i -1}(Q\setminus
\{\hat 0_Q\}). $$
If $P$ and $Q$ also have top elements then for all $r$,
$$\tilde H_r(\overline{P\times  Q } ) \cong_{G\times H}\,\,
\bigoplus_{i} \tilde H_i(\bar P) \otimes \tilde H_{r-i-2}(\bar Q). $$
\enth

\begin{xca} \label{prodcm} Let $P$ and $Q$ be  posets with  bottom
elements. Show that $P \times
Q
$ is Cohen-Macaulay if and only if
$P$ and
$Q$ are Cohen-Macaulay.  This result also holds in the nonpure case but is more difficult
to prove; see \cite{bww2}.
\end{xca}

The products in Theorem~\ref{prodrep} are known as {\em reduced products}.  There is a
similar formula for the homology of ordinary direct products which follows from the
K\"unneth theorem and another $(G\times H)$-homeomorphism  
$$P \times Q \cong \| P\| \times \|Q \|$$
of Quillen \cite{q} and Walker \cite{wa88}.

\begth \label{ordprod} Let $P$ be a $G$-poset and let $Q$ be an $H$ poset.  Then for all
$r$,
$$H_r(P \times Q) \cong_{G\times H} \bigoplus_{i} H_i(P) \otimes H_{r-i}(Q).$$
\enth

Next we consider the $n$-fold product $P^{\times n}$ of a $G$-poset $P$. 
 The group that acts
on
$P^{\times n}$ is the wreath product
$\s_n[G]$, defined in Section~\ref{symm}.
 The action on $P^{\times n}$ is given by 
$$(g_1,g_2,\dots,g_n;\sigma) (p_1,p_2 ,\dots,p_n) =
(g_1p_{\sigma^{-1}(1)},g_2p_{\sigma^{-1}(2)},\dots,g_n p_{\sigma^{-1}(n)}).$$ 

For each pair of positive  integers $i$ and  $m$,  define the
$\s_{m}[G]$-module 
$$ H(P,i,m): = \begin{cases} S^{(m)} [ H_{i-1}(P)]\quad  &\mbox{if $i$
is odd}\\
\\ S^{(1^{m})}[ H_{i-1}(P)]\quad  &\mbox{if $i$ is
even.}\end{cases}
$$
For a partition  $\lambda$  with $m_i$ parts equal to
$i$ for each
$i$, let 
$$\s_\lambda[G]:=\ltimes_{i: m_i > 0} \s_{m_i}[G]$$ and define the
$\s_\lambda[G]$-module
$$ H(P,\lambda) := \bigotimes_{i: m_i > 0}  H(P,i,m_i).$$

We can now state an analog of Theorem~\ref{ordprod}.

\begth[Sundaram and Welker \cite{suwe96}]\label{powerord}  Let $P$ be a $G$-poset.  Then
for all
$r$,
$$ H_r(P^{\times n}) \cong_{\s_n[G]}
\bigoplus_{\scriptsize\begin{array}{c}\lambda
\vdash r+n\\ l(\lambda) = n \end{array} } 
H(P,\lambda)\uparrow_{\s_\lambda[G]}^{\s_n[G]}.$$ 
\enth

For an analog of Theorem~\ref{prodrep}, define the 
$\s_{m}[G]$-modules 
$$\tilde H(P,i,m,1): = \begin{cases} S^{(m)} [\tilde H_{i-2}(P)]\quad  &\mbox{if $i$
is odd}\\
\\ S^{(1^{m})}[\tilde H_{i-2}(P)]\quad  &\mbox{if $i$ is
even}\end{cases}
$$ and
$$\tilde H(P,i,m,2): = \begin{cases} S^{(1^m)} [\tilde H_{i-2}(P)]\quad  &\mbox{if $i$
is odd}\\
\\ S^{(m)}[\tilde H_{i-2}(P)]\quad  &\mbox{if $i$ is
even}\end{cases}
$$ and  the $\s_\lambda[G]$-modules 
$$\tilde H(P,\lambda,1) := \bigotimes_{i: m_i > 0} \tilde H(P,i,m_i,1)$$ and 
$$\tilde H(P,\lambda,2) := \bigotimes_{i: m_i > 0} \tilde H(P,i,m_i,2).$$

\begth[Sundaram and Welker \cite{suwe96}] \label{powerrep} Let $P$ be a $G$-poset with a
bottom element
$\hat 0$. Then for all $r,$
$$\tilde H_r(P^{\times n}-\{(\hat 0, \dots, \hat 0)\}) \cong_{\s_n[G]}
\bigoplus_{\scriptsize\begin{array}{c}\lambda
\vdash r+n+1\\ l(\lambda) = n \end{array} } \tilde H(P-\{\hat
0\},\lambda,1)\uparrow_{\s_\lambda[G]}^{\s_n[G]}.$$ If $P$ also has a top element then for
all
$r$,
$$\tilde H_r(\overline{ P^{\times n}}) \cong_{\s_n[G]}
\bigoplus_{\scriptsize\begin{array}{c}\lambda
\vdash r +2\\ l(\lambda) = n \end{array} } \tilde H(\bar
P,\lambda,2)\uparrow_{\s_\lambda[G]}^{\s_n[G]}.$$
\enth 

We remark that we have stated Theorems~\ref{powerord} and ~\ref{powerrep}  in a form
different from the original given  in \cite{suwe96} but completely equivalent. 

\begin{example}[The Boolean algebra $B_n$] \label{boolprod}  We apply the second part of
Theorem~\ref{powerrep} to $ B_n = B_1^{\times n}$.  The trivial group $\s_1$ acts on
$B_1$ and this induces an action of $\s_n[\s_1] = \s_n$ on $B_n$, which is precisely the
 action given in (\ref{snaction}). Clearly $\tilde H_{i-2} (\bar B_1)  = 0$ unless
$i =1$, in which case 
$H_{-1}(\bar B_1)$ is the trivial representation of $S_1$.    It follows that 
$$\tilde H(B_1,1,m,2) = S^{(1^m)}[S^{(1)}] = S^{(1^m)},$$
and $\tilde H(B_1,i,m,2) = 0$ for $ i\ne 1.$  Hence
$$\tilde H(B_1,\lambda,2) = \begin{cases} S^{(1^n)} &\mbox{ if } \lambda =
(1^n)\\ 0 &\mbox{otherwise}.\end{cases}$$   It follows that the only nonzero term in the
decomposition given in Theorem~\ref{powerrep} is the term corresponding to $\lambda =
(1^n)$ and $r = n-2$.  Hence we recover the fact that
$\tilde H_r(\bar B_n) $ is the sign representation  for $r=n-2$ and is $0$ otherwise.
\end{example} 

\begin{xca}  The  face lattice  $C_n$ of the $n$-cross-polytope with its top element
removed can be expressed as a  product,
$$C_n \setminus {\hat 1} = (C_1 \setminus {\hat 1})^{\times n}.$$   Since
$\s_2$ acts on
$C_1 \setminus \{\hat 1\}$, the product induces  an action of the
hyperoctahedral group
$\s_n[\s_2]$ on
$C_n
\setminus {\hat 1}$.
\begin{enumerate}
\item[(a)] Show that this action  is  the  action given in
Example~\ref{crossaction}. 
\item[(b)] Use Theorem~\ref{powerrep} to prove that 
$$\tilde H_r(\bar C_n) \cong_{\s_n[\s_2]} \begin{cases} S^{(1^n)}[S^{(1^2)}] &\mbox{if } r
= n-1 \\ 0 &\mbox{otherwise}.
\end{cases}$$
\item[(c)] Use Theorem~\ref{powerrep} to  obtain a generalization of 
 the formula in (b) for the poset $X_m^{\times n}$, where $X_m$
is the poset  consisting of a bottom element
$\hat 0
$ and $m$  atoms. 
\end{enumerate}
\end{xca}

\begin{example}[\cite{suwe96}] \label{lowintpart} In this example we compute the  homology
of lower intervals in the partition lattice 
$\Pi_n$.  Each lower interval is isomorphic to  a product of smaller partition
lattices,
$$[\hat 0,x] \cong \ltimes_i \Pi_i^{\times m_i}, $$
where $m_i$ is equal to the number of blocks of $x$ of size $i$. The stabilizer
$(\s_n)_x$ of
$x$ under the action of $\s_n$ on $\Pi_n$ is isomorphic to $\ltimes_i
\s_{m_i}[\s_i]$.  The  $(\s_n)_x$-poset $[\hat 0,x]$ is isomorphic to the $\ltimes_i
\s_{m_i}[\s_i]$-poset $\ltimes_i \Pi_i^{\times m_i}$.  We use the second parts of
Theorems~\ref{prodrep} and~\ref{powerrep} to compute the representation of $\ltimes_i
\s_{m_i}[\s_i]$ on the homology of  $(\hat 0, x)$.   This yields the following formula of
Lehrer and Solomon \cite{ls} for the only nonvanishing  homology,
$$\tilde H_{r(x)-2}(\hat 0,x) \cong_{\ltimes_i
\s_{m_i}[\s_i]} \,\,S^{(m_1)}[\tilde H_{-2}(\bar \Pi_1) ]
\otimes  S^{(1^{m_2})}[\tilde H_{-1}(\bar \Pi_2) ] \otimes S^{({m_3})}[\tilde
H_{0}(\bar \Pi_3) ] \otimes \cdots.$$  By summing over all set partitions of rank $r$ and
  taking Frobenius characteristic  one gets,
$$\ch \bigoplus_{\scriptsize \begin{array}{c} x \in \Pi_n \\ r(x) =r \end{array}} \tilde
H_{r-2}(\hat 0,x) = 
\mbox{degree
$n$ term in } (-1)^{r} h_{n-r}
\Big{[}
\sum_{m
\ge 1} (-1)^{m-1}
\ch
\tilde H_{m-3}(\bar \Pi_m) \Big{]}. 
$$ 
\end{example}

\begin{xca}[\cite{su,su94}] \label{low2}  Use Theorems~\ref{prodrep} and~\ref{powerrep} to
compute the representation of  $\ltimes_i
\s_{m_i}[\s_i]$ on $H_i(\hat 0,x)$, where $m_i$ is equal to the number of blocks of $x$ of
size $i$ and $(\hat 0, x)$ is an interval in the
\begin{enumerate}
\item[(a)]   $d$-divisible partition lattice $\Pi_{nd}^{0 \bmod d}$,
\item[(b)] $1 \bmod d$ partition poset $\Pi_{nd+1}^{1 \bmod d}$.
\end{enumerate}
\end{xca}

Computations such as those of Example~\ref{lowintpart} and Exercise~\ref{low2}  are used in
applications of equivariant  versions of the  Orlik-Solomon formula and the
Goresky-MacPherson formula.  This is discussed further in Section~\ref{arrangsec}.

 The examples  above don't adequately  demonstrate the  power of
the product theorems because in these examples homology occurs  in only one dimension. A
demonstration of the full power   can  be found in \cite{su94},
\cite{suwa} and
\cite{wa2}, where  representations on  the homology of intervals of the
$k$-equal partition lattice $\Pi_{n,k}$, the at least $k$ partition lattice $\Pi_{n}^{\ge
k}$, the general $j \bmod d$ partition poset $\Pi_{nd+j}^{j\bmod d}$, and the other
restricted block size partition posets of Section~\ref{whit}, are computed.

\section{Quillen fiber lemma} \label{quilfibsec} In a seminal 1978 paper, Quillen
introduced several  poset fiber theorems.  In this section we shall present the most basic
of these, known as the  ``the Quillen fiber lemma'', which has proved to be one of the most
useful tools in  poset topology.

Given two posets 
$P$ and $Q$, a map $f:P \to Q$ is called a {\it poset
map} if it is order preserving, i.e., $x \le_P y $ implies $f(x) \le_Q f(y)$. 
Given two $G$-posets 
$P$ and $Q$, a  poset map $f:P \to Q$ is called a $G$-{\it poset
map} if $g f(x) = f(gx)$ for all $x \in P$.
Given two simplicial complexes $\Delta$ and $\Gamma$, a {\em simplicial map}  is a map 
$f:\Delta \to
\Gamma$ that takes $0$-faces of $\Delta$ to $0$-faces of $\Gamma$ and preserves
inclusions, i.e., $f$ is a poset map from $P(\Delta) $ to $P(\Gamma)$. 
If $\Delta$ and $\Gamma$ are $G$-simplicial complexes and $f$ commutes with the $G$-action
then $f$ is said to be a $G$-simplicial map.     
A $G$-continuous map $f:X \to Y$ from $G$-space
$X$ to $G$-space $Y$ is a continuous map that commutes with the
$G$-action.  
Clearly,  a $G$-poset map $f:P \to Q$ induces a $G$-simplicial  map $f:\Delta(P) \to
\Delta(Q)$, which in turn induces a $G$-continuous map $f:\|\Delta(P)\| \to
\|\Delta(Q)\|$. 

For any element $x$ 
of a poset $ P$, let
$$ P_{>x} := \{y\in P : y > x\} 
\,\,\,\mbox{ and }\,\,\, P_{\ge x} := \{y \in P : y \ge x\}.$$    
The subsets
$ P_{< x}$  and  $ P_{\le x}$ are defined similarly.  

The basic version of the Quillen fiber lemma pertains to homotopy type.  Quillen also gave
an integral homology version and a equivariant homology version. There is also an
equivariant homotopy version due to Th\'evenaz and Webb \cite{tw}.  We present only the
homotopy version and the equivariant homology version. Recall that 
$\simeq$  denotes homotopy equivalence.

\begth[Quillen fiber lemma \cite{q}]  Let $f: P \to Q$ be a poset map.  If the fiber
$f^{-1}(Q_{\le y})$ is contractible for all $y \in Q$ then $$P \simeq Q.$$
\enth

\begth[Equivariant homology version \cite{q}] Let $f: P \to Q$ be a $G$-poset map.  If the
fiber
$f^{-1}(Q_{\le y})$ is acyclic (over $\C$) for all $y \in Q$ then 
$$\tilde H_r(P) \cong_G \tilde H_r(Q)$$
for all $r$.
\enth

Let's look at  Quillen's original example \cite{q}.  For a finite group $G$, let
$\mathcal S_p(G)$ be the poset of nontrivial $p$-subgroups of $G $ ordered by inclusion
and let
$\mathcal A_p(G)$ be the induced subposet of nontrivial elementary abelian $p$-subgroups of
$G$.  The poset $\mathcal S_p(G)$ and its order complex
 were studied by Brown in his work on group cohomology.  Quillen
proposed the smaller poset
$\mathcal A_p(G)$ as a way of studying the homotopy invariants of the Brown complex.  He
used the Quillen fiber lemma  to establish homotopy equivalence between the Brown
complex and the Quillen complex.  Both posets are $G$-posets, where $G$ acts by
conjugation.

\begth[Quillen \cite{q}] For any finite group $G$ and prime $p$,
\bq\label{brown}\mathcal A_p(G) \simeq \mathcal S_p(G) ,\eq
and 
\bq\label{brownhom}\tilde H_j(\mathcal A_p(G)) \cong_G \tilde H_j(\mathcal S_p(G)
)\quad\forall j.\eq 
\enth
\begin{proof}

 The
inclusion map
$$i:\mathcal A_p(G) \to \mathcal S_p(G) $$ is a $G$-poset map.
  Let us show that the fiber $i^{-1}(S_p(G)_{\le H})$ is contractible for
all $H  \in \mathcal S_p(G)$ so that we can apply the Quillen fiber lemma.   Note that 
$$i^{-1}(S_p(G)_{\le H}) = \mathcal A_p(H).$$  Since $H$ is a nontrivial
$p$-group, it has a nontrivial center.  Let $B$ be the subgroup of the center consisting
of all elements of order $1$ or $p$. Consider the poset map
$f:  \mathcal A_p(H) \to \{BA: A\in \mathcal A_p(H)\}$ defined by
$f(A) = BA$.  The fibers of this map are contractible since they all have maximum
elements.  So by the Quillen fiber lemma  $$\mathcal A_p(H) \simeq \{BA: A\in\mathcal
A_p(H)\}.$$
Since  $\{BA: A\in\mathcal
A_p(H)\}$ has    minimum element $B$, it is contractible.  So $ \mathcal A_p(H)$ and thus 
$i^{-1}(S_p(G)_{\le H})$ is contractible.  Hence by the Quillen fiber lemma, (\ref{brown})
holds and by the equivariant homology version, (\ref{brownhom}) holds.
\end{proof}

For certain finite groups $G$, the homology of the Quillen complex $\Delta(\mathcal
A_p(G))$ is well-understood; namely  for groups of Lie type.  Quillen \cite{q} shows that
for   groups
$G$ of Lie type in characteristic $p$,  the Quillen complex
$\Delta(\mathcal A_p(G))$  is homotopic to a simplicial complex known as the building for
$G$, which has the homotopy type of a wedge of spheres of a single dimension.
Webb
\cite{web} shows that  the unique nonvanishing $G$-equivariant homology of  $\mathcal
A_p(G)$  is the same as that of the building, which is known as the Steinberg
representation.   For general finite groups,
 the Quillen complex is not nearly as well-behaved, nor  well-understood.  In
fact, it was only recently shown by Shareshian \cite{sh04} that for certain primes the
integral homology of the Quillen complex of the symmetric group   has torsion.

The following
long-standing conjecture of Quillen  imparts a sense of the significance  of poset
topology in group theory.
\begin{con}[Quillen Conjecture] For any finite group $G$ and prime $p$, the poset $\mathcal
A_p(G)$ is contractible if and only if $G$ has a nontrivial normal $p$-subgroup.
\end{con}
The necessity of contractibility was proved by Quillen;  the sufficiency is still open.
However significant progress has been made by Aschbacher and Smith 
\cite{as93}.

In order to gain understanding of the Quillen complex for the symmetric group, Bouc
\cite{bo}  considered the induced subposet $T_n$ of $\mathcal S_2(\s_n)$ consisting of
nontrivial 
$2$-subgroups of $\s_n$ that contain a transposition, and the induced subposet $T^\prime_n$
of
$\mathcal A_2(\s_n)$ consisting of nontrivial elementary abelian $2$-subgroups generated by
transpositions.  It is easy to see that $T^\prime_n$ is $\s_n$-homeomorphic to the
matching complex $M_n$ discussed in Section~\ref{matchex}.  Bouc used the Quillen fiber
lemma to prove
 
\begin{equation} \label{bouceq1} T_n \simeq T^\prime_n \simeq M_n\end{equation} and for
all
$i$,
\begin{equation} \label{bouceq2} \tilde H_i(T_n) \cong_{\s_n} \tilde
H_i(T^\prime_n) \cong_{\s_n} \tilde H_i(M_n).\end{equation}  In addition to computing the
representation of the symmetric group on the homology of the matching complex, Bouc
discovered torsion in  the integral homology of the matching complex.  See \cite{shwa2},
\cite{wa3} and
\cite{jon05} for further results on torsion in the matching complex.

\begin{xca} Prove  (\ref{bouceq1}) and (\ref{bouceq2}).
\end{xca}

The usefulness of 
the matching complex   in understanding  the topology of the
Quillen and Brown complexes
 was recently demonstrated by Ksontini \cite{ks,ks2}.  He  used simple connectivity of
 $M_n$ for $n \ge 8$, which was proved by Bouc, to establish 
simple connectivity of 
$\mathcal
S_2(\s_n)$ for $n \ge 8$.  It was also  shown by Ksontini
\cite{ks,ks2,ks3}, Shareshian
\cite{sh04}, and Shareshian and Wachs \cite{shwa3} that a hypergraph version of the 
matching complex  is
 useful in studying the topology of 
$\Delta(\mathcal S_p(\s_n))$ when $p \ge 3$.

\begin{xca} Let $\bar x$ be a closure operator  on $P$, i.e. $ x  \le \bar x$ and
$\bar{\bar x} =
\bar x$ for all $x \in P$.  Define $\mbox{cl}(P) := \{x \in P : \bar x = x\}$.  
\begin{itemize}
\item[(a)] Use   
the Quillen fiber lemma to prove: $$ P
\simeq \mbox{cl}(P).$$
\item[(b)] Use (a) to show that $$\bar{\mathcal W}(n,k) \simeq \bar{\mathcal
N}(n,k) .$$ (We already know
this and (c) from Exercise~\ref{injexer}.)
\item[(c)] Derive an equivariant homology version of the homotopy result in (a) and use it
to prove $$\tilde H_i(\bar{\mathcal W}(n,k)) \cong_{\s_n} \tilde H_i( \bar{\mathcal
N}(n,k) )\,\,\,\, \forall i.$$ 
\end{itemize}
\end{xca}

\begin{xca}[Homology version of Rota's crosscut theorem]   Let $L$ be a   lattice
 and let
$M$ be the subposet of $\bar L$ consisting of non-$\hat 0$ meets of coatoms.  Prove:
 \begin{enumerate}
\item[(a)] $\bar L \simeq M$.  
\item[(b)]If $G$ is a group
acting on 
$L$ then for all $i$,
 $\tilde H_i(\bar L) \cong_G \tilde H_i( M)$. 
\item[(c)] Let $\Gamma(L)$ be the simplicial complex whose vertex set is  the set of
coatoms of
$L$ and whose faces are sets of coatoms whose meet is not $\hat 0$ (this is known as the
cross-cut complex of $L$).  Show $\bar L \simeq \Gamma(L)$ and if  $G$ is a group
acting on 
$L$ then for all $i$,
 $\tilde H_i(\bar L) \cong_G \tilde H_i( \Gamma(L))$.
\end{enumerate} 
\end{xca}  

The next two examples are connected with the combinatorics of  knot spaces and arose  in
the  work of Vassiliev \cite{v93,v94,v99}. 
\begin{example} \label{discex} Let  $\mbox{NCG}_n$ be the poset of disconnected
graphs on node set $[n]$ ordered by inclusion of edge sets.
 Let 
$$f: \overline{\mbox{NCG}}_n \to \bar \Pi_n$$ be the poset map such that  $f(G)$ is the
partition of $[n]$ whose blocks are the node sets of the connected components of $G$.
The
fibers of $f$ are given by

$$f^{-1}((\hat 0, \pi]) = (\hat 0, G_\pi],$$ where $G_\pi$ is the graph whose
connected components are cliques on the blocks of $\pi$. Since the fibers have maximal
elements they are contractible.  Hence by the Quillen fiber lemma
$$ \overline {\rm NCG}_n \simeq \overline \Pi_{n},$$
which implies that  $\overline{\rm NCG}_n$ has the homotopy type of a wedge of $(n-1)!$
spheres of dimension $n-3$.   
By the
equivariant homology version of the Quillen fiber lemma the only nonvanishing homology of 
$\overline{\rm NCG}_n$ is given by 
\bq\label{homncg} \tilde H_{n-3}(\overline{\rm NCG}_n) \cong_{\mathfrak S_n} \tilde
H_{n-3} (\Pi_n).\eq

\begin{xca} Let $\mbox{CG}_n$ be the poset of connected
graphs on node set $[n]$ ordered by inclusion of edge sets.  From Alexander duality and
(\ref{homncg}), we have
that he only
nonvanishing homology of
$\overline{\rm CG}_n$ is given by 
$$\tilde H_{{n-1 \choose 2}-1}(\overline{\rm CG}_n) \cong_{\mathfrak S_n} \sgn^{\otimes
n}
\otimes \tilde
H_{n-3} (\Pi_n).$$
Show that  for $n \ge 3$, the poset  $\overline{\rm CG}_n$ has the homotopy type of a
wedge of
$(n-1)!$
spheres of dimension ${n-1 \choose 2}-1$.

\end{xca}
  \end{example}

\begin{example} A graph is said to be 
$k$-{\it connected} if removal of any set of at most $k-1$ vertices leaves  the graph
connected.  Note that, unless $k =1$, this is a different graph property from that of being
$k$-edge-connected, which was discussed in Section~\ref{moddsec}. Every triangle tree is
$2$-edge-connected, but not
$2$-connected if the number of nodes is greater than $3$.    Let 
$\mbox{NCG}^k_n$ be the poset of graphs on node set $[n]$ that are not
$k$-connected, ordered by inclusion of edge sets.  The following result solves a problem of
Vassiliev \cite{v99} which arose from his study of knots. 
\begth[Babson, Bj\"orner, Linusson, Shareshian, Welker \cite{bblsw}, Turchin \cite{tur}]
\label{bblswth}
\begin{itemize}
\item[]
\item[(a)]
$\overline{\mbox{NCG}}^2_n$ has the homotopy type of a wedge of $(n-2)!$ spheres of
dimension
$2n-5$. 

\item[(b)]$\tilde H_{n-3}(\overline{\mbox{NCG}}^2_n)
\downarrow^{\s_n}_{\s{n-1}}\,\,\cong_{\s_{n-1}}\,\, \tilde H_{n-4} (\Pi_{n-1}) \otimes
\sgn_{n-1}$

\item[(c)] $\tilde H_{n-3}(\overline{\mbox{NCG}}^2_n)
\,\,\cong_{\s_{n}}\,\, (\tilde H_{n-4} (\Pi_{n-1})\uparrow_{\s_{n-1}}^{\s_n} -
\tilde H_{n-3} (\Pi_{n})) \otimes \sgn_n$
\end{itemize}  
\enth  

Note that the representation on the right side of  (c) is the  Whitehouse
module discussed in Section~\ref{moddsec}. We  describe the proof method of Babson,
Bj\"orner, Linusson, Shareshian, and Welker. Define the poset map
$$f: \overline{\mbox{NCG}}^2_n \to \overline{B_{n-1} \times \Pi_{n-1}}$$  by letting
$$f(G)= (S,\pi),$$ where
$S$ is the set of nodes joined to node  $n$ by an edge, and $\pi$ is the partition whose
blocks correspond to the connected components of $G \setminus \{n\}$.  For example, if $G$
is the not $2$-connected  graph in Figure~\ref{fignot2} then 
$$f(G) = ( \{2,4,6\},123/45/6).$$

\begin{figure}\begin{center}
\includegraphics[width=3.5cm]{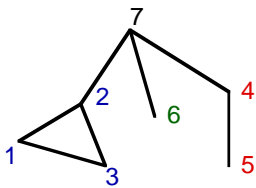}
\end{center}
\begin{center}\caption{Not 2-connected graph}\label{fignot2}
\end{center}
\end{figure}

Some of the fibers have maximum elements; so one sees immediately  that they are
contractible.  For example, 

$$f^{-1}((\hat 0,( \{2,4,6\},123/45/6)]) = (\hat 0, G],$$ where $G$ is the
graph in Figure~\ref{fignot2}.  But it is not so easy to see that the fibers without
maximum elements are contractible.  For example, the fiber $f^{-1}((\hat 0,(
\{1,2\},123456)])$ has no maximum.  To show that such fibers are contractible, Babson,
Bj\"orner, Linusson, Shareshian, and Welker used  discrete Morse theory, which had just
been introduced by Forman \cite{for}.  In fact, this  was the first of many applications of
discrete Morse theory in topological combinatorics.  

Now by the Quillen fiber lemma, we have $$\overline{\mbox{NCG}}^2_n \simeq
\overline{B_{n-1}
\times \Pi_{n-1}}.
$$
Since the product of shellable bounded posets is shellable (Theorem~\ref{prodshell}),
and the proper part of a shellable poset is shellable (Corollary~\ref{boundshell}),
$\overline{B_{n-1} \times
\Pi_{n-1}}$  has the homotopy type of a wedge of
$(2n-5)$-spheres.  The number of spheres is the absolute value of the M\"obius invariant,
which by Proposition~\ref{prodrule} and (\ref{parhom}) is $(n-2)!$.   Hence   Part (a) of
the theorem holds.

For Part (b), we use the equivariant homology version of the Quillen fiber lemma.  Clearly
the map
$f$ commutes with the permutations that fix $n$.
Hence, 
$$\tilde H_{2n-5}(\overline{\mbox{NCG}}^2_n)\downarrow_{\s_{n-1}}^{\s_n}\,\,
\cong_{\s_{n-1}}\,\,
\tilde H_{2n-5}(\overline{B_{n-1}
\times \Pi_{n-1}}).
$$
Theorem~\ref{prodrep} completes the proof of Part (b).

The proof of Part (c) uses Part (b) and  a computation of the fixed point M\"obius
invariant (Section~\ref{fixedp}).   

\begin{xca}\label{excg} Let ${\mbox{CG}}^k_n$ be the poset of $k$-connected graphs on node
set
$[n]$ ordered by inclusion of edge sets.  
\begin{enumerate}
\item [(a)] Use Theorem~\ref{bblswth} to determine the  homotopy type of
$\overline{\mbox{CG}}^2_n$.
\item [(b)](\cite{bblsw}) Show that $\overline{\mbox{CG}}^{n-2}_n$ is dual to the face
poset of the matching complex $M_n$.   Conclude that neither
$\overline{\mbox{CG}}^{n-2}_n$ nor $\overline{\mbox{NCG}}^{n-2}_n$ has the homotopy type 
of a wedge of spheres and that their integral homology has torsion.
\end{enumerate} 
\end{xca}

Shareshian \cite{sh01} uses discrete Morse theory to determine the homotopy type of
$\overline{\mbox{CG}}^2_n$ directly without resorting to Theorem~\ref{bblswth}, and to
obtain  bases for homology  of $\overline{\mbox{CG}}^2_n$, thereby solving another
problem of Vassiliev \cite{v99}.  The only  value of
$k$ other than
$k=1,2,n-2$, for which anything significant is known about the topology of the not
$k$-connected graph complex is $k=3$.  Discrete
Morse theory is used to obtain the following result.

\begth[Jonsson \cite{jon}]  $\overline{\mbox{NCG}}^3_n$ has the homotopy type of 
a wedge of $(n-3) á(n-2)!/2$ spheres of
dimension
$2n-4$.
\enth

\begin{prob} What can be said about the representation of $\s_n$ on the
homology of
$\overline{\mbox{NCG}}^3_n$.
\end{prob}

\begin{prob}(\cite{bblsw}) What can be said about the topology of
$\overline{\mbox{NCG}}^k_n$ for $3 < k <n-2$.  When does $\overline{\mbox{NCG}}^k_n$ fail
to be a wedge of spheres?  When does its integral homology
have torsion?  Recall that  the integral homology of
$\overline{\mbox{NCG}}^{n-2}_n$ has torsion (Exercise~\ref{excg} (b)). 
\end{prob}

\end{example}

\section{General poset fiber theorems} \label{genfibsec} In this section we present
homotopy and homology versions of a   fiber theorem of  Bj\"orner, Wachs and Welker,  which
generalizes  the Quillen fiber lemma and several other fiber theorems.  Recall the
definitions of $k$-connectivity and $k$-acyclicity, which appear after
Theorem~\ref{skelshell}.

\begth[Bj\"orner, Wachs, Welker \cite{bww1}]\label{main} Let $f:P \to Q$ be a poset map
such that  for all
$q \in Q$ the fiber
$\Delta(f^{-1}(Q_{\le q}))$ is
$l(f^{-1}(Q_{<q}))$-connected.   Then
\bq \label{eqmain}\Delta(P) \simeq \Delta(Q) \,\,\vee \,\,\left\{
\Delta(f^{-1}(Q_{\le q})) * \Delta(Q_{> q}) : q \in Q \right \},\eq   where $*$
denotes join defined in (\ref{joineq}) and  $\vee$ denotes the wedge (of $\Delta(Q)$ and
all
$\Delta(f^{-1}(Q_{\le q})) * \Delta(Q_{> q})$)  formed by  identifying the vertex $q$ 
in $\Delta(Q)$ with any
element of
$ f^{-1}(Q_{\le q})$,  for each $q \in Q$.  
\enth

For clarity,
let us remark that if $\Delta(Q)$ is connected then the space described on the 
right-hand side
of (\ref{eqmain}), which has $|Q|$ wedge-points, 
is homotopy equivalent to a {\it one-point} wedge, where arbitrarily chosen points of
$ f^{-1}(Q_{\le q})$, one for each $q\in Q$, are identified with some (arbitrarily 
chosen) point of $Q$. Thus (\ref{eqmain}) becomes
\bq \label{eqmaincon} \Delta(P) \simeq \Delta(Q) \vee \,\bigvee_{q \in
Q} 
\Delta(f^{-1}(Q_{\le q})) * \Delta(Q_{> q}).\eq

We will refer to a poset map $f:P \to Q$ such that for all $q \in Q$ the fiber 
$\Delta(f^{-1}(Q_{\le q}))$ is
$l(f^{-1}(Q_{<q}))$-connected as being {\it well-connected}.   
Note that the connectivity condition implies that each fiber
$f^{-1}(Q_{\le q})$ is nonempty.

\begin{figure}\begin{center}
\includegraphics[width=9cm]{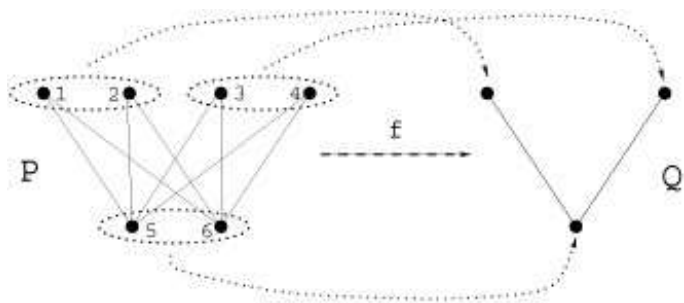}
\end{center}
\begin{center}\caption{A well-connected poset map}\label{figsample}
\end{center}
\end{figure}

\begin{example} \label{ex2} Let
$f : P
\rightarrow Q$ be the poset map depicted in Figure~\ref{figsample}.
For the two maximal elements of $Q$ the fiber
$\Delta (f^{-1}(Q_{\leq q}))$ is a $1$-sphere. For the bottom
element of $Q$ the fiber $\Delta (f^{-1}(Q_{\leq q}))$ is a $0$-sphere, and
$\Delta (Q_{> q})$ is a $0$-sphere too. So in either case 
$\Delta(f^{-1}(Q_{\le q})) * \Delta(Q_{> q})$ is homeomorphic to a
$1$-sphere.  
Hence the simplicial complex on the right side of (\ref{eqmain}) has a  $1$-sphere attached to each element of $Q$.  Thus
Theorem~\ref{main} determines
$\Delta(P)$  to have the homotopy type of  a wedge of three
$1$-spheres.   One can see this directly by observing that
$\Delta(P)$ is homeomorphic to two $1$-spheres intersecting in two points.
\end{example}

In \cite{bww1}, a version of the general fiber theorem for homology over the integers or
over any field is also given, as are equivariant homotopy and homology versions.  We 
state only the equivariant homology version here.

\begth[Bj\"orner, Wachs, Welker \cite{bww1}] \label{ghommain} Let $f:P
\to Q$ be a
$G$-poset map. If  $f^{-1}(Q_{\le q})$  is 
  $l(f^{-1}(Q_{< y}))$-acyclic
for all $ q \in Q$  then 
\bq \label{eqgmain}\\ \nonumber \tilde H_r(P) \,\,\cong_G \,\,\tilde H_r(Q)
\,\,\oplus 
 \bigoplus_{q\in Q/G} \,\,\bigoplus_{i=-1}^r
\left. \left (
\tilde H_i(f^{-1}(Q_{\le q}))
 \otimes
\tilde H_{r-i-1}(Q_{>q})\right)\right\uparrow_{G_q}^G. 
\eq
\enth

The general fiber theorems are proved using  techniques of \cite{wzz,zz}  involving
diagrams of spaces and spectral sequences, which were developed to study subspace
arrangements.  In Section~\ref{arrangsec}, we discuss the connection between the general
fiber theorems and subspace arrangements.

\begin{example}  Theorem~\ref{ghommain} can be applied to the well-connected poset map $f
: P
\rightarrow Q$ given in
Example~\ref{ex2}.  
Let $G$ be the cyclic group  $\s_2$ whose
non-identity element  acts by
$(1,2)(3,4)$ on $P$ and trivially on $Q$.   The map $f$ is clearly a $G$-poset map. For
each $q \in Q$, we have $G_q = \s_2$.  If $q$ is the bottom element of $Q$ then the
fiber
$f^{-1}(Q_{\le q})$ is $\s_2$-homeomorphic to a $0$-sphere on which $\s_2$ acts trivially. 
The same is true for  $Q_{>q}$.  It follows that the representation of $\s_2$ on $\tilde
H_0(f^{-1}(Q_{\le q}))
 \otimes
\tilde H_{0}(Q_{>q})$ is the trivial representation $S^{(2)}$. If $q$ is one of  the
maximal elements of
$Q$ then  the fiber
$f^{-1}(Q_{\le q})$ is
$\s_2$-homeomorphic to a circle with
$(1,2)(3,4)$ acting by reflecting the circle about the line spanned by a pair 
of antipodal points.  Hence  the representation of
 $\s_2$ on $\tilde H_1(f^{-1}(Q_{\le
q}))
$ is the sign representation $S^{(1^2)}$. Since $Q_{>q}$ is the empty
simplicial complex, the representation of
 $\s_2$ on $
\tilde H_{-1}(Q_{>q})$ is the trivial representation.   It follows that the representation
of
 $\s_2$ on $\tilde H_1(f^{-1}(Q_{\le
q}))
 \otimes
\tilde H_{-1}(Q_{>q})$ is $S^{(1^2)}$.   We conclude from (\ref{eqgmain}) that the
 $\mathfrak S_2$-module $\tilde H_1(P)$ 
decomposes into $S^{(2)} \oplus S^{(1^2)} \oplus  S^{(1^2)}$ and that   $\tilde
H_i(P) = 0$ for $i \ne 1$.
\end{example}

\begin{example} Consider the action of the hyperoctahedral group $\s_n[\Z_2]$ on the type B
partition lattice
$\Pi_n^B$. Define  the bar-erasing map 
\begin{equation} \label{barerase} f:
{\bar\Pi_n^B}
\to {\bar\Pi_{\{0,1,\dots,n\}}}\end{equation}
 by letting $f(\pi) $ be the partition obtained by erasing the bars from the 
barred partition $\pi$.   For
example
$$f(0 \,4 \,7\, / 1 \, \bar 5 \, 6 \, \bar 9\, / 2
\,\bar 3) = 0 \,4 \,7\, / 1 \,  5 \, 6 \,  9\, / 2
\, 3$$
It is clear that $f$ is a $\s_n[\Z_2]$-poset map. We establish well-connectedness in the
following exercise.

\begin{xca}[Wachs \cite{wa5}] \label{expn}
Let $P_n$ be the  
subposet of $\Pi_n^{ B}$ consisting of  barred partitions whose zero-block is  a
singleton, where the {\em zero-block} is the block that contains $0$. 
The poset $P_3$ with the zero-block suppressed is given in Figure~\ref{figpn}.  Show the
following.
\begin{itemize} 
\item[(a)]   The fibers of the bar-erasing map  have
the form
$$\Pi_{b_0}^{\small B} \times   P_{b_1} \times \cdots \times P_{b_k} -
\{\hat 0
\},$$ where $b_0 + b_1+\dots+b_k = n+1$.
\item[(b)]  $P_n$ is a geometric semilattice
\item[(c)]  The bar-erasing map (\ref{barerase}) is well-connected.
\item[(d)]  $|\mu(P_n \cup \{\hat 1\}| = (2(n-1) -1)!!$.
\end{itemize}
\end{xca}

\begin{figure}\begin{center}
\includegraphics[width=7cm]{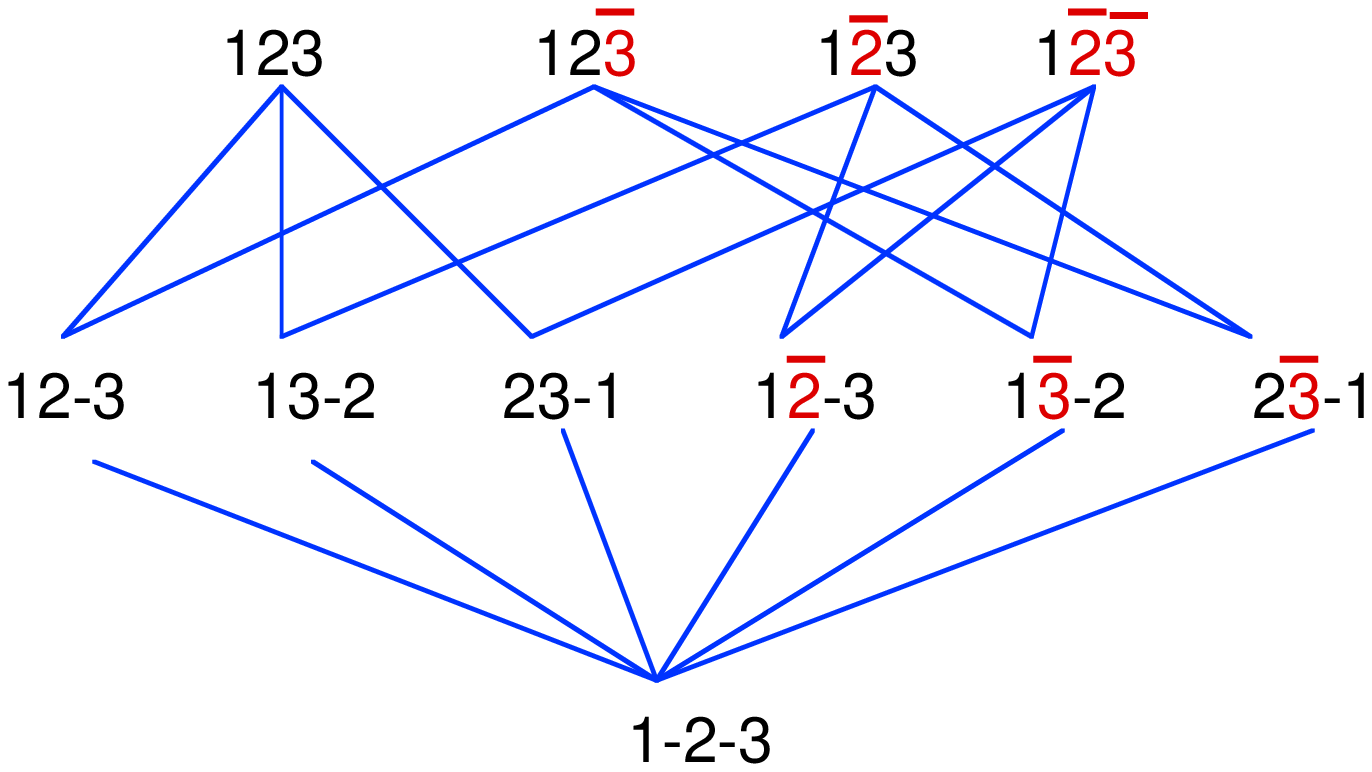}
\end{center}
\begin{center}\caption{$P_3$}\label{figpn}
\end{center}
\end{figure}

It follows from Exercise~\ref{expn} that the general fiber theorems can be used to
compute the homotopy type and $\s_n[\Z_2]$-equivariant homology  of $\Pi^B_n$.  Since 
$\overline{\Pi_n^B}$ has the homotopy type of a wedge of
$(2n-1)!!$ spheres of dimension  $n-2$ (see   Exercise~\ref{bnel}), the homotopy
equivalence (\ref{eqmaincon}) reduces to a combinatorial identity, which turns out to be 
easy to prove directly and  therefore not very interesting.  The equivariant homology
version of the general fiber theorem  does yield an interesting identity.  
Theorems~\ref{prodrep} and~\ref{powerrep} can be applied to the fibers given in Part (a)
of Exercise~\ref{expn}.  This and the homology version of the general fiber theorem are
used in \cite{wa5} to prove
\bq \label{dowling}\qquad \qquad \tilde
H_{n-2}(\overline{\Pi_n^{\small B}}) \,\,\cong_{\s_n[\Z_2]} \,\,
\bigoplus_{(\gamma_1,\dots,\gamma_k): \sum_i \gamma_i = n}  \tilde
H_{\gamma_1-2}(\overline{P_{\gamma_1}})  \bullet  \cdots \bullet \Black \tilde H_{\gamma_k
-2}(\overline{P_{\gamma_k}}) ,\eq
where $\bullet $ denotes induction product:
$$X \bullet Y := (X \otimes Y) \uparrow_{\s_n[\Z_2] \times \s_m[\Z_2]} ^{\s_{n+m}[\Z_2]}.$$
This example is considered in the more general setting of rank-selected Dowling
lattices in \cite{wa5}.
\end{example}

Nonpure versions of Cohen-Macaulay fiber
theorems of  Baclawski \cite{bac80}  and of Quillen \cite{q} follow from the general  fiber
theorems.  We state only the nonpure version of Baclawski's result.

\begth[Bj\"orner, Wachs, Welker \cite{bww1}] \label{bac} Let $P$ and $Q$ be semipure posets
and let 
$f:P \to Q$ be a surjective 
rank-preserving poset map.
  Assume that for all $q \in Q$, the fiber
$\Delta(f^{-1}(Q_{\le q}))$ is  Cohen-Macaulay over $\mathbf k$.  If $Q$ is sequentially
Cohen-Macaulay over $\mathbf k$, then so is
$P$.
\enth

\section{Fiber theorems and subspace arrangements} \label{arrangsec} 

Shortly after the  Goresky-MacPherson formula  (Theorem~\ref{GMform}) appeared,  Ziegler
and
\v Zivaljevi\'c obtained a  homotopy version,  and after that Sundaram
and Welker obtained an equivariant homology version.       In this section we show
that the Ziegler-\v Zivaljevi\'c  and   Sundaram-Welker formulas for linear subspace
arrangements  can  be easily derived from the general fiber theorems of
Section~\ref{genfibsec}. This connection should come as no surprise, since the method of
proof of the general fiber theorems is based on the  proofs of Ziegler-\v Zivaljevi\'c and
Sundaram-Welker; namely diagrams of spaces and spectral sequences. 

 Let $\mathcal A$ be a linear subspace arrangement   i.e., a finite
collection of  linear subspaces in Euclidean space
$\R^d$.    The {\it singularity link} $V_{\mathcal A}^o$ is defined as
$$V_{\mathcal A}^o = \mathbb S^{(d-1)} \cap \bigcup_{X \in \mathcal A} X ,$$  where
$\mathbb S^{(d-1)}$  is the unit $(d-1)$-sphere in $\R^d$.
Recall that  $L({\mathcal A})$ denotes the  intersection lattice of $\mathcal A$
(ordered by reverse inclusion).   Let
$\bar L({\mathcal A})$ denote 
$ L({\mathcal A})\setminus\{\hat 0,\hat 1\}$ if
$\mathcal A$ is essential (i.e., $\cap \mathcal A = \{0\}$), and 
$ L({\mathcal A})\setminus\{\hat 0\}$ otherwise.  

\begth[Ziegler \& \v Zivaljevi\'c \cite{zz}]\label{zzzz}  
Let  $\mathcal A$ be a linear subspace arrangement.  Then
\bq \label{zzeq} V_{\mathcal A}^o \,\,\simeq\,\, \Delta(\bar L({\mathcal A}))\,\, 
\vee \,\,\left\{ {\rm susp}^{\dim
x}(\Delta(\hat 0,x))
: x \in \bar L({\mathcal A})\right\},\eq where the wedge is formed by identifying 
each vertex $x$ in $\Delta(\bar L({\mathcal A}))$ with any  point
in
${\rm susp}^{\dim x}(\Delta(\hat 0,x))$. 
Consequently, for each $i$,
\bq \label{homzz} \tilde H_i(V_{\mathcal A}^o;\Z) \cong 
\bigoplus_{x \in  L(\mathcal A)\setminus \{\hat 0 \}} \tilde H_{i-\dim x}((\hat
0,x);\Z).\eq
\enth

\begin{proof}(given in \cite{bww1}). 
Suppose $\mathcal A = \{X_1,\dots,X_n\}$. 
Let $\mathcal H$ be an  essential hyperplane
arrangement in $\R^d$ such that each
$X_i$  is the intersection of a collection
of hyperplanes in $\mathcal H$.  The hyperplane arrangement $\mathcal H$ partitions
$\mathbb S^{(d-1)}$ into open cells of each dimension from $0$ to $d-1$.  Let $\mathcal
F(\mathcal H)$ be the poset of closures of the nonempty cells ordered by  inclusion.
Clearly
$\Delta(\mathcal F(\mathcal H))$ is a triangulation of $\mathbb S^{(d-1)}$.   For $X \in
\bar L({\mathcal A})$, let
$\mathcal F_{\mathcal H}(X)$ be the order ideal of  $\mathcal F(\mathcal H)$ consisting of
closed cells contained in $X$, and let
$\mathcal F_{\mathcal H} (\mathcal A)$ be the  order ideal  of
$\mathcal F(\mathcal H)$  consisting of closed cells contained in   $V_{\mathcal A}^o
$.  Note that   $\Delta(\mathcal F_{\mathcal H}( X))$ is a
triangulation of the $(\dim X -1)$-sphere $X \cap \mathbb S^{(d-1)}$ and $\Delta(\mathcal
F_{\mathcal H}(\mathcal A))$ is a triangulation of
$V_{\mathcal A}^o$.  

Now let
$$f :
\mathcal F_{\mathcal H}(\mathcal A) \to \bar L({\mathcal A})^*$$ be defined by
 $$f(\tau) = \bigcap_{i: X_i \supseteq \tau} X_i.$$ 
Clearly $f$ is order preserving. 
We claim that $f$ is well-connected. 
Observe that for all $X \in \bar L({\mathcal A})$, $$f^{-1}((\bar L({\mathcal A})^*)_{\le
X})=
\mathcal F_{\mathcal H}(X).$$ So $f^{-1}((\bar L({\mathcal A})^*)_{\le X})$ is a $(\dim X
-1)$-sphere, which is $(\dim X -2)$-connected. Since
$f^{-1}((\bar L({\mathcal A})^*)_{< X})$ has length at most $\dim X -2$, $f$ is indeed
well-connected,  
 and thus (\ref{zzeq}) follows from
Theorem~\ref{main}.  
\end{proof} 

The reason we refer to Theorem~\ref{zzzz} as  a homotopy version of the
Goresky-MacPherson formula is that    the Goresky-MacPherson formula
(for linear subspace arrangements) follows from  (\ref{homzz})  by
Alexander duality.  Indeed, one can simply  apply Theorem~\ref{poalex}  to   the
subposet
$\mathcal F_{\mathcal H}(\mathcal A)$ of
$\mathcal F(\mathcal H)$.  A homotopy version of the  Goresky-MacPherson formula
for general affine subspace arrangements is also given in \cite{zz}. 

Now let $G$ be a finite subgroup of the orthogonal group  $O_d$ that maps subspaces in
$\mathcal A$ to subspaces in $\mathcal A$. We say that $\mathcal A$ is a $G$-arrangement.
Clearly
$G$ acts as a group of poset maps on 
$L(\mathcal A)$ and as a group of 
 homomorphisms
on  $\tilde H_i(V_{\mathcal A}^o)$ and on $\tilde H_i( M_{\mathcal A})$.  For
each  $x \in L({\mathcal A})$, 
let $S_x$ be the $(\dim x -1)$-sphere  $x \cap
\mathbb S^{(d-1)}$.  Each $g$ in the  stabilizer
$G_x$ acts as an orientation preserving or reversing homeomorphism on $S_x$. 
Define $$\sgn_x(g) :=
\begin{cases} 1 \quad&\mbox{if $g$ is orientation preserving} \\
-1 \quad&\mbox{if $g$ is orientation reversing.}\end{cases}$$
By viewing 
$g$ as an element of $GL(x)$, we have
$\sgn_x(g) = \det(g) $. 
Note that
$\tilde H_{\dim x -1}(S_x)$ is a one dimensional representation of  $G_x$
whose character is given by $\sgn_x$.

\begth[Sundaram and Welker \cite{suwe}] \label{swthm} Let $\mathcal A$ be a $G$-arrangement
of linear subspaces in $\R^d$.  Then for all $i$,
\bq \label{homsw} \tilde H_i(V_{\mathcal A}^o) \cong_G 
\bigoplus_{x \in  L(\mathcal A)\setminus \{\hat 0 \}/G} (\tilde H_{i-\dim x}(\hat 0,x)
\otimes
\sgn_x)\uparrow_{G_x}^{G}.\eq
Consequently,
\bq \label{homsw} \tilde H^i(M_{\mathcal A}) \cong_G 
\bigoplus_{x \in  L(\mathcal A)\setminus \{\hat 0 \}/G} (\tilde H_{d-2-i-\dim x}(\hat 0,x)
\otimes
\sgn_{x}\otimes\sgn_{\hat 0})\uparrow_{G_x}^{G}.\eq
\enth

\begin{xca} Use Theorems~\ref{ghommain} and~\ref{poalex} to prove Theorem~\ref{swthm}. 
\end{xca}

In  \cite{suwe}  Theorem~\ref{swthm} is applied  to subspace arrangements whose
intersection lattices are  the $k$-equal
partition lattice discussed in 
 Section~\ref{kesec} and the $d$-divisible partition lattice discussed in 
Exercises~\ref{modd} and~\ref{equidivex}.   By computing the
multiplicity of the trivial representation in $\tilde H^{i}(M_{\mathcal A})$, Sundaram
and Welker obtain    the Betti number formula of Arnol'd given in (\ref{arn})  and an  
analogous  formula  for the space of monic polynomials of degree
$n$, with at least one root  multiplicity  not divisible by $d$.  For another recent
approach to studying these Betti numbers, see \cite{k02}.

In 
\cite{wa2},   Theorem~\ref{swthm} is applied  to the subspace arrangement  whose
intersection lattice  is the  restricted block size
partition lattice $\Pi_n^S \cup \{\hat 0\}$ of Theorem~\ref{genwac}.  Since each interval
$(\hat 0,x)$ of
 $\Pi_n^S$ is the  product of
smaller
$\Pi_m^S$,  one uses Theorem~\ref{genwac},
and the product formulas of Section~\ref{prodsec}  to compute the homology of $(\hat
0,x)$ (see Example~\ref{lowintpart} and Exercise~\ref{low2}).  This yields,

\begth[Wachs \cite{wa2}] \label{genwac2}
Suppose   $S\subseteq \{2,3,\dots\}$ is such that  $S$ and $\{s-\min S: s \in S\}$ are
closed under addition (eg., $S=\{k,k+1,\dots\}$ or more generally  $S=\{kd, (k+1)d,
\dots\})$. For
$n \in S$, let
$M_n^S$ be the manifold
$$\{{\mathbf z}:=(z_1,z_2,\dots,z_n) \in
\C^n: \mbox{ for some $i$ the number of occurrences of  $z_i$ is not in $S$}\}.$$
Then
\begin{eqnarray*}\quad \qquad\sum_{\scriptsize \begin{array}{c} m\in \Z\\ n \in
S\end{array}}
\ch\tilde H^{m}(M^{S}_{n})\,u^n\,(-t)^{2n-m-1} = \sum_{i \ge 1 }h_i t^i
\Big{[} \Big{(}\sum_{i \ge 1} h_{i}  \Big{)}^{[-1]}\Big{]}
\Big{[}\sum_{i \in S} h_{i}\, u^i\,\, t^ {\phi(i)} \Big{]},
\end{eqnarray*}
where $\phi(i) := \max\{j \in  \Z^+ : i-(j-1)\min S \in S\}$. (Our notation reflects the
fact that plethysm is associative.)
\enth

We leave it as an open problem to use this formula to compute the multiplicity of the
trivial representation in  $\tilde H^{m}(M^{S}_{n})$, and thereby obtain a formula
(which generalizes the above mentioned Sundaram-Welker formula) for the
Betti numbers of  the space of monic polynomials of degree
$n$  with at least one root  multiplicity  not in $S$.

\section{Inflations of simplicial complexes} \label{chessex} Let $\Delta$ be a simplicial
complex on vertex set
$[n]$ and let ${\mathbf m} = (m_1,\dots,m_n)$ be an $n$-tuple of positive
integers.   We form a new simplicial complex
$\Delta_{\mathbf m}$, called the ${\mathbf m}$-{\it inflation} of
$\Delta$, as follows.  The
 vertex set of $\Delta_{\mathbf m}$  is 
$\{(i,c) : i
\in [n], c \in [m_i] \}$ and the faces of $\Delta_{\mathbf m}$ are of the form
$\{(i_1,c_1),
\dots, (i_k,c_k)\}$ where  $\{i_1,\dots,i_k\}$ is a $(k-1)$-face of $\Delta$ and
$c_j
\in [m_{i_j}]$ for all $j = 1,\dots, k$.  We can think of
$c_j$ as a color assigned to vertex $i_j$ and of $\{(i_1,c_1), \dots, (i_k,c_k)\}$ as a coloring of the vertices of face
$\{i_1,\dots,i_k\}$.  A color for vertex
$i$ is chosen from $m_i$ colors. 

\begin{example} Let $P$ and $Q$ be the posets depicted in Figure~\ref{figsample}. 
Clearly $\Delta(P)$ is the  
$(2,2,2)$-inflation of $\Delta(Q)$. 
\end{example}

\begth[Bj\"orner, Wachs, and Welker \cite{bww1}] \label{color} Let $\Delta$ be a simplicial
complex on vertex set
$[n]$ and let
${\mathbf m}$ be an $n$-tuple positive integers. If $\Delta$ is connected then
$$\Delta_{\mathbf m} \simeq \bigvee_{F \in \Delta}  
(\mbox{susp}^{|F|}(\mbox{lk}_\Delta F))^{\vee\nu(F,{\mathbf m})},
$$ where  $\nu(F,{\mathbf m}) = \prod_{i \in F} (m_i -1)$, $\mbox{susp}^k(\Delta)$ denotes
the
$k$-fold suspension of $\Delta$, and $\Delta^{\lor k}$ denotes the $k$-fold wedge of
$\Delta$.  
\enth

We  prove this theorem in the following exercise.
\begin{xca} \label{inf} Let $f: \Delta_{{\mathbf m}} \to \Delta$ be the  simplicial map
that  sends each vertex $(i,c)$ of $\Delta_{{\mathbf m}}$ to 
vertex $i$ of $\Delta$. We call this the {\em deflating map}.  This can be viewed as  a
poset map 
$f:P(\Delta_{{\mathbf m}})
\to P(\Delta)$.
\begin{itemize}
\item[(a)]
 Show that the fiber 
$f^{-1}(P(\Delta)_{ \le F})$ is a geometric semilattice. 
\item[(b)] Show that each 
$f^{-1}(\overline{P(\Delta)}_{ \le F})$ is a wedge of
$\nu(F,\mathbf m)$ spheres of dimension $\dim F$. 
\item[(c)]  Use (b) and  the fact that the join operation is
 distributive over the wedge operation to prove Theorem~\ref{color}.
\item[(d)] Show that $\Delta_{\bf m}$ is sequentially Cohen-Macaulay if $\Delta$ is.
\end{itemize}

\end{xca}

\begth \label{infg} Let   $\Delta$ be a $G$-simplicial complex on vertex set $[n]$
 and let $\mathbf m$ be an $n$-tuple of positive integers. 
If
$G$ acts on the inflation $\Delta_{\mathbf m}$ and this action commutes
 with the deflating map,
then for all $r \in \Z$,
$$ \tilde H_r(\Delta_{\mathbf m}) \cong_G  \bigoplus_{F \in \Delta/G}  
\left(\tilde H_{|F|-1}(\langle F\rangle_{\mathbf m(F)}) \otimes
\tilde H_{r-|F|}(\mbox{lk}_\Delta F)\right)\uparrow_{G_F}^G,
$$ 
where $\langle F \rangle$ is the subcomplex generated by $F$ and
${\mathbf m}(F)$ is the subsequence 
$(m_{i_1},\dots,m_{i_t})$ of ${\mathbf m} = (m_1,\dots,m_n)$ for 
$F = \{i_1 < \dots < i_t\}$. 
\enth

Shareshian \cite{sh04} and Shareshian and Wachs \cite{shwa3} use Theorems~\ref{color} 
and~\ref{infg}, to derive information about the homology of the Quillen complex $\mathcal
A_p(\s_n)$ from a hypergraph version of the matching complex.  Another
application of Theorems~\ref{color} and~\ref{infg} can be found in 
\cite{py}.  Here we discuss the original example  that led to the general
fiber theorem, namely the colored chessboard complex \cite{wa4}.

Recall from Section~\ref{matchex} that the chessboard complex $M_{m,n}$ is the  
collection of rook placements on an $m\times n$ chessboard such that there is at most one 
rook in each row and each column.  
This complex
 arose in the work of Garst
\cite{gar} on Tits coset complexes, in the work of
\v Zivaljevi\'c and Vre\'cica \cite{zv}    on the  colored Tverberg problem,  in the
work of Reiner and Roberts \cite{rr} on Segre algebras, and in the work of Babson and
Reiner \cite{bare} on generalizations of Coxeter complexes.  
 Like
the matching complex, the chessboard complex  
 has been extensively studied in the literature; see the survey
article of Wachs
\cite{wa3}.
  We state just a few of the results on chessboard complexes here. 

\begth\label{chess} Let $1 \le m \le n$ and $\nu_{m,n} = \min \{m, \lfloor{m+n+1 \over
3}\rfloor\} -1$.
\begin{itemize}
\item[(a)]{\rm(Garst \cite{gar})} If $n \ge 2m-1$ then $M_{m,n}$ is Cohen-Macaulay.
\item[(b)]{\rm(Bj\"orner, Lov\'asz, Vre\'cica and \v Zivaljevi\'c \cite{blvz})} The complex
$M_{m,n}$ is $(\nu_{m,n}-1)$-connected.
\item[(c)] {\rm(Shareshian and Wachs \cite{shwa2})} For $n \le 2m-5$ and $m+n \equiv 1
\bmod 3$,
$$\tilde H_{\nu_{m,n}}(M_{m,n};\Z) \cong \Z_3.$$
\end{itemize}
\enth

 Part (a) follows from part (b), and Ziegler
\cite{z94} gives an alternative proof of part (b) by establishing shellability of the
$\nu_{m,n}$-skeleton of $M_{m,n}$.   
  Part (c) is one case of a more general torsion result, whose  proof 
uses  a chessboard complex analog of Bouc's representation theoretic
formula~(\ref{boucresult}) (due to Garst
\cite{gar} for top homology, and to Friedman and Hanlon
\cite{fh} in general) to construct a basis for the vector space
$\tilde H_{m-1}(M_{m,n})$.  This basis consists of fundamental cycles constructed
from pairs of standard Young tableaux via a classical combinatorial algorithm known as the
Knuth-Robinson-Schensted correspondence. We remark that Parts (b) and (c)  have analogs
for the matching complex \cite{blvz} \cite{bo} \cite{shwa2}. 

  Let $M_{m,n}^r$ be the
simplicial complex of   placements  of colored rooks on an $m\times n$
chessboard such that there is at most  one  rook in each row
and each column, and the colors come from the set $[r]$. The colored chessboard
complex  arose in Garst's work as a more general example of a Tits coset complex.  One
can easily see that
$M^r_{m,n}$ is the
$(r,r,\dots,r)$-inflation of
$M_{m,n}$.  

\begin{xca}[Wachs \cite{wa4}] Use Theorems~\ref{color}, \ref{infg}, and~\ref{chess} to
obtain a colored version  of Theorem~\ref{chess}, i.e. a generalization  to 
$M_{m,n}^r$, where
$r
\ge 1$.
\end{xca}

For open problems on the topology of chessboard complexes, matching
complexes and related complexes, see  \cite{wa3} and \cite{shwa2}.

%% file: biblio.tex
%

\bibliographystyle{amsalpha}